\documentclass{article}

\usepackage{longtable}
\usepackage{booktabs}
\def\tightlist{}

\usepackage{charter}
\usepackage{fancyhdr}
\usepackage[margin=1.6in]{geometry}

\usepackage{amsmath,amssymb}
\usepackage[colorlinks=true]{hyperref}
\usepackage{xcolor}
\definecolor{myurlcolor}{rgb}{0.6,0,0}
\definecolor{mycitecolor}{rgb}{0,0,0.8}
\definecolor{myrefcolor}{rgb}{0,0,0.8}
\hypersetup{linkcolor=myrefcolor,citecolor=mycitecolor,urlcolor=myurlcolor}

\usepackage[numbered]{bookmark}

\renewcommand{\texttt}[1]{%
  \begingroup
  \ttfamily
  \begingroup\lccode`~=`/\lowercase{\endgroup\def~}{/\discretionary{}{}{}}%
  \begingroup\lccode`~=`[\lowercase{\endgroup\def~}{[\discretionary{}{}{}}%
  \begingroup\lccode`~=`.\lowercase{\endgroup\def~}{.\discretionary{}{}{}}%
  \catcode`/=\active\catcode`[=\active\catcode`.=\active
  \scantokens{#1\noexpand}%
  \endgroup
}

\usepackage{titlesec}

\renewcommand{\thesection}{Week~\arabic{section}}
\titleformat{\section}[display]{\normalfont}{\Large\bfseries\thesection}{1em}{\large\normalfont}

\usepackage{titletoc}
\titlecontents{section}[0em]{\normalfont}{\bfseries\thecontentslabel\hspace{1em}\normalfont\small}{}{\titlerule*[0.3pc]{.}\small\itshape\thecontentspage}[\vspace{0.5em}]

\usepackage[toc]{multitoc}

\title{This Week's Finds in Mathematical Physics (1--50)}
\author{John Baez}
\date{January 19, 1993 to March 12, 1995}

\usepackage{bussproofs}

\usepackage{tikz}
\usetikzlibrary{knots}

\usepackage{etoolbox}
\AtBeginEnvironment{quote}{\itshape}

\usepackage{tikz-cd}

\usepackage{graphicx}
\usepackage[export]{adjustbox}

\usepackage[main]{embedall}

\begin{document}

\begin{titlepage}
  \begin{center}
    {\Huge\textbf{This Week's Finds in}}
  \\[0.7em]{\Huge\textbf{Mathematical Physics}}
  \\[1em]{\huge\textit{Weeks 1 to 50}}
  \\[4em]{\LARGE \textit{January 19, 1993} to \textit{March 12, 1995}}
  \\[4em]{\huge by John Baez}
  \\[0.5em]{\Large{typeset by Tim Hosgood}}
  \end{center}
  
\vskip 3em
  \[
  \begin{tikzpicture}[scale = 1.5]
    \begin{knot}[clip width=7]
      \strand[thick] (0,0)
        to [out=down,in=up] (1,-1)
        to [out=down,in=up] (2,-2)
        to (2,-3);
      \strand[thick] (1,0)
        to [out=down,in=up] (0,-1)
        to (0,-2)
        to [out=down,in=up] (1,-3);
      \strand[thick] (2,0)
        to (2,-1)
        to [out=down,in=up] (1,-2)
        to [out=down,in=up] (0,-3);
    \end{knot}
  \end{tikzpicture}
  \raisebox{6em}{\huge \quad =\quad}
  \begin{tikzpicture}[scale = 1.5]
    \begin{knot}[clip width=7]
      \strand[thick] (0,0)
        to (0,-1)
        to [out=down,in=up] (1,-2)
        to [out=down,in=up] (2,-3);
      \strand[thick] (1,0)
        to [out=down,in=up] (2,-1)
        to (2,-2)
        to [out=down,in=up] (1,-3);
      \strand[thick] (2,0)
        to [out=down,in=up] (1,-1)
        to [out=down,in=up] (0,-2)
        to (0,-3);
    \end{knot}
  \end{tikzpicture}
\] 
\end{titlepage}

\begin{center}
{\huge \bf Preface}
\end{center}

\vskip 1em

\noindent
These are the first 50 issues of \emph{This Week's Finds of Mathematical Physics}.
This series has sometimes been called the world's first blog, though it was originally 
posted on a ``usenet newsgroup'' called sci.physics.research --- a form of communication 
that predated the world-wide web.  I began writing this series 
as a way to talk about papers I was reading and writing, and in the first 50 issues 
I stuck closely to this format.   These issues focus rather tightly on quantum gravity,
topological quantum field theory, knot theory, and applications of $n$-categories to these
subjects.   However, there are also digressions into Lie algebras, elliptic curves, linear logic
and other subjects.

Tim Hosgood kindly typeset all 300 issues of \emph{This Week's Finds} in 2020. They
will be released in six installments of 50 issues each.  I have edited the issues here
to make the style a bit more uniform and also to change some references to preprints,
technical reports, etc.\ into more useful arXiv links.  This accounts for some anachronisms
where I discuss a paper that only appeared on the arXiv later.   

I thank Blake Stacey and Fridrich Valach for helping me fix many mistakes.  There are 
undoubtedly many still remaining.  If you find some, please contact me and I will try to fix them.

\tableofcontents

\hypertarget{week1}{%
\section{January 19, 1993}\label{week1}}

I thought I might try something that may become a regular feature on
sci.physics.research, if that group comes to be. The idea is
that I'll briefly describe the papers I have enjoyed this week in
mathematical physics. I am making no pretense at being exhaustive or
objective\ldots{} what I review is utterly a function of my own biases
and what I happen to have run into. I am not trying to ``rate'' papers
in any way, just to entertain some people and perhaps inform them of
some papers they hadn't yet run into. ``This week'' refers to when I
read the papers, not when they appeared (which may be much earlier or
also perhaps later, since some of these I am getting as preprints).

\begin{enumerate}
\def\labelenumi{\arabic{enumi})}
\item
  J.\ Scott Carter and Masahico Saito, ``Syzygies among elementary string 
  interactions in 2+1 dimensions'', \emph{Lett.\ Math.\ Phys.}
  \textbf{23} (1991), 287--300.

  J.\ Scott Carter and Masahico Saito, ``On formulations and solutions of simplex equations'', preprint.

  J.\ Scott Carter and Masahico Saito, ``A diagrammatic theory of knotted surfaces'', preprint.

 J.\ Scott Carter and Masahico Saito, ``Reidemeister moves for surface isotopies and their interpretations as moves to movies'', preprint.
\end{enumerate}

The idea here is to take what has been done for knots in
\(3\)-dimensional space and generalize it to ``knotted surfaces,'' that
is, embedded 2-manifolds in \(4\)-dimensional space. For knots it is
convenient to work with \(2\)-dimensional pictures that indicate over-
and under-crossings; there is a well-known small set of ``Reidemeister
moves'' that enable you to get between any two pictures of the same
knot. One way to visualize knotted surfaces is to project them down to
\(\mathbb{R}^3\); there are ``Roseman moves'' analogous to the
Reidemeister moves that enable to get you between any two projections of
the same knotted surface. Carter and Saito prefer to work with
``movies'' that display a knotted surface as the evolution of knots
(actually links) over time. Each step in such a movie consists of one of
the ``elementary string interactions.'' They have developed a set of
``movie moves'' that connect any two movies of the same knotted surface.
These papers contain a lot of fascinating pictures! And there does seem
to be more than a minor relation to string theory. For example, one of
the movie moves is very analogous to the 3rd Reidemeister move --- which
goes 
\[
  \begin{tikzpicture}
    \begin{knot}[clip width=7]
      \strand[thick] (0,0)
        to [out=down,in=up] (1,-1)
        to [out=down,in=up] (2,-2)
        to (2,-3);
      \strand[thick] (1,0)
        to [out=down,in=up] (0,-1)
        to (0,-2)
        to [out=down,in=up] (1,-3);
      \strand[thick] (2,0)
        to (2,-1)
        to [out=down,in=up] (1,-2)
        to [out=down,in=up] (0,-3);
    \end{knot}
  \end{tikzpicture}
  \raisebox{4em}{\quad=\quad}
  \begin{tikzpicture}
    \begin{knot}[clip width=7]
      \strand[thick] (0,0)
        to (0,-1)
        to [out=down,in=up] (1,-2)
        to [out=down,in=up] (2,-3);
      \strand[thick] (1,0)
        to [out=down,in=up] (2,-1)
        to (2,-2)
        to [out=down,in=up] (1,-3);
      \strand[thick] (2,0)
        to [out=down,in=up] (1,-1)
        to [out=down,in=up] (0,-2)
        to (0,-3);
    \end{knot}
  \end{tikzpicture}
\] 
I won't try to draw the corresponding movie move, but just as the 3rd
Reidemeister move is the basis for the Yang--Baxter equation
\(R_{23}R_{13}R_{12} = R_{12}R_{13}R_{23}\) (the subscripts indicate
which strand is crossing which), the corresponding movie move is the
basis for a variant of the ``Frenkel--Moore'' form of the ``Zamolodchikov
tetrahedron equation'' which first arose in string theory. This variant
goes like
\(S_{124}S_{135}S_{236}S_{456} = S_{456}S_{236}S_{135}S_{124}\), and
Carter and Saito draw pictures that make this equation almost as obvious
as the Yang--Baxter equations.

In any event, this is becoming a very hot subject, since topologists are
interested in generalizing the new results on knot theory to higher
dimensions, while some physicists (especially Louis Crane) are convinced
that this is the right way to tackle the ``problem of time'' in quantum
gravity (which, in the loop variables approach, amounts to studying the
relationship of knot theory to the 4th dimension, time). In particular,
Carter and Saito are investigating how to construct solutions of the
Zamolodchikov equations from solutions of the Yang--Baxter equation ---
the goal presumably being to find invariants of knotted surfaces that
are closely related to the link invariants coming from quantum groups.
This looks promising, since Crane and Yetter have just constructed a
4-dimensional topological quantum field theory from the quantum
\(\mathrm{SU}(2)\). But apparently nobody has yet done it.

Lovers of category theory will be pleased to learn that the correct
framework for this problem appears to be the theory of \(2\)-categories.
These are categories with objects, morphisms between objects, and also
``2-morphisms'' between objects. The idea is simply that tangles are
morphisms between sets of points (i.e., each of the tangles in the
picture above are morphisms from 3 points to 3 points), while surfaces
in \(\mathbb{R}^4\) are \(2\)-morphisms between tangles. The instigators
of the \(2\)-categorical approach here seem to be Kapranov and
Voevodsky, whose paper ``\(2\)-categories and Zamolodhikov tetrahedra
equations'', to appear in \emph{Proc.\ Symp.\ Pure Math.}, is something I will
have to get ahold of soon by any means possible (I can probably nab it
from Oleg Viro down the hall; he is currently hosting Kharmalov, who is
giving a series of talks on knotted surfaces at \(2\)-categories here at
UCR). But it seems to be Louis Crane who is most strongly proclaiming
the importance of \(2\)-categories in \emph{physics}.

\begin{enumerate}
\def\labelenumi{\arabic{enumi})}
\setcounter{enumi}{1}
\tightlist
\item
  Jorge Pullin, ``Knot theory and quantum gravity in loop space: a primer'', to appear in 
  \emph{Proceedings of the Vth Mexican School of
  Particles and Fields}, ed.~J. L. Lucio, World Scientific, Singapore,
  now available as
  \href{https://arxiv.org/abs/hep-th/9301028}{\texttt{hep-th/9301028}}.
\end{enumerate}

This is a review of the new work on knot theory and the loop
representation of quantum gravity. Pullin is among a group who has been
carefully studying the ``Chern--Simons state'' of quantum gravity, so his
presentation, which starts with a nice treatment of the basics, leads
towards the study of the Chern--Simons state. This is by far the
best-understood state of quantum gravity, and is defined by
\(\mathrm{SU}(2)\) Chern--Simons theory in terms of the connection
representation, or by the Kauffman bracket invariant of knots in the
loop representation. It is a state of Euclideanized quantum gravity with
nonzero cosmological constant, and is not invariant under CP. Ashtekar
has recently speculated that it is a kind of ``ground state'' for
gravity with cosmological constant (evidence for this has been given by
Kodama), and that its CP violation may be a ``reason'' for why the
cosmological constant is actually zero (this part is extremely
speculative). Louis Crane, on the other hand, seems convinced that the
Chern--Simons state (or more generally states arising from modular tensor
categories) is \textbf{the wavefunction of the universe}. In any event,
it's much nicer to have one state of quantum gravity to play with than
none, as was the case until recently.

\begin{enumerate}
\def\labelenumi{\arabic{enumi})}
\setcounter{enumi}{2}
\tightlist
\item
  Lee Smolin, ``Time, measurement and information loss in quantum cosmology'', 
  preprint available as
  \href{https://arxiv.org/abs/gr-qc/9301016}{\texttt{gr-qc/9301016}}.
\end{enumerate}
\noindent
This is, as usual for Smolin, a very ambitious paper. It attempts to
sketch a solution of some aspects of the problem of time in quantum
gravity (in terms of the loop representation). I might as well quote
from the introduction:

\begin{quote}

Thus, to return to the opening question, if we are, within a
nonperturbative framework, to ask what happens after a black hole
evaporates, we must be able to construct spacetime diffeomorphism
invariant operators that can give physical meaning to the notion of
``after the evaporation.'' Perhaps I can put it in the following way:
the questions about loss of information or breakdown of unitary
evolution rely, implicitly, on a notion of time. Without reference to
time it is impossible to say that something is being lost. In a quantum
theory of gravity, time is a problematic concept which makes it
difficult to even ask such questions at the nonperturbative level,
without reference to a fixed spacetime manifold. {[}I would prefer to
say ``fixed background metric'' --- JB{]} The main idea, which it is the
purpose of this paper to develop, is that the problem of time in the
nonperturbative framework is more than an obstacle that blocks any easy
approach to the problem of loss of information in black hole
evaporation. It may be the key to its solution.

As many people have argued, the problem of time is indeed the conceptual
core of the problem of quantum gravity. Time, as it is conceived in
quantum mechanics is a rather different thing than it is from the point
of view of general relativity. The problem of quantum gravity,
especially when put in the cosmological context, requires for its
solution that some single concept of time be invented that is compatible
with both diffeomorphism invariance and the principle of superposition.
However, looking beyond this, what is at stake in quantum gravity is
indeed no less and no more than the entire and ancient mystery: What is
time? For the theory that will emerge from the search for quantum
gravity is likely to be the background for future discussions about the
nature of time, as Newtonian physics has loomed over any discussion
about time from the seventeenth century to the present.

I certainly do not know the solution to the problem of time. Elsewhere I
have speculated about the direction in which we might search for its
ultimate resolution. In this paper I will take a rather different point
of view, which is based on a retreat to what both Einstein and Bohr
taught us to do when the meaning of a physical concept becomes confused:
reach for an operational definition. Thus, in this paper I will adopt
the point of view that time is precisely no more and no less than that
which is measured by physical clocks. From this point of view, if we
want to understand what time is in quantum gravity then we must
construct a description of a physical clock living inside a relativistic
quantum mechanical universe.

\end{quote}

Technically speaking, what Smolin does is roughly as follows. He
considers quantum gravity coupled to matter, modelled in such a way that
the Hilbert space is spanned by states labelled by isotopy classes of:
any number \(N\) loops in a compact 3-manifold \(M\) (``space'') and
\(N\) surfaces with boundary in \(M\). (This trick is something I hadn't
seen before, though Smolin gives references to it.) He then introduces a
``clock field,'' which is just a free scalar field coupled to the
gravity, and does gauge-fixing to see what evolution with respect to
this clock field looks like. I will have to read this a number of times!

\hypertarget{week2}{%
\section{January 24, 1993}\label{week2}}

Well, this week I have had guests and have not been keeping up with the
literature. So ``this week's finds'' are mostly papers that have been
sitting around in my office and that I am now filing away.

\begin{enumerate}
\def\labelenumi{\arabic{enumi})}
\tightlist
\item
  Daniel Armand Ugon,
  Rodolfo Gambini, and Pablo Mora, ``Link invariants for intersecting loops'',
  October 1992 preprint, available from
  Gambini, Instituto de F\'isica, Facultad de Ciencias, Trist\'an Narvaja
  1674, Montevideo, Uruguay.
\end{enumerate}
\noindent
The authors generalize the standard trick for getting link invariants
from solutions of the Yang--Baxter equations, and show how to get link
invariants applicable to generalized links with 4-valent or 6-valent
vertices, that is, transverse double points, like \[
  \begin{tikzpicture}
    \begin{knot}[clip width=0]
      \strand[thick] (0,0)
      to [out=down,in=up] (1,-1);
      \strand[thick] (1,0)
      to [out=down,in=up] (0,-1);
    \end{knot}
  \end{tikzpicture}
\] and transverse triple points. This involves working with a
generalization of the braid group that includes generators for these
vertices as well as the usual generators for crossings. In this case of
4-valent vertices, rigorously working out the generators and relations
was done by Joan Birman in the paper below, but various people had used
the answer already. The case of triple points is very important in
physics due to the connection with the loop representation of quantum
gravity (which is what Gambini is working on these days). In this
representation, states are invariants of (possibly generalized) links,
and only by considering links with triple points can one define
operators such as the ``total volume of the universe''.

It is thus quite interesting that the authors make progress on
determining the ``right'' extension of the HOMFLY polynomial invariant
of links to links with transverse triple points --- that is, the
extension that one gets by doing calculations in \(\mathrm{SU}(n)\)
Chern--Simons theory. In a special case, namely when the HOMFLY
polynomial reduces to the Kauffman bracket (which corresponds to the Lie
group \(\mathrm{SU}(2)\)), one gets a state of quantum gravity that has
been under extensive investigation these days. The authors compute the
Kaufmann bracket of links with triple points using first order
perturbation theory in Chern--Simons theory. A nonperturbative
calculation would be very good to have!

\begin{enumerate}
\def\labelenumi{\arabic{enumi})}
\setcounter{enumi}{1}
\tightlist
\item
  Joan Birman, ``New points of view in knot theory'', \emph{Bull. Amer. Math. Soc.}
  \textbf{28} (1993), 253--287.  Also available as
  \href{https://arxiv.org/abs/math/9304209}{\texttt{math/9304209}}.
\end{enumerate}
\noindent
This is a nice review of the recent work on Vassiliev invariants of
links. Given an invariant of oriented links, one can extend it to links
with arbitrarily many double points by setting the value of the
invariant on a link with a double point \[
  \begin{tikzpicture}
    \begin{knot}[clip width=0]
      \strand[thick] (0,0)
      to [out=down,in=up] (1,-1);
      \strand[thick] (1,0)
      to [out=down,in=up] (0,-1);
    \end{knot}
  \end{tikzpicture}
\] to be the invariant of the link with the double point changed to \[
  \begin{tikzpicture}
    \begin{knot}[clip width=7]
      \strand[thick] (0,0)
      to [out=down,in=up] (1,-1);
      \strand[thick] (1,0)
      to [out=down,in=up] (0,-1);
      \flipcrossings{1}
    \end{knot}
  \end{tikzpicture}
\] minus the invariant of the link with the double point changed to \[
  \begin{tikzpicture}
    \begin{knot}[clip width=7]
      \strand[thick] (0,0)
      to [out=down,in=up] (1,-1);
      \strand[thick] (1,0)
      to [out=down,in=up] (0,-1);
    \end{knot}
  \end{tikzpicture}
\] Note that the link has to be oriented for this rule to make sense,
and the strands shown in the pictures above should be pointing
\emph{downwards}. Now, having made this extension, we say a link
invariant is a Vassiliev invariant of degree \(n\) if it vanishes on all
links with \(n+1\) or more double points.

It is interesting that this rule for extending link invariants to links
with double points, when applied to the Kauffman bracket, does
\emph{not} give the extension computed by Gambini et al in the paper
above. (This is not surprising, actually, but it shows that some
interesting things are going on in the subject of invariants for links
with self-intersections and Chern--Simons theory.)

\begin{enumerate}
\def\labelenumi{\arabic{enumi})}
\setcounter{enumi}{2}
\tightlist
\item
  Louis Kauffman and
  P.\ Vogel, ``Link polynomials and a graphical calculus'', \emph{Jour.\ of Knot Theory and its Ramifications},
  \textbf{1} (1992), 59--104.
\end{enumerate}

This is another nice treatment of link invariants for generalized links
with self-intersections. It concentrates on the famous link invariants
coming from Chern--Simons theory --- the HOMFLY polynomial (from
\(\mathrm{SU}(n)\)) and the Kauffman polynomial (from
\(\mathrm{SO}(n)\)). Lots of good pictures.

And, switching back to the category theory, \(2\)-categories, and the
like, let me list these before filing them away:

\begin{enumerate}
\def\labelenumi{\arabic{enumi})}
\setcounter{enumi}{3}
\item
   Louis Crane, ``Categorical physics'', available as
  \href{https://arxiv.org/abs/hep-th/9301061}{\texttt{hep-th/9301061}}.

  Louis Crane and David Yetter, A categorical construction of 4d topological quantum field
  theories, available as
  \href{https://arxiv.org/abs/hep-th/9301062}{\texttt{hep-th/9301062}}.

  Louis Crane and Igor
  Frenkel, Hopf Categories and their representations, draft version.

  Louis Crane and Igor Frenkel, Categorification and the construction of topological quantum field
  theory, draft version.
\end{enumerate}

These outline Louis Crane's vision of an approach to generally covariant
\(4\)-dimensional quantum field theories (e.g.~quantum gravity or a
``theory of everything'') based on \(2\)-categories. ``Categorical
physics'' sketches the big picture, while the paper with Yetter provides
a juicy mathematical spinoff --- the first known four-dimensional TQFT,
based on the representations of quantum \(\mathrm{SU}(2)\) and very
similar in spirit to the Turaev--Viro construction of a 3d TQFT from
quantum \(\mathrm{SU}(2)\). The papers with Frenkel (apparently still
not in their final form) describe the game plan and hint at marvelous
things still to come. The conjecture is stated: ``a 4d TQFT can be
reconstructed from a tensor \(2\)-category''. This follows up on Crane's
earlier work on getting 3d TQFTs from modular tensor categories (big
example: the categories of representations of quantum groups at roots of
unity). And the authors define the notion of a Hopf category, show how
the category of module categories of a Hopf category is a tensor
\(2\)-category, and use ``categorification'' to turn the universal
enveloping algebra of a quantum group into a Hopf category. Sound
abstract? Indeed it is. But the aim is clear: to cook up 4d TQFTs from
quantum groups. Such quantum field theories might be physically
important; indeed, the one associated to \(\mathrm{SU}(2)\) is likely to
have a lot to do with quantum gravity.

I am currently perusing Kapranov and Voevodsky's massive paper on
\(2\)-categories, which seems to be the starting point for Crane's above
papers and also those of Carter/Saito that I mentioned last week. Next
week I should post an outline of what this paper does.

\begin{enumerate}
\def\labelenumi{\arabic{enumi})}
\setcounter{enumi}{4}
\tightlist
\item
   S.\ W.\ Hawking, R.\ Laflamme and G.\ W.\
  Lyons, ``The origin of time asymmetry'', preprint available as
  \href{https://arxiv.org/abs/gr-qc/9301017}{\texttt{gr-qc/9301017}}.
\end{enumerate}
\noindent
I haven't had a chance to read this one yet but it looks very ambitious
and is likely to be interesting. Let me just quote from the introduction
to get across the goal:

\begin{quote}
The laws of physics do not distinguish the future from the past
direction of time. More precisely, the famous CPT theorem says that the
laws are invariant under the combination of charge conjugation, space
inversion and time reversal. In fact effects that are not invariant
under the combination CP are very weak, so to a good approximation, the
laws are invariant under the time reversal operation T alone. Despite
this, there is a very obvious difference between the future and past
directions of time in the universe we live in. One only has to see a
film run backward to be aware of this.

There are several expressions of this difference. One is the
so-called psychological arrow, our subjective sense of time, the fact
that we remember events in one direction of time but not the other.
Another is the electromagnetic arrow, the fact that the universe is
described by retarded solutions of Maxwell's equations and not advanced
ones. Both of these arrows can be shown to be consequences of the
thermodynamic arrow, which says that entropy is increasing in one
direction of time. It is a non trivial feature of our universe that it
should have a well defined thermodynamic arrow which seems to point in
the same direction everywhere we can observe. Whether the direction of
the thermodynamic arrow is also constant in time is something we shall
discuss shortly.

There have been a number of attempts to explain why the universe should
have a thermodynamic arrow of time at all. Why shouldn't the universe be
in a state of maximum entropy at all times? And why should the direction
of the thermodynamic arrow agree with that of the cosmological arrow,
the direction in which the universe is expanding? Would the
thermodynamic arrow reverse, if the universe reached a maximum radius
and began to contract?

Some authors have tried to account for the arrow of time on the basis of
dynamic laws. The discovery that CP invariance is violated in the decay
of the K meson, inspired a number of such attempts but it is now
generally recognized that CP violation can explain why the universe
contains baryons rather than anti baryons, but it cannot explain the
arrow of time. Other authors have questioned whether quantum gravity
might not violate CPT, but no mechanism has been suggested. One would
not be satisfied with an ad hoc CPT violation that was put in by hand.

The lack of a dynamical explanation for the arrow of time suggests that
it arises from boundary conditions. The view has been expressed that the
boundary conditions for the universe are not a question for Science, but
for Metaphysics or Religion. However that objection does not apply if
there is a sense in which the universe has no boundary. We shall
therefore investigate the origin of the arrow of time in the context of
the no boundary proposal of Hartle \& Hawking. This was formulated in
terms of Einsteinian gravity which may be only a low energy effective
theory arising from some more fundamental theory such as superstrings.
Presumably it should be possible to express a no boundary condition in
purely string theory terms but we do not yet know how to do this.
However the recent COBE observations indicate that the perturbations
that lead to the arrow of time arise at a time during inflation when the
energy density is about \(10^{-12}\) of the Planck density. In this
regime, Einstein gravity should be a good approximation.

\end{quote}

\noindent
I'll skip some more technical stuff on the validity of perturbative
calculations in quantum gravity\ldots

\begin{quote}

One can estimate the wave functions for the perturbation modes by
considering complex metrics and scalar fields that are solutions of the
Einstein equations whose only boundary is the surface \(S\). When \(S\)
is a small three sphere, the complex metric can be close to that of part
of a Euclidean four sphere. In this case the wave functions for the
tensor and scalar modes correspond to them being in their ground state.
As the three sphere \(S\) becomes larger, these complex metrics change
continuously to become almost Lorentzian. They represent universes with
an initial period of inflation driven by the potential energy of the
scalar field. During the inflationary phase the perturbation modes
remain in their ground states until their wave lengths become longer
than the horizon size. The wave function of the perturbations then
remains frozen until the horizon size increases to be more than the wave
length again during the matter dominated era of expansion that follows
the inflation. After the wave lengths of the perturbations come back
within the horizon, they can be treated classically.

This behaviour of the perturbations can explain the existence and
direction of the thermodynamic arrow of time. The density perturbations
when they come within the horizon are not in a general state but in a
very special state with a small amplitude that is determined by the
parameters of the inflationary model, in this case, the mass of the
scalar field. The recent observations by COBE indicate this amplitude is
about \(10^{-5}\). After the density perturbations come within the
horizon, they will grow until they cause some regions to collapse as
proto-galaxies and clusters. The dynamics will become highly non linear
and chaotic and the coarse grained entropy will increase. There will be
a well defined thermodynamic arrow of time that points in the same
direction everywhere in the universe and agrees with the direction of
time in which the universe is expanding, at least during this phase.

The question then arises: If and when the universe reaches and maximum
size, will the thermodynamic arrow reverse? Will entropy decrease and
the universe become smoother and more homogeneous during the contracting
phase?

\end{quote}

\noindent
I'll skip some stuff on why Hawking originally thought entropy had to
decrease during the Big Crunch if the no-boundary proposal were
correct\ldots{} and why he no longer thinks so.

\begin{quote}

The thermodynamic arrow will agree with the cosmological arrow for half
the history of the universe, but not for the other half. So why is it
that we observe them to agree? Why is it that entropy increases in the
direction that the universe is expanding? This is really a situation in
which one can legitimately invoke the weak anthropic principle because
it is a question of where in the history of the universe conditions are
suitable for intelligent life. The inflation in the early universe
implies that the universe will expand for a very long time before it
contracts again. In fact, it is so long that the stars will have all
burnt out and the baryons will have all decayed. All that will be left
in the contracting phase will be a mixture of electrons, positrons,
neutrinos and gravitons. This is not a suitable basis for intelligent
life.

The conclusion of this paper is that the no boundary proposal can
explain the existence of a well defined thermodynamic arrow of time.
This arrow always points in the same direction. The reason we observe it
to point in the same direction as the cosmological arrow is that
conditions are suitable for intelligent life only at the low entropy end
of the universe's history.
\end{quote}

\hypertarget{week3}{%
\section{January 30, 1993}\label{week3}}

Here's this week's reading material. The first test will be in two
weeks. :-)

\begin{enumerate}
\def\labelenumi{\arabic{enumi})}
\tightlist
\item
   Dror Bar-Natan, ``On the Vassiliev knot invariants'', Harvard
  University ``pre-preprint''.
\end{enumerate}

I went to U.C.\ San Diego this week to give a talk, and the timing was
nice, because Dror Bar-Natan was there. He is a student of Witten who
has started from Witten's ideas relating knot theory and quantum field
theory and developed them into a beautiful picture that shows how knot
theory, the theory of classical Lie algebras, and abstract Feynman
diagrams are three faces of the same thing. To put it boldly, in a
deliberately exaggerated form, Bar-Natan has proposed a conjecture
saying that knot theory and the theory of classical Lie algebras are one
and the same!

This won't seem very exciting if you don't know what a classical Lie
algebra is. Let me give a brief and very sketchy introduction,
apologizing in advance to all the experts for the terrible sins I will
commit, such as failing to distinguish between complex and real Lie
algebras.

Well, remember that a Lie algebra is just a vector space equipped with a
``bracket'' such that the bracket \([x,y]\) of any two vectors \(x\) and
\(y\) is again a vector, and such that the following hold:

\begin{enumerate}
\def\labelenumi{\alph{enumi})}
\tightlist
\item
  skew-symmetry: \([x,y] = -[y,x]\).
\item
  bilinearity: \([x,ay] = a[x,y]\), \([x,y+z] = [x,y] + [x,z]\). (\(a\)
  is a number.)
\item
  Jacobi identity: \([x,[y,z]] + [y,[z,x]] + [z,[x,y]] = 0\).
\end{enumerate}

The best known example is good old \(\mathbb{R}^3\) with the cross
product as the bracket. But the real importance of Lie algebras is that
one can get one from any Lie group --- roughly speaking, a group that's
also a manifold, and such that the group operations are smooth maps. And
the importance of Lie groups is that they are what crop up as the groups
of symmetries in physics. The Lie algebra is essentially the
``infinitesimal version'' of the corresponding Lie group, as anyone has
seen who has taken physics and seen the relation between the group of
rotations in \(\mathbb{R}^3\) and the cross product. Here the group is
called \(\mathrm{SO}(3)\) and the Lie algebra is called
\(\mathfrak{so}(3)\). (So \(\mathbb{R}^3\) with its cross product is
called \(\mathfrak{so}(3)\).) One can generalize this to any number of
dimensions, letting \(\mathrm{SO}(n)\) denote the group of rotations in
\(\mathbb{R}^n\) and \(\mathfrak{so}(n)\) the corresponding Lie algebra.
(However, \(\mathfrak{so}(n)\) is not isomorphic to \(\mathbb{R}^n\)
except for \(n = 3\), so there is something very special about three
dimensions.)

Similarly, if one uses complex numbers instead of real numbers, one gets
a group \(\mathrm{SU}(n)\) and Lie algebra \(\mathfrak{su}(n)\). And if
one looks at the symmetries of a \(2n\)-dimensional classical phase
space --- so-called canonical transformations, or symplectic
transformations --- one gets the group \(\mathrm{Sp}(n)\) and Lie
algebra \(\mathfrak{sp}(n)\). To be precise, \(\mathrm{SO}(n)\) consists
of all \(n\times n\) orthogonal real matrices with determinant \(1\),
\(\mathrm{SU}(n)\) consists of all \(n\times n\) unitary complex
matrices with determinant \(1\), and \(\mathrm{Sp}(n)\) consists of all
\((2n)\times(2n)\) real matrices preserving a nondegenerate
skew-symmetric form.

These are all very important in physics. Indeed, all the ``gauge
groups'' of physics are Lie groups of a certain sort, so-called compact
Lie groups, and in the standard model all the forces are symmetrical
under some gauge group or other. Electromagnetism a la Maxwell is
symmetric under the group \(\mathrm{U}(1)\) of complex numbers of unit
magnitude, or ``phases''. The electroweak force (unified
electromagnetism and weak force) is symmetric under
\(\mathrm{U}(1) \times \mathrm{SU}(2)\), where one uses the fact that
one can build up bigger semisimple Lie groups as direct sums (also
called products) of smaller ones. The gauge group for the strong force
is \(\mathrm{SU}(3)\). And, finally, the gauge group of the whole
standard model is simply
\(\mathrm{U}(1) \times \mathrm{SU}(2) \times \mathrm{SU}(3)\), which
results from lumping the electroweak and strong gauge groups together.
This direct sum business also works for the Lie algebras, so the Lie
algebra relevant to the standard model is written
\(\mathfrak{u}(1) \times \mathfrak{su}(2) \times \mathfrak{su}(3)\).

There are certain very special Lie algebras called simple Lie algebras
which play the role of ``elementary building blocks'' in the world of
Lie algebras. They cannot be written as the direct sum of other Lie
algebras (and in fact there is an even stronger sense in which they
cannot be decomposed). On the other hand, the Lie algebra of any compact
Lie group is a direct sum of simple Lie algebras and copies of
\(\mathfrak{u}(1)\) --- the one-dimensional Lie algebra with zero Lie
bracket which, for technical reasons, people don't call ``simple''.

These simple Lie algebras were classified by the monumental work of
Killing, Cartan and others, and the classification is strikingly simple:
there are infinite series of ``classical'' Lie algebras of type
\(\mathfrak{su}(n)\), \(\mathfrak{so}(n)\), and \(\mathfrak{sp}(n)\),
and five ``exceptional'' Lie algebras called \(\mathrm{G}_2\),
\(\mathrm{F}_4\), \(\mathrm{E}_6\), \(\mathrm{E}_7\), and
\(\mathrm{E}_8\). Believe it or not, there is a deep connection between
the exceptional Lie algebras and the Platonic solids. But that is
another story, one I barely know\ldots.

Now, Witten showed how one could use quantum field theory to constuct an
invariant of knots, or even links, corresponding to any representation
of a compact Lie group. (You won't even need to know what a
representation is to understand what follows.) This had been done in a
different way, in terms of ``quantum groups,'' by Reshetikhin and Turaev
(following up on work by many other people). These invariants are
polynomials in a variable \(q\) (for ``quantum''), and if one writes
\(q\) as \(e^\hbar\) and expands a power series in \(\hbar\), the
coefficient of \(\hbar^n\) is a ``Vassiliev invariant of degree \(n\)''.
Recall from last week that given an invariant of oriented knots, one can
extend it to knot with arbitrarily many nice crossings by setting the
value of the invariant on a knot with a crossing like \[
  \begin{tikzpicture}
    \begin{knot}[clip width=0]
      \strand[thick] (0,0)
      to [out=down,in=up] (1,-1);
      \strand[thick] (1,0)
      to [out=down,in=up] (0,-1);
    \end{knot}
  \end{tikzpicture}
\] to be the invariant of the knot with the crossing changed to \[
  \begin{tikzpicture}
    \begin{knot}[clip width=7]
      \strand[thick] (0,0)
      to [out=down,in=up] (1,-1);
      \strand[thick] (1,0)
      to [out=down,in=up] (0,-1);
      \flipcrossings{1}
    \end{knot}
  \end{tikzpicture}
\] minus the invariant of the knot with the crossing changed to \[
  \begin{tikzpicture}
    \begin{knot}[clip width=7]
      \strand[thick] (0,0)
      to [out=down,in=up] (1,-1);
      \strand[thick] (1,0)
      to [out=down,in=up] (0,-1);
    \end{knot}
  \end{tikzpicture}
\] (Again, the knot has to be oriented for this rule to make sense, and
the strands shown in the pictures above should be pointing downwards.)
Having made this extension, one says a knot invariant is a Vassiliev
invariant of degree \(n\) if it vanishes on all knots with \(n+1\) or
more double points.

This is where Dror stepped in, roughly. First of all, he showed that the
Vassiliev invariant of degree \(n\) is just what you get when you do
Witten's quantum-field-theoretic calculations perturbatively using
Feynman diagrams and look at the terms of order \(n\) in Planck's
constant, \(\hbar\)! Secondly, and more surprisingly, he developed a
bunch of relationships between Feynman diagrams and pictures of knots!
The third and most amazing thing he did takes a bit longer to
explain\ldots{}

Roughly, he showed that any Vassiliev invariant of degree \(n\) is
determined by some combinatorial data called a ``weight system.'' He
showed that any representation of a Lie algebra determines a weight
system and hence a Vassiliev invariant. But the really interesting thing
he showed is that many of the things one can do for Lie algebras can be
done for arbitrary weight systems. This makes it plausible that
\textbf{every} weight system, hence every Vassiliev invariant, comes
from a representation of a simple Lie algebra. In fact, Dror conjectures
that every Vassiliev invariant comes from a representation of a
classical simple Lie algebra. Now there is another conjecture floating
around these days, namely that Vassiliev invariants almost form a
complete set --- that is, that if two knots cannot be distinguished by
any Vassiliev invariants, they must either be the same or differ simply
by reversing the orientation of all the strands. If \emph{both} these
conjectures are true, one has in some sense practically reduced the
theory of knots to the theory of the classical Lie algebras! This
wouldn't mean that all of sudden we know the answer to every question
about knots, but it would certainly help a lot, and more importantly, in
my opinion, it would show that the connection between topology and the
theory of Lie algebras is far more profound than we really understand.
The ramifications for physics, as I hope all my chatting about knots,
gauge theories and quantum gravity makes clear, might also be profound.

Well, we \emph{certainly} don't understand all this stuff yet, since we
don't know how to prove these conjectures! But Dror's conjecture ---
that all weight systems come from representations of simple Lie algebras
--- is tantalizingly close to being within grasp, since he has reduced
it to a fairly elementary combinatorial problem, which I will now state.
Note that ``elementary'' does not mean easy to solve! Just easy to
state.

Before I state the combinatorial problem, let me say something about the
evidence for the conjecture that all Vassiliev invariants come from
representations of classical Lie algebras. In addition to all sorts of
``technical'' evidence, Dror has shown the conjecture is true for
Vassiliev invariants of degree \(\leqslant 9\) by means of many hours of
computation using his Sparcstation. In fact, he said in his talk that he
felt guilty about having a Sparcstation unless it was always computing
something, and that even as he spoke his computer was busily verifying
the conjecture for higher degrees. (I suggested that it was the
Sparcstation that should feel guilty when it was not working, not him.)
He also advertised that his programs, a mixture of C and Mathematica
code, are available by anonymous ftp from \texttt{math.harvard}. Use
user name ``ftp'', go to the directory ``dror''. You folks with Crays
should feel \emph{very} guilty if they are just sitting there and not helping
Dror verify this important conjecture. (I suggest that you first read
his papers and the file README in his directory, then check out his
programs, and then ask him where he's at and what would be worth doing.
Please don't pester him unless you are a good enough mathematician to
discuss this stuff intelligently and have megaflops to burn. If you want
to make a fool of yourself, \emph{don't} say I sent you.)

Okay, with no further ado, here's the conjecture in its elementary
combinatorial form. Let \(B\) be the vector space spanned by finite
graphs with univalent and ``oriented'' trivalent vertices, modulo some
relations\ldots{} first of all, a trivalent vertex is ``oriented'' if
there is a cyclic ordering of the three incident edges. That is, we
``orient'' the vertex \[
  \begin{tikzpicture}
    \draw[thick] (50:1.5) to (0:0);
    \draw[thick] (130:1.5) to (0,0);
    \draw[thick] (270:1.5) to (0,0);
    \draw[ultra thick,white] (0,0) circle (6pt);
    \draw[ultra thick,white] (0,0) circle (7pt);
    \node[scale=2,rotate=150] at (0,0) {$\circlearrowleft$};
  \end{tikzpicture}
\] by drawing a little clockwise or counterclockwise-pointing circle at
the vertex. (Or, for those of an algebraic bent, label the edges by
\(1\),\(2\),\(3\) but then mod out by cyclic permutations.) The
relations are:

\begin{enumerate}
\def\labelenumi{\arabic{enumi})}
\item
  if we reverse the orientation of a trivalent vertex, that's equivalent
  to multiplying the graph by \(-1\). (Remember we're in a vector space
  spanned by graphs.)
\item
  \[
      \vcenter{\hbox{\begin{tikzpicture}
        \draw[thick] (-1,1) to (1,1);
        \draw[thick] (0,-1) to (0,1);
        \draw[thick] (-1,-1) to (1,-1);
      \end{tikzpicture}}}
      \quad=\quad
      \vcenter{\hbox{\begin{tikzpicture}[rotate=90]
        \draw[thick] (-1,0.85) to (1,0.85);
        \draw[thick] (0,-0.85) to (0,0.85);
        \draw[thick] (-1,-0.85) to (1,-0.85);
      \end{tikzpicture}}}
      \quad-
      \vcenter{\hbox{\begin{tikzpicture}
        \draw[thick] (-0.7,1) to (0.7,-1);
        \draw[thick] (0.7,1) to (-0.7,-1);
      \end{tikzpicture}}}
    \] (That is, we can make this substitution anywhere we want; these
  pictures might be part of a bigger graph. Note that the ``X'' is not a
  vertex, since there aren't quadrivalent vertices; it's just one edge
  going over or under another. It doesn't matter whether it goes over or
  under since these are abstract graphs, not graphs embedded in space.)
\end{enumerate}

Now, let \(B_m\) be the vector space spanned by ``labelled'' finite
graphs with univalent and oriented trivalent vertices, modulo some
relations\ldots{} but first I have to say what ``labelled'' means. It
means that each edge is labelled with a \(1\) or \(-1\). The relations
are:

\begin{enumerate}
\def\labelenumi{\arabic{enumi})}
\item
  if we reverse the orientation of a trivalent vertex, it's the same as
  multiplying the labellings of all three incident edges by \(-1\).
\item
  \[
      \vcenter{\hbox{\begin{tikzpicture}
        \draw[thick] (-1,1) to (1,1);
        \draw[thick] (0,-1) to (0,1);
        \draw[thick] (-1,-1) to (1,-1);
      \end{tikzpicture}}}
      \quad=\quad
      \vcenter{\hbox{\begin{tikzpicture}[rotate=90]
        \draw[thick] (-1,0.85) to (1,0.85);
        \draw[thick] (0,-0.85) to (0,0.85);
        \draw[thick] (-1,-0.85) to (1,-0.85);
      \end{tikzpicture}}}
    \] \emph{if} the internal edge is labelled with a 1. (Here the 4
  external edges can have any labellings and we don't mess with that.)
\end{enumerate}

Now, define a linear map from \(B\) to \(B_m\) by mapping any graph to
the signed sum of the \(2^\text{number of edges}\) ways of labelling the
edges with \(1\) or \(-1\). Symbolically, \[
  \begin{tikzpicture}
    \draw[thick] (-1,0) to (1,0);
  \end{tikzpicture}
  \quad\mapsto\quad
  \begin{tikzpicture}
    \draw[thick] (-1,0) to (1,0);
    \node at (0,0.3) {$1$};
  \end{tikzpicture}
  \quad-\quad
  \begin{tikzpicture}
    \draw[thick] (-1,0) to (1,0);
    \node at (0,0.3) {$-1$};
  \end{tikzpicture}
\] Of course, one must work a bit to show this map is well-defined.
(This just takes a paragraph --- see Proposition 6.5 of Dror's paper.)

Okay, the conjecture is:

\begin{center}
            THIS MAP IS ONE-TO-ONE.
\end{center}
\noindent
If you can solve it, you've made great progress in showing that knots
and classical Lie groups are just two aspects of the same branch of
mathematics. Don't work on it, though, until you get Dror's paper and
make sure I stated it exactly right!!!!!

\begin{enumerate}
\def\labelenumi{\arabic{enumi})}
\setcounter{enumi}{1}
\tightlist
\item
  Abhay Ashtekar, ``Mathematical problems of non-perturbative quantum general relativity'',
  lectures delivered at the 1992 Les Houches summer
  school on Gravitation and Quantization, December 2, 1992, available as
  \href{https://arxiv.org/abs/gr-qc/9302024}{\texttt{gr-qc/9302024}}.
\end{enumerate}

This is a good overview of the loop variables approach to quantizing
general relativity as it currently stands. It begins with a review of
the basic difficulties with quantizing gravity, as viewed from three
perspectives: the particle physicist, the mathematical physicist, and
the general relativist. Technically, a main problem is that general
relativity consists of both evolution equations and constraint equations
on the initial data (which are roughly the metric of space at a given
time and its first time derivative, or really ``extrinsic curvature'').
So Ashtekar reviews Dirac's ideas on quantizing constrained systems
before sketching how this program is carried out for general relativity.

Then he considers a ``toy model'' --- quantum gravity in 2+1 dimensions.
This is a funny theory because \emph{classically} Einstein's equations
in 2+1 dimensions simply say that spacetime is flat (in a vacuum)! No
gravitational waves exist as in 3+1 dimensions, and one can say that the
information in the gravitational field is ``purely global'' --- locally,
everywhere looks the same as everywhere else (like Iowa), but there may
be global ``twists'' that you notice when going around a noncontractible
loop. There has been a lot of work on 2+1 gravity recently --- in a
sense this problem has been solved, by a number of methods --- and this
allows one to understand \emph{some} of the conceptual difficulties of
honest 3+1-dimensional quantum gravity without getting caught in an
endless net of technical complications.

Then Ashtekar jumps back to 3+1 dimensions and gives a more thorough
introduction to the loop variables approach. He ends by going through
some of the many open problems and possible ways to attack them.

I have worn myself out trying to do justice to Bar-Natan's work, so I
will postpone until next week a review of Kapranov and Voevodsky's paper
on \(2\)-categories.

\hypertarget{week4}{%
\section{February 8, 1993}\label{week4}}

I will begin with a couple of small things and then talk about the work
of Kapranov and Voevodsky.

\begin{enumerate}
\def\labelenumi{\arabic{enumi})}
\tightlist
\item
  R.\ Sole, D.\ Lopez, M.\ Ginovart and J.\ Valls, ``Self-organized criticality in Monte Carlo 
  simulated ecosystems'', \emph{Phys. Lett.}
  \textbf{A172} (1992), p.~56.
\end{enumerate}

This is mainly of interest to me thanks to a reference to some earlier
work on Conway's game of Life. At MIT, Tom Toffoli, Norm Margolus, and
grad students in the Physics of Computation group build special-purpose
computers for simulating cellular automata, the so-called CAM machines.
I have spent many enjoyable hours watching beautiful patterns do their
thing on a big-screen color TV while CAM 6 busily simulates them on a
\(256 \times 256\) lattice at the rate of many generations a second.
(The CAM 8 chip was still being debugged when I last checked.) More
recently, Jim Gilliam, a grad student here at UCR, found a very nice
program for the game of Life on Xwindows, called xlife. On my
Sparcstation it is even bigger and faster than CAM6. One can zoom in and
out, and, zooming all the way out, one sees something vaguely
reminiscent of nebulae of distant stars twinkling in the night
sky\ldots{} My computer science pal, Nate Osgood, muttered something
about the author, Chuck Silvers, using
cleverly optimized loops. It apparently can be found using the program
archie. Please don't ask \emph{me} for a copy, since it involves many
files.

The game of Life is actually one of the less fun cellular automata to
watch, since, contrary to its name, if one starts with a random
configuration it almost always eventually lapses into an essentially
static configuration (perhaps with some blinkers executing simple
periodic motions). I am pleased to find that this seemingly dull final
state might be fairly interesting in the study of self-organized
criticality! Recall that this is the phenomenon whereby a physical
system naturally works its way into a state such that the slightest
disturbance can have an arbitrarily large effect. The classic example is
a sand dune, which apparently works its way towards slopes close to the
critical one at which an avalanche occurs. Drop one extra grain of sand
on it and you can get a surprisingly wide distribution of possible sizes
for the resulting sandslide! Similar but more formidable effects may be
at work in earthquakes. The above paper cites

\begin{enumerate}
\def\labelenumi{\arabic{enumi})}
\setcounter{enumi}{1}
\tightlist
\item
Per Bak, Kan Chen, and Michael Creutz, ``Self-organized criticality in the `Game of Life'",
\emph{Nature} \textbf{342} (1989), 780--782.
\end{enumerate}
\noindent
which claims that in the final state of the game of Life, the density of
clusters \(D(s)\) of size \(s\) scales as about \(s^{-1.4}\), and that
the probability that a small perturbation will cause a flurry of
activity lasting a time \(t\) scales as about \(t^{-1.6}\). I'm no
expert, but I guess that the fact that the latter is a power law rather
than an exponential would be a signal of self-organized criticality. But
the paper also cites

\begin{enumerate}
\def\labelenumi{\arabic{enumi})}
\setcounter{enumi}{2}
\tightlist
\item
Charles Bennett and Marc S.\ Bourzutschy,```Life' not critical?'', 
\emph{Nature} \textbf{350} (1991), 468.
\end{enumerate}
\noindent
who claim that the work of Bac, Chen and Creutz is wrong. I haven't
gotten to read these papers; if anyone wants to report on them I'd be
interested.

The paper I read considers fancier variations on this theme,
investigating the possibility that ecosystems are also examples of
self-organized criticality. It's hard to know how to make solid science
out of this kind of thing, but I think there would be important
consequences if it turned out that the ``balance of nature,'' far from
being a stable equilibrium, was typically teetering on the brink of
drastic change.

\begin{enumerate}
\def\labelenumi{\arabic{enumi})}
\setcounter{enumi}{3}
\tightlist
\item
  H.\ D.\ Zeh, ``There are no quantum jumps, nor are there particles!'', \emph{Phys.\ Lett.} \textbf{A173} (1993), 189--192.
\end{enumerate}
\noindent
Having greatly enjoyed Zeh's book \emph{The Physical Basis for the Direction
of Time} --- perhaps the clearest account of a famously murky subject ---
I naturally took a look at this particle despite its overheated title.
(Certainly exclamation marks in titles should add to one's crackpot
index.) It is a nice little discussion of ``quantum jumps'' and the
``collapse of the wavefunction'' that takes roughly the viewpoint I
espouse, namely, that all one needs is Schr\"odinger's equation (and lots
of hard work) to understand what's going on in quantum theory --- no
extra dynamical mechanisms. It's not likely to convince anyone who
thinks otherwise, but it has references that might be useful no matter
which side of the debate one is on.

Also, this just in --- what you've all been waiting for --- another
interpretation of quantum mechanics! It's a book by David Bohm and Basil
Hiley, entitled \emph{The Undivided Universe --- An Ontological
Interpretation of Quantum Theory}. I have only seen an advertisement so
far; it's published by Routledge. The contents include such curious
things as ``the ontological interpretation of boson fields.'' Read it at
your own risk.

\begin{enumerate}
\def\labelenumi{\arabic{enumi})}
\setcounter{enumi}{4}
\tightlist
\item
  M.\ M.\ Kapranov and V.\ A.\ Voevodsky, 
   ``2-categories and Zamolodchikov tetrahedra equations'' in \emph{Algebraic 
   Groups and their Generalization: Quantum and Infinite-Dimensional Methods, 
   University Park, PA (1991)}, eds.\ W.\ J.\ Haboush and B.\ J.\ Parshall, 
   Proc.\ Sympos.\ Pure Math.\ \textbf{56}, American Mathematical Society,
   Providence, Rhode Island, 1994, pp.\ 177--259.
\end{enumerate}

This serious and rather dry paper is the basis for a lot of physicists
are just beginning to try to do: burst from the confines of 3
dimensions, where knots and topological quantum field theories like
Chern--Simons theory live, to 3+1 dimensions, where we live. The
``incomplete version'' I have now is 220 pages long, mostly commutative
diagrams, and doesn't have much to say about physics. But I have a hunch
that it will become required reading for many people fairly soon, so I'd
like to describe the main ideas in fairly simple terms.

I will start from scratch and then gradually accelerate. First, what's a
category? A category consists of a set of `objects' and a set of
`morphisms'. Every morphism has a `source' object and a `target' object.
(The easiest example is the category in which the objects are sets and
the morphisms are functions. If \(f\colon X\to Y\), we call \(X\) the
source and \(Y\) the target.) Given objects \(X\) and \(Y\), we write
\(\mathrm{Hom}(X,Y)\) for the set of morphisms `from' \(X\) `to' \(Y\)
(i.e., having \(X\) as source and \(Y\) as target).

The axioms for a category are that it consist of a set of objects and
for any 2 objects \(X\) and \(Y\) a set \(\mathrm{Hom}(X,Y)\) of
morphisms from \(X\) to \(Y\), and

\begin{enumerate}
\def\labelenumi{\alph{enumi})}
\item
  Given a morphism \(g\) in \(\mathrm{Hom}(X,Y)\) and a morphism \(f\)
  in \(\mathrm{Hom}(Y,Z)\), there is morphism which we call \(f\circ g\)
  in \(\mathrm{Hom}(X,Z)\). (This binary operation \(\circ\) is called
  `composition'.)
\item
  Composition is associative: \((f\circ g)\circ h = f\circ (g\circ h)\).
\item
  For each object \(X\) there is a morphism \(\mathrm{id}_X\) from \(X\)
  to \(X\), called the `identity' on \(X\).
\item
  Given any \(f\) in \(\mathrm{Hom}(X,Y)\),
  \(f \circ \mathrm{id}_X = f\) and \(\mathrm{id}_Y \circ f = f\).
\end{enumerate}

Again, the classic example is \(\mathsf{Set}\), the category with sets
as objects and functions as morphisms, and the usual composition as
composition! But lots of the time in mathematics one is some category or
other, e.g.:

\begin{itemize}
\tightlist
\item
  \(\mathsf{Vect}\) --- vector spaces as objects, linear maps as
  morphisms
\item
  \(\mathsf{Group}\) --- groups as objects, homomorphisms as morphisms
\item
  \(\mathsf{Top}\) --- topological spaces as objects, continuous
  functions as morphisms
\item
  \(\mathsf{Diff}\) --- smooth manifolds as objects, smooth maps as
  morphisms
\item
  \(\mathsf{Ring}\) --- rings as objects, ring homomorphisms as
  morphisms
\end{itemize}
\noindent
or in physics:

\begin{itemize}
\tightlist
\item
  \(\mathsf{Symp}\) --- symplectic manifolds as objects,
  symplectomorphisms as morphisms
\item
  \(\mathsf{Poiss}\) --- Poisson manifolds as objects, Poisson maps as
  morphisms
\item
  \(\mathsf{Hilb}\) --- Hilbert spaces as objects, unitary operators as
  morphisms
\end{itemize}
\noindent
(The first two are categories in which one can do classical physics. The
third is a category in which one can do quantum physics.)

Now, what's a \(2\)-category? This has all the structure of a category
but now there are also ``2-morphisms,'' that is, morphisms between
morphisms! This is rather dizzying at first. Indeed, much of category
theory is rather dizzying until one has some good examples to lean on
(at least for down-to-earth people such as myself), so let us get some
examples right away, and leave the definition to Kapranov and Voevodsky!
My favorite example comes from homotopy theory. Take a topological space
\(X\) and let the objects of our category be points of \(X\). Given
\(x\) and \(y\) in \(X\), let \(\mathrm{Hom}(x,y)\) be the set of all
unparametrized paths from \(x\) to \(y\). We compose such paths simply
by sticking a path from \(x\) to \(y\) and a path from \(y\) to \(z\) to
get a path from \(x\) to \(z\), and we need unparametrized paths to make
composition associative. Now given two paths from \(x\) to \(y\), say
\(f\) and \(g\), let \(\mathrm{Hom}(f,g)\), the set of \(2\)-morphisms
from \(f\) to \(g\), be the set of unparametrized homotopies from \(f\)
to \(g\) --- that is, ways of deforming the path \(f\) continuously to
get the path \(g\), while leaving the endpoints fixed.

This is a very enlightening example since homotopies of paths are really
just ``paths of paths,'' making the name \(2\)-morphism quite
appropriate. (Some of you will already be pondering \(3\)-morphisms,
\(4\)-morphisms, but it's too late, they've already been invented! I
won't discuss them here.) The notation for \(2\)-morphisms is quite
cute: given \(f,g\) in \(\mathrm{Hom}(x,y)\), we write \(F\) in
\(\mathrm{Hom}(f,g)\) as the following diagram:
\[\includegraphics[scale=0.3]{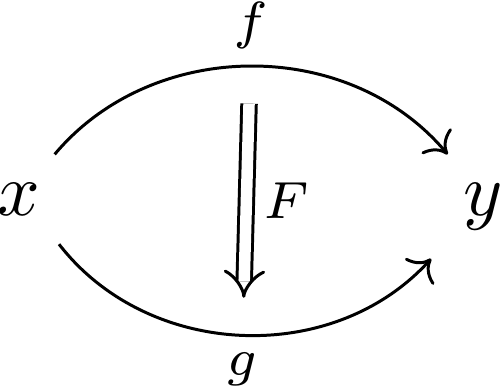}\] 
In other words,
while ordinary morphisms are \(1\)-dimensional objects (arrows),
\(2\)-morphisms are \(2\)-dimensional ``cells'' filling in the space
between two ordinary morphisms. We thus see that going up to ``morphisms
between morphisms'' is closely related to going up to higher dimensions.
And this is really why ``braided monoidal \(2\)-categories'' may play as
big a role in four-dimensional field theory as ``braided monoidal
categories'' do in three-dimensional field theory!

Rather than write down the axioms for a \(2\)-category, which are in
Kapranov and Voevodsky, let me note the key new thing about
\(2\)-morphisms: there are two ways to compose them, ``horizontally''
and ``vertically''. First of all, given the following situation:
\[\includegraphics[scale=0.3]{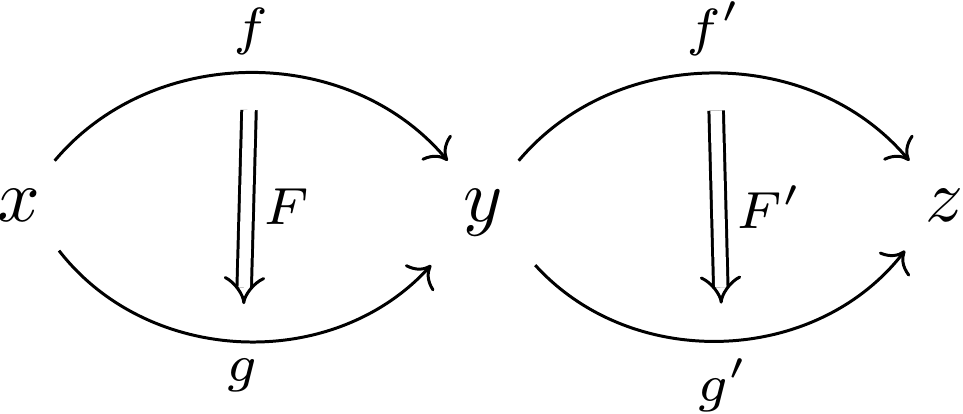}\] we can compose
\(F\) and \(F'\) horizontally to get a \(2\)-morphism from \(f'\circ f\)
to \(g \circ g'\). (Check this out in the example of homotopies!) But
also, given the following situation:
\[\includegraphics[scale=0.3]{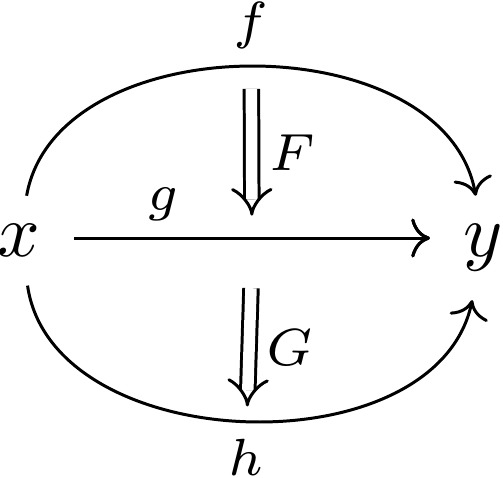}\] (\(f,g,h\) in
\(\mathrm{Hom}(x,y)\), \(F\) in \(\mathrm{Hom}(f,g)\), and \(G\) in
\(\mathrm{Hom}(g,h)\)), we can compose \(F\) and \(G\) vertically to get
a \(2\)-morphism from \(f\) to \(g\).

As Kapranov and Voevodsky note: ``Thus \(2\)-categories can be seen as
belonging to the realm of a new mathematical discipline which may be
called \textbf{2-dimensional algebra} and contrasted with usual
\(1\)-dimensional algebra dealing with formulas which are written in
lines.'' This is actually very important because already in the theory
of braided monoidal categories we began witnessing the rise of
mathematics that incorporated aspects of geometry into the notation
itself.

The theory of \(2\)-categories is not new; it was apparently invented by
Ehresmann, Benabou and Grothendieck in an effort to formalize the
structure possessed by the category of all categories. (If this notion
seems dangerously close to Russell's paradox, you are right --- but I
will not worry about such issues in what follows.) This category has as
its objects categories and as its morphisms ``functors'' between
categories. It is, in fact, a \(2\)-category, taking as the
\(2\)-morphisms ``natural transformations'' between functors. (For a
brief intro to functors and natural transformations, try my webpage
\href{http://math.ucr.edu/home/baez/categories.html}{``categories"}.) 
What is new to Kapranov and Voevodsky is the notion of a
monoidal \(2\)-category --- where one can take tensor products of
objects, morphisms, and 2-morphisms --- and ``braided'' monoidal
\(2\)-category --- where one has ``braidings'' that switch around the
two factors in a tensor product.

Let me turn to the possible relevance of all this to mathematical
physics. Here there is a nice \(2\)-category, namely the category of
``2-tangles.'' First recall the category of tangles:

The objects are simply the natural numbers \(\{0,1,2,3,\ldots\}\). We
think of the object \(n\) as a horizontal row of \(n\) points. The
morphisms in \(\mathrm{Hom}(n,m)\) are tangles connecting a row of \(n\)
points above to a row of \(m\) points below. Rather than define
``tangles'' I will simply draw pictures of some examples. Here is an
element of \(\mathrm{Hom}(2,4)\): \[
  \begin{tikzpicture}
    \begin{knot}[clip width=7]
      \strand[thick] (0,0)
        to [out=down,in=up] (1,-1)
        to [out=down,in=up] (2,-2);
      \strand[thick] (1,0)
        to [out=down,in=up] (0,-1)
        to (0,-2);
      \strand[thick] (1,-2)
        to [out=up,in=up,looseness=1.5] (3,-2);
    \end{knot}
  \end{tikzpicture}
\] and here is an element of \(\mathrm{Hom}(4,0)\): \[
  \begin{tikzpicture}
    \begin{knot}[clip width=7]
      \strand[thick] (0,0)
        to (0,-0.5)
        to [out=down,in=down,looseness=1.5] (1,-0.5)
        to (1,0);
      \strand[thick] (2,0)
        to (2,-0.5)
        to [out=down,in=down,looseness=1.5] (3,-0.5)
        to (3,0);
      \strand[thick] (2,-1.5)
        to [out=up,in=left] (2.5,-0.5);
      \strand[thick] (2.5,-0.5)
        to [out=right,in=up] (3,-1.5);
      \strand[thick] (2,-1.5)
        to [out=down,in=down] (3,-1.5);
      \flipcrossings{2}
    \end{knot}
  \end{tikzpicture}
\] Note that we can ``compose'' these tangles to get one in
\(\mathrm{Hom}(2,0)\): \[
  \begin{tikzpicture}
    \begin{knot}[clip width=7]
      \strand[thick] (0,0)
        to [out=down,in=up] (1,-1)
        to [out=down,in=up] (2,-2);
      \strand[thick] (1,0)
        to [out=down,in=up] (0,-1)
        to (0,-2);
      \strand[thick] (1,-2)
        to [out=up,in=up,looseness=1.5] (3,-2);
    \end{knot}
    \begin{knot}[clip width=7]
      \strand[thick] (0,-2)
        to (0,-2.5)
        to [out=down,in=down,looseness=1.5] (1,-2.5)
        to (1,-2);
      \strand[thick] (2,-2)
        to (2,-2.5)
        to [out=down,in=down,looseness=1.5] (3,-2.5)
        to (3,-2);
      \strand[thick] (2,-3.5)
        to [out=up,in=left] (2.5,-2.5);
      \strand[thick] (2.5,-2.5)
        to [out=right,in=up] (3,-3.5);
      \strand[thick] (2,-3.5)
        to [out=down,in=down] (3,-3.5);
      \flipcrossings{2}
    \end{knot}
  \end{tikzpicture}
\] Now, given tangles \(f,g\) in \(\mathrm{Hom}(m,n)\), a \(2\)-morphism
from \(f\) to \(g\) is a ``2-tangle.'' I won't define these either, but
we may think of a 2-tangle from \(f\) to \(g\) roughly as a ``movie''
whose first frame is the tangle \(f\) and last frame is the tangle
\(g\), and each of whose intermediate frames is a tangle except at
certain times when a catastrophe occurs. For example, here's a 2-tangle
shown as a movie\ldots{} \[
  \begin{array}{p{8em}p{8em}p{8em}}
    \begin{tikzpicture}
      \begin{knot}
        \strand[thick] (0,0)
          to (0,-2);
        \strand[thick] (1,0)
          to (1,-2);
      \end{knot}
      \node at (0.5,-2.5) {Frame 1};
    \end{tikzpicture}
    &
    \begin{tikzpicture}
      \begin{knot}
        \strand[thick] (0,0)
          to[out=right,in=right,looseness=0.4] (0,-2);
        \strand[thick] (1,0)
          to[out=left,in=left,looseness=0.4] (1,-2);
      \end{knot}
      \node at (0.5,-2.5) {Frame 2};
    \end{tikzpicture}
    &
    \begin{tikzpicture}
      \begin{knot}
        \strand[thick] (0,0)
          to[out=right,in=right,looseness=0.7] (0,-2);
        \strand[thick] (1,0)
          to[out=left,in=left,looseness=0.7] (1,-2);
      \end{knot}
      \node at (0.5,-2.5) {Frame 3};
    \end{tikzpicture}
    \\[1em]
    \hspace*{-1.8em}
    \raisebox{-1.6em}{
    \begin{tikzpicture}
      \begin{knot}[clip width=0]
        \strand[thick] (0,0)
          to (1,-2);
        \strand[thick] (1,0)
          to (0,-2);
      \end{knot}
      \node at (0.5,-2.5) {Frame 4};
      \node at (0.5,-3) {(catastrophe!)};
    \end{tikzpicture}
    }
    &
    \begin{tikzpicture}
      \begin{knot}
        \strand[thick] (0,0)
          to[out=down,in=down,looseness=3] (1,0);
        \strand[thick] (0,-2)
          to[out=up,in=up,looseness=3] (1,-2);
      \end{knot}
      \node at (0.5,-2.5) {Frame 5};
    \end{tikzpicture}
    &
    \begin{tikzpicture}
      \begin{knot}
        \strand[thick] (0,0)
          to[out=down,in=down,looseness=1.5] (1,0);
        \strand[thick] (0,-2)
          to[out=up,in=up,looseness=1.5] (1,-2);
      \end{knot}
      \node at (0.5,-2.5) {Frame 6};
    \end{tikzpicture}
  \end{array}
\] Well, it'll never win an Academy Award, but this movie is pretty
important. It's a picture of the \(3\)-dimensional slices of a
2-dimensional surface in (3+1)-dimensional spacetime, and this surface
is perfectly smooth but has a saddle point which we are seeing in frame
4. It is one of what Carter and Saito (see
\protect\hyperlink{week2}{``Week 2''}) call the ``elementary string
interactions.'' The relevance to string theory is pretty obvious: we are
seeing a movie of part of a string worldsheet, which is a surface in
(3+1)-dimensional spacetime. My interest in 2-tangles and
\(2\)-categories is precisely because they offer a bridge between string
theory and the loop variables approach to quantum gravity, which may
actually be the \textbf{same thing} in two different disguises. You
heard it here first, folks!

The reader may have fun figuring out what the two ways of composing
2-morphisms amount to in the category of 2-tangles.

There are, in fact, many clues as to the relation between string theory
and \(2\)-categories, one being the Zamolodchikov equation. This is the
analog of the Yang--Baxter equation --- an equation important in the
theory of braids --- one dimension up. It was discovered by
Zamolodchikov in 1980; a 1981 paper that might be a bit easier to get is

\begin{enumerate}
\def\labelenumi{\arabic{enumi})}
\setcounter{enumi}{5}
\tightlist
\item
A.\ B.\ Zamolodchikov, ``Tetrahedron equations and relativistic S-matrix for 
straight strings in 2+1 dimensions,'' \textsl{Commun.\ Math.\ Phys.\ }\textbf{79} 
(1981), 489--505. 
\end{enumerate}
\noindent
(It plays a
different role in the (3+1)-dimensional context, though.) Just as
braided monoidal categories are a good way to systematically find
solutions of the Yang--Baxter equation, braided monoidal
\(2\)-categories, as defined by Kapranov and Voevodsky, seem to be a
good way for finding solutions for the Zamolodchikov equation. (I will
post in a while about a new paper by Soibelman and Kazhdan that does
this. Also see the paper by Crane and Frenkel in
\protect\hyperlink{week2}{``Week 2''}.)

There are also lots of tantalizing ties between the loop variables
approach to quantum gravity and \(2\)-categories; one can see some of
these if one reads the work of Carter and Saito in conjunction with my
paper ``Quantum gravity and the algebra of tangles'' 
(\href{https://arxiv.org/abs/hep-th/9205007}{\texttt{hep-th/9205007}}).   
I hope to make these a lot clearer as time goes by.

\hypertarget{week5}{%
\section{February 13, 1993}\label{week5}}

I think I'll start out this week's list of finds with an elementary
introduction to Lie algebras, so that people who aren't ``experts'' can
get the drift of what these are about. Then I'll gradually pick up
speed\ldots{}

\begin{enumerate}
\def\labelenumi{\arabic{enumi})}
\tightlist
\item
  Vyjayanathi Chari and Alexander Premet, ``Indecomposable restricted
  representations of quantum \(\mathfrak{sl}_2\)'', 
  \emph{Publications of the Research Institute for Mathematical Sciences}
  \textbf{30} (1994), 335--352.
\end{enumerate}

\noindent
Vyjayanathi is our resident expert on quantum groups, and Sasha, who's
visiting, is an expert on Lie algebras in characteristic \(p\). They
have been talking endlessly across the hall from me and now I see that
it has born fruit. This is a pretty technical paper and I am afraid I'll
never really understand it, but I can see why it's important, so I'll
try to explain that!

Let me start with the prehistory, which is the sort of thing everyone
should learn. Recall what a Lie algebra is\ldots{} a vector space with a
``bracket'' operation such that the bracket \([x,y]\) of any two vectors
\(x\) and \(y\) is again a vector, and such that the following hold:

\begin{enumerate}
\def\labelenumi{\alph{enumi})}
\tightlist
\item
  skew-symmetry: \([x,y] = -[y,x]\).
\item
  bilinearity: \([x,ay] = a[x,y]\), \([x,y+z] = [x,y] + [x,z]\). (\(a\)
  is any number)
\item
  Jacobi identity: \([x,[y,z]] + [y,[z,x]] + [z,[x,y]] = 0\).
\end{enumerate}

These conditions, especially the third, may look sort of weird if you
are not used to them, but the examples make it all clear. The easiest
example of a Lie algebra is \(\mathfrak{gl}(n,\mathbb{C})\), which just
means all \(n\times n\) complex matrices with the bracket defined as the
``commutator'': \[[x,y] = xy-yx.\] Then straightforward calculations
show a)--c) hold\ldots{} so these conditions are really encoding the
essence of the commutator.

Now recall that the trace of a matrix, the sum of its diagonal entries,
satisfies \(\operatorname{tr}(xy)\) = \(\operatorname{tr}(yx)\). So the
trace of any commutator is zero, and if we let
\(\mathfrak{sl}(n,\mathbb{C})\) denote the matrices with zero trace, we
see that it's a sub-Lie algebra of \(\mathfrak{gl}(n,\mathbb{C})\) ---
that is, if \(x\) and \(y\) are in \(\mathfrak{sl}(n,\mathbb{C})\) so is \([x,y]\),
so we can think of \(\mathfrak{sl}(n,\mathbb{C})\) as a Lie algebra in
its own right. Going from \(\mathfrak{sl}(n,\mathbb{C})\) to
\(\mathfrak{gl}(n,\mathbb{C})\) is essentially a trick for booting out
the identity matrix, which commutes with everything else (hence has
vanishing commutators).  Multiples of identity matrix are the only ones with this
property, so they're sort of weird, and it simplifies things to get rid of
them here.

The simplest of the \(\mathfrak{sl}(n,\mathbb{C})\)'s is the Lie algebra
\(\mathfrak{sl}(2,\mathbb{C})\), affectionately known simply as
\(\mathfrak{sl}(2)\), which is a \(3\)-dimensional Lie algebra with a
basis given by matrices people call \(E\), \(F\), and \(H\) for
mysterious reasons:
\[E=\left(\begin{array}{cc}0&1\\0&0\end{array}\right) \quad F=\left(\begin{array}{cc}0&0\\1&0\end{array}\right) \quad H=\left(\begin{array}{cc}1&0\\0&-1\end{array}\right)\]
You will never be an expert on Lie algebras until you know by heart that
\[[H,E] = 2E, \quad [H,F] = -2F, \quad [E,F]  = H.\] Typically that's
the sort of remark I make before screwing up by a factor of two or
something, so you'd better check! This is a cute little multiplication
table\ldots{} but very important, since \(\mathfrak{sl}(2)\) is the
primordial Lie algebra from which the whole theory of ``simple'' Lie
algebras unfolds.

Physicists are probably more familiar with a different basis of
\(\mathfrak{sl}(2)\), the Pauli matrices:
\[\sigma_1=\left(\begin{array}{cc}0&1\\1&0\end{array}\right) \quad \sigma_2=\left(\begin{array}{cc}0&-i\\i&0\end{array}\right) \quad \sigma_3=\left(\begin{array}{cc}1&0\\0&-1\end{array}\right)\]
For purposes of Lie algebra theory it's actually better to divide each
of these matrices by \(i\) and call the resulting matrices \(I\), \(J\),
and \(K\), respectively. We then have
\[IJ = -JI = K, \quad JK = -KJ = I, \quad KI = -IK = J, \quad I^2 = J^2 = K^2 = -1\]
which is just the multiplication table of the quaternions! From the
point of view of Lie algebras, though, all that matters is
\[[I,J] = 2K, \quad [J,K] = 2I, \quad [K,I] = 2J.\] Given the relation
of these things and cross products, it should be no surprise that the
Pauli matrices have a lot to with angular momentum around the \(x\),
\(y\), and \(z\) axes in quantum mechanics.

If we take all \emph{real} linear combinations of \(E\),\(F\),\(H\) we
get a Lie algebra over the \emph{real} numbers called
\(\mathfrak{sl}(2,\mathbb{R})\), and if we take all real linear
combinations of \(I\),\(J\),\(K\) we get a Lie algebra over the reals
called \(\mathfrak{su}(2)\). These two Lie algebras are two different
``real forms'' of \(\mathfrak{sl}(2)\).

Now, people know just about everything about \(\mathfrak{sl}(2)\) that
they might want to. Well, there's always something more, but I'm
certainly personally satisfied! I recall when as an impressionable
student I saw a book by Serge Lang titled simply
``\(\mathrm{SL}(2,\mathbb{R})\)'', big and fat and scary inside. I knew
what \(\mathrm{SL}(2,\mathbb{R})\) was, but not how one could think of a
whole book's worth of things to write about it! A whole book on \(2 \times 2\)
matrices??

Part of how one gets so much to say about a puny little Lie algebra like
\(\mathfrak{sl}(2)\) is by talking about its representations. What's a
representation? Well, first you have to temporarily shelve the idea that
\(\mathfrak{sl}(2)\) consists of \(2\times2\) matrices, and think of it
more abstractly simply as a 3-dimensional vector space with basis
\(E\),\(F\),\(H\), equipped with a Lie algebra structure given by the
multiplication table \([H,E] = 2E\), \([H,F] = -2F\), \([E,F] = H\). If
this is how I'd originally defined it, it would then be a little
\emph{theorem} that this Lie algebra has a ``representation'' as
\(2\times2\) matrices. And it would turn out to have other
representations too. For example, there's a representation as
\(3\times3\) matrices given by sending
\[E\mapsto\left(\begin{array}{ccc}0&1&0\\0&0&2\\0&0&0\end{array}\right) \quad F\mapsto\left(\begin{array}{ccc}0&0&0\\2&0&0\\0&1&0\end{array}\right) \quad H\mapsto\left(\begin{array}{ccc}2&0&0\\0&0&0\\0&0&-2\end{array}\right)\]
In other words, these matrices satisfy the same commutation relations as
\(E\),\(F\), and \(H\) do.

More generally, and more precisely, we say an \(n\)-dimensional
representation of a Lie algebra \(L\) (over the complex numbers) is a
linear function \(R\) from \(L\) to \(n\times n\) matrices such that
\[R([x,y]) = [R(x),R(y)]\] for all \(x,y\) in \(L\). Note that on the
left the brackets are the brackets in \(L\), while on the right they
denote the commutator of \(n\times n\) matrices.

One good way to understand the essence of a Lie algebra is to figure out
what representations it has. And in quantum physics, Lie algebra
representations are where it's at: the symmetries of the world are
typically Lie groups, each Lie group has a corresponding Lie algebra,
the states of a quantum system are unit vectors in a Hilbert space, and
if the system has a certain Lie group of symmetries there will be a
representation of the Lie algebra on the Hilbert space. As any particle
physicist can tell you, you can learn a lot just by knowing which
representation of your symmetry group a given particle has.

So the name of the game is classifying Lie algebra
representations\ldots{} and many tomes have been written on this by now.
To keep things from becoming too much of a mess it's crucial to make two
observations. First, there's an easy way to get new representations by
taking the ``direct sum'' of old ones: the sum of an \(n\)-dimensional
representation and an \(m\)-dimensional one is an \((n+m)\)-dimensional
one, for example. Another way, not so easy, to get new representations
from an old one is to look for ``subrepresentations'' of the given
representation. In particular, a direct sum of two representations has
them as subrepresentations. (I won't define ``direct sum'' and
``subrepresentation'' here\ldots{} hopefully those who don't know will
be tempted to look it up.)

So rather than classifying \emph{all} representations, it's good to
start by classifying ``irreducible'' representations --- those that have
no suprepresentations (other than themselves and the trivial
0-dimensional representation). This is sort of like finding prime
numbers\ldots{} they are ``building blocks'' in representation theory.
But things are a little bit messier, alas. We say a representation is
``completely reducible'' if it is a direct sum of irreducible
representations. Unfortunately, not all representations need be
completely reducible!

Let's consider the representations of \(\mathfrak{sl}(2,\mathbb{C})\).
(The more sophisticated reader should note that I am implicitly only
considering finite- dimensional complex representations!) Here life is
as nice as could be: all representations are completely reducible, and
there is just one irreducible \(n\)-dimensional representation for each
\(n\), with the 2-dimensional and \(3\)-dimensional representations as
above. (By the way, I really mean that there is only one irreducible
\(n\)-dimensional representation up to a certain equivalence relation!)
Physicists --- who more often work with the real form
\(\mathfrak{su}(2)\) --- call these the spin-\(0\), spin-\(\frac{1}{2}\),
spin-\(1\), etc. representations. The ``spin'' of a particle is, in
mathematical terms, just the thing that tells you which representation
of \(\mathfrak{su}(2)\) it corresponds to!

Now let me jump up several levels of sophistication. In the last few
years people have realized that Lie groups are just a special case of
something called ``quantum groups''\ldots{} nobody talks about ``quantum
Lie algebras'' but that's essentially a historical accident: quantum
groups are \emph{not} groups, they're a generalization of them, and they 
\emph{don't} have Lie algebras, but they have a generalization of them --- so-called
quantized enveloping algebras.

Quantum groups can be formed from simple Lie algebras, and they depend
on a parameter \(q\), a nonzero complex parameter. This parameter ---
\(q\) is for quantum, naturally --- can be thought of as \[e^\hbar\] The
exponential of Planck's constant! When we set \(\hbar = 0\) we get
\(q = 1\), and we get back to the ``classical case'' of plain
old-fashioned Lie algebras and groups. Every representation of a quantum
group gives an invariant of links (actually even tangles), and these
link invariants are functions of \(q\). If we take the \(n\)th
derivative of one of these invariants with respect to \(\hbar\) and
evaluate it at \(\hbar = 0\) we get a ``Vassiliev invariant of degree
\(n\)'' (see \protect\hyperlink{week3}{``Week 3''} for the definition).
Better than that, when \(q\) is a root of unity each quantum group gives
us a \(3\)-dimensional ``topological quantum field theory,'' or TQFT
known as Chern--Simons theory. In particular, we get an invariant of
compact oriented 3-manifolds. So there is a hefty bunch of mathematical
payoffs from quantum groups. And there are good reasons to think of them
as the right generalization of groups for dealing with symmetries in the
physics of 2 and 3 dimensions. If string theory \emph{or} the loop
variables approach to quantum gravity have any truth to them, quantum
groups play a sneaky role in honest \(4\)-dimensional physics too.

In particular, there is a quantum version of \(\mathfrak{sl}(2)\) called
\(\mathfrak{sl}_q(2)\). When \(q = 1\) we essentially have the good old
\(\mathfrak{sl}(2)\). Chari and Premet have just worked out a lot of the
representation theory of \(\mathfrak{sl}_q(2)\). First of all, it's been
known for some time that as long as \(q\) is not a root of unity ---
that is, as long as we don't have \[q^n = 1\] for some integer \(n\) ---
the story is almost like that for ordinary \(\mathfrak{sl}(2)\). Namely,
there is one irreducible representation of each dimension, and all
representations are completely reducible. This fails at roots of unity
--- which turns out to be the reason why one can cook up TQFTs in this
case. It turns out that if \(q\) is an \(n\)th root of unity one can
still define representations of dimension 0,1,2,3, etc., more or less
just like the classical case, but only those of dimension
\(\leqslant n\) are irreducible. There are, in fact, exactly \(n\)
irreducible representations, and the fact that there are only finitely
many is what makes all sorts of neat things happen. The
\(k\)-dimensional representations with \(k \geqslant n\) are not
completely reducible. And, besides the representations that are
analogous to the classical case, there are a bunch more. They have not
been completely classified --- they are, according to Chari, a mess! But
she and Premet have classified a large batch of the ``indecomposable''
ones, that is, the ones that aren't direct sums of other ones. I guess
I'll leave it at that.

\begin{enumerate}
\def\labelenumi{\arabic{enumi})}
\setcounter{enumi}{1}
\tightlist
\item
  David Kazhdan and Iakov Soibelman, ``Representations of the quantized
  function algebras, \(2\)-categories and Zamolodchikov tetrahedra
  equations'', \emph{The Gelfand Mathematical Seminars 1990--1992}, Springer,
  Berlin, 1993, pp.\ 163--71.
\end{enumerate}
\noindent
In this terse paper, Kazhdan and Soibelman construct a braided monoidal
\(2\)-category using quantum groups at roots of unity. As I've said a
few times, people expect braided monoidal \(2\)-categories to play a
role in generally covariant 4d physics analogous to what braided
monoidal categories do in 3d physics. In particular, one might hope to
get invariants of \(4\)-dimensional manifolds, or of surfaces embedded
in 4-manifolds, this way. (See last week's post for a little bit about
the details.) I don't feel I understand this construction well enough
yet to want to say much about it, but it is clearly related to the
construction of a braided monoidal \(2\)-category from the category of
quantum group representations given by Crane and Frenkel (see
\protect\hyperlink{week2}{``Week 2''}).

\begin{enumerate}
\def\labelenumi{\arabic{enumi})}
\setcounter{enumi}{2}
\tightlist
\item
  Adrian Ocneanu, ``A note on simplicial dimension shifting'', available
  as \href{http://arxiv.org/abs/hep-th/9302028}{\texttt{hep-th/9302028}}.
\end{enumerate}
\noindent
Ouch! This paper claims to show that the very charming 4d TQFT
constructed by Crane and Yetter in ``A categorical construction of 4d
topological quantum field theories''
(\href{http://arxiv.org/abs/hep-th/9301062}{\texttt{hep-th/9301062}}) is
trivial! In particular, he says the resulting invariant of compact
oriented 4-manifolds is identically equal to 1. If so, it's back to the
drawing board. Crane and Yetter took the 3d TQFT coming from
\(\mathfrak{sl}_q(2)\) at roots of unity and then used a clever trick to
get 3-manifolds from a simplicial decomposition of a 4-manifold to get a
4d TQFT. Ocneanu claims this trick, which he calls ``simplicial
dimension shifting,'' only gives trivial 4-manifold invariants.

I am not yet in a position to pass judgement on this, since both
Crane/Yetter and Ocneanu are rather sketchy in key places. If indeed
Ocneanu is right, I think people are going to have to get serious about
facing up to the need for \(2\)-categorical thinking in
\(4\)-dimensional generally covariant physics. I had asked Crane, a big
proponent of \(2\)-categories, why they played no role in his 4d TQFT,
and he said that indeed he felt like the kid who took apart a watch, put
it back together, and found it still worked even though there was a
piece left over. So maybe the watch didn't really work without that
extra piece after all. In late March I will go to the Conference on
Quantum Topology thrown by Crane and Yetter (at Kansas State U. at
Manhattan), and I'm sure everyone will try to thrash this stuff out.

\begin{enumerate}
\def\labelenumi{\arabic{enumi})}
\setcounter{enumi}{3}
\tightlist
\item
  Abhay Ashtekar and Jerzy Lewandowski, ``Representation theory of
  analytic holonomy \(C^*\)-algebras'', available as
  \href{http://arxiv.org/abs/gr-qc/9311010}{\texttt{gr-qc/9311010}}.
\end{enumerate}
\noindent
This paper is a follow-up of the paper

\begin{enumerate}
\def\labelenumi{\arabic{enumi})}
\setcounter{enumi}{4}
\tightlist
\item
  Abhay Ashtekar and Chris Isham, ``Representations of the holonomy
  algebras of gravity and non-Abelian gauge theories'', \emph{Journal of
  Classical and Quantum Gravity} \textbf{9} (1992), 1069--1100. Also
  available as
  \href{http://arxiv.org/abs/hep-th/9202053}{\texttt{hep-th/9202053}}.
\end{enumerate}
\noindent
and sort of complements another,

\begin{enumerate}
\def\labelenumi{\arabic{enumi})}
\setcounter{enumi}{5}
\tightlist
\item
  John Baez, ``Link invariants, holonomy algebras and functional
  integration'', available as
  \href{http://arxiv.org/hep-th/9301063}{\texttt{hep-th/9301063}}.
\end{enumerate}

The idea here is to provide a firm mathematical foundation for the loop
variables representation of gauge theories, particularly quantum
gravity. Ashtekar and Lewandowski consider an algebra of gauge-invariant
observables on the space of \(\mathfrak{su}(2)\) connections on any
real-analytic manifold, namely that generated by piecewise analytic
Wilson loops. This is the sort of thing meant by a ``holonomy algebra''.
They manage to construct an explicit diffeomorphism-invariant state on
this algebra. They also relate this algebra to a similar algebra for
\(\mathfrak{sl}(2)\) connections --- the latter being what really comes
up in quantum gravity. And they do a number of other interesting things,
all quite rigorously. My paper dealt instead with an algebra generated
by ``regularized'' or ``smeared'' Wilson loops, and showed that there
was a 1-1 map from diffeomorphism-invariant states on this algebra to
invariants of framed links --- thus showing that the loop variables
picture, in which states are given by link invariants, doesn't really
lose any of the physics present in traditional approaches to gauge
theories. I am busy at work trying to combine Ashtekar and Lewandowski's
ideas with my own and push this program further --- my own personal goal
being to make the Chern--Simons path integral rigorous --- it being one
of those mysterious ``measures on the space of all connections mod gauge
transformations'' that physicists like, which unfortunately aren't
really measures, but some kind of generalization thereof. What it
\emph{should} be is a state (or continuous linear functional) on some
kind of holonomy algebra.

\hypertarget{week6}{%
\section{February 20, 1993}\label{week6}}

\begin{enumerate}
\def\labelenumi{\arabic{enumi})}
\tightlist
\item
  Alexander Vilenkin, ``Quantum cosmology'', talk given at
  Texas/Pascos 1992 at Berkeley, available as
  \href{https://arxiv.org/abs/gr-qc/9302016}{\texttt{gr-qc/9302016}}.
\end{enumerate}
\noindent
This is, as Vilenkin notes, an elementary review of quantum cosmology.
It won't be news to anyone who has kept up on that subject (except
perhaps for a few speculations at the end), but for those who haven't
been following this stuff, like myself, it might be a good way to get
started.

Let's get warmed up\ldots.

Quantizing gravity is mighty hard. For one thing, there's the ``problem
of time'' --- the lack of a distinguished time parameter in
\emph{classical} general relativity means that the usual recipe for
quantizing a dynamical system --- ``represent time evolution by the
unitary operators \(\exp(-iHt)\) on the Hilbert space of states, where \(t\) is
the time and \(H\), the Hamiltonian, is a self-adjoint operator'' --- breaks
down! As Wheeler so picturesquely put it, in general relativity we have
``many-fingered time''; there are lots of ways of pushing a spacelike
surface forwards in time.

But if we simplify the heck out of the problem, we might make a little
progress. (This is a standard method in physics, and whether or not it's
really justified, it's often the only thing one can do!) For one thing,
note that in the Big Bang cosmology there is a distinguished ``rest
frame'' (or more precisely, field of timelike vectors) given by the
galaxies, if we discount their small random motions. In reality these
are maybe not so small, and maybe not so random --- such things as the
``Virgo flow'' show this --- but we're talking strictly theory here,
okay? --- so don't bother us with facts! So, if we imagine that things
go the way the simplest Big Bang models predict, the galaxies just sit
there like dots on a balloon that is being inflated, defining a notion
of ``rest'' at each point in spacetime. This gives a corresponding
notion of time, since one can measure time using clocks that are at rest
relative to the galaxies. Then, since we are pretending the universe is
completely homogeneous and isotropic --- and let's say it's a closed
universe in the shape of a 3-sphere, to be specific --- the metric is
given by
\[dt^2-r(t)^2[(d\psi)^2 + (\sin \psi)^2{(d\theta)^2 + (\sin \theta)^2 (d\varphi)^2}]\]
What does all this mean? Here \(r(t)\) is the radius of the universe as
a function of time, the following stuff is just the usual metric on the
unit 3-sphere with hyperspherical coordinates \(\psi\), \(\theta\),
\(\varphi\) generalizing the standard coordinates on the 2-sphere we all
learn in college:
\[(d\psi)^2 + (\sin \psi)^2{(d\theta)^2 + (\sin \theta)^2 (d\varphi)^2}\]
and the fact that the metric on spacetime is \(dt^2\) minus a bunch of
stuff reflects the fact that spacetime geometry is ``Lorentzian,'' just
as in flat Minkowski space the metric is 
\[dt^2-dx^2-dy^2-dz^2.\] 
The name of the game in this simple sort of Big Bang cosmology is thus
finding the function \(r(t)\)! To do this, of course, we need to see
what Einstein's equations reduce to in this special case, and since
Einstein's equations tell us how spacetime curves in response to the
stress-energy tensor, this will depend on what sort of matter we have
around. We are assuming that it's homogeneous and isotropic, whatever it
is, so it turns out that all we need to know is its density \(\rho\) and
pressure \(P\) (which are functions of time). We get the equations
\[\begin{aligned}r''/r &= -\frac{4\pi}{3}(\rho+3P) \\ (r')^2 &= \frac{8\pi}{3}\rho r^2-1\end{aligned}\]
Here primes denote differentiation with respect to \(t\), and I'm using
units in which the gravitational constant and speed of light are equal
to 1.

Let's simplify this even more. Let's assume our matter is ``dust,''
which is the technical term for zero pressure. We get two equations:
\[\begin{aligned}r''/r &= -\frac{4\pi}{3}\rho \\ (r')^2 &= \frac{8\pi}{3}\rho r^2-1\end{aligned}\tag{1}\]
Now let's take the second one, differentiate with respect to \(t\),
\[2r''r' = \frac{8\pi}{3}(\rho'r^2 + 2 \rho r r')\] plug in what the
first equation said about \(r''\),
\[-\frac{8\pi}{3} \rho r r' = \frac{8\pi}{3}(\rho' r^2 + 2 \rho r r')\]
clear out the crud, and lo: \[3 \rho r' = -\rho' r\] or, more
enlighteningly, \[\frac{d(\rho r^3)}{dt} = 0.\] This is just
``conservation of dust'' --- the dust density times the volume of the
universe is staying constant. This, by the way, is a special case of the
fact that Einstein's equations \emph{automatically imply} local
conservation of energy (i.e., that the stress-energy tensor is
divergence-free).

Okay, so let's say \(\rho r^3 = D\), with \(D\) being the total amount
of dust. Then we can eliminate \(\rho\) from equations (1) and get:
\[\begin{aligned}r'' &= -\frac{4\pi D}{3r^2} \\ (r')^2-\frac{8\pi}{3}\frac{D}{r} &= -1 \end{aligned}\tag{2}\]
What does this mean? Well, the first one looks like it's saying there's
a force trying to make the universe collapse, and that the strength of
this force is proportional to \(1/r^2\). Sound vaguely familiar? It's
actually misleadingly simple --- if we had put in something besides dust
it wouldn't work quite this way --- but as long as we don't take it too
seriously, we can just think of this as gravity trying to get the
universe to collapse. And the second one looks like it's saying that the
``kinetic'' energy proportional to \((r')^2\), plus the ``potential''
energy proportional to \(-1/r\), is constant! In other words, we have a
nice analogy between the Big Bang cosmology and a very old-fashioned
system, a classical particle in one dimension attracted to the origin by
a \(1/r^2\) force!

It's easy enough to solve this equation, and easier still to figure it
out qualitatively. The key thing is that since the total ``energy'' in
the second equation of (2) is negative, there won't be enough ``energy''
for \(r\) to go to infinity, that is, there'll be a Big Bang and then a
big crunch. Here's \(r\) as a function of \(t\), roughly: \[
  \begin{tikzpicture}
    \draw[->] (0,0) to (5,0) node[label=below left:{$t$}]{};
    \draw[->] (0,0) to (0,2.5) node[label=below left:{$r$}]{};
    \draw[thick] plot [smooth,tension=1] coordinates {(0,0) (2,1.5) (4,0)};
  \end{tikzpicture}
\] 
\vskip -1em
\centerline{Figure 1}
\vskip 1em
\noindent
What goes up, must come down! This curve, which I haven't drawn too
well, is just a cycloid, which is the curve traced out by a point on the
rim of rolling wheel. So, succumbing to romanticism momentarily we could
call this picture \emph{one turn of the great wheel of time}\ldots.
But there is \emph{no} reason to expect further turns, because the
differential equation simply becomes singular when \(r = 0\). We may
either say it doesn't make sense to speak of ``before the Big Bang'' or
``after the big crunch'' --- or we can look for improved laws that avoid
these singularities. (I should repeat that we are dealing with
unrealistic models here, since for example there is no evidence that
there is enough matter around to ``close the universe'' and make this
solution qualitatively valid --- it may well be that there's a Big Bang
but no big crunch. In this case, there's only one singularity to worry
about, not two.)

People have certainly not been too ashamed to study the \emph{quantum}
theory of this system (and souped-up variants) in an effort to get a
little insight into quantum gravity. We would expect that quantum
effects wouldn't matter much until the radius of the universe is very
small, but when it \emph{is} very small they would matter a lot, and
maybe --- one might hope --- they would save the day, preventing the
nasty singularities. I'm not saying they \emph{do} --- this is hotly debated
--- but certainly some people hope they do. Of course, serious quantum
gravity should take into account the fact that geometry of spacetime has
all sorts of wiggles in it --- it isn't just a symmetrical sphere. This may
make a vast difference in how things work out. (For example, the big
crunch would be a lot more exciting if there were lots of black holes
around by then.) The technical term for the space of all metrics on
space is ``superspace'' (sigh), and the toy models one gets by ignoring
all but finitely many degrees of freedom are called ``minisuperspace''
models.

Let's look at a simple minisuperspace model. The simplest thing to try
is to take the classical equations of motion (2) and try to quantize
them just like one would a particle in a potential. This is a delicate
business, by the way, because one can't just take some classical
equations of motion and quantize them in any routine way. There are lots
of methods of quantization, but all of them require a certain amount of
case-by-case finesse.

The idea of ``canonical quantization'' of a classical system with one
degree of freedom --- like our Big Bang model above, where the one
degree of freedom is \(r\) --- is to turn the ``position'' (that's
\(r\)) into a multiplication operator and the ``momentum'' (often that's
something like \(r'\), but watch out!) into a differentation operator,
say \(-i \hbar \frac{d}{dr}\), so that we get the ``canonical
commutation relations'' \[[-i \hbar \frac{d}{dr}, r] = -i \hbar.\] We
then take the formula for the energy, or Hamiltonian, in terms of
position and momentum, and plug in these operators, so that the
Hamiltonian becomes an operator. (Here various ``operator-ordering''
problems can arise, because the position and momentum commuted in the
original classical system but not anymore!) To explain what I mean, why
don't I just do it!

So: I said that the formula
\[(r')^2-\frac{8\pi}{3} \frac{D}{r} = -1 \tag{3}\] looks a lot like a
formula of the form ``kinetic energy plus potential energy is
constant''. Of course, we could multiply the whole equation by anything
and get a valid equation, so it's not obvious that the ``right''
Hamiltonian is \[(r')^2-\frac{8\pi}{3} \frac{D}{r}\] or (adding 1
doesn't hurt) \[(r')^2-\frac{8\pi}{3} \frac{D}{r} + 1\] In fact, note
that multiplying the Hamiltonian by some function of \(r\) just amounts
to reparametrizing time, which is perfectly fine in general relativity.
In fact, Vilenkin and others before him have decided it's better to
multiply the Hamiltonian above by \(r^2\). Why? Well, it has to do with
figuring out what the right notion of ``momentum'' is corresponding to
the ``position'' \(r\). Let's do that. We use the old formula
\[p = \frac{dL}{dq'}\] relating momentum to the Lagrangian, where for us
the position, usually called \(q\), is really \(r\).

The Lagrangian of general relativity is the ``Ricci scalar'' \(R\) --- a
measure of curvature of the metric --- and in the present problem it
turns out to be
\[R = 6 \left(\frac{r''}{r} + \frac{(r')^2}{r^2}\right)\] But we are
reducing the full field theory problem down to a problem with one degree
of freedom, so our Lagrangian should be the above integrated over the
3-sphere, which has volume \(16 \pi r^3/3\), giving us
\[32\pi (r''r^2 + (r')^2 r)\] However, the \(r''\) is a nuisance, and we
only use the integral of the Lagrangian with respect to time (that's the
action, which classically is extremized to get the equations of motion),
so let's do an integration by parts, or in other words add a total
divergence, to get the Lagrangian \[L = -32\pi (r')^2 r.\]
Differentiating with respect to \(r'\) we get the momentum ``conjugate
to \(r\)'', \[p = -64\pi r'r.\] Now I notice that Vilenkin uses as the
momentum simply \(-r'r\), somehow sweeping the monstrous \(64\pi\) under
the rug. I have the feeling that this amounts to pushing this factor
into the definition of \(\hbar\) in the canonical commutation relations.
Since I was going to set \(\hbar\) to 1 in a minute anyway, this is okay
(honest). So let's keep life simple and use \[p = -r'r.\] Okay! Now
here's the point, we want to exploit the analogy with good old quantum
mechanics, which typically has Hamiltonians containing something like
\(p^2\). So let's take our preliminary Hamiltonian
\[(r')^2-\frac{8\pi}{3} \frac{D}{r} + 1\] and multiply it by \(r^2\),
getting \[H = p^2-\frac{8\pi D}{3}r + r^2.\] Hey, what's this? A
harmonic oscillator! (Slightly shifted by the term proportional to
\(r\).) So the universe is just a harmonic oscillator\ldots{} I guess
that's why they stressed that so much in all my classes!

Actually, despite the fact that we are working with a very simple model
of quantum cosmology, it's not quite \emph{that} simple. First of all,
recall our original classical equation, (3). This constrained the energy
to have a certain value. I.e., we are dealing not with a Hamiltonian in
the ordinary sense, but a ``Hamiltonian constraint'' --- typical of
systems with time reparametrization invariance. So our quantized
equation says that the ``wavefunction of the universe,'' \(\psi(r)\),
must satisfy \[H \psi = 0.\] Also, unlike the ordinary harmonic
oscillator we have the requirement that \(r\geqslant0\). In other word,
we're working with a problem that's like a harmonic oscillator and a
``wall'' that keeps \(r\geqslant0\). Think of a particle in a potential
like this: 
\newpage
\[
  \begin{tikzpicture}
    \draw[->] (0,0) to (5,0) node [label=below left:{$r$}]{};
    \draw[->] (0,-0.7) to (0,2.5) node[label=below left:{$V(r)$}]{};
    \draw[thick] plot [smooth,tension=1] coordinates {(0,0) (1.8,-0.5) (4,2.5)};
  \end{tikzpicture} 
\] 
\vskip -1em
\centerline{Figure 2}
\vskip 1em
\noindent
Here \(V(r) = -(8\pi D/3)r + r^2\). The minimum of \(V\) is at
\(r = 4 \pi D/3\) and the zeroes are at \(r = 0\) and \(8 \pi D/3\).
Classically, a particle with zero energy starting at \(r = 0\) will roll
to the right and make it out to \(r = 4\pi D/3\) before rolling back to
\(r = 0\). This is basically the picture we had in Figure 1, except that
we've reparametrized time so we have simple harmonic motion instead of
cycloid.

Quantum mechanically, however one must pick boundary conditions at
\(r = 0\) to make the problem well-defined!

This is where the fur begins to fly!! Hawking and Vilenkin have very
different ideas about what the right boundary conditions are. And note
that this is not a mere technical issue, since they determine the
wavefunction of the universe in this approach! I will not discuss this
since Vilenkin does so quite clearly, and if you understand what I have
written above you'll be in a decent position to understand him. I will
just note that Vilenkin, rather than working with a universe full of
``dust,'' considers a universe in which the dominant contribution to the
stress-energy tensor is the cosmological constant, that is, the negative
energy density of a ``false vacuum'', which believers in inflation (such
as Vilenkin) think powered the exponential growth of the universe at an
early stage. So his equations are slightly different from those above
(and are only meant to apply to the early history of the universe).

{[}Let me just interject a question to the experts if I may --- since
I've written this long article primarily to educate myself. It would
seem to me that the equation \(H \psi = 0\) above would only have a
normalizable solution if the boundary conditions were fine-tuned! I.e.,
maybe the equation \(H \psi = 0\) itself determines the boundary
conditions! This would be very nice; has anyone thought of this? It
seems reasonable because, with typical boundary conditions, the operator
\(H\) above will have pure point spectrum (only eigenvalues) and it
would be rather special for one of them to be \(0\), allowing a
normalizable solution of \(H \psi = 0\). Also, corrections and education
of any sort are welcomed. I would love to discuss this with some
experts.{]}

Anyway, suppose we find some boundary conditions and calculate \(\psi\),
the ``wavefunction of the universe.'' (I like repeating that phrase
because it sounds so momentous, despite the fact that we are working
with a laughably oversimplified toy model.) What then? What are the
implications for the man in the street?

Let me get quite vague at this point. Think of the radius of the
universe as analogous to a particle moving in the potential of Figure 2.
In the current state of affairs classical mechanics is an excellent
approximation, so it seems to trace out a classical trajectory. Of
course it is really obeying the laws of quantum mechanics, so the
trajectory is really a ``wave packet'' --- technically, we use the WKB
approximation to see how the wave packet can seem like a classical
trajectory. But near the Big Bang or big crunch, quantum mechanics
matters a lot: there the potential is rapidly varying (in our simple
model it just becomes a ``wall'') and the wave packet may smear out
noticeably. (Think of how when you shoot an electron at a nucleus it
bounces off in an unpredictable direction --- it's wavefunction just
tells you the \emph{probability} that it'll go this way or that!) So
some quantum cosmologists have suggested that if there is a big crunch,
the universe will pop back out in a highly unpredictable, random kind of
way!

I should note that Vilenkin has a very different picture. Since this
stuff makes large numbers of assumptions with very little supporting
evidence, it is science that's just on the brink of being mythology.
Still, it's very interesting.

\begin{enumerate}
\def\labelenumi{\arabic{enumi})}
\setcounter{enumi}{1}
\tightlist
\item
  Lee Smolin, ``Finite, diffeomorphism invariant observables in quantum
  gravity'', available as
  \href{https://arxiv.org/abs/gr-qc/9302011}{\texttt{gr-qc/9302011}}.
\end{enumerate}

The big problem in canonical quantization of gravity, once one gets
beyond ``mini-'' and ``midisuperspace'' models, is to find enough
diffeomorphism-invariant observables. There is a certain amount of
argument about this stuff, and various approaches, but one common
viewpoint is that the ``physical'' observables, that is, the really
observable observables, in general relativity are those that are
invariant under all diffeomorphisms of spacetime. I.e., those that are
independent of any choice of coordinates. For example, saying ``My
position is \((242,2361,12,-17)\)'' is not diffeomorphism-invariant, but
saying ``I'm having the time of my life'' is. It's hard to find lots of
(tractable) diffeomorphism invariant observables --- or even any! Try
figuring out how you would precisely describe the shape of a rock
without introducing any coordinates, and you'll begin to see the
problem. (The quantum mechanical aspects make it harder.)

Rovelli came up a while back with a very clever angle on this problem.
It's rather artificial but still a big start. Using a ``field of
clocks'' he was able to come up with interesting diffeomorphism
invariant observables. The idea is simply that if you had clocks all
around you could say ``when the bells rang 2 a.m.\ I was having the time
of my life'' --- and this would be a diffeomorphism-invariant statement,
since rather than referring to an abstract coordinate system it
expresses the coincidence of two physical occurences, just like ``the
baseball broke through the window''. Then he pushed this idea to define
``evolving constants of motion'' --- a deliberate oxymoron --- to deal
with the famous ``problem of time'' in general relativity: how to treat
time evolution in a coordinate-free manner on a spacetime that's not
flat and, worse, whose geometry is ``uncertain'' a la Heisenberg? This
is treated, by the way, in

\begin{enumerate}
\def\labelenumi{\arabic{enumi})}
\setcounter{enumi}{2}
\tightlist
\item
  Carlo Rovelli, ``Time in quantum gravity: an hypothesis'', \emph{Phys.\
  Rev.} \textbf{D43} (1991), 442--456.
\end{enumerate}
\noindent
Also, an excellent and very thorough review of the problem of time and
various proposed solutions, including Rovelli's, is given in

\begin{enumerate}
\def\labelenumi{\arabic{enumi})}
\setcounter{enumi}{3}
\tightlist
\item
  Chris J. Isham, ``Canonical quantum gravity and the problem of time'',
  125 pages, available as
  \href{https://arxiv.org/abs/gr-qc/9210011}{\texttt{gr-qc/9210011}}.
\end{enumerate}
Anyway, in a paper I very briefly described in
\protect\hyperlink{week1}{``Week 1''}:
\begin{enumerate}
\def\labelenumi{\arabic{enumi})}
\setcounter{enumi}{4}
\tightlist
\item
  Lee Smolin, ``Time, measurement and information loss in quantum
  cosmology'', available as
  \href{https://arxiv.org/abs/gr-qc/9301016}{\texttt{gr-qc/9301016}}.
\end{enumerate}
Smolin showed, how, using a clever trick sketched in the present paper
to get ``observables'' invariant under spatial but not temporal
observables, together with Rovelli's idea, one could define lots of \emph{real}
observables, invariant under spacetime diffeomorphisms that is, thus
making a serious bite into this problem.

I warn the reader that there is a fair amount that is not too realistic
about these methods. First there's the ``clock field'' --- this can
actually be taken as a free massless scalar field, but in so doing there
is the likelihood of serious technical problems. Some of these are
discussed in

\begin{enumerate}
\def\labelenumi{\arabic{enumi})}
\setcounter{enumi}{5}
\tightlist
\item
  P. Hajicek, ``Comment on `Time in quantum gravity --- an
  hypothesis''', \emph{Phys. Rev.} \textbf{D44} (1991), 1337--1338.
\end{enumerate}
(But I haven't actually read this, just Isham's description.) Also, the
clever trick of the present paper is to couple gravity to an
antisymmetric tensor gauge field so that in addition to having loops as
part of one's ``loop representation,'' one has surfaces --- a ``surface
representation''. But this antisymmetric tensor gauge field is not the
sort of thing that actually seems to arise in physics (unless I'm
missing something). Still, it's a start. I think I'll finish by quoting
Smolin's abstract:

\begin{quote}
Two sets of spatially diffeomorphism invariant operators are constructed
in the loop representation formulation of quantum gravity. This is done
by coupling general relativity to an anti- symmetric tensor gauge field
and using that field to pick out sets of surfaces, with boundaries, in
the spatial three manifold. The two sets of observables then measure the
areas of these surfaces and the Wilson loops for the self-dual
connection around their boundaries. The operators that represent these
observables are finite and background independent when constructed
through a proper regularization procedure. Furthermore, the spectra of
the area operators are discrete so that the possible values that one can
obtain by a measurement of the area of a physical surface in quantum
gravity are valued in a discrete set that includes integral multiples of
half the Planck area. These results make possible the construction of a
correspondence between any three geometry whose curvature is small in
Planck units and a diffeomorphism invariant state of the gravitational
and matter fields. This correspondence relies on the approximation of
the classical geometry by a piecewise flat Regge manifold, which is then
put in correspondence with a diffeomorphism invariant state of the
gravity-matter system in which the matter fields specify the faces of
the triangulation and the gravitational field is in an eigenstate of the
operators that measure their areas.
\end{quote}

\begin{center}\rule{0.5\linewidth}{0.5pt}\end{center}

\begin{quote}
\emph{In the Space and Time marriage we have the greatest Boy meets Girl
story of the age. To our great-grandchildren this will be as poetical a
union as the ancient Greek marriage of Cupid and Psyche seems to us.}

--- Lawrence Durrell
\end{quote}

\hypertarget{week7}{%
\section{March 1, 1993}\label{week7}}

\begin{enumerate}
\def\labelenumi{\arabic{enumi})}
\tightlist
\item
  Abhay Ashtekar, ``Mathematical problems of non-perturbative quantum general relativity''
  (lectures delivered at the 1992 Les Houches summer school on
  Gravitation and Quantization), 87 pages, Plain TeX,
  available as
  \href{https://arxiv.org/abs/gr-qc/9302024}{\texttt{gr-qc/9302024}}.
\end{enumerate}

I described this paper in \protect\hyperlink{week3}{``Week 3''}, but now
it's available from gr-qc. It's a good quick introduction to the loop
representation of quantum gravity.

\begin{enumerate}
\def\labelenumi{\arabic{enumi})}
\setcounter{enumi}{1}
\tightlist
\item
  Abhay Ashtekar, \emph{Lectures on Non-perturbative Canonical Gravity}, World Scientific Press, 1991.
\end{enumerate}

This book, which I finally obtained, is \emph{the} introduction to the
loop representation of quantum gravity. What's the loop representation?
Well, this is a long story, so you really should read the book. But just
to get you going, let me describe Ashtekar's ``new variables,'' which
form the basis for Rovelli and Smolin's construction of the loop
representation.

First, recall that general relativity is usually thought of as a theory
about a metric on spacetime --- more precisely, a Lorentzian metric.
Here spacetime is a \(4\)-dimensional manifold, and a Lorentzian metric
allows you to calculate the ``dot product'' of any two tangent vectors
at a point. This is in quotes because, while a normal dot product might
look like
\[(v_0,v_1,v_2,v_3)\cdot(w_0,w_1,w_2,w_3) = v_0w_0 + v_1w_1 + v_2w_2 + v_3w_3\]
relative to some basis, for a Lorentzian metric we can always find a
basis of the tangent space such that
\[(v_0,v_1,v_2,v_3)\cdot(w_0,w_1,w_2,w_3) = v_0w_0 -v_1w_1 -v_2w_2 -v_3w_3\]
Now the metric in general relativity defines a ``connection,'' which
tells you a tangent vector might ``twist around'' as you parallel
translate it, that is, move it along while trying to keep it from
rotating unnecessarily. Here ``twist around'' is in quotes because,
since you are parallel translating the vector, it's not really
``twisting around'' in the usual sense, but it might seem that way
relative to some coordinate system. For example, if you used polar
coordinates to describe parallel translation on the plane, it might seem
that the unit vector in the \(r\) direction ``twisted around'' towards the
\(\theta\) direction as you dragged it along. But in another coordinate
system --- say the usual \(x\)-\(y\) system --- it would not appear to
be ``twisting around''. This fact is expressed by saying ``the
connection is not a tensor''.

But from the connection we can cook up a big fat tensor, the ``Riemann
tensor'' \(R^i_{jkl}\), which says how much the vector in the \(l\)th
direction (here the indices range from 0 to 3) twists towards the
\(i\)th direction when you move it around a teeny little square in the
\(j\)-\(k\) plane. The Lagrangian in ordinary GR is just the integral of
the ``Ricci scalar curvature,'' \(R\), which is gotten from the Riemann
tensor by ``contraction'', i.e.~summing over the indices in a certain
way: \[R = R^i_{ji}{}^j\] where we are raising indices using the metric
in a manner beloved by physicists and feared by many mathematicians. If
you integrate the Lagrangian over a region of spacetime you get the
``action'', and in classical general relativity (in a vacuum, for
simplicity) one can formulate the laws of motion simply by saying: any
teeny change in the metric that vanishes on the boundary of the region
should leave the action constant to first order. In other words, the
solutions of the equations of general relativity are the *stationary
points* of the action. If you know how to do variational calculus you
can derive Einstein's equations from this variational principle, as it's
called. Mathematicians will be pleased to know that Hilbert beat
Einstein to the punch here, so the integral of \(R\) is called the
``Einstein--Hilbert'' action for general relativity.

But there's another formulation of general relativity in terms of an
action principle. This is called the ``Palatini'' action --- and
actually I'm going to describe a slight variation on it, that is
conceptually simpler, and apparently appears for the first time in
Ashtekar's book. The Palatini approach turns out to be more elegant and
is a nice stepping-stone to the Ashtekar approach. In the Palatini
approach one thinks of general relativity not as being a theory of a
metric, but of a ``tetrad'' and an ``\(\mathfrak{so}(3,1)\)
connection''. To explain what these are, I will cut corners and assume
all the fiber bundles lurking around are trivial; the experts will
easily be able to figure out the general case. So: an (orthonormal)
tetrad, or ``vierbein,'' is a just a kind of field on spacetime which at
each point consists of an (ordered) orthonormal basis of the tangent
space. If we express the metric in terms of a tetrad, it looks just like
the formula for the standard ``inner product''
\[(v_0,v_1,v_2,v_3)\cdot(w_0,w_1,w_2,w_3) = v_0w_0 -v_1w_1 -v_2w_2 -v_3w_3\]
This allows us to identify the group of linear transformations of the
tangent space that preserve the metric with the group of linear
transfomations preserving the standard ``inner product,'' which is
called \(\mathrm{SO}(1,3)\) since there's one plus sign and three
minuses. And from the connection mentioned above one gets an
\(\mathrm{SO}(1,3)\) connection, or, what's more or less the same thing,
an \(\mathfrak{so}(1,3)\)-valued \(1\)-form, that is, a kind of field
that can eat a tangent vector at any point and spits out element of the
Lie algebra \(\mathfrak{so}(1,3)\).

What's \(\mathfrak{so}(1,3)\)? Well, elements of \(\mathfrak{so}(1,3)\)
include ``infinitesimal'' rotations and Lorentz transformations, since
\(\mathrm{SO}(1,3)\) is generated by rotations and Lorentz
transformations. More precisely, \(\mathfrak{so}(1,3)\) is a
6-dimensional Lie algebra having as a basis the three infinitesimal
rotations \(J_1\), \(J_2\), and \(J_3\) around the three axes, and the
three infinitesimal Lorentz transformations or ``boosts'' \(K_1\),
\(K_2\), \(K_3\). The bracket in this most important Lie algebra is
given by
\[\begin{aligned}\,[J_i,J_j] & = J_k \\ [K_i,K_j] &= -J_k \\ [J_i,K_j] &= K_k\end{aligned}\]
where \((i,j,k)\) is a cyclic permutation of \((1,2,3)\). (I hope I
haven't screwed up the signs.) Note that the \(J\)'s by themselves form
a Lie subalgebra called \(\mathfrak{so}(3)\), the Lie algebra of the
rotation group \(\mathrm{SO}(3)\). Note that \(\mathfrak{so}(3)\) is
isomorphic to the cute little Lie algebra \(\mathfrak{su}(2)\) I
described in my post \protect\hyperlink{week5}{``Week 5''}; \(J_1\),
\(J_2\), and \(J_3\) correspond to the guys \(I\), \(J\), and \(K\)
divided by two.

The \(\mathfrak{so}(1,3)\) connection has a curvature, and using the
tetrads again we can identify this with the Riemann curvature tensor. So
the Palatini trick is to rewrite the Einstein--Hilbert action in terms of
the curvature of the \(\mathfrak{so}(1,3)\) connection and the tetrad
field. This is called the Palatini action. Charmingly, even though the
tetrad field is utterly unphysical, we can treat it and the
\(\mathfrak{so}(1,3)\) connection as independent fields and, doing
calculus of variations to find stationary points of the action, we get
equations equivalent to Einstein's equations.

Ashtekar's ``new variables'' --- from this point of view --- rely on a
curious and profound fact about \(\mathfrak{so}(1,3)\). Note that
\(\mathfrak{so}(1,3)\) is a Lie algebra over the real numbers. But if we
allow ourselves to form \emph{complex} linear combinations of the
\(J\)'s and \(K\)'s, thus:
\[\begin{aligned}M_i &= (J_i + iK_i)/2 \\ N_i &= (J_i -iK_i)/2\end{aligned}\]
(please don't mix up the subscript \(i = 1,2,3\) with the other \(i\),
the square root of minus one) we get the following brackets:
\[\begin{aligned}\,[M_i,M_j] &= M_k \\ [N_i,N_j] &= N_k \\ [M_i,N_j] &= 0\end{aligned}\]
I think the signs all work but I wouldn't trust me if I were you. The
wonderful thing here is that the \(M\)'s and \(N\)'s commute with each
other, and each set has commutation relations just like the \(J\)'s! The
\(J\)'s, recall, are infinitesimal rotations, and the Lie algebra they
span is \(\mathfrak{so}(3)\). So in a sense the Lie algebra of the
Lorentz group can be ``split'' into ``left-handed'' and ``right-handed''
copies of \(\mathfrak{so}(3)\), also known as ``self-dual'' and
``anti-self-dual'' copies. This is, in fact, what lies behind the
handedness of neutrinos, and many other wonderful things.

But let me phrase this result more precisely. Since we allowed ourselves
complex linear combinations of the \(J\)'s and \(K\)'s, we are now
working in the ``complexification'' of the Lie algebra
\(\mathfrak{so}(3,1)\), and we have shown that this Lie algebra over the
complex numbers splits into two copies of
\(\mathfrak{so}(3,\mathbb{C})\), the complexification of
\(\mathfrak{so}(3)\).

Ashtekar came up with some ``new variables'' for general relativity in
the context of the Hamiltonian approach. Here we are working in the
Lagrangian approach, where things are simpler because they are
``generally covariant,'' not requiring a split of spacetime into space
and time. The Lagrangian approach to the new variables is due to Samuel,
Jacobson and Smolin, and in this approach all they amount to is this:
take the \(\mathfrak{so}(1,3)\) connection of the Palatini approach, think of 
\(\mathfrak{so}(1,3)\) as sitting inside the complexification thereof,
and consider only the ``right-handed'' part! Thus, from an
\(\mathfrak{so}(1,3)\) connection, we get a
\(\mathfrak{so}(3,\mathbb{C})\) connection. The ``new variables'' are
just the tetrad field and this \(\mathfrak{so}(3,\mathbb{C})\)
connection.

I have tried to keep down the indices but I think I will write down the
Palatini Lagrangian and then the ``new variables'' Lagrangian, without
explaining exactly what they mean, just to show how amazingly
similar-looking they are. In the Palatini approach we have a tetrad
field, which now we write in its full glory as \(e_I^i\), and the
curvature of the \(\mathfrak{so}(1,3)\) connection, which now we write
as \(\Omega_{ij}^{IJ}\). The Lagrangian is then
\[e_I^i e_J^j \Omega_{ij}^{IJ}\] (which we integrate against the usual
volume form to get the action). In the new variables approach we have a
tetrad field again, and we write the curvature of the
\(\mathfrak{so}(3,\mathbb{C})\) connection as \(F_{ij}^{IJ}\). (This
turns out to be just the ``right-handed'' part of \(\Omega_{ij}^{IJ}\).)
The Lagrangian is \[e_I^i e_J^j F_{ij}^{IJ} !\] Miraculously, this also
gives Einstein's equations.

What's utterly unclear from what I've said so far is why this helps so
much in trying to quantize gravity. I may eventually get around to
writing about that, but for now, read the book!

\begin{enumerate}
\def\labelenumi{\arabic{enumi})}
\setcounter{enumi}{2}
\tightlist
\item
  David Yetter and Louis Crane, ``We are not stuck with gluing'',
  available as
  \href{https://arxiv.org/abs/hep-th/9302118}{\texttt{hep-th/9302118}}.
\end{enumerate}
\noindent
Well, in \protect\hyperlink{week2}{``Week 2''} I mentioned Crane and
Yetter's marvelous construction of a 4d topological quantum field theory
using the representations of the quantum group \(SU_q(2)\) --- and in
\protect\hyperlink{week5}{``Week 5''} I mentioned Ocneanu's ``proof''
that the resulting 4-manifold invariants were utterly trivial (equal to
1 for all 4-manifolds). Now Crane and Yetter have replied, saying that
their 4-manifold invariants are not trivial and that Ocneanu interpreted
their paper incorrectly. I look forward to the conference on quantum
topology in Kansas at the end of May, where the full story will
doubtless come out.

\begin{enumerate}
\def\labelenumi{\arabic{enumi})}
\setcounter{enumi}{3}
\tightlist
\item
   R.\ Capovilla, J.\ Dell and T.\ Jacobson, 
   ``The initial value problem in light of Ashtekar's variables'', available as
  \href{https://arxiv.org/abs/gr-qc/9302020}{\texttt{gr-qc/9302020}}.
\end{enumerate}
\noindent
The advantage of Ashtekar's new variables is that they simplify the form
of the constraint equations one gets in the initial-value problem for
general relativity. This is true both of the classical and quantum
theories. Rovelli and Smolin used this to find, for the first time, lots
of states of quantum gravity defined by link invariants. Here the above
authors are trying to apply the new variables to the \emph{classical}
theory.

\begin{enumerate}
\def\labelenumi{\arabic{enumi})}
\setcounter{enumi}{4}
\tightlist
\item
  Sergey Piunikhin, ``Combinatorial expression for universal Vassiliev link invariant'', 
  available as
  \href{https://arxiv.org/abs/hep-th/9302084}{\texttt{hep-th/9302084}}.
\end{enumerate}

Somebody ought to teach those Russians how to use the word ``the'' now
that the cold war is over. Anyway, this paper defines a kind of
universal object for Vassiliev invariants, which is sort of similar to
what I was trying to do in

\begin{enumerate}
\def\labelenumi{\arabic{enumi})}
\setcounter{enumi}{5}
\tightlist
\item
John Baez, ``Link invariants of finite type and perturbation theory'', \emph{Lett.\ Math.\ Phys.} \textbf{26} (1992) 43--51.  Also available as \href{https://arxiv.org/abs/hep-th/9207041}
{\texttt{hep-th/9207041}}.
\end{enumerate}
\noindent
but more concrete, and (supposedly) simpler than Kontsevich's approach.
My parenthesis simply indicates that I haven't had time to figure out
what's going on here.

\hypertarget{week8}{%
\section{March 5, 1993}\label{week8}}

I was delighted to find that Louis Kauffman wants to speak at the
workshop at UCR on knots and quantum gravity; he'll be talking on
``Temperley--Lieb recoupling theory and quantum invariants of links and
manifolds''. His books

\begin{enumerate}
\def\labelenumi{\arabic{enumi})}
\item
   Louis Kauffman, 
   \emph{On Knots}, Annals of Mathematics Studies \textbf{115}, Princeton
   U.\ Press, Princeton, 1987.
\end{enumerate}
\noindent
and more recently

\begin{enumerate}
\def\labelenumi{\arabic{enumi})}
\setcounter{enumi}{1}
\item
    Louis Kauffman, 
  \emph{Knots and Physics}, World Scientific Press, Singapore, 1991.
\end{enumerate}
\noindent
are a lot of fun to read, and convinced me to turn my energies towards
the intersection of knot theory and mathematical physics. As you can see
by the title of the series he's editing, he is a true believer the deep
significance of knot theory. This was true even before the Jones
polynomial hit the mathematical physics scene, so he was well-placed to
discover the relationship between the Jones polynomial (and other new
knot invariants) and statistical mechanics, which seems to be what won
him his fame. He is now the editor of a journal, \emph{Journal of Knot
Theory and its Ramifications}.

He sent me a packet of articles and preprints which I will briefly
discuss. If you read \emph{any} of the stuff below, \emph{please} read
the delightful reformulation of the 4-color theorem in terms of cross
products that he discovered! I am strongly tempted to assign it to my
linear algebra class for homework\ldots.

\begin{enumerate}
\def\labelenumi{\arabic{enumi})}
\setcounter{enumi}{2}
\item
   Louis Kauffman, ``Map coloring and the vector cross product'',
  \emph{J.\ Comb.\ Theory} \textbf{48} (1990), 145--154.

   Louis Kauffman, ``Map coloring, $q$-deformed spin networks, and Turaev--Viro invariants
  for 3-manifolds'', \emph{Int.\ Jour.\  Mod.\ Phys.}
  \textbf{6} (1992), 1765--1794.

   Louis Kauffman, ``An algebraic approach to the planar colouring problem'', 
   \emph{Commun.\ Math.\ Phys.} \textbf{152} (1993), 565--590.
\end{enumerate}
\noindent
As we all know, the usual cross product of vectors in \(\mathbb{R}^3\)
is not associative, so the following theorem is slightly interesting:

\begin{quote}
{\rm 
\textbf{Theorem.} Consider any two bracketings of a product of any
finite number of vectors, e.g.:
\[L = a \times (b \times ((c \times d) \times e) \quad\text{and}\quad  R = ((a \times b) \times c) \times (d \times e).\]
Let \(i\), \(j\), and \(k\) be the usual canonical basis for
\(\mathbb{R}^3\): \[i = (1,0,0), \quad j = (0,1,0), \quad k = (0,0,1).\]
Then we may assign \(a,b,c,\ldots\) values taken from \(\{i,j,k\}\) in
such a way that \(L = R\) and both are nonzero.
} 
\end{quote}

But what's really interesting is:

\begin{quote}
{\rm 
\textbf{Meta-Theorem.} The above proposition is equivalent to the
4-color theorem. Recall that this theorem says that any map on the plane
may be colored with 4 colors in such a way that no two regions with the
same color share a border (an edge).
}
\end{quote}

What I mean here is that the only way known to prove this Theorem is to
deduce it from the 4-color theorem, and conversely, any proof of this
Theorem would easily give a proof of the 4-color theorem! As you all
probably know, the 4-color theorem was a difficult conjecture that
resisted proof for about a century before succumbing to a computer-based
proof that required the consideration of many, many special cases:

\begin{enumerate}
\def\labelenumi{\arabic{enumi})}
\setcounter{enumi}{3}
\item
  Kenneth Appel and Wolfgang Haken, \emph{Every Planar Map is Four Colorable}
  Contemporary Mathematics \textbf{98} American Mathematical Society, Providence
  Rhode Island, 1989.
\end{enumerate}

So the Theorem above may be regarded as a \emph{profoundly} subtle
result about the ``associativity'' of the cross product!

Of course, I hope you all rush out now and find out how this Theorem is
equivalent to the 4-color theorem. For starters, let me note that it
uses a result of Tait: first, to prove the 4-color theorem it's enough
to prove it for maps where only 3 countries meet at each vertex (since
one can stick in a little new country at each vertex), and second,
4-coloring such a map is equivalent to coloring the \emph{edges} with 3
colors in such a way that each vertex has edges of all 3 colors
adjoining it. The 3 colors correspond to \(i\), \(j\), and \(k\)!

Kauffman and Saleur (the latter a physicist) come up with another
algebraic formulation of the 4-color theorem in terms of the
Temperley--Lieb algebra. The Temperley--Lieb algebra \(TL_n\) is a cute
algebra with generators \(e_1, \ldots, e_{n-1}\) and relations that
depend on a constant \(d\) called the ``loop value'':
\[\begin{aligned}e_i^2 &= de_i \\ e_i e_{i+1} e_i &= e_i \\ e_i e_{i-1} e_i &= e_i \\ e_i e_j &= e_j e_i \quad\text{for } |i-j| > 1.\end{aligned}\]
The point of it becomes clear if we draw the \(e_i\) as tangles on \(n\)
strands. Let's take \(n = 3\) to keep life simple. Then \(e_1\) is \[
  \begin{tikzpicture}
    \begin{knot}[clip width=7]
      \strand[thick] (0,0)
        to [out=down,in=down,looseness=2] (1,0);
      \strand[thick] (0,-2)
        to [out=up,in=up,looseness=2] (1,-2);
      \strand[thick] (2,0)
        to (2,-2);
    \end{knot}
  \end{tikzpicture}
\] while \(e_2\) is \[
  \begin{tikzpicture}
    \begin{knot}[clip width=7]
      \strand[thick] (1,0)
        to [out=down,in=down,looseness=2] (2,0);
      \strand[thick] (1,-2)
        to [out=up,in=up,looseness=2] (2,-2);
      \strand[thick] (0,0)
        to (0,-2);
    \end{knot}
  \end{tikzpicture}
\] In general, \(e_i\) ``folds over'' the \(i\)th and \((i+1)\)st
strands. Note that if we square \(e_i\) we get a loop --- e.g., \(e_1\)
squared is \[
  \begin{tikzpicture}
    \begin{knot}[clip width=7]
      \strand[thick] (0,0)
        to [out=down,in=down,looseness=2] (1,0);
      \strand[thick] (0,-2)
        to [out=up,in=up,looseness=2] (1,-2);
      \strand[thick] (0,-2)
        to [out=down,in=down,looseness=2] (1,-2);
      \strand[thick] (0,-4)
        to [out=up,in=up,looseness=2] (1,-4);
      \strand[thick] (2,0)
        to (2,-4);
    \end{knot}
  \end{tikzpicture}
\] 
Here we are using the usual product of tangles (see my webpage  
\href{http://math.ucr.edu/home/baez/tangles.html}{``tangles''}). 
Now the rule in Temperley--Lieb land
is that we can get rid of a loop if we multiply by the loop value \(d\);
that is, the loop ``equals'' \(d\). So \(e_1\) squared is just \(d\)
times \[
  \begin{tikzpicture}
    \begin{knot}[clip width=7]
      \strand[thick] (0,0)
        to [out=down,in=down,looseness=2] (1,0);
      \strand[thick] (0,-4)
        to [out=up,in=up,looseness=2] (1,-4);
      \strand[thick] (2,0)
        to (2,-4);
    \end{knot}
  \end{tikzpicture}
\] which --- since we are doing topology --- is the same as \(e_1\).
That's why \(e_i^2 = de_i\).

The other relations are even more obvious. For example, \(e_1 e_2 e_1\)
is just \[
  \begin{tikzpicture}
    \begin{knot}[clip width=7]
      \strand[thick] (0,0)
        to [out=down,in=down,looseness=2] (1,0);
      \strand[thick] (0,-4)
        to [out=up,in=up,looseness=2] (1,-4);
      \strand[thick] (2,0)
        to (2,-1.5)
        to [out=down,in=down,looseness=1.4] (1,-1.5)
        to [out=up,in=up,looseness=1.4] (0,-1.5)
        to (0,-2.5)
        to [out=down,in=down,looseness=1.4] (1,-2.5)
        to [out=up,in=up,looseness=1.4] (2,-2.5)
        to (2,-4);
    \end{knot}
  \end{tikzpicture}
\] which, since we are doing topology, is just \(e_1\)! Similarly,
\(e_2 e_1 e_2 = e_1\), and \(e_i\) and \(e_j\) commute if they are far
enough away to keep from running into each other.

As an exercise for combinatorists: figure out the dimension of \(TL_n\).

Okay, very cute, one might say, but so what? Well, this algebra was
actually first discovered in statistical mechanics, when Temperley and
Lieb were solving a \(2\)-dimensional problem:

\begin{enumerate}
\def\labelenumi{\arabic{enumi})}
\setcounter{enumi}{4}
\item
H. N. V. Temperley and E. H. Lieb, ``Relations between the `percolation' 
and `coloring' problem and other
graph-theoretical problems associated with regular planar lattices:
some exact results on the `percolation' problem'', \emph{Proc. Roy. Soc. Lond.}
\textbf{322} (1971), 251--280.
\end{enumerate}
\noindent
It gained a lot more fame when it appeared as the explanation for the
Jones polynomial invariant of knots --- although Jones had been using it
not for knot theory, but in the study of von Neumann algebras, and the
Jones polynomial was just an unexpected spinoff. Its importance in knot
theory comes from the fact that it is a quotient of the group algebra of
the braid group (as explained in \emph{Knots and Physics}). It has also
found a lot of other applications; for example, I've used it in my paper
on quantum gravity and the algebra of tangles. So it is nice to see that
there is also a formulation of the 4-color theorem in terms of the
Temperley--Lieb algebra (which I won't present here).

\begin{enumerate}
\def\labelenumi{\arabic{enumi})}
\setcounter{enumi}{5}
\item
  Louis Kauffman, ``Knots and physics'', \emph{Proc. Symp. Appl.
  Math.} \textbf{45} (1992), 131--246.

   Louis Kauffman, ``Spin networks, topology and discrete physics'',
  University of Illinois at Chicago preprint.

   Louis Kauffman,``Vassiliev invariants and the Jones polynomial'',
  University of Illinois at Chicago preprint.

  Louis Kauffman, ``Gauss codes and quantum groups'', University of
  Illinois at Chicago preprint.

   Louis Kauffman and H. Saleur,``Fermions and link invariants'',
  Yale University preprint YCTP-P21-91, July 5, 1991.

  Louis Kauffman, ``State models for link polynomials'',
  \emph{L'Enseignement Mathematique}, \textbf{36} (1990), 1--37.

  F. Jaeger, Louis Kauffman and H.
  Saleur, ``The Conway polynomial in \(\mathbb{R}^3\) and in thickened surfaces: a
  new determinant formulation'', preprint.
\end{enumerate}

These are a variety of papers on knots, physics and everything\ldots.
The more free-wheeling among you might enjoy the comments at the end of
the first paper on ``knot epistemology.''

I am going to a conference on gravity at U.C.\ Santa Barbara on Friday and
Saturday, which I why I am posting this early, and why I have no time to
describe the above papers. I'll talk about my usual obsessions, and hear
what other people are up to, perhaps bringing back some words of wisdom
for next week's ``This Week's Finds''.

\hypertarget{week9}{%
\section{March 12, 1993}\label{week9}}

\begin{enumerate}
\def\labelenumi{\arabic{enumi}.}
\tightlist
\item
  Boguslaw Broda,``Surgical invariants of four-manifolds'', preprint
  available as
  \href{https://arxiv.org/abs/hep-th/9302092}{\texttt{hep-th/9302092}}.
\end{enumerate}

There are a number of attempts underway to get invariants of
four-dimensional manifolds (and 4d topological quantum field theories)
by techniques analogous to those that worked in three dimensions. The
3-manifold invariants and 3d topological quantum field theories got
going with the work of Witten on Chern--Simons theory, but since this was
not rigorous a number of ways were devised to make it so. These seem
different at first glance but all give the same answer. Two approaches
that use a lot of category theory are the Heegard splitting approach
(due to Crane, Kohno and Kontsevich, in which one writes a 3-manifold as
two solid \(n\)-holed tori glued together by a diffeomorphism of their
boundaries), and the surgery on links approach (due to Reshetikhin and
Turaev, in which one builds up 3-manifolds by starting with the
3-sphere, cutting out thickened links and gluing them back in a
different way, allowing one to define invariants of 3-manifolds from
link invariants). In the case of 3 dimensions a nice paper relating the
Heegard splitting and the surgery on links approaches is

\begin{enumerate}
\def\labelenumi{\arabic{enumi})}
\setcounter{enumi}{1}
\item
 Sergey Piunikhin, ``Reshetikhin-Turaev and Crane-Kohno-Kontsevich 3-manifold invariants
coincide'', \emph{Journal of Knot Theory and its Ramifications} \textbf{2} (1993), 65--95.
\end{enumerate}

People are now trying to generalize all these ideas to 4-manifolds.
There is already an interesting bunch of 4-manifold invariants out
there, the Donaldson invariants, which are hard to compute, but were
shown (heuristically) by Witten to be related to a quantum field theory.
Lately people have been trying to define invariants using category
theory; these may or may not turn out to be the same.

I've already been trying to keep you all updated on the story about
Crane and Yetter's 4d TQFT. This week I'll discuss another approach,
with a vast amount of help from Daniel Ruberman, a topologist at
Brandeis. Any errors in what I write on this are likely to be due to my
misunderstandings of what he said --- caveat emptor! Broda's paper is
quite terse --- probably due to the race that is going on --- and is
based on:

\begin{enumerate}
\def\labelenumi{\arabic{enumi})}
\setcounter{enumi}{2}
\item
E. Cesar de Sa, ``A link calculus for 4-manifolds'', in \emph{Topology
of Low-Dimensional Manifolds, Proc. Second Sussex Conf.}, Springer Lecture
Notes in Mathematics \textbf{722}, Springer, Berlin, 1979, pp.~16--30.
\end{enumerate}
\noindent
so I should start by describing what little I understand of de Sa's
work.

One can describe (compact, smooth) 4-manifolds in terms of handlebody
decomposition. This allows one to actually draw pictures representing
4-manifolds. A lot of times when people first hear about topology they
get they impression that it's all about rubber doughnuts, M\"obius strips,
and other Dali-esque wiggly objects in hyperspace. Then, when they take
courses in it, they are confronted with nasty separation axioms and
cohomology theories! This is just to scare away outsiders! Handlebody
theory really \emph{is} about wiggly objects in hyperspace, and it's
lots of fun --- though to be good in it you need to know your point set
topology and your algebraic topology, I'm afraid --- and much better
than I do!

Recall:
\[\begin{aligned}D^n &= \mbox{unit ball in $\mathbb{R}^n$} \\ S^n &= \mbox{unit sphere in $\mathbb{R}^{n+1}$}\end{aligned}\]
In particular note that \(S^0\) is just two points. Note that:

\begin{itemize}
\tightlist
\item
  the boundary of \(D^4\) is \(S^3\)
\item
  the boundary of \(D^3 \times D^1\) is
  \(D^3 \times S^0 \cup S^2 \times D^1\)
\item
  the boundary of \(D^2 \times D^2\) is
  \(D^2 \times S^1 \cup S^1 \times D^2\)
\item
  the boundary of \(D^1 \times D^3\) is
  \(D^1 \times S^2 \cup S^0 \times D^3\)
\item
  the boundary of \(D^4\) is \(S^3\)
\end{itemize}

I have written this rather redundant chart in a way that makes the
pattern very clear and will come in handy below for those who aren't
used to this stuff.

To build up a 4-manifold we can start with a ``0-handle,'' \(D^4\),
which has as boundary \(S^3\).

Then we glue on ``1-handles,'' that is, copies of \(D^3 \times D^1\).
Note that part of the boundary of \(D^1 \times D^3\) is
\(D^3 \times S^0\), which is two \(D^3\)'s; when we glue on a 1-handle
we simply attach these two \(D^3\)'s to the \(S^3\) by a diffeomorphism.
The resulting space is not really a smooth manifold, but it can be
smoothed. It then becomes a smooth 4-manifold with boundary.

Then we glue on ``2-handles'' by attaching copies of \(D^2 \times D^2\)
along the part of their boundary that is \(D^2 \times S^1\). Then we
smooth things out.

Then we glue on ``3-handles'' by attaching copies of \(D^1 \times D^3\)
along the part of their boundary that is \(D^1 \times S^2\). Then we
smooth things out.

Then we glue on ``4-handles'' by attaching copies of \(D^4\) along their
boundary, i.e.~\(S^3\).

We can get any compact oriented 4-manifold this way using attaching maps
that are compatible with the orientations. The reader who is new to this
may enjoy constructing 2-manifolds in an analogous way. Compact oriented
2-manifolds with boundary are just \(n\)-holed tori.

What's cool is that with some tricks one can still \emph{draw} what's
going in the case of 3-manifolds and 4-manifolds. Here I'll just
describe how it goes for 4-manifolds, since that's what Cesar de Sa and
Broda are thinking about. By the way, a good introduction to this stuff
is

\begin{enumerate}
\def\labelenumi{\arabic{enumi})}
\setcounter{enumi}{3}
\item
\emph{The Topology of 4-manifolds}, by Robion C. Kirby, Springer
Lecture Notes in Mathematics \textbf{1374}, 1989.
\end{enumerate}
\noindent
So --- here is how we \emph{draw} what's going on. I apologize for being
somewhat sketchy here (sorry for the pun, too). I am a bit rushed since
I'm heading off somewhere else next weekend\ldots{} and I am not as
familiar with this stuff as I should be.

So, when we start with our 0-handle, or \(D^4\), we ``draw'' its
boundary, \(S^3\). Think of \(S^3\) as \(\mathbb{R}^3\) and a point at
infinity. Since we use perspective when drawing pictures of 3-d objects,
this boils down to pretending that our blackboard is a picture of
\(S^3\)!

As we add handles we continue to ``draw'' what's happening at the
boundary of the 4-manifold we have at each stage of the game. 1-handles
are attached by gluing a \(D^3 \times D^1\) onto the boundary along two
\(D^3\)'s --- or balls --- so we can just draw the two balls.

2-handles are attached by gluing a \(D^2 \times D^2\) onto the boundary
of the 4-manifold we have so far along a \(D^2 \times S^1\) --- or solid
torus, so we just need to figure out how to draw an embedded solid
torus. Well, for this we just need to draw a knot (that is, an embedded
circle), and write an integer next to it saying how many times the
embedded solid torus ``twists'' --- plus or minus depending on clockwise
or counterclockwise --- as we go around the circle. In other words, an
embedded solid torus is (up to diffeomorphism) essentially the same as a
framed knot. If we are attaching a bunch of 2-handles we need to draw a
framed link.

Things get a bit hairy in the case when one of the framed links goes
through one of the 1-handles that we've already added. It's easier to
draw this situation if we resort to another method of drawing the
1-handles. It's a bit more subtle, and took me quite a while to be able
to visualize (unfortunately I seem to have to visualize this stuff to
believe it). So let's go back to the situation where we have \(D^4\),
with \(S^3\) as its boundary, and we are adding 1-handles. Instead of
drawing two balls, we draw an unknotted circle with a dot on it! The dot
is just to distinguish this kind of circle from the framed links we
already have. But what the circle \emph{means} is this. The circle is
the boundary of an obvious \(D^2\), and we can push the interior of this
\(D^2\) (which is sitting in the \(S^3\)) into the interior of \(D^4\).
If we then remove a neighborhood of the \(D^2\), what we have left is
\(S^1 \times D^3\), which is just the result of adding a 1-handle to
\(D^4\).

This is probably easier to visualize one dimension down: if we have a
good old unit ball, \(D^3\), and slap an interval, or \(D^1\), onto its
boundary, and then push the interior of the interval into the interior
of the ball, and remove a neighborhood of the interval, what we have
left is just an \(S^1 \times D^2\).

So in short, we can draw all the 1-handles by drawing unlinked,
unknotted circles with dots on them, and then draw all the 2-handles by
drawing framed links that don't intersect these circles.

At this point, if you have never seen this before, you are probably
dreading the 3-handles and 4-handles. Luckily a theorem comes to our
rescue! If we start at the other end of our handlebody decomposition, as
it were, we start with 4-handles and glue on 3-handles. If you ponder
the chart and see what the pattern of what we're doing is, you'll see
that a single 4-handle with some 3-handles stuck on is just the same as
a 0-handle with some 1-handles stuck on. So when we now glue this thing
(or things) onto the stuff we've built out of 0-, 1-, and 2-handles, we
are doing so using a diffeomorphism of its boundary. But a theorem of
Laudenbach and Poenaru,

\begin{enumerate}
\def\labelenumi{\arabic{enumi})}
\setcounter{enumi}{4}
\item 
F. Laudenbach and V. Poenaru, 
``A note on \(4\)-dimensional handlebodies'', \emph{Bull. Math. Soc. France} \textbf{100} (1972),
337--344.
\end{enumerate}
\noindent
says that any such diffeomorphism extends to one of the interior. This
means that it doesn't make a darn bit of difference which diffeomorphism
we use to glue it on. In short, all the information is contained in the
1- and 2-handles, so we can \emph{draw} 4-manifolds by first drawing a
batch of unknotted unlinked circles with dots on them and then drawing a
framed link in the complement.

(A question for the experts, since I'm just learning this stuff: in
the above we seem to be assuming that there's only one 0-handle. Is this
an okay assumption or do we need something fancier if there's more?)

Now a given 4-manifold may have lots of different handlebody
decompositions. So, as usual, we would like to have a finite set of
``moves'' that allow us to get between any pair of handlebody
decompositions of the same 4-manifold. Then we can construct a
4-manifold invariant by cooking up a number from a handlebody
decomposition --- presented as a picture as above, if we want --- and
showing that it doesn't change under these ``moves''.

So, what de Sa did was precisely to find such a set of moves. (There,
that's what I understand of his work!)

And what Broda did was precisely to use the Kauffman bracket invariant
of framed links to cook up an invariant of 4-manifolds from the
handlebody decomposition --- which, note, involves lots of links. Recall
that the Kauffman bracket assigns to each link a polynomial in one
variable, \(q\). Here ``\(q\)'' is just the same \(q\) that appears in
the quantum group \(SU_q(2)\). As I mentioned in
\protect\hyperlink{week5}{``Week 5''}, this acts quite differently when
\(q\) is a root of unity, and the 3d topological quantum field theories
coming from quantum groups, as well as Crane and Yetter's 4d topological
quantum field theory, come from considering this root-of-unity case. So
it's no surprise that Broda requires \(q\) to be a root of unity.

Ruberman had some other remarks about Broda's invariant, but I think I
would prefer to wait until I understand them\ldots.

\begin{enumerate}
\def\labelenumi{\arabic{enumi})}
\setcounter{enumi}{5}
\item
  Abhay Ashtekar,  Ranjeet S.\ Tate and Claes Uggla, ``Minisuperspaces: symmetries and quantization'', available as
  \href{https://arxiv.org/abs/gr-qc/9302026}{\texttt{gr-qc/9302026}}.

  Abhay Ashtekar, Ranjeet S.\ Tate and Claes Uggla, ``Minisuperspaces: observables and quantization'', available as
  \href{https://arxiv.org/abs/gr-qc/9302027}{\texttt{gr-qc/9302027}}.
\end{enumerate}

I was just at the Pacific Coast Gravity Meeting last weekend and heard
Ranjeet Tate talk on this work. Recall first of all that minisuperspaces
are finite-dimensional approximations to the phase space of general
relativity, and are used to get some insight into quantum gravity. I
went through an example in \protect\hyperlink{week6}{``Week 6''}. In
these papers, the authors quantize various ``Bianchi type''
minisuperspace models. The ``Bianchi type'' business comes from a
standard classification of homogeneous (but not necessarily isotropic)
cosmologies and having a lot of symmetry. It is based in part on
Bianchi's classification of 3-dimensional Lie algebras into nine types.
The second paper gives a pretty good review of this stuff before diving
into the quantization, and I should learn it!

The most exciting aspect of these papers, at least to the dilettante
such as myself, is that one can quantize these models and show that
quantization does NOT typically remove the singularities (``Big Bang''
and/or ``big crunch''). Of course, these models have only finitely many
degrees of freedom, and are only a caricature of full-fledged quantum
gravity, so one can still argue that \emph{real} quantum gravity will
get rid of the singularities. But a number of general relativists are
arguing that this is not the case, and we simply have to learn to live
with singularities. So it's good to look at models, however simple,
where one can work things out in detail, and not just argue about
generalities.

\begin{enumerate}
\def\labelenumi{\arabic{enumi})}
\setcounter{enumi}{6}
\tightlist
\item
  Alan D.\ Rendall, ``Unique determination of an inner product by adjointness relations in
  the algebra of quantum observables'',
  Max-Planck-Institut f\"ur Astrophysik preprint.
\end{enumerate}

I had known Rendall from his work on the perturbative expansion of the
time evolution operators in classical general relativity. He became
interested in quantum gravity a while ago and visited Ashtekar and
Smolin at Syracuse University, since (as he said) the best way to learn
is by doing. There he wrote this paper on Ashtekar's approach to finding
the right inner product for the space of states of quantum gravity. I
had heard about this paper, but hadn't seen it until I met Rendall at
the gravity meeting last weekend. He gave me a copy and explained it. It
is a simple and beautiful paper --- such nice mathematical results that
I am afraid someone else may have found them earlier somewhere.

Ashtekar's idea is to fix the inner product by requiring that the
physical observables, which are operators on the space of states, be
self-adjoint. Rendall shows the following. Let \(A\) be a *-algebra
acting on a vector space \(V\). Let us say that an inner product on
\(V\) is ``strongly admissable'' if 1) the representation is a
*-representation with respect to this inner product, 2) for each element
of \(A\), the corresponding linear transformation on \(V\) is bounded
relative to the norm given by this inner product, and 3) the completion
of \(V\) in the inner product is a topologically irreducible
representation of \(A\). Rendall shows the uniqueness of a strongly
admissable inner product on any representation \(V\) of \(A\) (up to a
constant multiple). Of course, such an inner product need not exist, but
when it does, it is unique. This is as nice a result along these lines
as one could hope for. He also has a more complicated result that
applies to unbounded operators. A good piece of work on the foundations
of quantum theory!

\begin{enumerate}
\def\labelenumi{\arabic{enumi})}
\setcounter{enumi}{7}
\tightlist
\item
  Arlen Anderson, ``Thawing the frozen formalism: the difference between observables and
  what we observe'', preprint available
  as \href{https://arxiv.org/abs/gr-qc/9211028}{\texttt{gr-qc/9211028}}.
\end{enumerate}

There were a number of youngish folks giving talks at the gravity
meeting who have clearly been keeping up with the recent work on the
problem of time and other conceptual problems in quantum gravity. In
very brief terms, the problem of time is that in general relativity, we
have not a Hamiltonian in the traditional sense, but a ``Hamiltonian
constraint'' \(H = 0\), so when we quantize it superficially appears
that there are no dynamics whatsoever (as it seems like we have a zero
Hamiltonian!). That's the reason for the term ``frozen formalism'' ---
and the desire to ``thaw'' it, or find the dynamics lurking in it. In
fact, the Hamiltonian constraint is just a reflection of the fact that
general relativity has no preferred time coordinate, and we are just
learning how to deal with the quantum theory of such systems. For a good
survey of the problem and some new proposed solutions, I again refer
everyone to Isham's paper:

\begin{enumerate}
\def\labelenumi{\arabic{enumi})}
\setcounter{enumi}{8}
\tightlist
\item
  Chris J. Isham, ``Canonical quantum gravity and the problem of time'', 
  125 pages, available as
  \href{https://arxiv.org/abs/gr-qc/9210011}{\texttt{gr-qc/9210011}}.
\end{enumerate}

In particular, one interesting approach is due to Rovelli, and is called
``evolving constants of motion'' (a deliberate and very accurate
oxymoron). While there are serious technical problems with this
approach, it's very natural from a physical point of view --- at least
once you get used to it. I have the feeling that the younger physicists
are, as usual, getting used to it a lot more quickly than the older
folks who have been pondering the problem of time for many years.
Anderson is one of these younger folks, and his paper develops Rovelli's
approach in a toy model, namely the case of two free
particles satisfying the Schr\"odinger equation.

\begin{enumerate}
\def\labelenumi{\arabic{enumi})}
\setcounter{enumi}{9}
\tightlist
\item
  Cayetano Di Bartolo, Rodolfo Gambini and
  Jorge Griego, ``The extended loop group: an infinite dimensional manifold associated
  with the loop space'', 42 pages, available as
  \href{https://arxiv.org/abs/gr-qc/9303010}{\texttt{gr-qc/9303010}}.
\end{enumerate}

Unfortunately I don't have the time now to give this paper the
discussion it deserves. Gambini is one of the original inventors of the
loop representation of gauge theories, so his work is especially worth
paying attention to. He explained the idea of this paper to me a while
back. Its aim is to provide a workable ``calculus'' for the loop
representation by enlarging the ordinary loop group to a larger group
which is actually an infinite-dimensional Lie group --- the point being
that the usual loop group doesn't have a Lie algebra, but this one does.
As one might expect, the Lie algebra of this group is closely related to
the theory of Vassiliev invariants. The paper considers some
applications to quantum gravity and knot theory.

\hypertarget{week10}{%
\section{March 20, 1993}\label{week10}}

The most substantial part of this issue is some remarks by Daniel
Ruberman on the paper I was talking about last time by Boguslaw Broda. They
apparently show that Broda's invariant is not as new as it might appear.
But they're rather technical, so I'll put them near the end, and start
off with something on the light side, and then note some interesting
progress on the Vassiliev invariant scene.

\begin{enumerate}
\def\labelenumi{\arabic{enumi})}
\tightlist
\item
  Marcia Bartusiak, ``Beyond Einstein --- is space loopy?'',
  \emph{Discover}, April 1993.
\end{enumerate}

In the airport in Montreal I ran into this article, which was the cover
story, with an upside-down picture of Einstein worked into a bunch of
linked key-rings. I bought it --- how could I resist? --- since it is
perhaps the most ``pop'' exposition of the loop representation of
quantum gravity so far. Those interested in the popularization of modern
physics might want to compare

\begin{enumerate}
\def\labelenumi{\arabic{enumi})}
\tightlist
\setcounter{enumi}{1}
\item
John Horgan, ``Gravity quantized? A radical theory of gravity weaves space from tiny
loops'', \emph{Scientific American}, September 1992.
\end{enumerate}

Given the incredible hype concerning superstring theory, which seems to
have faded out by now, I sort of dread the same thing happening to the
loop representation of quantum gravity. It is intrinsically less
hype-able, since it does not purport to be a theory of everything, and it
comes right after superstrings were supposed to have solved all the
mysteries of the universe. Also, its proponents are (so far) a more
cautious breed than the string theorists --- note the question marks in
both titles! But we will see\ldots.

Marcia Bartusiak is a contributing editor of Discover and the author of
a book on current topics in astronomy and astrophysics, \emph{Thursday's
Universe}, which I haven't read. She'll be coming out with a book in
June, \emph{Through a Universe Darkly}, that's supposed to be about how
theories of cosmology have changed down through the ages. She does a
decent job of sketching vaguely the outlines of the loop representation
to an audience who must be presumed ignorant of quantum theory and
general relativity. Of course, there is also a certain amount of
human-interest stuff, with Ashtekar, Rovelli and Smolin (quite rightly)
coming off as the heroes of the story. There are, as usual, little boxes
with gee-whiz remarks like

\begin{verbatim}
            WITH REAMS OF PAPER
                         SPREAD OUT
               OVER THE KITCHEN TABLE
                 THEY FOUND
               SOLUTION AFTER SOLUTION
               FOR EQUATIONS
                     THOUGHT IMPOSSIBLE TO SOLVE
\end{verbatim}
(which is, after all, true --- nobody had previously found solutions to
the constraint equations in canonical quantum gravity, and all of a
sudden here were lots of 'em!). And there are some amusing discussions
of personality: ``Affable, creative, and easy-going, Rovelli quickly
settled into the role of go-between, helping mesh the analytic powers of
the quiet, contemplative Ashtekar with the creativity of the brash,
impetuous Smolin.'' And discussions of how much messier Smolin's office
is than Ashtekar's.

In any event, it's a fun read, and I recommend it. Of course, I'm
biased, so don't trust me.

\begin{enumerate}
\def\labelenumi{\arabic{enumi})}
\setcounter{enumi}{2}
\item
   Sergey Piunikhin,``Vassiliev invariants contain more information than all knot
  polynomials'', preprint. 

  Sergey Piunikhin, ``Turaev--Viro and Kauffman-Lins invariants for 3-manifolds coincide'',
  \emph{Journal of Knot Theory and its
  Ramifications} \textbf{1} (1992), 105--135.

  Sergey Piunikhin, ``Different presentations of 3-manifold invariants arising in rational
  conformal field theory'', preprint.

   Sergey Piunikhin, ``Weights of Feynman diagrams, link polynomials and Vassiliev knot
  invariants'', preprint.

  Sergey Piunikhin ``Reshetikhin--Turaev and Crane--Kohno--Kontsevich 3-manifold invariants
  coincide'', preprint.
\end{enumerate}

I received a packet of papers by Piunikhin a while ago. The most new and
interesting thing is the first paper listed above. In
\protect\hyperlink{week3}{``Week 3''} I noted a conjecture of Bar-Natan
that all Vassiliev invariants come from quantum group knot invariants
(or in other words, from Lie algebra representations.) Piunikhin claims
to refute this by showing that there is a Vassiliev invariant of degree
6 that does not. (However, other people have told me his claim is
misleading!) I have been too busy to read this paper yet.

\begin{enumerate}
\def\labelenumi{\arabic{enumi})}
\setcounter{enumi}{3}
\tightlist
\item
  Bernd Br\"ugmann, Bibliography of publications related to classical and quantum
  gravity in terms of the Ashtekar variables, available as
  \href{https://arxiv.org/abs/gr-qc/9303015}{\texttt{gr-qc/9303015}}.
\end{enumerate}

Let me just quote the abstract; this should be a handy thing:

\begin{quote}
This bibliography attempts to give a comprehensive overview of all the
literature related to the Ashtekar variables. The original version was
compiled by Peter Huebner in 1989, and it has been subsequently updated
by Gabriela Gonzalez and Bernd Br\"ugmann. Information about additional
literature, new preprints, and especially corrections are always
welcome.
\end{quote}

\begin{enumerate}
\def\labelenumi{\arabic{enumi})}
\setcounter{enumi}{4}
\tightlist
\item
  Boguslaw Broda,``Surgical invariants of four-manifolds'', preprint
  available as
  \href{https://arxiv.org/ps/hep-th/9302092}{\texttt{hep-th/9302092}}.
  (Revisited --- see \protect\hyperlink{week9}{``Week 9''})
\end{enumerate}

Let me briefly recall the setup: we describe a compact 4-manifold by a
handlebody decomposition, and represent this decomposition using a link
in \(S^3\). The 2-handles are represented by framed knots, while the
1-handles are represented by copies of the unknot (which we may think of
as having the zero framing). The 1-handles and 2-handles play quite a
different role in constructing the 4-manifold --- which is why one
normally draws the former as copies of the unknot with a \emph{dot} on
them --- but Broda's construction does NOT care about this. Broda simply
takes the link, forgetting the dots, and cooks up a number from it,
using cabling and the Kauffman bracket at an root of unity. Let's call
Broda's invariant by \(b(M)\) --- actually for each primitive \(r\)th
root of unity, we have \(b_r(M)\).

Broda shows that this is a 4-manifold invariant by showing it doesn't
change under the de Sa moves. One of these consists of adding or
deleting a Hopf link \[
  \begin{tikzpicture}
    \begin{knot}[clip width=7pt]
      \strand[thick] (0,0)
        to [out=up,in=up,looseness=2] (1,0)
        to [out=down,in=down,looseness=2] (0,0);
      \strand[thick] (0.5,0)
        to [out=up,in=up,looseness=2] (1.5,0)
        to [out=down,in=down,looseness=2] (0.5,0);
      \flipcrossings{2}
    \end{knot}
  \end{tikzpicture}
\] in which both components have the zero framing and one represents a
1-handle and the other a 2-handle. This move depends on the fact that we
can ``cancel'' a 1-handle and 2-handle pair, a special case of a general
result in Morse theory.

But since Broda's invariant doesn't care which circles represent
1-handles and which represent 2-handles, Broda's invariant is also
invariant under adding two 2-handles that go this way. This amounts to
taking a connected sum with \(S^2 \times S^2\). I.e.,
\(b(M) = b(M\# S^2 \times S^2)\).

Now, Ruberman told me a while back that we must also have
\(b(M) = b(M\#\mathbb{CP}^2\#-\mathbb{CP}^2)\), that is, the invariant
doesn't change under taking a connected sum with a copy of
\(\mathbb{CP}^2\) (complex projective 2-space) and an
orientation-reversed copy of \(\mathbb{CP}^2\). This amounts to adding
or deleting a Hopf link in which one component has the zero framing and
the other has framing 1. I didn't understand this, so I pestered
Ruberman some more, and this is what he says (modulo minor edits). I
have not had time to digest it yet:

\begin{quote}
The first question you asked was about the different framings on a
2-handle which goes geometrically once over a 1-handle, i.e.~makes a
Hopf link in which one of the circles is special (i.e is really a
1-handle, i.e.~in Akbulut--Kirby's notation is drawn with a dot.) The
answer is that the framing doesn't matter, since the handles cancel.
This is explained well (in the PL case) in Rourke-Sanderson's book.
(Milnor's book on the h-cobordism theorem explains it in terms of Morse
functions, in the smooth case.)

From this, it follows that
\(b(M) = b(M\#S^2 \times S^2) = b(M\#\mathbb{CP}^2\#-\mathbb{CP}^2)\).
For \(M\) is unchanged if you add a cancelling 1,2 pair, independent of
the framing on the 2-handle. If you change the special circle to an
ordinary one, \(b(M)\) doesn't change. On the other hand, \(M\) has been
replaced by its sum with either \(S^2 \times S^2\) or
\(\mathbb{CP}^2 \# -\mathbb{CP}^2\), depending on whether the framing on
the 2-handle is even or odd. (Exercise: why is only the parity
relevant?)

Now as I pointed out before, if one replaces all of the 1-handles
(special circles) of a 4-manifold with 2-handles, the invariant doesn't
change. This operation corresponds to doing surgery on the 4-manifold,
along the cores of the 1-handles. In particular, the manifold has
changed by a cobordism. (This is a basic construction; when you do
surgery you produce a cobordism, in this case it's \(M \times I\) with
2-handles attached to it along the circles which you surgered.)

From this, I will now show that Broda's invariant is determined by the
signature. (This is in the orientable case. Actually it seems that his
invariant is really an invariant of an oriented manifold.) The argument
above says that for any \(M\), there is an \(M'\), with
\(b(M) = b(M')\), where \(M'\) has no 1-handles, and where \(M\) and
\(M'\) are cobordant. In particular, \(M'\) is simply connected. So it
suffices to show that \(b(N) = b(N')\) if \(N\) and \(N'\) are simply
connected.

So now you can assume you have two simply connected manifolds \(N,N'\)
which are cobordant via a \(5\)-dimensional cobordism W, which you can
also assume simply connnected. By high-dimensional handlebody theory,
you can get rid of the 1-handles and 4-handles of \(W\), and assume that
all the 2-handles are added, then all of the 3-handles. If you add all
the 2-handles to N, you get
\(N\;\#\;k(S^2 \times S^2)\;\#\;l(\mathbb{CP}^2\#(-\mathbb{CP}^2))\) for some
\(k\) and \(l\). (Here is where simple connectivity is relevant; the
attaching circle of a 2-handle is null-homotopic, and therefore isotopic
to an unknotted circle. It's a simple exercise to see what happens when
you do surgery on a trivial circle, ie you add on \(S^2 \times S^2\) or
\(\mathbb{CP}^2 \# -\mathbb{CP}^2\). On the other hand you get the same
manifold as the result of adding 2-handles to \(N'\). So
\[N\;\#\;k(S^2 \times S^2)\;\#\;l(\mathbb{CP}^2\#(-\mathbb{CP}^2)) = 
N'\;\#\;k'(S^2 x S^2)\;\# \;l'(\mathbb{CP}^2\#(-\mathbb{CP}^2)),\]
so by previous remarks \(b(N) = b(N')\), i.e \(b\) is a cobordism
invariant.

Now: \(b\) is also multiplicative under connected sum, because connected
sum just takes the union of the link diagrams. The cobordism group is
\(\mathbb{Z}\), detected by the signature, so \(b\) must be a multiple
of the signature, modulo some number. (Maybe at this point I realize
\(b\) should be \(b_r\) or some such). If you compute (as a grad student
Tian-jin Li did for me) \(b_r(\mathbb{CP}^2)\), you find that \(b_r\)
lives in the group of \(r\)th (or maybe \(4r\)th; I'm at home and don't
have my note) roots of unity.

My conclusion: this invariant is a rather complicated way to compute the
signature of a 4-manifold (modulo \(r\) or \(4r\)) from a link diagram
of the manifold.

There is an important moral of the story, which is perhaps not obvious
to someone outside of 4-manifolds. Any invariant which purports to go
beyond classical ones (i.e.\ invariants of the intersection form) must
treat \(\mathbb{CP}^2\) and \(-\mathbb{CP}^2\) very differently. It
seems to be the case that many manifolds which are different (i.e.\
nondiffeomorphic) become diffeomorphic after you add on
\(\mathbb{CP}^2\). Thus any useful invariant should get rather
obliterated by adding \(\mathbb{CP}^2\). On the other hand,
non-diffeomorphic manifolds seem to stay non-diffeomorphic, no matter
how many \(-\mathbb{CP}^2\)'s you add on. This phenomenon doesn't seem
to be exhibited by any of the quantum-group type constructions for
3-manifolds; as it shouldn't, since (from the 3-manifold point of view)
an unknot with framing +1 or -1 doesn't change the 3-manifold. So if
you're looking for a combinatorial invariant, it seems critical that you
try to build in the asymmetry wrt orientation which 4-manifolds seem to
possess.

Exercise: do the nonorientable case. The answer should be that \(b\) is
determined by the Euler characteristic, mod 2.
\end{quote}

\hypertarget{week11}{%
\section{March 23, 1993}\label{week11}}

I'm hitting the road again tomorrow and will be going to the Quantum
Topology conference in Kansas until Sunday, so I thought I'd post this
week's finds early. As a result they'll be pretty brief. Let me start
with one that I mentioned in \protect\hyperlink{week9}{``Week 9''} but
is now easier to get:

\begin{enumerate}
\def\labelenumi{\arabic{enumi})}
\tightlist
\item
  Alan D.\ Rendall, ``Unique determination of an inner product by adjointness relations in
  the algebra of quantum observables'',
  available as
  \href{https://arxiv.org/abs/gr-qc/9303026}{\texttt{gr-qc/9303026}}.
\end{enumerate}
\noindent
and then mention another thing I've gotten as a spinoff from the gravity
conference at UCSB:

\begin{enumerate}
\def\labelenumi{\arabic{enumi})}
\setcounter{enumi}{1}
\tightlist
\item
  Ranjeet S.\ Tate, \emph{An Algebraic Approach to the Quantization of Constrained Systems:
  Finite Dimensional Examples}, Ph.D.\ thesis, Physics Department,
  Syracuse University, August 1992,
  \texttt{SU-GP-92/8-1}.
\end{enumerate}

Both the technical problems of ``canonical'' quantum gravity and one of
the main conceptual problems --- the problem of time --- stem from the
fact that general relativity is a system in which the initial data have
constraints. So improving our understanding of quantizing constrained
classical systems is important in understanding quantum gravity.

Let me say a few words about these constraints and what I mean by
``canonical'' quantum gravity.

First consider the wave equation in 2 dimensions. This is an equation
for a function from \(\mathbb{R}^2\) to \(\mathbb{R}\), say
\(\varphi(t,x)\), where \(t\) is a timelike and \(x\) is a spacelike
coordinate. The equation is simply
\[\frac{d^2\varphi}{dt^2} -\frac{d^2\varphi}{dx^2} = 0.\] Now this
equation can be rewritten as an evolutionary equation for initial data
as follows. We consider pairs of functions \((Q,P)\) on \(\mathbb{R}\)
--- which we think of as the functions \(\varphi\) and \(d\varphi/dt\) on ``space'', that
is, on a surface \(t = \text{constant}\). And we rewrite the
second-order equation above as a first-order equation:
\[\frac{d}{dt}(Q,P) = \left(P,\frac{d^2Q}{dx^2}\right).\tag{1}\] This is
a standard trick. We call the space of pairs \((Q,P)\) the ``phase
space'' of the theory. In canonical quantization, we treat this a lot
like the space \(\mathbb{R}^2\) of pairs \((q,p)\) describing the
initial position and momentum of a particle. Note that for a harmonic
oscillator we have an equation a whole lot like (1):
\[\frac{d}{dt}(q,p) = (p,-q).\] This is why when we quantize the wave
equation it's a whole lot like the harmonic oscillator.

Now in general relativity things are similar but more complicated. The
analog of the pairs \((\varphi, d\varphi/dt)\) are pairs \((Q,P)\) where
\(Q\) is the metric on spacetime restricted to a spacelike hypersurface
--- that is, the ``metric on space at a given time'' --- and \(P\) is
concocted from the extrinsic curvature of that hypersurface as it sits
in spacetime. Now the name of the game is to turn Einstein's equation
for the metric into a first-order equation sort of like (1). The problem
is, in general relativity there is no god-given notion of time. So we
need to \emph{pick} a ``lapse function'' on our hypersurface, and a
``shift vector field'' on our hypersurface, which say how we want to
push our hypersurface forwards in time. The lapse function says at each
point how much we push it in the normal direction, while the shift
vector field says at each point how much we push it in some tangential
direction. These are utterly arbitrary and give us complete flexibility
in how we want to push the hypersurface forwards. Even if spacetime was
flat, we could push the hypersurface forwards in a dull way like: \[
  \begin{tikzpicture}
    \draw[thick] (0,0) to (3,0) node[label=right:{new}]{};
    \draw[thick] (0,-0.5) to (3,-0.5) node[label=right:{old}]{};
  \end{tikzpicture}
\] or in a screwy way like \[
  \begin{tikzpicture}
    \draw[thick] plot [smooth] coordinates {(0,0) (0.7,1) (1.5,0) (2.5,0.5) (2.8,0) (3,0)} node[label=right:{new}]{};
    \draw[thick] (0,-0.5) to (3,-0.5) node[label=right:{old}]{};
  \end{tikzpicture}
\] Of course, in general relativity spacetime is usually not flat, which
makes it ultimately impossible to decide what counts as a ``dull way''
and what counts as a ``screwy way,'' which is why we simply allow all
possible ways.

Anyway, having \emph{chosen} a lapse function and shift vector field, we
can rewrite Einstein's equations as an evolutionary equation. This is a
bit of a mess, and it's called the ADM (Arnowitt--Deser--Misner)
formalism. Schematically, it goes like
\[\frac{d}{dt}(Q,P) = (\text{stuff}_1,\text{stuff}_2).\tag{2}\] 
where both ``\(\text{stuff}_1\)'' and ``\(\text{stuff}_2\)'' depend on 
both \(Q\) and \(P\) in a pretty
complex way.

But there is a catch. While the evolutionary equations are equivalent to
6 of Einstein's equations (Einstein's equation for general relativity is
really 10 scalar equations packed into one tensor equation), there are 4
more of Einstein's equations which turn into \emph{constraints} on \(Q\)
and \(P\). 1 of these constraints is called the Hamiltonian constraint
and is closely related to the lapse function; the other 3 are called the
momentum or diffeomorphism constraints and are closely related to the
shift vector field.

For those of you who know Hamiltonian mechanics, the reason why the
Hamiltonian constraint is called what it is is that we can write it as
\[H(Q,P) = 0\] for some combination of \(Q\) and \(P\), and this
\(H(Q,P)\) acts a lot like a Hamiltonian for general relativity in that
we can rewrite (2) using the Poisson brackets on the ``phase space'' of
all \((Q,P)\) pairs as
\[\begin{aligned}\frac{d}{dt}Q &= \{P,H(Q,P)\} \\ \frac{d}{dt}P &= \{Q,H(Q,P)\}.\end{aligned}\]
The funny thing is that \(H\) is not zero on the space of all \((Q,P)\)
pairs, so the equations above are nontrivial, but it does vanish on the
submanifold of pairs satisfying the constraints, so that, in a sense,
``the Hamiltonian of general relativity is zero''. But one must be
careful in saying this because it can be confusing! It has confused lots
of people worrying about the problem of time in quantum gravity, where
they naively think ``What --- the Hamiltonian is zero? That means
there's no dynamics at all!''

The problem in quantizing general relativity in the ``canonical''
approach is largely figuring out what to do with the constraints. It was
Dirac who first seriously tackled such problems, but the constraints in
general relativity always seemed intractible (when quantizing) until
Ashtekar invented his ``new variables'' for quantum gravity, that all of
a sudden make the constraints look a lot simpler. Ashtekar also has
certain generalizations of Dirac's general approach to quantizing
systems with constraints, and part of what Tate (who was a student of
Ashtekar) is doing is to study a number of toy models to see how
Ashtekar's ideas work.

I should note that there are lots of other ways to handle problems with
constraints, like BRST quantization, that aren't mentioned here at all.

Well, I'm off to Kansas and I hope to return with a bunch of goodies and
some gossip about 4-manifold invariants, topological quantum field
theories and the like. Lee Smolin will be talking there too so I will
try to extract some information about quantum gravity from him.

\hypertarget{week12}{%
\section{April 10, 1993}\label{week12}}

I had a lot of fun at the ``Quantum Topology'' conference at Kansas
State University, in Manhattan (yes, that's right, Manhattan, Kansas,
the so-called ``Little Apple''), and then spent a week recovering. Now
I'm back, ready for the next quarter\ldots{}

The most novel idea I ran into at the conference was due to Oleg Viro,
who, ironically, is right here at U. C. Riverside. He spoke on work with
Turaev on generalizing the Alexander module (a classical knot invariant)
to get a similar sort of module from any \(3\)-dimensional topological
quantum field theory. A ``topological quantum field theory,'' or TQFT
for short, is (in the language of physics) basically just a generally
covariant quantum field theory, one that thinks all coordinate systems
are equally good, just as general relativity is a generally covariant
classical field theory. For a more precise definition of TQFTs (which
even mathematicians who know nothing of physics can probably follow),
see my article
\href{http://math.ucr.edu/home/baez/symmetries.html}{``symmetries''}. In
any event, I don't think Viro's work exists in printed form yet; I'll
let you all know when something appears.

The most lively talk was one by Louis Crane and David Yetter, the
organizers of the conference. As I noted a while back, they claimed to
have constructed a FOUR-dimensional TQFT based on some ideas of Ooguri,
who was working on \(4\)-dimensional quantum gravity. This would be very
exciting as long as it isn't ``trivial'' in some sense, because all the
TQFTs developed so far only work in \(3\)-dimensional spacetime. A
rigorous 4-dimensional TQFT might bring us within striking distance of a
theory of quantum gravity --- this is certainly Crane's goal. Ocneanu,
however, had fired off a note claiming to prove that the Crane--Yetter
TQFT was trivial, in the sense that the partition function of the field
theory for any compact oriented 4-manifold equalled 1! In a TQFT, the
partition function of the field theory on a compact manifold is a
invariant of the manifold, and if it equalled 1 for all manifolds, it
would be an extremely dull invariant, hence a rather trivial TQFT.

So, on popular demand, Crane and Yetter had a special talk at 8 pm in
which they described their TQFT and presented results of calculations
that showed the invariant did NOT equal 1 for all compact oriented
4-manifolds. So far they have only calculated it in some special cases:
\(S^4\), \(S^3 \times S^1\), and \(S^2 \times S^2\). Amusingly, Yetter
ran through the calculation in the simplest case, \(S^4\), in which the
invariant \emph{does} happen to equal 1. But he persuaded most of us (me
at least) that the invariant really is an invariant and that he can
calculate it. I say ``persuade'' rather than prove because he didn't
present a proof that it's an invariant; the current proof is grungy and
computational, but Viro and Kauffman (who were there) pointed out some
ways that it could be made more slick, so we should see a comprehensible
proof one of these days. However, it's still up in the air whether this
invariant might be ``trivial'' in some more sophisticated sense, e.g.,
maybe it's a function of well-known invariants like the signature and
Euler number. Unfortunately, Ocneanu decided at the last minute not to
attend. Nor did Broda (inventor of another 4-manifold invariant that
Ruberman seems to have shown ``trivial'' in previous This Week's Finds)
show up, though he had been going to.

On a slightly more technical note, Crane and Yetter's TQFT depends on
chopping up the 4-manifold into simplices (roughly speaking,
4-dimensional versions of tetrahedra). Their calculation involves
drawing projections of these beasts into the plane and applying various
rules; it was quite fun to watch Yetter do it on the blackboard. Turaev
and Viro had constructed such a ``simplicial'' TQFT in 3 dimensions, and
Ooguri had been working on simplicial quantum gravity. As I note below,
Lee Smolin has a new scheme for doing \(4\)-dimensional quantum gravity
using simplices. During the conference he was busy trying to figure out
the relation of his ideas to Crane and Yetter's.

Also while at the conference, I found a terrible error in
\protect\hyperlink{week10}{``Week 10''} in my description of

\begin{quote}
{\rm Sergey Piunikhin,``Vassiliev invariants contain more information than all knot
polynomials'', preprint.}
\end{quote}

I had said that Piunikhin had discovered a Vassiliev invariant that
could distinguish knots from their orientation-reversed versions. No!
The problem was a very hasty reading on my part, together with the
following typo in the paper, that tricked my eyes:

\begin{quote}
Above constructed Vassiliev knot invariant \(w\) of order six does knot
detect orientation of knots.
\end{quote}
\noindent
Ugh! Also, people at the conference said that Piunikhin's claim in this
paper to have found a Vassiliev invariant not coming from quantum group
knot polynomials is misleading. I don't understand that yet.

Here are some papers that have recently shown up\ldots{}

\begin{enumerate}
\def\labelenumi{\arabic{enumi})}
\tightlist
\item
  Karel Kuchar, ``Canonical quantum gravity'', available as
  \href{https://arxiv.org/abs/gr-qc/9304012}{\texttt{gr-qc/9304012}}.
\end{enumerate}
\noindent
Kuchar (pronounced Koo-kahsh, by the way) is one of the grand old men of
quantum gravity, one of the people who stuck with the subject for the
many years when it seemed absolutely hopeless, who now deserves some of
the credit for the field's current resurgence. He has always been very
interested in the problem of time, and for anyone who knows a little
general relativity and quantum field theory, this is a very readable
introduction to some of the key problems in canonical quantum gravity. I
should warn the naive reader, however, that Kuchar's views about the
problem of time expressed in this paper go strongly against those of
many other experts! It is a controversial problem.

Briefly, many people believe that physical observables in quantum
gravity should commute with the Hamiltonian constraint (cf.\
\protect\hyperlink{week11}{``Week 11''}); this means that they are
time-independent, or constants of motion, and this makes the dynamics of
quantum gravity hard to ferrett out. Kuchar calls such quantities
``perennials.'' But Rovelli has made a proposal for how to recover
dynamics from perennials, basically by considering 1-parameter families
\(A_t\) of perennials, ironically called ``evolving constants of
motion.'' Kuchar argues against this proposal on two grounds: first, he
does not think physical observables need to commute with the Hamiltonian
constraint, and second, he argues that there may be very few if any
perennials. The latter point is much more convincing to me than the
former, at least at the \emph{classical} level, where the presence of
enough perennials would be close to the complete integrability of
general relativity, which is most unlikely. But on the quantum level
things are likely to be quite different, and Smolin has recently been at
work attempting to construct perennials in quantum gravity (cf.\
\protect\hyperlink{week1}{``Week 1''}). As for Kuchar's former point,
that observables in quantum gravity need not be perennials, his
arguments seem rather weak. In any event, read and enjoy, but realize
that the subject is a tricky one!

\begin{enumerate}
\def\labelenumi{\arabic{enumi})}
\setcounter{enumi}{1}
\tightlist
\item
  John Fischer, ``2-categories and 2-knots'', preprint, last
  revised Feb.~6 1993. 
\end{enumerate}
\noindent
This is the easiest way to learn about the \(2\)-category of 2-knots.
Recall (from \protect\hyperlink{week1}{``Week 1''} and
\protect\hyperlink{week4}{``Week 4''}) that a 2-knot is a surface
embedded in \(\mathbb{R}^4\), which may visualized as a ``movie'' of
knots evolving in time. Fischer shows that the algebraic structure of
2-knots is captured by a braided monoidal \(2\)-category, and he
describes this \(2\)-category.

\begin{enumerate}
\def\labelenumi{\arabic{enumi})}
\setcounter{enumi}{2}
\tightlist
\item
  Mark Miller and Lee Smolin,``A new discretization of classical and quantum general relativity'',
   available as
  \href{https://arxiv.org/abs/gr-qc/9304005}{\texttt{gr-qc/9304005}}.
\end{enumerate}
\noindent
Here Smolin proposes a new simplicial approach to general relativity
(there is an older one known as the Regge calculus) which uses
Ashtekar's ``new variables,'' and works in terms of the
Capovilla--Dell--Jacobson version of the Lagrangian. Let me just quote the
abstract, I'm getting tired:

\begin{quote}

We propose a new discrete approximation to the Einstein equations, based
on the Capovilla--Dell--Jacobson form of the action for the Ashtekar
variables. This formulation is analogous to the Regge calculus in that
it results from the application of the exact equations to a restricted
class of geometries. Both a Lagrangian and Hamiltonian formulation are
proposed and we report partial results about the constraint algebra of
the Hamiltonian formulation. We find that in the limit that the
\(\mathrm{SO}(3)\) gauge symmetry of frame rotations is reduced to the
abelian \(\mathrm{U}(1)^3\), the discrete versions of the diffeomorphism
constraints commute with each other and with the Hamiltonian constraint.

\end{quote}

\begin{enumerate}
\def\labelenumi{\arabic{enumi})}
\setcounter{enumi}{3}
\tightlist
\item
  ``Higher algebraic structures and quantization'', by Dan Freed, available as
  \href{https://arxiv.org/abs/hep-th/9212115}{\texttt{hep-th/9212115}}.
\end{enumerate}
\noindent
This is about TQFTs and the high-powered algebra needed to do justice to
the ``ladder of field theories'' that one can obtain starting with a
\(d\)-dimensional TQFT --- gerbs, torsors, \(n\)-categories and other such
scary things. I am too beat to do this justice.

\hypertarget{week13}{%
\section{April 20, 1993}\label{week13}}

Well, folks, this'll be the last ``This Week's Finds'' for a while,
since I'm getting rather busy preparing for my conference on knots and
quantum gravity, and I have a paper that seems to be taking forever to
finish.

\begin{enumerate}
\def\labelenumi{\arabic{enumi})}
\tightlist
\item
  Anthony W. Knapp, \emph{Elliptic Curves}, Mathematical Notes,
  Princeton University Press, 1992.
\end{enumerate}

This is a shockingly user-friendly introduction to a subject that can
all too easily seem intimidating. I'm certainly no expert but maybe just
for that reason I should sketch a brief ``introduction to the
introduction'' that may lure some of you into studying this beautiful
subject.

What I will say will perhaps appeal to people who like complex analysis
or mathematical physics, but Knapp concentrates on the aspects related
to number theory. For other approaches one might try

\begin{enumerate}
\def\labelenumi{\arabic{enumi})}
\setcounter{enumi}{1}
\item
  Serge Lang, \emph{Elliptic Functions}, 2nd edition, Springer, Berlin
  1987.
\item
  Dale Husemoeller, \emph{Elliptic Curves}, Springer, Berlin, 1987.
\end{enumerate}

Okay, where to start? Well, how about this: the sine function is an
analytic function on the complex plane with the property that
\[\sin(z + 2\pi) = \sin z.\] It also satisfies a nice differential
equation \[(\sin' z)^2 = 1 -(\sin z)^2\] and for this reason, we could,
if we hadn't noticed the sine function otherwise, have run into it when
we tried to integrate
 \[\frac{1}{\sqrt{1 -u^2}}\] 
The differential equation above implies that the integral is nice to do by the
substitution \(u = \sin z\), and we get the answer \(\arcsin u\). If the
sine function --- or more generally, trig functions --- didn't exist
yet, we would have invented them when we tried to do integrals involving
square roots of quadratic polynomials.

Elliptic functions are a beautiful generalization of all of this stuff.
Say we wanted, just for the heck of it, an analytic function that was
periodic not just in one direction on the complex plane, like the sine
function, but in \emph{two} directions. For example, we might want some
function \(P(z)\) with \[P(z + 2\pi) = P(z)\] and also
\[P(z + 2\pi i) = P(z).\] This function would look just the same on each
\(2\pi\)-by-\(2\pi\) square: \[
  \begin{tikzpicture}
    \draw[->] (-0.4,0) to (4,0) node [label=below:{$\Re(z)$}]{};
    \draw[->] (0,-0.4) to (0,4) node[label=left:{$\Im(z)$}]{};
    \foreach \x in {1,2,3} {
      \draw[thick,dotted] (\x,-0.4) to (\x,-0.1);
      \draw[thick] (\x,-0.1) to (\x,3.5);
      \draw[thick,dotted] (\x,3.5) to (\x,4);
    }
    \foreach \y in {1,2,3} {
      \draw[thick,dotted] (-0.4,\y) to (-0.1,\y);
      \draw[thick] (-0.1,\y) to (3.5,\y);
      \draw[thick,dotted] (3.5,\y) to (4,\y);
    }
    \draw[<->] (0.1,-0.2) to node[label=below:{$2\pi$}]{} (0.9,-0.2);
    \draw[<->] (-0.2,0.1) to node[label=left:{$2\pi$}]{} (-0.2,0.9);
  \end{tikzpicture}
\] so if we wanted, we could think of it as being a function on the
torus formed by taking one of these squares and identifying its top side
with its bottom side, and its left side with its right side.

More generally --- while we're fantasizing about this wonderful
doubly-periodic function --- we could ask for one that was periodic in
any old two directions. That is, fixing two numbers \(\omega_1\) and
\(\omega_2\) that aren't just real-valued multiples of each other, we
could hope to find an analytic function on the complex plane with
\(\omega_1\) and \(\omega_2\) as periods:
\[\begin{aligned}P(z + \omega_1) &= P(z) \\ P(z + \omega_2) &= P(z).\end{aligned}\]
Then \(P(z)\) would be the same at all points on the ``lattice'' of
points \(n\omega_1 + m \omega_2\), which might look like the squares
above or might be like \[
  \begin{tikzpicture}
    \draw[->] (-0.4,0) to (4,0) node [label=below:{$\Re(z)$}]{};
    \draw[->] (0,-0.4) to (0,4) node[label=left:{$\Im(z)$}]{};
    \foreach \x in {-0.3,0.7,1.7,2.7} {
      \draw[thick,dotted] (\x,-0.5) to ({\x+0.1},-0.2);
      \draw[thick] ({\x+0.1},-0.2) to ({\x+1.4},3.7);
      \draw[thick,dotted] ({\x+1.4},3.7) to ({\x+1.5},4);
    }
    \foreach \y in {0.6} {
      \draw[thick,dotted] (-0.6,\y) to (-0.3,{\y-0.05});
      \draw[thick] (-0.3,{\y-0.05}) to (3.6,{\y-0.7});
      \draw[thick,dotted] (3.6,{\y-0.7}) to (3.9,{\y-0.75});
    }
    \foreach \y in {1.6} {
      \draw[thick,dotted] (-0.4,\y) to (-0.1,{\y-0.05});
      \draw[thick] (-0.1,{\y-0.05}) to (3.8,{\y-0.7});
      \draw[thick,dotted] (3.8,{\y-0.7}) to (4.1,{\y-0.75});
    }
    \foreach \y in {2.6} {
      \draw[thick,dotted] (-0.2,\y) to (0.1,{\y-0.05});
      \draw[thick] (0.1,{\y-0.05}) to (4,{\y-0.7});
      \draw[thick,dotted] (4,{\y-0.7}) to (4.3,{\y-0.75});
    }
    \foreach \y in {3.6} {
      \draw[thick,dotted] (0,\y) to (0.3,{\y-0.05});
      \draw[thick] (0.3,{\y-0.05}) to (4.2,{\y-0.7});
      \draw[thick,dotted] (4.2,{\y-0.7}) to (4.5,{\y-0.75});
    }
  \end{tikzpicture}
\] or some such thing.

Let's think about this nice function \(P(z)\) we are fantasizing about.
Alas, if it were analytic on the whole plane (no poles), it would be
bounded on each little parallelogram, and since it's doubly periodic, it
would be a bounded analytic function on the complex plane, hence
\textbf{constant} by Liouville's theorem. Well, so a constant function
has all the wonderful properties we want --- but that's too boring!

So let's allow it to have poles! But let's keep it as nice as possible,
so let's have the only poles occur at the lattice points
\[L = \{n \omega_1 + m \omega_2\}.\] And let's make the poles as nice as
possible. Can we have each pole be of order one? That is, can we make
\(P(z)\) blow up like \(1/(z -\omega)\) at each lattice point \(\omega\)
in \(L\)? No, because if it did, the integral of \(P\) around a nicely
chosen parallelogram around the pole would be zero, because the
contributions from opposite sides of the parallelogram would cancel by
symmetry. (A fun exercise.) But by the Cauchy residue formula this means
that the residue of the pole vanishes, so it can't be of order one.

Okay, try again. Let's try to make the pole at each lattice point be of
order two. How can we cook up such a function? We might try something
obvious: just sum up, for all \(\omega\) in the lattice \(L\), the
functions \[\frac{1}{(z -\omega)^2}\] We get something periodic with
poles like \(1/(z -\omega)^2\) at each lattice point \(\omega\). But
there's a big problem --- the sum doesn't converge! (Another fun
exercise.)

Oh well, try again. Let's act like physicists and \textbf{renormalize}
the sum by subtracting off an infinite constant! Just subtract the sum
over all \(\omega\) in \(L\) of \(1/\omega^2\). Well, all \(\omega\)
except zero, anyway. This turns out to work, but we really should be
careful about the order of summation here: really, we should let
\(P(z)\) be \(1/z^2\) plus the sum for all nonzero \(\omega\) in the
lattice \(L\) of \(1/(z -\omega)^2 -1/\omega^2\). This sum does converge
and the limit is a function \(P(z)\) that's analytic except for poles of
order two at the lattice points. This is none other than the Weierstrass
elliptic function, usually written with a fancy Gothic \(\mathfrak{P}\)
to intimidate people. Note that it really depends on the two periods
\(\omega_1\) and \(\omega_2\), not just \(z\).

Now, it turns out that \(P(z)\) really \emph{is} a cool generalization
of the sine function. Namely, it satisfies a differential equation like
the one the sine does, but fancier:
\[P'(z)^2 = 4 P(z)^3 - g_2 P(z) - g_3\] 
where \(g_2\) and \(g_3\)
are some constants that depend on the periods \(\omega_1\) and
\(\omega_2\). Just as with the sine function we can use the
\emph{inverse} of Weierstrass \(\mathfrak{P}\) function to do some
integrals, but this time we can do integrals involving square roots of
cubic polynomials! If you look in big nasty books of special functions
or tables of integrals, you will see that there's a big theory of this
kind of thing that was developed in the 1800's --- back when heavy-duty
calculus was hip.

There are, however, some other cool ways of thinking about what's going
on here. First of all, remember that we can think of \(P(z)\) as a
function on the torus. We can think of this torus as being
``coordinatized'' --- I use the word loosely --- by \(P(z)\) and its
first derivative \(P'(z)\). I.e., if we know \(x = P(z)\) and
\(y = P'(z)\) we can figure out where the point \(z\) is on the torus.
But of course \(x\) and \(y\) can't be any old thing; the differential
equation above says they have to satisfy \[y^2 = 4x^3 -g_2 x -g_3.\]
Here \(x\) and \(y\) are complex numbers of course. But look what this
means: it means that if we look at the pairs of complex numbers
\((x,y)\) satisfying the above cubic equation, we get something that
looks just like a torus! This is called an elliptic curve, since for
algebraic geometers a ``curve'' is the set of solutions \((x,y)\) of
some polynomial in two \emph{complex} variables --- not two real
variables.

So, an ``elliptic curve'' is basically just the solutions of a cubic
equation in two variables. Actually, we want to rule out curves that
have singularities, that is, places where there's no unique tangent line
to the curve, as in \(y^2 = x^3\) or \(y^2 = x^2(x+1)\) --- draw these
in the real plane and you'll see what I mean. Anyway, all elliptic
curves can, by change of variables, be made to look like our favorite
one, \[y^2 = 4x^3 -g_2 x -g_3.\] There are lots of more fancy ways of
thinking about elliptic curves, and one is to think of the fact that
they look like a torus as the key part. In a book on algebraic geometry
you might see an elliptic curve as a curve with genus one (i.e., with
one ``handle,'' like a torus has). One nice thing about a torus is that
is a group. That is, we know how to add complex numbers, and we can add
modulo elements of the lattice \(L\), so the torus becomes a group with
addition mod \(L\) as the group operation. This is simple enough, but
it means that when we look at the solutions of
\[y^2 = 4x^3 -g_2 x -g_3\] they must form a group somehow, and viewed
this way it's not at all obvious! Nonetheless, there is a beautiful
geometric description of the group operation in these terms --- I'll
leave this for Knapp to explain.

Let me wrap this up --- the story goes on and on, but I'm getting tired
--- with a bit about what it has to do with number theory. It has a lot
to do with Diophantine equations, where one wants integer, or rational
solutions to a polynomial equation. Suppose that \(g_2\) and \(g_3\) are
rational, and one has some solutions to the equation
\[y^2 = 4x^3 - g_2 x - g_3.\] 
Then it turns out that one can use the
group operation on the elliptic curve to get new solutions! Actually, it
seems as if Diophantus knew this way back when in some special cases.
For example, for the problem \[y(6 -y) = x^3 -x\] Diophantus could start
with the trivial solution \((x,y) = (-1,0)\), do some mysterious stuff,
and get the solution \((17/9,26/27)\). Knapp explains how this works in
the Overview section, but then more deeply later. Basically, it uses the
fact that this curve is an elliptic curve, and uses the group structure.

In fact, one can solve mighty hard-seeming Diophantine problems using
these ideas. Knapp talks a bit about a problem Fermat gave to Mersenne
in 1643 --- this increased my respect for Fermat a bit. He asked, find a
Pythagorean triple \((X,Y,Z)\), that is: 
\[X^2 + Y^2 = Z^2,\]
such that \(Z\) is a square number and \(X + Y\) is too! One can solve
this using elliptic curves. I don't know if Mersenne got it --- the
answer is at the end of this post, but heavy-duty number theorists out
there might enjoy trying this one if they don't know it already.

Some more stuff:

\begin{enumerate}
\def\labelenumi{\arabic{enumi})}
\setcounter{enumi}{3}
\tightlist
\item
  Jim Stasheff, ``Closed string field theory, strong homotopy Lie algebras and the
  operad actions of moduli spaces'', available
  as
  \href{https://arxiv.org/abs/hep-th/9304061}{\texttt{hep-th/9304061}}.
\end{enumerate}
\noindent
One conceptually pleasing approach to string theory is closed string
field theory, where one takes as the basic object unparametrized maps
from circle into a manifold \(M\) representing ``space'', i.e., elements
of \[\mathsf{Maps}(S^1,M)/\mathsf{Diff}^+(S^1).\] A state of closed
string field theory would be roughly a function on the above set. Then
one tries to define all sorts of operations on these states, in order to
write down ways the strings can interact. For example, there is a
``convolution product'' on these functions which almost defines a Lie
algebra structure. However, the Jacobi identity only holds ``up to
homotopy,'' so we have an algebraic structure called a homotopy Lie
algebra. Physicists would say that the Jacobi identity holds modulo a
BRST exact term. This is just the beginning of quite a big bunch of
mathematics being developed by Stasheff, Zwiebach, Getzler, Kapranov and
many others. My main complaint with the physics is that all these
structures seem to depend on choosing a Riemannian metric on \(M\) --- a
so-called ``background metric.'' Since string theory is supposed to
include a theory of quantum gravity it is annoying to have this
God-given background metric stuck in at the very start. Perhaps I just
don't understand this stuff. I am looking around for stuff on
background-independent closed string field theory, since I have lots of
reason to believe that it's related to the loop representation of
quantum gravity. Unfortunately, I scarcely know the subject --- I had
hoped Stasheff's work would help me, but it seems that this metric
always enters.

\begin{enumerate}
\def\labelenumi{\arabic{enumi})}
\setcounter{enumi}{4}
\tightlist
\item
  Sergey Matveev and Michael Polyak, ``A geometrical presentation of the surface mapping 
  class group and surgery'', preprint.
\end{enumerate}
\noindent
This paper shows how to express the mapping class group of a surface in
terms of tangles. This gives a nice relationship between two approaches
to 3d TQFTs (topological quantum field theories): the Heegard
decomposition approach, and the surgery on links approach.

\begin{enumerate}
\def\labelenumi{\arabic{enumi})}
\setcounter{enumi}{5}
\tightlist
\item
  Michael  Polyak, ``Invariants of 3-manifolds and conformal field theories'', preprint.
\end{enumerate}
\noindent
The main good thing about this paper in my opinion is that it simplifies
the definition of a modular tensor category. Recall that Moore and
Seiberg showed how any string theory (more precisely, any rational
conformal field theory) gave rise to a modular tensor category, and then
Crane showed that any modular tensor category gave rise to a 3d TQFT.
Unfortunately a modular tensor category seems initially to be a rather
baroque mathematical object. In this paper Polyak shows how to get lots
of the structure of a modular tensor category from just the ``fusion''
and ``braiding'' operators, subject to some mild conditions. I have a
conjecture that all nonnegative link invariants (in the sense of my
paper on tangles and quantum gravity) give rise to modular tensor
categories, and this simplifies things to the point where maybe I might
eventually be able to prove it. There are lots of nice pictures here,
too, by the way.

\begin{center}\rule{0.5\linewidth}{0.5pt}\end{center}

\textbf{Answer to puzzle:} \[X = 1061652293520\] \[Y = 4565486027761\]
\[Z = 4687298610289\]

\hypertarget{week14}{%
\section{May 8, 1993}\label{week14}}

Things are moving very fast in the quantum gravity/4d topology game, so
I feel I should break my vow not to continue this series until after
next weekend's conference on Knots and Quantum Gravity.

Maybe I should recall where things were when I left off. The physics
problem motivating a lot of work in theoretical physics today is
reconciling general relativity and quantum theory. The key feature of
general relativity is that time and space do not appear as a
``background structure,'' but rather are dynamical variables. In
mathematical terms, this just means that there is not a fixed metric;
instead gravity \emph{is} the metric, and the metric evolves with time
like any other physical field, satisfying some field equations called
the Einstein equations.

But it is worth stepping back from the mathematics and trying to put
into simple words why this makes general relativity so special. Of
course, it's very hard to put this sort of thing into words. But
roughly, we can say this: in Newtonian mechanics, there is a universal
notion of time, the ``\(t\)'' coordinate that appears in all the
equations of physics, and we assume that anyone with a decent watch will
be able to keep in synch with everyone else, so there is no confusion
about what this ``\(t\)'' is (apart from choosing when to call
\(t = 0\), which is a small sort of arbitrariness one has to live with).
In special relativity this is no longer true; watches moving relative to
each other will no longer stay in synch, so we need to pick an
``inertial frame,'' a notion of rest, in order to have a ``\(t\)''
coordinate to play with. Once we pick this inertial frame, we can write
the laws of physics as equations involving ``\(t\)''. This is not too
bad, because there is only a finite-parameter family of inertial frames,
and simple recipes to translate between them, and also because nothing
going on will screw up the functioning of our (idealized) clocks: that
is, the ``\(t\)'' coordinate doesn't give a damn about the \emph{state}
of the universe. That's what is meant by saying a ``background
structure'' --- it's some aspect of the universe that is unaffected by
everything else that's going on.

In general relativity, things get much more interesting: there is no
such thing as an inertial frame that defines coordinates on spacetime,
because there is no way you can get a lot of things at different places
to remain at rest with each other --- this is what is meant by saying
that spacetime is curved. You can measure time with your watch,
so-called ``proper time,'' but this applies only near you. More
interestingly still, to compare what your watch is doing to what someone
else's is doing, you actually need to know a lot about the state of the
universe, e.g., whether there are any heavy masses around that are
curving spacetime. The ``metric,'' whereby one measures distances and
proper time, depends on the state of the universe --- or more properly,
it is part of the state of the universe.

Trying to do \emph{quantum} theory in this context has always been too
hard for people. Part of the reason why is that built into the heart of
traditional quantum theory is the ``Hamiltonian,'' which describes the
evolution of the state of the system relative to a God-given
``background'' notion of ``\(t\)''. Anyone who has taken quantum
mechanics will know that the star of the show is the Schr\"odinger
equation: \[i\frac{d\psi}{dt} = H\psi\] saying how the wavefunction
\(\psi\) changes with time in a way depending on the Hamiltonian \(H\).
No ``\(t\),'' no ``\(H\)'' --- this is one basic problem with trying to
reconcile quantum theory with general relativity.

Actually, it turns out that the analog to Schr\"odinger's equation for
quantum gravity is the Wheeler--DeWitt equation. The Hamiltonian is
replaced by an operator called the ``Hamiltonian constraint'' and we
have \[H\psi = 0.\] Note how this cleverly avoids mentioning ``\(t\)''!
The problem is, people still aren't quite sure what to do with the
solutions to this equation --- we're so used to working with
Schr\"odinger's equation.

Now in 1988 Witten wrote a paper in which he coined the term
``topological quantum field theory,'' or TQFT, for short. This was meant
to capture in a rigorous way what field theories like quantum gravity
should be like. Actually, Witten was working on a different theory
called Donaldson theory, which also has the property of having no
background structures. Shortly thereafter the mathematician Atiyah came
up with a formal definition of a TQFT. To get an idea of this
definition, try my notes on
\href{http://math.ucr.edu/home/baez/symmetries.html}{``symmetries''} and
(if you don't know what categories are)
\href{http://math.ucr.edu/home/baez/categories.html}{``categories''}.
For a serious tour of TQFTs and the like, try his book:

\begin{enumerate}
\def\labelenumi{\arabic{enumi})}
\tightlist
\item
Michael Atiyah, \emph{The Geometry and Physics of Knots}, Cambridge
U. Press, Cambridge, 1990.
\end{enumerate}
\noindent
One can think of a TQFT as a framework in which a Wheeler--DeWitt-like
equation governs the dynamics of a quantum field theory. Experts may
snicker here, but it is true, if not as enlightening as other things one
can say.

I won't bother to define TQFTs here, but I think Smolin put it very well
when he said the idea of TQFTs really helped us break out of our
traditional idea of fields as being something defined at every point of
spacetime, wiggling around, and allowed us to see field theory from many
new angles. For example, TQFTs let us wiggle out of the old conundrum of
whether spacetime is continuous or discrete, because many TQFTs can be
\emph{equivalently} described in either of two ways: via a continuum
model of spacetime, or via a discrete one in which spacetime is given a
``simplicial structure,'' like a big tetrahedral tinkertoy lattice kind
of thing. The latter idea appears to be due to Turaev and Viro, although
certainly physicists have had similar ideas for years, going back to
Ponzano and Regge, who worked on simplicial quantum gravity.

Now the odd thing is that while interesting 3d TQFTs have been found,
the most notable being Chern--Simons theory, nobody has quite been able
to make 4d TQFTs rigorous. Witten's original work on Donaldson theory
has led to many interesting things, but not yet a full-fledged TQFT in
the rigorous sense of Atiyah. And quantum gravity still resists being
formulated as a TQFT.

A while back I noted that Crane and Yetter had invented a 4d TQFT using
the simplicial approach. There has been a lot of argument over whether
this TQFT is interesting or ``trivial.'' Of course, trivial is not a
precise concept. For a while Ocneanu claimed that the partition function
of every compact 4-manifold equalled 1 in this TQFT, which counts as
very trivial. But this appears not to be the case. Broda invented
another 4d TQFT and here on ``This Week's Finds'' Ruberman showed it was
trivial in the sense that the partition function of any compact
4-manifold was a function of the ``signature'' of the 4-manifold. This
is trivial because the signature is a well-understood invariant and if
we are trying to do something new and interesting that just isn't good
enough.

In the following paper:

\begin{enumerate}
\def\labelenumi{\arabic{enumi})}
\setcounter{enumi}{1}
\tightlist
\item
Justin Roberts, ``Skein theory and Turaev--Viro invariants'',
\emph{Topology}, \textbf{34}, 771--787.  Available as
\href{https://citeseerx.ist.psu.edu/viewdoc/summary?doi=10.1.1.138.8587}{https://citeseerx.ist.psu.edu/viewdoc/summary?doi=10.1.1.138.8587}.
\end{enumerate}
\noindent
Roberts \emph{almost} claims to show that the Crane--Yetter invariant is
trivial in the same sense, namely that the partition function of any
compact 4-manifold is an exponential of the signature. Now if Crane and
Yetter's own computations are correct, this cannot be the case, but it
\emph{could} be an exponential of a linear combination of the signature
and the Euler characteristic, as far as I know. The catch is that
Roberts does not normalize his version of the Crane--Yetter invariant in
the same way that Crane and Yetter do, so it is hard to compare results.
But Roberts says: ``The normalisations here do not agree with those in
Crane and Yetter, and I have not checked the relationship. However, when
dealing with the {[}3d TQFT{]} invariants, different normalisations of
the initial data change the invariants by factors depending on standard
topological invariants (for example Betti numbers), so there is every
reason to belive that these {[}4d TQFT{]} invariants are trivial (that
is, they differ from 1 only by standard invariant factors) in all
normalisations.''

This is a bit of a disappointment, because Crane at least had hoped that
their TQFT might actually turn out to \emph{be} quantum gravity. This
was not idle dreaming; it was because the Crane--Yetter construction was
a rigorous analog of some work by Ooguri on simplicial quantum gravity.

Then, about a week ago, Rovelli put a paper onto the net:

\begin{enumerate}
\def\labelenumi{\arabic{enumi})}
\setcounter{enumi}{2}
\tightlist
\item
  Carlo Rovelli, ``The basis of the Ponzano--Regge--Turaev--Viro--Ooguri model is the loop
  representation basis'',  available as
  \href{https://arxiv.org/abs/hep-th/9304164}{\texttt{hep-th/9304164}}.
\end{enumerate}

This is a remarkable paper that I have not been able to absorb yet.
First it goes over 3d quantum gravity --- which \emph{has} been made
into a rigorous TQFT. It works with the simplicial formulation of the
theory. That is, we consider our (3-dimensional) spacetime as being
chopped up into tetrahedra, and assign to each edge a length, which is
required to be \(0,\frac{1}{2},1,\frac{3}{2},\ldots\). This idea of quantized
edge-lengths goes back to 4d work of Ponzano and Regge, but recently
Ooguri showed that in 3d this assumption gives the same answers as
Witten's continuum approach to 3d quantum gravity. The ``half-integers''
\(0,\frac{1}{2},1,\frac{3}{2},\ldots\) should remind physicists of spin, which is
quantized in the same way, and mathematically this is exactly what is
going on: we are really labelling edges with representations of the
group \(\mathrm{SU}(2)\), that is, spins. What Rovelli shows is that if
one starts with the loop representation of 3d quantum gravity (yet
another approach), one can prove it equivalent to Ooguri's approach, and
what's more, using the loop representation one can \emph{calculate} the
lengths of edges of triangles in a given state of space (space here is a
\(2\)-dimensional triangulated surface) and \emph{show} that lengths are
quantized in units of the Planck length over 2. (Here the Planck length
\(L\) is the fundamental length scale in quantum gravity, about
\(1.6 \times 10^{-33}\) meters.)

And, most tantalizing of all, he sketches a generalization of the above
to 4d. In 4d it is known that in the loop representation of quantum
gravity it is areas of surfaces that are quantized in units of
\(L^2/2\), rather than lengths. Rovelli considers an approach where one
chops 4-dimensional spacetime up into simplices and assigns to each
2-dimensional face a half-integer area. He uses this to write down a
formula for the inner product in the Hilbert space of quantum gravity
--- thus, at least formally, answering the long-standing ``inner product
problem'' in quantum gravity. The problem is that, unlike in 3d quantum
gravity, here one must sum over ways of dividing spacetime into
simplices, so the formula for the inner product involves a sum that does
not obviously converge. This is however sort of what one might expect,
since in 4d quantum gravity, unlike 3d, there are ``local excitations''
--- local wigglings of the metric, if you will --- and this makes the
Hilbert space be infinite-dimensional, whereas the 3d TQFTs have
finite-dimensional Hilbert spaces.

I think I'll quote him here. It's a bit technical in patches, but worth
it\ldots{}

\begin{quote}
We conclude with a consideration on the formal structure of 4-d quantum
gravity, which is important to understand the above construction.
Standard quantum field theories, as QED and QCD, as well as their
generalizations like quantum field theories on curved spaces and
perturbative string theory, are defined on metric spaces. Witten's
introduction of the topological quantum field theories has shown that
one can construct quantum field theories defined on a manifold which has
only its differential structure, and no fixed metric structure. The
theories introduced by Witten and axiomatized by Atiyah have the
following peculiar feature: they have a finite number of degrees of
freedom, or, equivalently, their quantum mechanical Hilbert spaces are
finite dimensional; classically this follows from the fact that the
number of fields is equal to the number of gauge transformations.
However, not any diff-invariant field theory on a manifold has a finite
number of degrees of freedom. Witten's gravity in 3-d is given by the
action
\[S[A,E] = \int F\wedge E \]
which has finite number of degrees of freedom. Consider the action
\[S[A,E] = \int F\wedge e\wedge e \]
in 3+1 dimensions, for a (self dual) \(\mathrm{SO}(3,1)\) connection
\(A\) and a (real) one form \(e\) with values in the vector
representation of \(\mathrm{SO}(3,1)\). This theory has a strong
resemblance with its 2+1 dimensional analog: it is still defined on a
differential manifold without any fixed metric structure, and is
diffeomorphism invariant. We expect that a consistent quantization of
such a theory should be found along lines which are more similar to the
quantization of the \(\int(F\wedge E)\) theory than to the quantization
of theories on flat space, based on the Wightman axioms namely on
\(n\)-points functions and related objects. Still, the theory
\(\int(F\wedge e\wedge e)\) has genuine field degrees of freedom: its
physical phase space is infinite dimensional, and we expect that its
Hilbert state space will also be infinite dimensional. There is a
popular belief that a theory defined on a differential manifold without
metric and diffeomorphism invariant has necessarily a finite number of
degrees of freedom (``because thanks to general covariance we can gauge
away any local excitation''). This belief is of course wrong. A theory
as the one defined by the action \(\int(F\wedge e\wedge e)\) is a theory
that shares many features with the topological theories, in particular,
no quantity defined ``in a specific point'' is gauge invariant; but at
the same time it has genuinely infinite degrees of freedom. Indeed, this
theory is of course nothing but (Ashtekar's form of) standard general
relativity.

The fact that ``local'' quantities like the \(n\)-point functions are
not appropriate to describe quantum gravity non-perturbatively has been
repeatedly noted in the literature. As a consequence, the issue of what
are the quantities in terms of which a quantum theory of gravity can be
constructed is a much debated issue. The above discussion indicates a
way to face the problem: The topological quantum field theories studied
by Witten and Atiyah provide a framework in terms of which quantum
gravity itself may be framed, in spite of the infinite degrees of
freedom. In particular, Atiyah's axiomatization of the topological field
theories provides us with a clean way of formulating the problem. Of
course, we have to relax the requirement that the theory has a finite
number of degrees of freedom. These considerations leads us to propose
that the correct general axiomatic scheme for a physical quantum theory
of gravity is simply Atiyah's set of axioms up to finite dimensionality
of the Hilbert state space. We denote a structure that satisfies all
Atiyah's axioms, except the finite dimensionality of the state space, as
a \textbf{generalized topological theory}.

The theory we have sketched is an example of such a generalized
topological theory. We associate to the connected components of the
boundary of \(M\) the infinite dimensional state space of the Loop
Representation of quantum gravity. Eq. 5 {[}the magic formula I alluded
to --- jb{]}, then, provides a map, in Atiyah's sense, between the state
spaces constructed on two of these boundary components. Equivalently, it
provides the definition of the Hilbert product in the state space.

One could argue that the framework we have described cannot be
consistent, because it cannot allow us to recover the ``broken phase of
gravity'' in which we have a nondegenerate background metric: in the
proposed framework one has only non-local observables on the boundaries,
while in the broken phase a \emph{local} field in \(M\) has
non-vanishing vacuum expectation value. We think that this argument is
weak because it disregards the diffeomorphism invariance of the theory:
in classical general relativity no experiment can distinguish a
Minkowskian spacetime metric from a non-Minkowkian flat metric. The two
are physically equivalent, as two gauge-related Maxwell potentials. For
the same reason, no experiment could detect the absolute \emph{position}
of, say, a gravitational wave, (while of course the position of an e.m.\
wave is observable in Maxwell theory). Physical locality in general
relativity is only defined as coincidence of some physical variable with
some other physical variable, while in non general relativistic physics
locality is defined with respect to a fixed metric structure. In
classical general relativity, there is no gauge-invariant obervable
which is local in the coordinates. Thus, any observation can be
described by means of the value of the fields on arbitrary boundaries of
spacetime. This is the correct consequence of the fact that ``thanks to
general covariance we can gauge away any local excitation'', and this is
the reason for which one can have the ADM ``frozen time'' formalism. The
spacetime manifold of general relativity is, in a sense, a much weaker
physical object than the spacetime metric manifold of ordinary theories.
All the general relativistic physics can be read from the boundaries of
this manifold. At the same time, however, these boundaries still carry
an infinite dimensional number of degrees of freedom.
\end{quote}

Next, let me take the liberty of describing some work of my own:

\begin{enumerate}
\def\labelenumi{\arabic{enumi})}
\setcounter{enumi}{3}
\tightlist
\item
  John Baez, ``Diffeomorphism-invariant generalized measures on the space of
  connections modulo gauge transformations'', \emph{Proceedings of the Conference 
  on Quantum Topology}, ed.\ David N.\ Yetter, World Scientific Press, Singapore, 
  1994, pp.\ 21--43.  Also available as \href{https://arxiv.org/abs/hep-th/9305045}
  {\texttt{hep-th/9305045}}.
\end{enumerate}
\noindent
This is an extremely interesting paper by a very good mathematician.
Whoops! Let's be objective here. In the loop representation of quantum
gravity, states of quantum gravity are given naively by certain
``measures'' on a space \(A/G\) of connections modulo gauge
transformations. The Chern--Simons path integral is also such a
``measure''. In both cases, the ``measure'' in question is invariant
under diffeomorphisms of space. And in both cases, the loop transform
allows one to think of these measures as instead being functions of
multiloops (collections of loops in space). This is the origin of the
relationship to knot theory.

The problem, as always in quantum field theory, is that the notion of
``measure'' must be taken with a big grain of salt --- it's not the sort
of measure they taught you about in real analysis. Instead, these
measures are a kind of ``generalized measure'' that allows you to
integrate not all continuous functions on \(A/G\) but only those lying
in an algebra called the ``holonomy algebra,'' defined by Ashtekar,
Isham and Lewandowski. To be precise and technical, this is the closure
in the \(L^\infty\) norm of the algebra of functions on \(A/G\)
generated by ``Wilson loops,'' or traced holonomies around loops. So
what we are really interested in is not diffeomorphism-invariant
measures on \(A/G\), but diffeomorphism invariant elements of the dual
of the holonomy algebra. I begin with a review of generalized measures,
introduce the holonomy algebra, and then do a bunch of new work in which
I show how to rigorously construct lots of diffeomorphism-invariant
elements of the dual of the holonomy algebra by doing lattice gauge
theory on graphs embedded in space. Again, as with the work discussed
above, we see that the discrete and continuum approaches to space go
hand-in-hand! And we see that there are some interesting connections
between singularity theory and group representation theory showing up
when we try to understand ``measures'' on the space \(A/G\).

The following is a part of a paper discussed in
\protect\hyperlink{week5}{``Week 5''}, now available from gr-qc:

\begin{enumerate}
\def\labelenumi{\arabic{enumi})}
\setcounter{enumi}{4}
\tightlist
\item
  Abhay Ashtekar and Jerzy Lewandowski,``Completeness of Wilson loop functionals 
  on the moduli space of
  \(\mathrm{SL}(2,\mathbb{C})\) and \(\mathrm{SU}(1,1)\)-connections'',
   available as
  \href{https://arxiv.org/abs/gr-qc/9304044}{\texttt{gr-qc/9304044}}.
\end{enumerate}
\noindent
I didn't discuss this aspect of the paper, so let me quote the abstract:

\begin{quote}
The structure of the moduli spaces \(M := A/G\) of (all, not just flat)
\(\mathrm{SL}(2,\mathbb{C})\) and \(\mathrm{SU}(1,1)\) connections on a
\(n\)-manifold is analysed. For any topology on the corresponding spaces
\(A\) of all connections which satisfies the weak requirement of
compatibility with the affine structure of \(A\), the moduli space \(M\)
is shown to be non-Hausdorff. It is then shown that the Wilson loop
functionals --- i.e., the traces of holonomies of connections around
closed loops --- are complete in the sense that they suffice to separate
all separable points of \(M\). The methods are general enough to allow
the underlying \(n\)-manifold to be topologically non-trivial and for
connections to be defined on non-trivial bundles. The results have
implications for canonical quantum general relativity in 4 and 3
dimensions.
\end{quote}

By the way, someone should extend this result to more general noncompact
semisimple Lie groups, and also show that for all compact semisimple Lie
groups the Wilson loop functionals in any faithful representation
\emph{do} separate points (this is known for the fundamental
representation of \(\mathrm{SU}(n)\)). If I had a bunch of grad students
I would get one to do so.

The following was discussed in an earlier edition of this series,
\protect\hyperlink{week11}{``Week 11''}, but is now available from
gr-qc:

\begin{enumerate}
\def\labelenumi{\arabic{enumi})}
\setcounter{enumi}{5}
\tightlist
\item
   Ranjeet S.\ Tate, \emph{An Algebraic Approach to the Quantization of 
   Constrained Systems: Finite Dimensional Examples}, 
  Ph.D.~Thesis, Department of Physics, Syracuse University, 124 pages, 
  available as
  \href{https://arxiv.org/abs/gr-qc/9304043}{\texttt{gr-qc/9304043}}.
\end{enumerate}
\noindent
I haven't read the following one but it seems like an interesting
application of loop variables to more down-to-earth physics; Gambini was
one of the originators of the loop representation, and intended it for
use in QCD:

\begin{enumerate}
\def\labelenumi{\arabic{enumi})}
\setcounter{enumi}{5}
\tightlist
\item
   Rodolfo  Gambini and Leonardo Setaro, ``\(\mathrm{SU}(2)\) QCD in the 
   path representation'', available as
  \href{https://arxiv.org/abs/hep-lat/9305001}{\texttt{hep-lat/9305001}}.
  (``\texttt{hep-lat}'' is the computational and lattice physics
  preprint list.)
\end{enumerate}
\noindent
Let me quote the abstract:

\begin{quote}
We introduce a path-dependent hamiltonian representation (the path
representation) for \(\mathrm{SU}(2)\) with fermions in 3 + 1
dimensions. The gauge-invariant operators and hamiltonian are realized
in a Hilbert space of open path and loop functionals. We obtain a new
type of relation, analogous to the Mandelstam identity of second kind,
that connects open path operators with loop operators. Also, we describe
the cluster approximation that permits to accomplish explicit
calculations of the vacuum energy density and the mass gap.
\end{quote}

\hypertarget{week15}{%
\section{May 23, 1993}\label{week15}}

Last weekend we had a conference on Knots and Quantum Gravity here at
Riverside. I will briefly describe the talks, many of which will
eventually appear in a conference proceedings volume. I think that to be
nice I will list these talks in order of how technical my descriptions
will be, rather than chronologically.

Oleg Viro spoke on ``simplicial topological quantum field theories''.
There has been a lot of interest recently in constructing topological
quantum field theories using triangulations of manifolds. The most
famous of these is due to Turaev and Viro. Witten showed that this one
is the same as quantum gravity in 2+1 dimensions. The nice thing is that
this gives an alternate description of the Turaev--Viro theory in terms
of more traditional ideas from field theory, so the same theory has a
``discrete'' and a ``continuum'' formulation --- some evidence for my
notion that quantum gravity will resolve the old ``is space continuous
or discrete'' argument by saying ``both, and neither,'' just as quantum
mechanics resolved the old ``is light a wave or a particle'' dispute!
(Hegel would've loved it.)

When constructing simplicial topological quantum field theories, one has
to prove that the answer you get is independent of the triangulation.
Viro reviewed a couple sets of ``moves'' whereby one can get between any
two triangulations of the same manifold --- the Alexander moves, and the
Pachner moves. He also discussed an alternate, and more convenient, way
of describing manifolds by ``special spines''. Here the idea is as
follows. Pick a bunch of points in the manifold. From each one, start
blowing a little bubble, which grows bigger and bigger until it bumps
into the other bubbles. The result is something very much like a foam of
soap bubbles, which generically have polyhedral faces, with edges and
vertices of a special sort. Look at a mess of foam sometime if you don't
know what I mean: in 3 dimensions, three bubbles meet at an edge, and
four at a vertex. One can describe this situation purely
combinatorially, and it contains all the information about the manifold.
There are a certain set of moves, the Matveev moves, relating any two
such ``special spines'' for the same manifold. One can figure out what
these moves might be by staring some foam and watching how the bubbles
move.

Louis Kaufmann changed his talk from the announced subject and instead
talked about Vassiliev invariants of knots and their relation to
perturbative Chern--Simons theory. Let me just recall what all this is
about. Chern--Simons theory is a TQFT (topological quantum field theory)
in 3 dimensions in which the field is a connection. In physical terms, a
connection is just a generalization of the vector potential in
electromagnetism. Recall that in 3 dimensional space, the vector
potential is a vector field \(A\) whose curl is the magnetic field. In
quantum theory, the significance of the vector potential is as follows.
If we take a particle and carry it around a loop, its wavefunction gets
multiplied by a phase, that is, a complex number of absolute value 1.
These ``phases'' form a group, since the product of two phases is a
phase. This group is called \(\mathrm{U}(1)\), since we can think of
phases as \(1\times1\) unitary matrices. A key idea in modern physics is
to generalize the heck out of electromagnetism by allowing other groups
to play the role of phases. The group we choose is called the ``gauge
group.'' The second simplest choice after \(\mathrm{U}(1)\) is
\(\mathrm{SU}(2)\), the \(2\times2\) unitary matrices with determinant =
1.

You can do Chern--Simons theory with any gauge group but it's especially
simple with gauge group \(\mathrm{SU}(2)\). An \(\mathrm{SU}(2)\)
connection is just a kind of field that lets one do ``parallel
translation'' around a loop in space and get an element of
\(\mathrm{SU}(2)\). Mathematicians call this the holonomy of the
connection around the loop. Physicists typically take the trace of the
group element (in this case, just the sum of the diagonal entries of the
\(2\times2\) matrix), and call that the ``Wilson loop observable,'' a
function of the connection that depends on the loop.

Now the great thing about Chern--Simons theory is that the theory is
independent of any choice of coordinates or background metric on
spacetime. This is part of what we mean by saying the theory is a TQFT.
Another aspect of a TQFT is that there is a ``vacuum state'' and we can
calculate the expectation value of a Wilson loop in the vacuum state.
The idea one should have is that the connection is undergoing all sorts
of ``quantum fluctuations'' in the vacuum state, but that we can ask for
the \emph{average} value for the trace of the connection of the holonomy
around a loop in the vacuum state. Given a knot \(K\), we write this
expectation value as \(\langle K \rangle\). Now the great thing about
Chern--Simons theory is that the vacuum state does not care what
coordinates you use to describe it. Thus \(\langle K \rangle\) does not
depend on the \emph{geometry} of the knot \(K\) (which would take
coordinates or a metric to describe), but only on its topology. In other
words, \(\langle K \rangle\) is a knot invariant. In fact we can define
\(\langle K \rangle\) not just for knots, but also for links (bunches of
knots, possibly intertangled), by taking the expectation value of the
product of the Wilson loops, one for each knot. So Chern--Simons theory
really gives a link invariant. Witten showed that this link invariant is
just the Kauffman bracket, which is an invariant easily calculated using
the rules:

\textbf{Rule 1:} If \(K\) is the ``empty link,'' the link with NO
components whatsoever --- i.e., just the empty set --- we have
\[\langle K \rangle = 1.\] This is sort of a normalization rule.

\textbf{Rule 2:} If \(K'\) is obtained from \(K\) by adding an unlinked
copy of the unknot (an unknotted circle) to \(K\),
\[\langle K' \rangle = -(a^2 + a^{-2})\langle K \rangle.\] Here \(a\) is
an adjustable parameter that appears in Chern--Simons theory --- a
function of the coupling constant.

\textbf{Rule 3:} Suppose \(K\), \(L\), and \(L'\) are 3 knots or links
differing at just one crossing (we're supposing them to be drawn as
pictures in 2 dimensions). And suppose at this crossing they look as
follows: \[
  \begin{tikzpicture}
    \begin{knot}[clip width=7pt]
      \strand[thick] (0,0)
        to [out=down,in=up] (1,-1);
      \strand[thick] (1,0)
        to [out=down,in=up] (0,-1);
      \flipcrossings{1}
    \end{knot}
    \node at (0.5,-1.5) {$K$};
  \end{tikzpicture}
  \qquad\qquad
  \begin{tikzpicture}
    \begin{knot}
      \strand[thick] (0,0)
        to (0,-1);
      \strand[thick] (1,0)
        to (1,-1);
    \end{knot}
    \node at (0.5,-1.5) {$L$};
  \end{tikzpicture}
  \qquad\qquad
  \begin{tikzpicture}
    \begin{knot}
      \strand[thick] (0,0)
        to [out=down,in=down,looseness=1.3] (1,0);
      \strand[thick] (1,-1)
        to [out=up,in=up,looseness=1.3] (0,-1);
      \flipcrossings{1}
    \end{knot}
    \node at (0.5,-1.5) {$L'$};
  \end{tikzpicture}
\] Any rotated version of this picture is fine too.

Then
\[\langle K \rangle = a\langle L\rangle + a^{-1}\langle L'\rangle.\]
That's it! I leave as an exercise for the reader to calculate
\(\langle K \rangle\) for the trefoil knot, \[
  \begin{tikzpicture}
    \begin{knot}[clip width=7pt]
      \strand[thick] (0,0)
        to [out=down,in=up] (1,-1)
        to [out=down,in=up] (0,-2)
        to [out=down,in=up] (1,-3);
      \strand[thick] (1,0)
        to [out=down,in=up] (0,-1)
        to [out=down,in=up] (1,-2)
        to [out=down,in=up] (0,-3);
      \flipcrossings{1,3}
      \strand[thick] (1,-3)
        to [out=down,in=down] (2,-3)
        to (2,0)
        to [out=up,in=up] (1,0);
      \strand[thick] (0,-3)
        to [out=down,in=down] (3,-3)
        to (3,0)
        to [out=up,in=up] (0,0);
    \end{knot}
  \end{tikzpicture}
\] and get \(-a^5 -a^{-3} -a^{-7}\). Then try the mirror-image trefoil,
or ``left-handed trefoil,'' and see what you get.

Now in quantum field theory people like doing perturbative calculations,
and that's interesting here even though we know the exact answer.
Namely, there is a coupling constant \(c\) in Chern--Simons theory such
that \(a = \exp(c)\), and if one uses Feynman diagrams and the rest of
the usual machinery for quantum field theory and does a perturbative
calculation of the vacuum expectation value of Wilson loops, one gets
the same answer, but as a power series in \(c\). The coefficient of the
\(c^n\) is a special sort of link invariant called a ``Vassiliev
invariant'' of degree \(n\). I discussed these a lot in
\protect\hyperlink{week3}{``Week 3''} (see below for how to get that
article), so I won't repeat myself here. In any event, Kauffman gave a
nice discussion of this sort of thing.

Viktor Ginzburg talked about his work with Milgram on Vassiliev
invariants. He had hoped to show that these were an almost complete set
of knot invariants, able to distinguish between any two knots that
weren't just orientation-reversed versions of each other (here we equip
the knots with an orientation, or field of arrows running along the
knot). He came to the conference, as he said, sadder and wiser. He
presented a nice result on Vassiliev invariants that might be a step
towards proving completeness.

Dana Fine spoke on ``Chern--Simons theory and the Wess--Zumino--Witten
model''. There is a very interesting ``ladder of field theories'' that
contains ``topological quantum gravity'' in 4 dimensions, Chern--Simons
theory in 3 dimensions, and the Wess--Zumino--Witten model in 2
dimensions. Dana Fine spoke on the bottom 2 rungs of this ladder. He
described a very explicit, although still formal, reduction of the
Chern--Simons path integral (the integral one does to compute the
expectation values I mentioned above) to the path integral in the
Wess--Zumino--Witten model. The relation between CS theory and the WZW
model is what Witten used in his original argument that CS theory gives
interesting link invariants, so this is of interest in knot theory as
well as physics.

On Saturday morning we had a nice trio of talks from the Syracuse gang,
Syracuse University being a hotbed of new work on quantum gravity.
Ashtekar and Smolin are there, as are a bunch of good grad students
(Bernd Br\"ugmann was there until very recently) and postdocs, including
Jerzy Lewandowski and Renate Loll. The whole gang is moving down to Penn
State this summer, and they will be hiring Jorge Pullin, now at Utah
State. There are not many people working on quantum gravity --- and not
many jobs in the field --- so Penn State will become arguably \emph{the}
center, at least in the US (let us not forget Penrose, Hawking, Isham,
et al!), for work on this subject.

Renate Loll spoke on the ``Loop representation of gauge theory and
gravity.'' This was an introduction to the ideas of Gambini, Trias,
Rovelli, Smolin, et al on doing quantum field theory solely in terms of
Wilson loops. It is this approach that makes the connection between
quantum gravity and knot theory.

Abhay Ashtekar spoke on ``Loop transforms.'' The process of encoding a
connection in terms of all its Wilson loops is called the loop
transform, and it can be regarded as a nonlinear generalization of the
Fourier transform. Ashtekar has led an effort to make this transform
thoroughly rigorous and mathematically respectable, and he discussed
this work.

Jorge Pullin spoke on ``The quantum Einstein equations and the Jones
polynomial''. He outlined his work with Gambini and Br\"ugmann in which
they show perturbatively that the Kaufmann bracket (or, alternatively,
Chern--Simons theory) gives a state of 4-d quantum gravity. This is
perhaps the most exciting aspect of the ``ladder of field theories''
mentioned above, since it touches upon the real world --- or at least
comes darn close.

On Sunday, Gerald Johnson started things off with an introduction to his
work on making the Feynman path integral rigorous. This is relevant
because a main problem with Witten's otherwise marvelous work is that
the path integral in Chern--Simons theory has not been made rigorous.
Dana Fine's talk offered one approach, and my own talk offered another
(based on my recent paper).

Perhaps the most novel talk was by Paolo Cotta-Ramusino (``4d quantum
gravity and knot theory'') describing his work with Maurizio Martellini
on \(4\)-dimensional TQFTs and invariants of 2-knots, that is, embedded
surfaces in \(\mathbb{R}^4\) (or more general 4-manifolds). This is an
attempt to push the Wilson loop story up one dimension, in an effort to
make it applicable to theories similar to quantum gravity. These
theories are the so-called ``BF theories,'' whose Lagrangian is of the
form \(\operatorname{tr}(B \wedge F)\), where \(B\) is a Lie algebra
valued \(2\)-form and \(F\) is the curvature of a connection. Martellini
and Cotta-Ramusino's work on this is still in a preliminary stage but it
seems rather promising.

Perhaps the most controversial talk was by Louis Crane, entitled
``Quantum gravity, spin geometry and categorical physics.'' This was
about his ideas on using category theory to construct \(4\)-dimensional
TQFTs. He also emphasized the importance of TQFTs that use
triangulations but wind up being independent of the triangulation, thus
slipping through the discrete/continuous distinction. Many of his
assertions provoked violent reactions from the physicists present.

Finally, I spoke on ``Strings, tangles and gauge fields,'' beginning by
pointing out some relationships between closed string field theory and
the loop representation of quantum gravity, and then retreating to safer
ground and describing my work on trying to make the Chern--Simons path
integral and the loop representation more mathematically rigorous. I
will write a paper on this subject this summer and try to further build
up my case for the conjecture that string theory and gauge field theory
are in a sense dual descriptions of certain TQFTs.

A rather technical introduction to currently interesting topics in
closed string field theory has just appeared:

\begin{enumerate}
\def\labelenumi{\arabic{enumi})}
\tightlist
\item
  Barton
  Zwiebach,``Closed string field theory --- an introduction'', available as
  \href{https://arxiv.org/abs/hep-th/9305026}{\texttt{hep-th/9305026}}.
\end{enumerate}

Unfortunately for me, he mainly treats theories of strings moving around
on a manifold with a background metric, while it seems that the loop
representation of quantum gravity is very like a ``background-free''
string field theory. A paper that recently came out and appears to
support my notions is the following:

\begin{enumerate}
\def\labelenumi{\arabic{enumi})}
\setcounter{enumi}{1}
\tightlist
\item
  S.\ G.\
  Naculich, H.\ A.\ Riggs, and H.\ J.\ Schnitzer, ``Two-dimensional Yang--Mills 
  theories are string theories'', available as
  \href{https://arxiv.org/abs/hep-th/9305097}{\texttt{hep-th/9305097}}.
\end{enumerate}
\noindent
Apparently this builds on work by Gross, Taylor, and Minahan which
treated \(\mathrm{SU}(n)\) Yang--Mills theories in 2 dimensions as string
theories, and does something similar for the gauge groups
\(\mathrm{SO}(n)\) and \(\mathrm{Sp}(n)\).

I have a pack of interesting papers to describe but I am already worn
out, so I will put that off until next week, except for the following
paper by Smolin that I seem to have neglected:

\begin{enumerate}
\def\labelenumi{\arabic{enumi})}
\setcounter{enumi}{2}
\tightlist
\item
  Lee Smolin, ``What can we learn from the study of non-perturbative quantum general
  relativity?'', available as
  \href{https://arxiv.org/abs/gr-qc/9211019}{\texttt{gr-qc/9211019}}.
\end{enumerate}
\noindent
This is a nice introduction to current issues associated to the loop
representation of quantum gravity and nonperturbative quantum gravity in
general. As should be evident from my weekly reviews, the subject seems
to be moving faster and faster, and it is best to read some of the
review papers like this one by Smolin if one wants to figure out where
things are now and where they might be heading.

\hypertarget{week16}{%
\section{May 30, 1993}\label{week16}}

A nice crop of papers has built up while I've been taking a
break\ldots{} In \protect\hyperlink{week15}{``Week 15''} I talked a bit
about constructing topological quantum field theories starting with a
triangulation of spacetime, and how this seems to sneak around the old
``is spacetime continuous or discrete'' argument. Let me describe a bit
about one of the more mathematically elegant physics papers I've run
across in a while, which treats exactly this issue. Then I'll describe
two review articles, one on gravity in 2+1 dimensions (which is closely
related to the lattice business), and one on Lagrangians for quantum
gravity.

\begin{enumerate}
\def\labelenumi{\arabic{enumi})}
\tightlist
\item
  Stephen-wei Chung, Masafumi Fukuma and Alfred
  Shapere, ``Structure of topological lattice field theories in three
  dimensions'', available as
  \href{https://arxiv.org/abs/hep-th/9305080}{\texttt{hep-th/9305080}}.
\end{enumerate}

What's a \(2\)-dimensional ``topological lattice field theory''?
According to the definition used in this paper, it goes like this. First
take a compact oriented 2-manifold without boundary \(M\), that is, an
\(n\)-holed torus. (One could also discuss the case when there is a
boundary, but to keep life simple we won't here.) We want to calculate a
number \(Z(M)\), the partition function of \(M\), since the partition
function is a basic ingredient in Feynman's approach to quantum field
theory. We first triangulate \(M\)\ldots{} so a patch might look like:
\[
  \begin{tikzpicture}
    \draw[thick] (0,0) to (1,1.6) to (2,0) to cycle;
    \draw[thick] (2,0) to (1,1.6) to (3,1.6) to cycle;
    \draw[thick] (0,0) to (2,0) to (1,-1.6) to cycle;
  \end{tikzpicture}
\] Then ``disassemble'' \(M\) into separate triangles, like this: \[
  \begin{tikzpicture}
    \draw[thick] (0,0) to (1,1.6) to (2,0) to cycle;
    \draw[thick] (2.5,0) to (1.5,1.6) to (3.5,1.6) to cycle;
    \draw[thick] (0,-0.5) to (2,-0.5) to (1,-2.1) to cycle;
  \end{tikzpicture}
\] Now assign to each edge of the disassembled version of \(M\) a
``color'' taken from a fixed finite set \(S\). Note that there are twice
as many edges in the disassembled version of \(M\) as in the original
triangulation of \(M\). Any way of assigning a color to each edge of the
disassembled \(M\) will be called a ``coloring''. We think of a coloring
as a ``history of the world'' and we will compute \(Z(M)\) by summing a
certain quantity over all colorings.

To compute this quantity, we need two pieces of data that determine our
theory. First, for each \(i,j\) in \(S\), we fix a complex number
\(g^{ij}\). We require that the matrix \(g^{ij}\) be invertible. We
define \(g_{ij}\) to be the matrix inverse of \(g^{ij}\). We can raise
and lower indices with \(g\) as if it were a metric. The matrix \(g\)
will be used when we glue two edges of the disassembled \(M\) together
in the process of rebuilding \(M\). Second, for each \(i,j,k\) in \(S\),
we fix a number \(c_{ijk}\). This number comes in because each triangle
has three edges.

Here's how we calculate \(Z(M)\). Write down one index next to each edge
of the disassembled \(M\) --- by ``index'' I mean something like
\(i,j,k\) running over \(S\). Then write down the obvious factor of
\(g\) for each pair of edges that get glued together when we form \(M\),
and write down the obvious factor of \(c\) for each triangle in \(M\).
Finally, sum over all colorings to get \(Z(M)\).

For example, if \(M\) were a torus that we triangulated with two
triangles like this \[
  \begin{tikzpicture}
    \draw[thick] (0,0) to (1,1.6) to (2,0) to cycle;
    \draw[thick] (0,0) to (2,0) to (1,-1.6) to cycle;
  \end{tikzpicture}
\] --- with opposite edges of the parallelogram identified --- we would
dissasemble \(M\) and label the edges like this, say: \[
  \begin{tikzpicture}
    \draw[thick] (0,0) to node[fill=white]{$l$} (1,1.6) to node[fill=white]{$m$} (2,0) to node[fill=white]{$j$} cycle;
    \draw[thick] (0,-0.5) to node[fill=white]{$i$} (2,-0.5) to node[fill=white]{$k$} (1,-2.1) to node[fill=white]{$n$} cycle;
  \end{tikzpicture}
\] To form \(M\) we glue \(i\) to \(j\), \(k\) to \(l\), and \(m\) to
\(n\). So we write down \[g^{ij}g^{kl}g^{mn}c_{jml}c_{ink}\] and then
sum over \(i,j,k,l,m,n\) to get \(Z(M)\). Notice that for this procedure
to be well defined it had better not matter whether we write \(g^{ij}\)
or \(g^{ji}\), since we have no way of knowing which to use. So \(g\)
had better be symmetric. Similarly, we had better have
\(c_{ijk} = c_{jki}\) --- invariance under cyclic permutations. Note
that since \(M\) is oriented we can (and will) require that we go around
each triangle counterclockwise when writing down things like
\(c_{ink}\), as we have done above.

Okay, this is a pretty scheme, but the real point is that it should be
independent of the triangulation of \(M\) we chose, for us to have
something that deserves to be called ``topological.'' This imposes extra
conditions on \(g\) and \(c\). Here it is handy to know that we can get
between any two triangulations of \(M\) using a sequence of two moves
and their inverses. The first move is called the ``(2,2) move.'' It
looks like this: \[
  \begin{tikzpicture}
    \draw[thick] (0,0) node{$\bullet$} to (1,1.5) node{$\bullet$} to (2,0) node{$\bullet$} to (1,-1.5) node{$\bullet$} to cycle;
    \draw[thick] (1,1.5) to (1,-1.5);
  \end{tikzpicture}
  \raisebox{4.5em}{$\qquad\longleftrightarrow\qquad$}
  \begin{tikzpicture}
    \draw[thick] (0,0) node{$\bullet$} to (1,1.5) node{$\bullet$} to (2,0) node{$\bullet$} to (1,-1.5) node{$\bullet$} to cycle;
    \draw[thick] (0,0) to (2,0);
  \end{tikzpicture}
\] It is called the (2,2) move since it really amounts to taking 2 faces
of a tetrahedron and replacing them with the other 2 faces! There is a
similar (3,1) move that takes 3 faces of a tetrahedron and replaces them
with the other 1, as follows: \[
  \begin{tikzpicture}
    \draw[thick] (0,0) node{$\bullet$} to (1.8,3) node{$\bullet$} to (3.6,0) node{$\bullet$} to cycle;
    \node at (1.8,1.1) {$\bullet$};
    \draw[thick] (0,0) to (1.8,1.1);
    \draw[thick] (1.8,3) to (1.8,1.1);
    \draw[thick] (3.6,0) to (1.8,1.1);
  \end{tikzpicture}
  \raisebox{4.5em}{$\qquad\longleftrightarrow\qquad$}
  \begin{tikzpicture}
    \draw[thick] (0,0) node{$\bullet$} to (1.8,3) node{$\bullet$} to (3.6,0) node{$\bullet$} to cycle;
  \end{tikzpicture}
\] (This drawing done by my friend Bruce Smith in a fit of insomnia!)
These are examples of the ``Pachner moves,'' and the same idea works in
any dimension. But in 2 dimensions we can use a move called the ``bubble
move'' instead of the (3,1) move. Here is where drawing vertices as
\(\bullet\)'s is crucial: \[
  \begin{tikzpicture}
    \draw[thick] (1,1.5) to [bend left=70] (1,-1.5) to [bend left=70] cycle;
    \draw[thick] (1,1.5) node{$\bullet$} to node{$\bullet$} (1,-1.5) node{$\bullet$};
  \end{tikzpicture}
  \raisebox{4.5em}{$\qquad\longleftrightarrow\qquad$}
  \begin{tikzpicture}
    \draw[thick] (1,1.5) node{$\bullet$} to (1,-1.5) node{$\bullet$};
  \end{tikzpicture}
\] On the left, we have two hideously deformed triangles (remember, this
is topology!) that are attached along TWO edges, leaving two edges
exposed, and in the right we have collapsed them down to a single edge.
We leave it as a fun exercise to show that you can do anything with the
(2,2) move and the bubble move that you can do with the (2,2) move and
the (3,1) move.

Requiring that \(Z(M)\) be invariant under the (2,2) moves amounts to
the following equation --- if you check it, you will make sure you
understand what's going on: \[c_{xy}^u c_{uz}^w = c_{xu}^w c_{yz}^u.\]
Here I have raised indices using the ``metric'' \(g\). This equation
looks sort of hairy, but it's actually something very nice in disguise.
We need to tease out its inner essence! Suppose we take a vector space
\(A\) having the colors in \(S\) as a basis, and use the tensor
\(c_{ij}^k\) to define a bilinear map from \(A \times A\) to \(A\). Then
the equation above says this map is an associative product! If you
ponder the picture of the (2,2) move for a while, this should become
obvious to you. Think of each triangle as being a gadget that you can
feed vectors into from two sides and have the ``product'' pop out on the
third side. Then the equation \[
  \begin{tikzpicture}
    \draw[thick] (0,0) node{$\bullet$} to (1,1.5) node{$\bullet$} to (2,0) node{$\bullet$} to (1,-1.5) node{$\bullet$} to cycle;
    \draw[thick] (1,1.5) to (1,-1.5);
  \end{tikzpicture}
  \raisebox{4.5em}{$\qquad=\qquad$}
  \begin{tikzpicture}
    \draw[thick] (0,0) node{$\bullet$} to (1,1.5) node{$\bullet$} to (2,0) node{$\bullet$} to (1,-1.5) node{$\bullet$} to cycle;
    \draw[thick] (0,0) to (2,0);
  \end{tikzpicture}
\] really is just associativity! To understand this in a deeper way,
read Kapranov and Voevodsky's paper (reviewed in
\protect\hyperlink{week4}{``Week 4''}), especially the section on the
``associahedron''.

Requiring that \(Z(M)\) be invariant under the bubble move amounts to
the following: \[c_{xu}^v c_{yv}^u = g_{xy}\] Here \(g_{xy}\) is the
matrix inverse of \(g^{xy}\). Again, I leave it as an exercise to show
this is the right equation. It is a formula expressing the ``metric''
\(g\) on \(A\) in terms of the product on \(A\)! In fact, it has a
beautiful algebraic interpretation: it says that the algebra \(A\) is
``semisimple.'' A semisimple algebra is just a direct sum of matrix
algebras, and in such algebras the inner product \(g(a,b)\) of any two
elements is just equal to \(\operatorname{tr}(a^T b)\), where \(a^T\) is
the transpose of \(a\).

So we discover a charming fact: there is a one-to-one correspondence
between topological lattice field theories in 2 dimensions and
finite-dimensional semisimple algebras over the complex numbers!

Actually this was apparently already shown by

\begin{enumerate}
\def\labelenumi{\arabic{enumi})}
\setcounter{enumi}{1}
\tightlist
\item
  C.\ Bachas and P.\ M.\ S.\ Petropoulos, ``Topological models on the lattice 
  and a remark on string theory cloning'', \emph{Commun.\ Math.\ Phys.}
  \textbf{152} (1993) 191.
\end{enumerate}

and

\begin{enumerate}
\def\labelenumi{\arabic{enumi})}
\setcounter{enumi}{2}
\tightlist
\item
   M.\ Fukuma, S.\ Hosono and H.\ Kawai, ``Lattice topological field theory in
    two-dimensions'', available as
  \href{https://arxiv.org/abs/hep-th/921254}{\texttt{hep-th/921254}}.
\end{enumerate}

The big result of the present paper is to generalize this to 3
dimensions. The authors consider a specific definition of 3d topological
lattice field theories in which one chops a 3d manifold up into
tetrahedra and assigns colors to edges. They claim to get a one-to-one
correspondence between these and finite-dimensional Hopf algebras for
which the antipode squared is the identity! If you don't know what a
Hopf algebra is, let me simply say it is a very beautiful sort of thing
that has both a product and a ``coproduct,'' and they come up all the
time in group theory, knot theory, and the study of quantum groups. So
we are seeing that there is a profound correspondence between topology
and algebra, with higher-dimensional topology giving more subtle
algebra.

(In fact, I am a little worried that the authors haven't stated the
theorem quite precisely enough to have it be quite true, but it's
basically right --- I am afraid only only gets a particular class of
Hopf algebras, those which are semisimple and cosemisimple. I may be
missing something.)

Let me conclude with a few very exciting open problems.

A. One could instead consider theories in which colors are only assigned
to faces. This turns out not to broaden the class of examples: any Hopf
algebra has a dual Hopf algebra, and one just gets the theory associated
to the dual Hopf algebra this way! But if one considers theories in
which colors are assigned BOTH to edges AND faces one apparently gets a
larger class of 3d examples. What algebraic structure do these
correspond to? B. The Turaev--Viro theory of quantum gravity ---
described below --- is a 3d topological lattice field theory of some
sort. Where does it fit into this picture? The authors ask this question
but don't answer it. Also, a more difficult problem --- where does
Chern--Simons gauge theory fit into this picture? C. The 64,000 dollar
question: how does all this generalize to 4 dimensions? What sort of
algebraic structure corresponds to a 4d topological lattice field
theory? It is becoming increasingly clear that 4d field theories will
involve some kind of ``higher algebra'' that we are only beginning to
understand.

\begin{center}\rule{0.5\linewidth}{0.5pt}\end{center}

\begin{enumerate}
\def\labelenumi{\arabic{enumi})}
\setcounter{enumi}{3}
\tightlist
\item
   Steven Carlip, ``Six ways to quantize (2+1)-dimensional gravity'', 
  available as \href{https://arxiv.org/abs/gr-qc/9305020}{\texttt{gr-qc/9305020}}.
\end{enumerate}

While we have no real way to quantize gravity in 3+1 dimensions ---
although lots of good ideas --- we have six, count 'em, six, ways to do
it in 2+1 dimensions! Sometimes this sort of thing makes one yearn to be
a physicist in some other, lower-dimensional universe. However, lest one
make such wish prematurely to a genie passing by, one should note that
life in 2+1 dimensions is boring compared to our 3+1-dimensional world.
The reason can be seen from the following count of the number of
independent components of the Riemann tensor \(R_{ijkl}\) which vanishes
when spacetime is flat, and the Einstein tensor, which vanishes when the
vacuum Einstein equations hold:

\begin{longtable}[]{@{}cll@{}}
\toprule
dimension & Riemann & Einstein\tabularnewline
\midrule
\endhead
1 & 0 & 0\tabularnewline
2 & 1 & 1\tabularnewline
3 & 6 & 6\tabularnewline
4 & 20 & 10\tabularnewline
\bottomrule
\end{longtable}

What this means is that, until one gets up to dimension 4, the vacuum
Einstein equations imply that spacetime is flat. That means that there
are no gravitational waves in empty space; there are only global,
topological effects. Typically this means that if space is compact there
are only finitely many degrees of freedom. This means that 2+1 quantum
gravity is really quantum mechanics, not full-fledged quantum field
theory (which deals with \emph{local} excitations --- wiggles in the
metric and such --- and infinitely many degrees of freedom). The good
news is, this means that 2+1 gravity is somewhat understandable --- no
nasty infinities or ill-defined integrals needed, etc. The bad news is,
it means 2+1 gravity is not too much like 3+1 gravity. But still, many
of the most puzzling qualitative features of quantum gravity are present
in the 2+1 case. In particular, one has a testing ground in which to
look at the interlocking triad of problems that stump us in the 3+1
case: the problem of time, the problem of observables, and the inner
product problem. In brief these are: what is time evolution in quantum
gravity, what are the observables in quantum gravity, and what is the
inner product on the space of states of quantum gravity? As you can see,
we are overwhelmingly ignorant about quantum gravity! I think that work
on 2+1 gravity has given us some interesting clues about these problems.

Carlip describes 6 approaches to 2+1 gravity. I'll list them and comment
on them briefly below. But one point to make is that these approaches
have \emph{not} all been shown to be equivalent; on the contrary, they
seem to give different answers. Part of the problem in my opinion is
that we do not have enough criteria for a ``good'' theory of 2+1 quantum
gravity. Certainly one would like to see that in the \(\hbar\to 0\)
limit the theory reduces to classical gravity in some sense or other
(but this is a bit vague). Perhaps another thing one could hope for is
that the theory be a 2+1-dimensional TQFT. I am not sure which of the
approaches below give a TQFT (although \#6 definitely does and probably
so does \#2):

\hypertarget{reduced-adm-phase-space-quantization}{%
\subparagraph{\#1 Reduced ADM phase space
quantization}\label{reduced-adm-phase-space-quantization}}

The ``ADM'' or Arnowitt--Deser--Misner formalism amounts to what people
would typically call canonical quantization: one writes down a
description of the phase space of quantum gravity in terms of initial
data, figures out the Poisson brackets of functions on this phase space,
and then tries to quantize by turning them into commutators. In gravity
the roles of ``position'' and ``momentum'' variables are played by the
metric on space at a given time, and the extrinsic curvature (or more
precisely, something cooked up from it).

\hypertarget{chern-simons-theoryconnection-representation}{%
\subparagraph{\#2 Chern--Simons theory/Connection
representation}\label{chern-simons-theoryconnection-representation}}

This is essentially the 2+1 analog of Ashtekar's approach in 3+1
dimensions, in that a connection and triad field play the main role,
rather than the metric. However, in 2+1 dimensions we can lump the triad
field and the connection together to get an ``\(I\mathrm{SO}(2,1)\)
connection'' --- where \(I\mathrm{SO}(2,1)\) is mildly terrifying
notation for the Poincar\'e group in 2+1 dimensions (or ``inhomogeneous
Lorentz group,'' hence the ``I''). The action for the theory then
becomes the Chern--Simons action, as noted by Witten.

\hypertarget{covariant-canonical-quantization}{%
\subparagraph{\#3 Covariant canonical
quantization}\label{covariant-canonical-quantization}}

This might sound oxymoronic to some, but what it means is that the phase
space of solutions is described in a manifestly covariant way, rather
than in terms of initial data, and then one tries to turn Poisson
brackets into commutators.

\hypertarget{loop-representation}{%
\subparagraph{\#4 Loop representation}\label{loop-representation}}

The loop representation of quantum gravity starts with the connection
representation and then takes traces of holonomies around loops ---
so-called Wilson loops --- as the basic variables to quantize. This
suffers irritating technical problems in 2+1 dimensions, as noted in the
following recent paper:

\begin{enumerate}
\def\labelenumi{\arabic{enumi})}
\setcounter{enumi}{4}
\tightlist
\item
  Donald Marolf, An illustration of 2+1 gravity loop transform troubles, 
  available as \href{https://arxiv.org/abs/gr-qc/9305015}{\texttt{gr-qc/9305015}}.
\end{enumerate}
\noindent
I know that Ashtekar and Loll are attacking these problems right
now; Loll discussed this a bit in a lecture she gave in my seminar.

\hypertarget{the-wheeler-dewitt-equation}{%
\subparagraph{\#5 The Wheeler--DeWitt
equation}\label{the-wheeler-dewitt-equation}}

Here we proceed as in approach \#1 but attempt to impose the Hamiltonian
and diffeomorphism constraints after quantizing. That is, we start with
an overly large phase space of initial data for general relativity ---
overly large because a given solution of Einstein's equations will have
many different initial data on different spacelike slices --- quantize
by turning Poisson brackets into commutators, and \emph{then} try to take care
of the mistake we made by defining the ``physical'' states to be those
annihilated by certain operators, the Hamiltonian and diffeomorphism
constraints. I gave a brief intro to this in
\protect\hyperlink{week11}{``Week 11''}. This is the most traditional
approach in 3+1 gravity.

\hypertarget{lattice-approaches}{%
\subparagraph{\#6 Lattice approaches}\label{lattice-approaches}}

These are closely related to the topological lattice field theories
described above. Here we treat spacetime as discrete, that is, as a kind
of lattice. One approach here is due to Regge and Ponzano, and recently
worked out rigorously by Turaev and Viro. To get going in this theory,
you ``triangulate'' your \(3\)-dimensional spacetime, that is, chop it
into tetrahedra. All we need to work with is this ``simplicial complex''
consisting of tetrahedra, their triangular faces, their line-segment
edges, and the vertex points. We assume for simplicity that spacetime is
compact, so we can use finitely many tetrahedra. Thus everything in
sight is finite and discrete. A ``history of the world'' in this theory
amounts to labelling each edge with a length, or ``spin'', that must be
\(0,\frac{1}{2},1,\frac{3}{2},\ldots\) or \(j/2\). There are thus finitely many possible
histories. To do calculations in this theory, we follow Feynman's
procedure and ``sum over histories'' --- write down a formula for the
quantity we are interested in, and add up its value for all histories,
weighted by a quantity depending on the history, the exponential of the
action of that history, to obtain the vacuum expectation value of the
quantity. The formula for the action is very familiar to folks
knowledgeable about quantum theory. Each tetrahedron has 6 edges
labelled by spins, and we calculate a quantity called the ``\(6j\)
symbol'' from these spins and then add it up for all tetrahedra. In the
Turaev--Viro version, we have replaced the gauge group \(\mathrm{SU}(2)\)
by the corresponding quantum group, with the quantum parameter \(q\) a
root of unity, so there are only finitely many irreducible
representations, or spins, to sum over. (See
\protect\hyperlink{week5}{``Week 5''} for the vaguest of introductions
to quantum groups and their representations!) The beauty of this theory
is that the answer one gets is independent of the triangulation one has
chosen.

While I'm at it, let me list some key references to the subject of
lattice 2+1 gravity, a subject I'm fascinated by these days.

The grandaddy of them all, the Ponzano--Regge paper, is:

\begin{enumerate}
\def\labelenumi{\arabic{enumi})}
\setcounter{enumi}{5}
\tightlist
\item
  G.\ Ponzano and T.\ Regge, ``Semiclassical limits of Racach coefficients'', 
  in F.\ Bloch (ed.), \emph{Spectroscopic and
  Group Theoretical Methods in Physics}, Amsterdam: North-Holland 1968.
\end{enumerate}

Then there are:

\begin{enumerate}
\def\labelenumi{\arabic{enumi})}
\setcounter{enumi}{6}
\item
  Edward Witten, ``(2+1)-dimensional gravity as an exactly soluble system'', 
  \emph{Nucl.\ Phys.} \textbf{B311} (1988), 46--78.
\item
  V.\ G.\ Turaev and O.\ Y.\ Viro,  ``State sum invariants of 3-manifolds and 
  quantum 6j-symbols'', \emph{Topology} \textbf{31} (1992), 865--902.
\end{enumerate}

Also Ooguri wrote a paper on 3+1 lattice gravity that has been quite
influential:

\begin{enumerate}
\def\labelenumi{\arabic{enumi})}
\setcounter{enumi}{8}
\tightlist
\item
  Hirosi Ooguri, ``Topological lattice models in four dimensions'', 
  \emph{Mod.\ Phys.\ Lett.} \textbf{A7} (1992), 2799--2810.
  Available as
  \href{https://arxiv.org/abs/hep-th/9205090}{\texttt{hep-th/9205090}}.  
\end{enumerate}
\noindent
And there is also the recent paper by Rovelli, which I discussed in
\protect\hyperlink{week14}{``Week 14''}. This is very readable (once you
know what's going on!) and conceptual.

\begin{enumerate}
\def\labelenumi{\arabic{enumi})}
\setcounter{enumi}{9}
\tightlist
\item
  Peter Peldan, ``Actions for gravity, with generalizations: a review'', 61 pages, 
  available as \href{https://arxiv.org/abs/gr-qc/9305011}{\texttt{gr-qc/9305011}}.
\end{enumerate}

The classic action principle for general relativity is the
Einstein--Hilbert action: the Ricci scalar times the volume form
associated to the metric. An important modification, often called the
Palatini action, takes a connection and tetrad (aka vierbein or frame
field) as basic. More recently, Plebanski invented an action using the
self-dual part of the connection and a tetrad field; this turns out to
be closely related to an action naturally associated with the Ashtekar
``new variables'' (a self-dual connection and tetrad field), although
this was realized only subsequently by Capovilla, Dell, and Jacobson.
More recently still, there is the Capovilla-Dell-Jacobson action. These
new action principles shed a very interesting new light on gravity,
particularly when it comes to quantizing it. Of course it must be
remembered that actions that give the \emph{same} classical dynamics can
(and typically DO) give \emph{different} quantum theories. So a
traditionalist might question whether these new actions give the
``right'' quantum theory of gravity. Of course, the correct response to
such a traditionalist is ``well, you come up with the `right' quantum
theory of gravity and then we can compare!'' The point is that the good
old Einstein--Hilbert action is extremely intractable when it comes to
quantization --- so perhaps it is not the ``right'' one, and \emph{any}
quantization is more enlightening than none at this stage.

Peldan presents a grand tour of the various Lagrangian formulations of
gravity, and on page 3 of this large manuscript there is a large diagram
of the main Lagrangian and Hamiltonian approaches to gravity in 3+1
dimensions, while on page 35 there is a somewhat smaller chart for 2+1
gravity. (A very brief preliminary warmup on some of these formulations
appears in my earlier article, \protect\hyperlink{week7}{``Week 7''}.) I
plan on going through this carefully in order to be able to make up for
years of neglect on my part of this sort of thing.

\hypertarget{week17}{%
\section{June 13, 1993}\label{week17}}

This'll be the last ``This Week's Finds'' for a few weeks, as I am going
up to disappear until July. I've gotten some requests for introductory
material on gauge theory, knot theory, general relativity, TQFTs and
such recently, so I just made a list of
\href{http://math.ucr.edu/home/baez/books.html}{some of my favorite
books} on this kind of thing --- with an emphasis on the readable ones.

Also, just as a little plug here, a graduate student here at UCR (Javier
Muniain) and I are turning my course notes from this year into a book
called ``Gauge Fields, Knots and Gravity,'' meant to be an elementary
introduction to these subjects. This will eventually be published by
World Scientific if all goes well. It will gently remind the reader
about manifolds, differential forms, Lagrangians, etc., develop a little
gauge theory, knot theory, and general relativity, and at the very end
it'll get to the relationship between knot theory and quantum gravity
--- at which point one could read more serious stuff on the subject.

A while back Lee Rudolph asked my opinion of the following article:

\begin{enumerate}
\def\labelenumi{\arabic{enumi})}
\tightlist
\item
   Arthur Jaffe and Frank Quinn, ``Theoretical mathematics: toward a cultural synthesis 
   of mathematics and theoretical physics", \emph{Bull.\ Amer.\ Math.\ Soc.} \textbf{29} (1993),     
   1--13.  Available as \href{https://arxiv.org/abs/math/9307227}{\texttt{arXiv:math/9307227}}.
\end{enumerate}

People who are seriously into mathematical physics will know that with
string theory the interaction between mathematicians and physicists,
especially mathematicians who haven't traditionally been close to
physics (e.g.~algebraic geometers), has strengthened steadily for the
last 10 years or so. Physicists are coming up with lots of exciting
mathematical ``results'' --- often \emph{not} rigorously proved! --- and
mathematicians are getting very interested. Let me quote the abstract:

\begin{quote}
Is speculative mathematics dangerous? Recent interactions between
physics and mathematics pose the question with some force: traditional
mathematical norms discourage speculation; but it is the fabric of
theoretical physics. In practice there can be benefits, but there can
also be unpleasant and destructive consequences. Serious caution is
required, and the issue should be considered before, rather than after
obvious damage occurs. With the hazards carefully in mind, we propose a
framework that should allow a healthy and a positive role for
speculation.
\end{quote}

Replies have been solicited, so there may be a debate on this timely
subject in the AMS Bulletin. This subject has a great potential for
flame wars --- or, as Greeks and academics refer to them, ``polemics.''
Luckily Jaffe and Quinn take a rather careful and balanced tone. I think
anyone interested in the culture of mathematics and physics should take
a look at this.

Now for two books:

\begin{enumerate}
\def\labelenumi{\arabic{enumi})}
\setcounter{enumi}{1}
\tightlist
\item
  DeWitt L.\ Sumner, \emph{New Scientific Applications of Geometry and Topology},
   \emph{Proc.\ Symp.\ Appl.\ Math.} \textbf{45}, AMS.
\end{enumerate}

This volume has a variety of introductory papers on applications of knot
theory; the titles are roughly ``Evolution of DNA topology,'' ``Geometry
and topology of DNA and DNA-protein interactions,'' ``Knot theory and
DNA,'' ``Topology of polymers,'' ``Knots and Chemistry,'' and ``Knots
and physics.''

\begin{enumerate}
\def\labelenumi{\arabic{enumi})}
\setcounter{enumi}{2}
\tightlist
\item
  Louis Kauffman and Sostenes Lins, \emph{Temperley--Lieb Recoupling Theory and 
  Invariants of 3-Manifolds}, Princeton: Princeton U.\ Press, 1994.
\end{enumerate}

This is an elegant exposition of the 3-manifold invariants obtained from
the quantum group \(SU_q(2)\) --- or in other words, from Chern--Simons
theory. In part this is a polishing of existing work, but it also
contains some interesting new ideas.

And now for some papers:

\begin{enumerate}
\def\labelenumi{\arabic{enumi})}
\setcounter{enumi}{3}
\tightlist
\item
   M.\  Carfora, M.\ Martellini and A.\
  Marzuoli, ``\(12j\)-symbols and four-dimensional quantum gravity'', 
  Dipartimento di Fisica, Universita di Roma ``La Sapienza'' preprint.
\end{enumerate}
\noindent
This is an attempt to do for 4d quantum gravity what Regge, Ponzano and
company so nicely did for 3d quantum gravity (see
\protect\hyperlink{week16}{``Week 16''}) --- describe it using
triangulated manifolds and angular momentum theory.

\begin{enumerate}
\def\labelenumi{\arabic{enumi})}
\setcounter{enumi}{4}
\tightlist
\item
   Y.\ S.\ Soibelman, ``Selected topics in quantum groups'', 
   Lectures for the European School of Group Theory, Harvard University preprint.
\end{enumerate}
\noindent
This is a nice review of Soibelman's work on quantum groups, quantum
spheres, and other aspects of ``quantum geometry.''

\begin{enumerate}
\def\labelenumi{\arabic{enumi})}
\setcounter{enumi}{5}
\tightlist
\item
 J.\ Scott Carter and Masahico Saito , ``Braids and movies", 
 \emph{Journal of Knot Theory and Its Ramifications} \textbf{5} (1996), 589--608.
\end{enumerate}
\noindent
Just as every knot or link is given as the closure of a braid --- for
example, the trefoil knot \[
  \begin{tikzpicture}
    \begin{knot}[clip width=7pt]
      \strand[thick] (0,0)
        to [out=down,in=up] (1,-1)
        to [out=down,in=up] (0,-2)
        to [out=down,in=up] (1,-3);
      \strand[thick] (1,0)
        to [out=down,in=up] (0,-1)
        to [out=down,in=up] (1,-2)
        to [out=down,in=up] (0,-3);
      \flipcrossings{1,3}
      \strand[thick] (1,-3)
        to [out=down,in=down] (2,-3)
        to (2,0)
        to [out=up,in=up] (1,0);
      \strand[thick] (0,-3)
        to [out=down,in=down] (3,-3)
        to (3,0)
        to [out=up,in=up] (0,0);
    \end{knot}
  \end{tikzpicture}
\] is the closure of \[
  \begin{tikzpicture}
    \begin{knot}[clip width=7pt]
      \strand[thick] (0,0)
        to [out=down,in=up] (1,-1)
        to [out=down,in=up] (0,-2)
        to [out=down,in=up] (1,-3);
      \strand[thick] (1,0)
        to [out=down,in=up] (0,-1)
        to [out=down,in=up] (1,-2)
        to [out=down,in=up] (0,-3);
      \flipcrossings{1,3}
    \end{knot}
  \end{tikzpicture}
\] --- every ``2-knot'' or ``2-link'' --- that is, a surface embedded in
\(\mathbb{R}^4\), is the closure of a ``2-braid''. Just as there are
``Markov moves'' that say when two links come from the same braid, there
are moves for 2-braids --- as discussed here.

\begin{enumerate}
\def\labelenumi{\arabic{enumi})}
\setcounter{enumi}{6}
\tightlist
\item
  Danny Birmingham and Mark Rakowski, ``Combinatorial invariants from four 
  dimensional lattice  models: II'', available as
  \href{https://arxiv.org/abs/hep-th/9305022}{\texttt{hep-th/9305022}}.
\end{enumerate}
\noindent
The previous paper obtains some invariants of 4-manifolds by
triangulating them and doing a kind of ``state sum'' much like those I
described in \protect\hyperlink{week16}{``Week 16''}. This paper shows
those invariants are trivial --- at least for compact manifolds, where
one just gets the answer ``1''. This seems to be happening a lot lately.

\begin{enumerate}
\def\labelenumi{\arabic{enumi})}
\setcounter{enumi}{7}
\tightlist
\item
  Boguslaw Broda, ``A note on the four-dimensional Kirby calculus'', 
  available as
  \href{https://arxiv.org/abs/hep-th/9305101}{\texttt{hep-th/9305101}}.
\end{enumerate}
\noindent
An earlier attempt by Broda to construct 4-manifold invariants along
similar lines was discussed here in \protect\hyperlink{week9}{``Week
9''} and \protect\hyperlink{week10}{``Week 10''} --- the upshot being
that the invariant was trivial. This is a new attempt and Broda has told
me that the arguments for the earlier invariant being trivial do not
apply. Here's hoping!

\begin{enumerate}
\def\labelenumi{\arabic{enumi})}
\setcounter{enumi}{8}
\tightlist
\item
  H.-J.\ Matschull, ``Solutions to the Wheeler DeWitt constraint of canonical gravity
  coupled to scalar matter fields'', available as
  \href{https://arxiv.org/abs/gr-qc/9305025}{\texttt{gr-qc/9305025}}.
\end{enumerate}
\noindent
One very important technical issue in the loop representation of quantum
gravity is how to introduce matter fields into the picture. Let me quote:

\begin{quote}
It is shown that the Wheeler DeWitt constraint of canonical gravity
coupled to Klein Gordon scalar fields and expressed in terms of
Ashtekar's variables admits formal solutions which are parametrized by
loops in the three dimensional hypersurface and which are extensions of
the well known Wilson loop solutions found by Jacobson, Rovelli and
Smolin.
\end{quote}

\begin{enumerate}
\def\labelenumi{\arabic{enumi})}
\setcounter{enumi}{9}
\tightlist
\item
   Luis J.\ Garay, ``Hilbert space of wormholes'', available as
  \href{https://arxiv.org/abs/gr-qc/9306002}{\texttt{gr-qc/9306002}}.
\end{enumerate}
\noindent
I think I'll just quote the abstract on this one:

\begin{quote}
Wormhole boundary conditions for the Wheeler--DeWitt equation can be
derived from the path integral formulation. It is proposed that the
wormhole wave function must be square integrable in the maximal analytic
extension of minisuperspace. Quantum wormholes can be invested with a
Hilbert space structure, the inner product being naturally induced by
the minisuperspace metric, in which the Wheeler--DeWitt operator is
essentially self--adjoint. This provides us with a kind of probabilistic
interpretation. In particular, giant wormholes will give extremely small
contributions to any wormhole state. We also study the whole spectrum of
the Wheeler--DeWitt operator and its role in the calculation of Green's
functions and effective low energy interactions.
\end{quote}

\begin{enumerate}
\def\labelenumi{\arabic{enumi})}
\setcounter{enumi}{10}
\tightlist
\item
   Vipul Periwal, ``Chern--Simons theory as topological closed string'',
   available as
  \href{https://arxiv.org/abs/hep-th/9305115}{\texttt{hep-th/9305115}}.
\end{enumerate}
\noindent
Lately people have been getting interested in gauge theories that can be
interpreted as closed string field theories. I mentioned one recent
paper along these lines in \protect\hyperlink{week15}{``Week 15''},
which considers Yang--Mills in 2 dimensions. (This was not the first
paper to do so, I should emphasize.) A while back Witten wrote a paper
on Chern--Simons gauge theory in 3 dimensions as a background-free open
string field theory, but I was unable to understand it. This paper seems
conceptually simpler, although it uses some serious mathematics. I think
I might be able to understand it. It starts:

\begin{quote}
The perturbative expansion of any quantum field theory (qft) with fields
transforming in the adjoint representation of \(\mathrm{SU}(N)\) is a
topological expansion in surfaces, with \(N^{-2}\) playing the role of a
handle-counting parameter. For \(N\) large, one hopes that the dynamics
of the qft is approximated by the sum (albeit largely intractable) of
all planar diagrams. The topological classification of diagrams has
nothing a priori to do with approximating the dynamics with a theory of
strings evolving in spacetime.

Gross (see also refs\ldots) has shown recently that the large \(N\)
expansion does actually provide a way of associating a theory of strings
in QCD. Maps of two-dimensional string worldsheets into two-dimensional
spacetimes are necessarily somewhat constricted. What one would like is
a qft with fields transforming in the adjoint representation in
\(d > 2\), which is at the same time exactly solvable. One could then,
in principle, attempt to associate a theory of strings with such a qft
by exhibiting a `sum over connected surfaces' interpretation for the
free energy of the qft. There is no guaranty that such an association
will exist.
\end{quote}

The author argues that Chern--Simons theory is a ``rara avis among QFTs''
for which such an association exists. He takes the free energy for
\(\mathrm{SU}(N)\) Chern--Simons theory on \(S^3\), does a large-\(N\)
expansion on it, and shows that the coefficient of the \(N^{2-2g}\) term
is the (virtual) Euler characteristics of the moduli space of surfaces
with g handles. I wish I understood this better at a very pedestrian
level! E.g., is there some string-theoretic reason why one might expect
the free energy to be of this form? Anyway, then he considers \(T^3\),
and gets something related to surfaces with a single puncture in them.

\hypertarget{week18}{%
\section{September 11, 1993}\label{week18}}

I will be resuming this series of articles this fall, though perhaps not
at a rate of one \protect\hyperlink{week}{``Week''} per week, as I'll be
pretty busy. For those of you who haven't seen this series before, let
me explain. It's meant to be a guide to some papers, mostly in preprint
form, that I have found interesting. I should emphasize that it's an
utterly personal and biased selection --- if more people did this sort
of thing, we might get a fairer sample, but I'll be unashamed in
focussing on my own obsessions, which these days lean towards quantum
gravity, topological quantum field theories, knot theory, and the like.

Quite a pile of papers has built up over the summer, but I will start by
describing what I did over my summer vacation:

\begin{enumerate}
\def\labelenumi{\arabic{enumi})}
\tightlist
\item
  John Baez, ``Strings, loops, knots, and gauge fields'', 
  available as
	  \href{https://arxiv.org/abs/hep-th/9309067}{\texttt{hep-th/9309067}}.
\end{enumerate}
\noindent
When I tell layfolk that I'm working on the loop representation of
quantum gravity, and try to describe its relation to knot theory, I
usually say that in this approach one thinks of space, not as a smoothly
curved manifold (well, I try not to say ``manifold''), but as a bunch of
knots linked up with each other. If they have been exposed to physics
popularizations they will usually ask me at this point if I'm talking
about superstring theory. To which I used to respond, somewhat annoyed,
that no, it was quite different. Superstring theory, I explained, is a
grandiose ``theory of everything'' that tries to describe all known
forces and particles, and lots more, too, as being vibrating loops of
string hurling around in 349-dimensional space. (Well, maybe just 10, or
26.) It is a complicated mishmash of all previous failed approaches to
unifying gravity with the other forces: Yang--Mills theory, Kaluza-Klein
models, strings, and supersymmetry. (The last is a symmetry principle
that postulates for every particle another one, a mysterious
``superpartner,'' despite the fact that no such superpartners have been
seen.) And it has made no testable predictions as of yet. The loop
representation of quantum gravity, on the other hand, is a much more
conservative project. It simply attempts to use some new mathematics to
reconcile two theories which both seem true, but up to now have been as
immiscible as oil and water: quantum field theory, and general
relativity. If it works, it will still be only the first step towards
unifying gravity with the other forces. If the questioner has the gall
to ask if \emph{it} has made any testable predictions, I say that so far
it is essentially a mathematics project. On the one hand, here are
Einstein's equations; on the other hand, here are the rules of thumb for
``quantizing'' some equations. Is there a consistent and elegant way of
applying those rules to those equations? People have tried for 40 years
or so without real success, but quite possibly they just weren't being
clever enough, since the rules of thumb leave a lot of scope for
creativity. Then a physicist named Ashtekar came along and reformulated
Einstein's equations using some new variables (usually known by experts
as the ``new variables''). This made the equations look much more like
those that describe the other forces in physics. This led to renewed
hope that Einstein's equations might be consistently quantized after
all. Then physicists named Rovelli and Smolin , working with Ashtekar,
made yet \emph{another} change of variables, based on the new variables.
Rovelli and Smolin's variables were labelled by loops in space, so they
are called the loop variables. These loops are quite unlike strings,
since they are merely mathematical artifacts for playing with Einstein's
equations, not actual little \emph{objects} whizzing about. But using
them, Rovelli and Smolin were able to quantize Einstein's equations and
actually find a lot of solutions! However, they were making up a lot of
new mathematics as they went along, and, as usual in theoretical
physics, it wasn't 100\% rigorous (which, as we know, is like the
woman who could trace her descent from William the Conqueror with only
two gaps). So I, as a mathematician, got interested in this and am
trying to help out and see how much of this apparently wonderful
development is for real.

The odd thing is that there are a lot of mathematical connections
between string theory and the loop representation. Gradually, as time
went on, I became more and more convinced that maybe the layfolk were
right --- maybe the loop representation of quantum gravity really \emph{was}
string theory in disguise, or vice versa. This made me a little embarrassed
by how much I had been making fun of string theory. Still, it could be a
very good thing. On the one hand, the loop representation of quantum
gravity is much more well-motivated from basic physical principles than
string theory --- it's not as baroque --- a point I still adhere to. So
maybe one could use it to understand string theory a lot more clearly.
On the other hand, string theory really attempts to explain, not just
gravity, but a whole lot more --- so maybe it might help people see what
the loop representation of quantum gravity has to do with the other
forces and particles (if in fact it actually works).

I decided to write a paper about this, and as I did some research I was
intrigued to find more and more connections between the two approaches,
to the point where it is clear that while they are presently very
distinct, they come from the same root, historically speaking.

Here's what I wound up saying:

\begin{quote}
The notion of a deep relationship between string theories and gauge
theories is far from new. String theory first arose as a model of hadron
interactions. Unfortunately this theory had a number of undesirable
features; in particular, it predicted massless spin-2 particles. It was
soon supplanted by quantum chromodynamics (QCD), which models the strong
force by an \(\mathrm{SU}(3)\) Yang--Mills field. However, string models
continued to be popular as an approximation of the confining phase of
QCD. Two quarks in a meson, for example, can be thought of as connected
by a string-like flux tube in which the gauge field is concentrated,
while an excitation of the gauge field alone can be thought of as a
looped flux tube. This is essentially a modern reincarnation of
Faraday's notion of ``field lines,'' but it can be formalized using the
notion of Wilson loops. If A denotes a classical gauge field, or
connection, a Wilson loop is simply the trace of the holonomy of A
around a loop in space. If instead A denotes a quantized gauge field,
the Wilson loop may be reinterpreted as an operator on the Hilbert space
of states, and applying this operator to the vacuum state one obtains a
state in which the Yang--Mills analog of the electric field flows around
the loop.

In the late 1970's, Makeenko and Migdal, Nambu, Polyakov, and others
attempted to derive equations of string dynamics as an approximation to
the Yang--Mills equation, using Wilson loops. More recently, D. Gross and
others have been able to \emph{exactly} reformulate Yang--Mills theory in
\(2\)-dimensional spacetime as a string theory by writing an asymptotic
series for the vacuum expectation values of Wilson loops as a sum over
maps from surfaces (the string worldsheet) to spacetime. This
development raises the hope that other gauge theories might also be
isomorphic to string theories. For example, recent work by Witten and
Periwal suggests that Chern--Simons theory in 3 dimensions is also
equivalent to a string theory.

String theory eventually became popular as a theory of everything
because the massless spin-2 particles it predicted could be interpreted
as the gravitons one obtains by quantizing the spacetime metric
perturbatively about a fixed ``background'' metric. Since string theory
appears to avoid the renormalization problems in perturbative quantum
gravity, it is a strong candidate for a theory unifying gravity with the
other forces. However, while classical general relativity is an elegant
geometrical theory relying on no background structure for its
formulation, it has proved difficult to describe string theory along
these lines. Typically one begins with a fixed background structure and
writes down a string field theory in terms of this; only afterwards can
one investigate its background independence. The clarity of a manifestly
background-free approach to string theory would be highly desirable.

On the other hand, attempts to formulate Yang--Mills theory in terms of
Wilson loops eventually led to a full-fledged ``loop representation'' of
gauge theories, thanks to the work of Gambini, Trias, and others. After
Ashtekar formulated quantum gravity as a sort of gauge theory using the
``new variables,'' Rovelli and Smolin were able to use the loop
representation to study quantum gravity nonperturbatively in a
manifestly background-free formalism. While superficially quite
different from modern string theory, this approach to quantum gravity
has many points of similarity, thanks to its common origin. In
particular, it uses the device of Wilson loops to construct a space of
states consisting of ``multiloop invariants,'' which assign an amplitude
to any collection of loops in space. The resemblance of these states to
wavefunctions of a string field theory is striking. It is natural,
therefore, to ask whether the loop representation of quantum gravity
might be a string theory in disguise --- or vice versa.

The present paper does not attempt a definitive answer to this question.
Rather, we begin by describing a general framework relating gauge
theories and string theories, and then consider a variety of examples.
Our treatment of examples is also meant to serve as a review of
Yang--Mills theory in 2 dimensions and quantum gravity in 3 and 4
dimensions.
\end{quote}

I should add that the sort of string theory I talk about in this paper
is fairly crude compared to that which afficionados of the subject
usually concern themselves with. It treats strings only as maps from a
surface (the string worldsheet) into spacetime, and only cares about
such maps up to diffeomorphism, i.e., smooth change of coordinates. In
most modern string theory the string worldsheet is equipped with more
geometrical structure (a conformal structure) --- it looks locally like
the complex plane, so one can talk about holomorphic functions on it and
the like. This is why string theorists are always muttering about
conformal field theory. But the sort of string theory that Gross and
others (Taylor, Minahan, and Polychronakos, particularly) have been
using to describe 2d Yang--Mills theory does not require a conformal
structure on the string worldsheet, so it's at least \emph{possible}
that more interesting theories like 4d quantum gravity can be formulated
as string theories without reference to conformal structures. (Of
course, if one integrates over all conformal structures, that's a way of
referring to conformal structures without actually picking one.) I guess
I'm rambling on here a bit, but this is really the most mysterious point
as far as I'm concerned.

One hint of what might be going on is as follows. And here, I'm afraid,
I will be quite technical. As noted by Witten and formalized by Moore,
Seiberg, and Crane, a rational conformal field theory gives rise to a
particularly beautiful sort of category called a modular tensor
category. This contains, as it were, the barest essence of the theory.
Any modular tensor category gives rise in turn to a 3d topological
quantum field theory --- examples of which are Chern--Simons theory and
quantum gravity in 3 dimensions. And Crane and Frenkel have shown (or
perhaps it's fairer to say that if they ever finish their paper they
\emph{will} have shown) that the nicest modular tensor categories give
rise to braided tensor \(2\)-categories, which should, if there be
justice, give 4d topological quantum field theories. (For more
information on all these wonderful things --- which no doubt seem
utterly intimidating to the uninitiated --- check out previous ``This
Week's Finds.'') Quantum gravity in 4 dimensions is presumably something
roughly of this sort, if it exists. So what I'm hinting at, in brief, is
that a bunch of category theory may provide the links between
\emph{modern} string theory with its conformal fields and the loop
representation of quantum gravity. This is not as outre as it may
appear. The categories being discussed here are really just ways of
talking about \emph{symmetries} (see my stuff on
\href{http://math.ucr.edu/home/baez/categories.html}{categories} and
\href{http://math.ucr.edu/home/baez/symmetries.html}{symmetries} for
more on this). As usual in physics, the clearest way to grasp the
connection between two seemingly disparate problems is often by
recognizing that they have the same symmetries.

\hypertarget{week19}{%
\section{September 27, 1993}\label{week19}}

I will now start catching up on some of the papers that have accumulated
over the summer. This time I'll say a bit about recent developments in
quantum field theory and \(4\)-dimensional topology.

The quantum field theories that describe three of the forces of nature
(electromagnetic, strong and weak) depend for their formulation on a
fixed metric on spacetime --- that is, a way of measuring distance and
time. Indeed, it seems pretty close to being true that spacetime is
\(\mathbb{R}^4\), and that the ``interval'' between any two points in
\(4\)-dimensional space is given by the Minkowski metric
\[dt^2 -dx^2 -dy^2 -dz^2\] where \(dt\) is the change in the time, or
\(t\), coordinate, \(dx\) is the change in the spatial \(x\) coordinate,
and so on. However, it's apparently not quite true. In fact, the
presence of matter or energy distorts this metric a little, and the
effect of the resulting ``curvature of spacetime'' is perceived as
gravity. This is the basic idea of general relativity, which is nicely
illustrated by the way in which the presence of the sun bends starlight
that passes nearby.

Gravity is thus quite different from the other forces, at least to our
limited understanding. The other forces we have quantum theories of, and
these theories depend on a \emph{fixed} (that is, pre-given) metric. We
have no quantum theory of gravity yet, only a classical theory, and this
theory is precisely a set of equations describing a \emph{variable}
metric, that is, one dependent upon the state of the universe. These
are, of course, Einstein's equations.

In fact it is no coincidence that we have no quantum theory of gravity.
For most of the last 50 years or so physicists have been working very
hard at inventing and understanding quantum field theories that rely for
their formulation on a fixed metric. Indeed, physicists spent huge
amounts of effort trying to make a theory of quantum gravity along
essentially these lines! This is what one calls ``perturbative'' quantum
gravity. Here one says, ``Well, we know the metric isn't quite the
Minkowski metric, but it's awfully close, so we'll write it as the
Minkowski metric plus a small perturbation, derive equations for this
perturbation from Einstein's equations, and make a quantum field theory
based on \emph{those} equations.'' That way we could use the good old
Minkowski metric as a ``background metric'' and thus use all the methods
that work for other quantum field theories. This was awfully fishy from
the standpoint of \emph{elegance}, but if it had worked it might have
been a very good thing, and indeed we learned a lot from its failure to
work. Mainly, though, we learned that we need to bite the bullet and
figure out how to do quantum field theory without any background metric.

A recent big step was made when people (in particular Witten and Atiyah)
formulated the notion of a ``topological quantum field theory.'' This is
a precise list of properties one would like a quantum field theory
independent of any background metric to satisfy. A wish list, as it
were. One of the best-understood examples of such a ``TQFT'' is
Chern--Simons theory. This is a quantum field theory that makes sense in
3-dimensional spacetime, not 4d spacetime, so in a sense it has no shot
at being ``true.'' However, it connects up to honest 4d physics in some
very interesting ways, it serves as warmup for more serious physics yet
to come, AND it has done wonders for the study of topology.

It is also worth noting that one particular case of Chern--Simons theory
is equivalent to quantum gravity in 3d spacetime. Here I am being a bit
sloppy; there are various ways of doing quantum gravity in 3 dimensions
and they are not all equivalent, but the approach that relates to
Chern--Simons theory is, in my opinion, the nicest. This approach to 3d
quantum gravity was the advantage that it can also be described using a
``triangulation'' of spacetime. In other words, if we prefer the
discrete to the continuum, we can ``triangulate'' it, or cut it up into
tetrahedra, and formulate the theory solely in terms of this
triangulation. Of course, it's pretty common in numerical simulations to
approximate spacetime by a lattice or grid like this. What's amazing
here is that one gets \emph{exact} answers that are \emph{independent}
of the triangulation one picks. The idea for doing this goes back to
Ponzano and Regge, but it was all done quite rigorously for 3d quantum
gravity by Turaev and Viro just a few years ago. In particular, they
were able to show the 3d quantum gravity is a TQFT using only
triangulations, no ``continuum'' stuff.

It is tempting to try to do something like this for 4 dimensions. But it
is unlikely to be so simple. A number of people have recently tried to
construct 4d TQFTs copying tricks that worked in 3d. Some papers along
these lines that I have mentioned before are:

\begin{quote}
{\rm
Louis Crane and David Yetter,
``A categorical construction of 4d topological quantum field theories'', available as
\href{https://arxiv.org/abs/hep-th/9301062}{\texttt{hep-th/9301062}}. 
(\protect\hyperlink{week2}{Week 2})}
\end{quote}

\begin{quote}
{\rm   Boguslaw Broda, ``Surgical invariants of four-manifolds'',
available as
\href{https://arxiv.org/abs/hep-th/9302092}{\texttt{hep-th/9302092}}.
(\protect\hyperlink{week9}{Week 9} and \protect\hyperlink{week10}{Week 10})}
\end{quote}
\noindent
(I have listed which ``Week'' I discussed
these in case anyone wants to go back and check out some of the details.)

These papers ran into stiff opposition as soon as they came out! First
Ocneanu claimed that the Crane--Yetter construction was trivial, in the
sense that the number it associated to any compact \(4\)-dimensional
spacetime manifold was 1. (This number is called the partition function
of the quantum field theory, and having it be 1 for all spacetimes means
the theory is deadly dull.)

\begin{quote}
{\rm
Adrian Ocneanu, ``A note on simplicial dimension shifting'',
available as
\href{https://arxiv.org/abs/hep-th/9302028}{\texttt{hep-th/9302028}}.
(\protect\hyperlink{week5}{``Week 5''}) }
\end{quote}
\noindent
Crane and Yetter wrote a rebuttal noting that Ocneanu was not dealing
with quite the same theory:

\begin{quote}
{\rm David Yetter and Louis Crane, ``We are not stuck with gluing'', 
available as
\href{https://arxiv.org/abs/hep-th/9302118}{\texttt{hep-th/9302118}}. 
(\protect\hyperlink{week7}{Week 7})}
\end{quote}
\noindent
They also presented, at their conference this spring, calculations
showing that their partition function was not equal to 1 for certain
examples.

In my discussions of Broda's work I extensively quoted some
correspondence with Dan Ruberman, who showed that in Broda's original
construction, the partition function of a \(4\)-dimensional manifold was
just a function of its signature and possibly some Betti numbers ---
these being well-known invariants, it's not especially exciting from the
point of view of topology. This was also shown by Justin Roberts:

\begin{quote}
{\rm 
Justin Roberts, ``Skein theory and Turaev--Viro invariants'',
\emph{Topology}, \textbf{34}, 771--787.   Available as
\href{https://citeseerx.ist.psu.edu/viewdoc/summary?doi=10.1.1.138.8587}{https://citeseerx.ist.psu.edu/viewdoc/summary?doi=} \break
\href{https://citeseerx.ist.psu.edu/viewdoc/summary?doi=10.1.1.138.8587}{10.1.1.138.8587}.
 (\protect\hyperlink{week14}{Week 14})}
\end{quote}
\noindent
He suggested that the Crane--Yetter partition function was also a
function of the signature and Betti numbers, but did not check their
precise normalization conventions, and so did not quite \emph{prove}
this. However, more recently Crane and Yetter, together with Kauffman,
have shown this themselves:

\begin{enumerate}
\def\labelenumi{\arabic{enumi})}
\tightlist
\item
  Louis Crane, Louis H.\  Kauffman and David N.\ Yetter,  
  ``Evaluating the Crane--Yetter Invariant'', available as
  \href{https://arxiv.org/abs/hep-th/9309063}{\texttt{hep-th/9309063}}.
\end{enumerate}

\begin{quote}
Abstract: We provide an explicit formula for the invariant of
4-manifolds introduced by Crane and Yetter (in 
\href{https://arxiv.org/abs/hep-th/9301062}{\texttt{hep-th/9301062}}).
A consequence of our result is the existence of a combinatorial formula
for the signature of a 4-manifold in terms of local data from a
triangulation. Potential physical applications of our result exist in
light of the fact that the Crane--Yetter invariant is a rigorous version
of ideas of Ooguri on \(B \wedge F\) theory.
\end{quote}
\noindent
They also have shown that Broda's original construction, and also a
souped-up construction of his, give a partition function that depends
only on the signature:

\begin{enumerate}
\def\labelenumi{\arabic{enumi})}
\setcounter{enumi}{1}
\tightlist
\item
Louis Crane, Louis H. Kauffman and David N. Yetter, 
  ``On the Classicality of Broda's \(\mathrm{SU}(2)\) Invariants of
  4-Manifolds'', available as
  \href{https://arxiv.org/abs/hep-th/9309102}{\texttt{hep-th/9309102}}.
\end{enumerate}

\begin{quote}
Abstract: Recent work of Roberts has shown that the first surgical
4-manifold invariant of Broda and (up to an unspecified normalization
factor) the state-sum invariant arising from the TQFT of Crane--Yetter
are equivalent to the signature of the 4-manifold. Subsequently Broda
defined another surgical invariant in which the 1- and 2- handles are
treated differently. We use a refinement of Roberts' techniques
developed by the authors in \texttt{hep-th/9309063} to show that the
``improved'' surgical invariant of Broda also depends only on the
signature and Euler character.
\end{quote}

Now let me say just a little bit about what this episode might mean for
physics as well as mathematics. The key is the ``\(B \wedge F\)'' theory
alluded to above. This is a quantum field theory that makes sense in 4
dimensions. I have found that the nicest place to read about it is:

\begin{enumerate}
\def\labelenumi{\arabic{enumi})}
\setcounter{enumi}{2}
\tightlist
\item
 Gary Horowitz,  ``Exactly soluble diffeomorphism-invariant theories'', 
  \emph{Commun. Math. Phys.} \textbf{125} (1989), 417--437.
\end{enumerate}
\noindent
This theory is a kind of simplified version of 4d quantum gravity that
is a lot closer in \emph{character} to Chern--Simons theory. Like
Chern--Simons theory, there are no ``local degrees of freedom'' --- every
solution looks pretty much like every other one as long as we don't take
a big tour of space and notice that funny things happen when we go
around a noncontractible loop, which is the sort of thing that can only
exist if space has a nontrivial topology. 4d quantum gravity, on the
other hand, should have loads of local degrees of freedom --- the local
curving of spacetime!

What Crane and Yetter were dreaming of doing was constructing 4d quantum
gravity as a TQFT using triangulations of spacetime. What they really
did, it turns out, was to construct \(B\wedge F\) theory as a TQFT using
triangulations. (Broda constructed it another way.) On the one hand, the
simplicity of \(B\wedge F\) theory compared to honest-to-goodness 4d
quantum gravity makes it possible to understand it a lot better, and
calculate it out explicitly. On the other hand, \(B\wedge F\) theory is
so simple that it doesn't tell us much new about topology, at least not
the topology of \(4\)-dimensional manifolds per se. Via Donaldson theory
and the work of Kronheimer and Mrowka it's probably telling us a lot
about the topology of \(2\)-dimensional surfaces embedded in
\(4\)-dimensional manifolds --- but alas, I don't understand this stuff
very well yet!

Getting our hands on 4d quantum gravity as a TQFT along these lines is
still, therefore, an unfinished business. But we are, at last, able to
study some examples of 4d TQFTs and ponder more deeply what it means to
do quantum field theory without any background metric. The real thing
missing is local degrees of freedom. Without them, any model is really
just a ``toy model'' not much like physics as we know it. The loop
representation of quantum gravity has these local degrees of freedom (to
the extent that we understand the loop representation!), and so the
challenge (well, one challenge!) is to better relate it to what we know
about TQFTs.

\hypertarget{week20}{%
\section{October 2, 1993}\label{week20}}

I think I'll depart from my usual concerns this week and talk about a
book I'd been meaning to get my hands on for ages:

\begin{enumerate}
\def\labelenumi{\arabic{enumi})}
\tightlist
\item
  John H. Conway and Neil J. A. Sloane, \emph{Sphere Packings, Lattices
  and Groups}, second edition, Grundlehren der mathematischen
  Wissenschaften \textbf{290}, Springer, Berlin, 1993.
\end{enumerate}

This is a mind-boggling book. I have always regarded research in
combinatorics as a pleasure I must deny myself, for the study of the
\emph{discrete} presents endless beautiful tapestries in which one could
easily lose oneself for life, and I regard this as a kind of
self-indulgence when there is, after all, the physical universe out
there waiting to be understood. Of course, it's good that \emph{someone}
does combinatorics, since even the most obscure corners have a strong
tendency to become useful eventually --- and if they write books about
it, I can have my cake and eat it too, by \emph{reading} about it.
Conway and Sloane are two masters of combinatorics, and this book is
like a dessert tray piled so high with delicacies that it's hard to know
where to begin. Rather than attempt to describe it, let me simply show
you a few things I found in it. The book is 679 pages long, so what I'll
say is only the most minute sample!

Let's begin with sphere packings. Say one is stacking cannonballs on
one's lawn, which is a quaint custom now that cannons have been replaced
by far more horrible weapons. A nice way to do it is to lay out a
triangle of cannonballs thus:
\[\includegraphics[max width=0.65\linewidth]{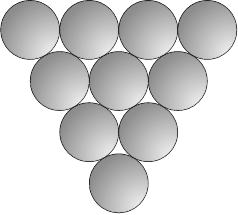}\]
and then set down another triangle of cannonballs on top of the first
layer, one in every other hole, thus:
\[\includegraphics[max width=0.65\linewidth]{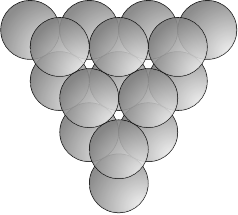}\]
and another:
\[\includegraphics[max width=0.65\linewidth]{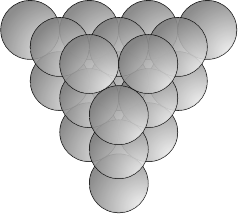}\]
and finally one on top --- which I won't attempt to draw.  
(Actually all the figures here were made by Simon Burton.)  Very nice
stack! There are many mathematical reflections it could lead one to, one
of which is: is this the \emph{densest} possible way one can pack
spheres of equal radius? The density is \(\pi/\sqrt{18} = .7405\ldots\);
can one do better? Apparently it was Kepler who first conjectured that
the answer is ``no.'' This is a famous hard problem. In 1991 W.-Y.
Hsiang announced that he had a proof, but I am not sure how many experts
have read it and been convinced.

This packing is a very regular sort of packing, what people call a
lattice packing. A lattice is a discrete subset of \(n\)-dimensional
Euclidean space closed under addition, and the packing above corresponds
to a lattice called the ``face-centered cubic'' or fcc lattice, which is
the set of all points \((x,y,z)\) where \(x\), \(y\), and \(z\) are
integers adding up to an \emph{even} number. (One might have a bit of
fun drawing some of those points and seeing why they do the job.)
Naturally, because of their regularity lattice packings are easier to
study than non-lattice ones. In fact, Gauss showed in 1831 that the fcc
lattice is the densest of all \emph{lattice} packings of spheres in 3
dimensions. This justifies the practice of stacking grapefruit this way
in the supermarket.

Let's take a bit closer look at what's going on here. Imagine an fcc
lattice going off infinitely in all directions\ldots{} each sphere is
touching 12 others: 6 in its own layer, as it were, 3 in the layer above
and 3 in the layer below. If we only pay attention to a given sphere and
those touching it we see something like:
\[\includegraphics[max width=0.65\linewidth]{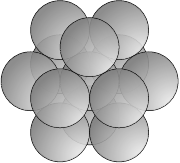}\]
If we remove the central sphere to clean up the picture a bit, the
centers of the remaining are the vertices of a shape called the
cuboctahedron:
\[\includegraphics[max width=0.65\linewidth]{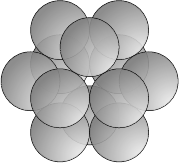}\]
It has 12 vertices, 8 triangular faces, and 6 square faces. One can get
it either by cutting off the corners of a cube just right, or by cutting
off the corners of an octahedron!

Let's return to thinking about the 12 spheres touching the central one.
This raises another question: is 12 the largest number of equal-radius
spheres able to touch a given one? This is the ``kissing spheres''
problem --- not quite the same as the packing problem! In 1694, Newton
conjectured that 12 was the largest number one could achieve. His
correspondent, David Gregory, thought 13 might be possible. It turns out
that Newton was right (some proofs appeared in the late 1800's), but it
is important to realize that Gregory's guess was not as dumb as it might
sound!

Why? Well, there's another way to get 12 spheres to touch a central one.
Namely, locate them at the vertices of a (regular) icosahedron. If one
wants to get to know the icosahedron a bit better one might read my
article \href{http://math.ucr.edu/home/baez/six.html}{Some Thoughts on
the Number Six}. But I hope you have an icosahedron available, or at
least a good mental image of one. It's easy to describe the vertices of
the icosahedron mathematically in terms of the golden ratio,
\[\frac{\sqrt{5} + 1}{2} = 1.61803398874989484820458683437\ldots\] This
usually goes by the name of ``\(\Phi\)'', or sometimes ``\(\tau\)''. The magic properties
of \(\Phi\) are too numerous to list here, but what counts is that one gets
the 12 vertices of a (regular) icosahedron by taking the points \[
  \begin{gathered}
    (\pm \Phi, \pm1, 0),
  \\(\pm1, 0, \pm \Phi),
  \\(0, \pm \Phi, \pm1).
  \end{gathered}
\] This is easy to check.

Now the \emph{interesting} thing is that when one gets 12 spheres to
touch a central one using the icosahedron, the 12 sphere don't touch
each other! There's room to move 'em around a bit, and perhaps (thought
Gregory) even enough room to stick in another one! Well, Gregory was
wrong, but one can do something pretty cool with this wiggle room.
First, though, let's check that there really \emph{is} a little space
between those outer spheres. First, compute the distance between
neighboring vertices of the icosahedron by taking two and working it
out: 
\[\|(\Phi,1,0) -(\Phi,-1,0)\| = 2\] 
Then, compute the distance from any
of the vertices to the origin, which is the center of the central
sphere: 
\[\|(\Phi,1,0)\| = \sqrt{\Phi^2 + 1}\] 
Using the charming fact that
\(\Phi^2 = \Phi + 1\) this simplifies to \(\sqrt{\Phi + 2}\), but then I guess we
just need to grind it out: 
\[\sqrt{\Phi + 2} =  1.902\ldots\] 
The point is
that this is less than 2, so two neighboring spheres surrounding the
central one are farther from each other than from the central one, i.e.,
they don't touch.

Now here's an interesting question: say we labelled the 12 spheres
touching the central one with numbers 1-12. Is there enough room to roll
these spheres around, always touching the central one, and permute the
spheres 1-12 in an interesting way? Notice that it's trivially easy to
roll them around in a way that amounts to rotating the icosahedron, but
are there more interesting permutations one can get?

Recall that the group of permutations of 12 things is called \(S_{12}\),
and it has \[12! = 479001600\] elements. The rotational symmetry group
of the icosahedron is much smaller. One can count its elements as
follows: pick a vertex and say which new vertex it gets rotated to ---
there are 12 possibilities --- and then note that there are 5 ways that
could have happened, for a total of 60. In fact (as I show in the
article \href{http://math.ucr.edu/home/baez/six.html}{``Six''}) the
rotational symmetry group of the icosahedron is \(A_5\), the group of
\emph{even} permutations of 5 things. So we have a nice embedding of the
group \(A_5\) into \(S_{12}\), but we are hoping that by some clever
wiggling we can get a bigger subgroup of \(S_{12}\) as the group of
permutations we can achieve by rolling spheres around.

In fact, one can achieve all \emph{even} permutations of the 12 spheres this
way! In other words, we get the group \(A_{12}\), with
\[\frac{12!}{2} = 239500800\] elements. How? Well, I wish I could draw
an icosahedron for you, but I can't really do so well on this medium.
The best way is to draw a top view:
\[\includegraphics[max width=0.65\linewidth]{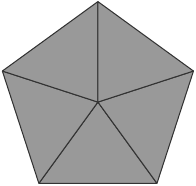}\]
and a bottom view:
\[\includegraphics[max width=0.65\linewidth]{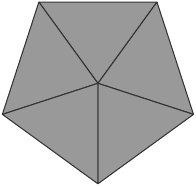}\]
Now, if we take the top 6 spheres and bunch them up so they are all
touching, and the bottom 6 spheres and bunch them up so \emph{they} all
touch, we can, in fact, twist the top 6 counterclockwise around \(1/5\)
of a turn. That is, we can map
\[\includegraphics[max width=0.65\linewidth]{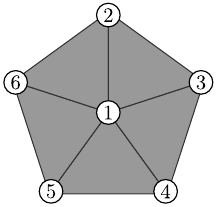}\]
to
\[\includegraphics[max width=0.65\linewidth]{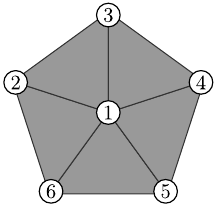}\]
I haven't actually \emph{done} this, or even proved I can, but Conway
and Sloane say so. And then the point is that all of \(A_{12}\) is
generated by ``twists'' of this form. Conway and Sloane give a
sophisticated and quick proof of this fact, which I can't resist
mentioning. Readers who don't know much group theory can skip the next
paragraph!

First, let \(t(i)\) be the element of \(S_{12}\) corresponding to the
clockwise twist about the \(i\)th vertex of the icosahedron (so that
what's drawn above is \(t(1)\)). This is a 5-cycle, and we need to show
these 12 5-cycles generate \(A_{12}\). Consider the subgroup generated
by elements of the form \(t(i)t(j)^{-1}\) --- a clockwise twist followed
by a counterclockwise twist. This is the Mathieu group \(M_{12}\), a
most remarkable group! In particular, its action on the vertices of the
icosahedron is able to map any 5 vertices to any 5 others (we say it's
quintuply transitive), so by conjugating \(t(i)\) with elements of the
Mathieu group we can get \emph{any} 5-cycle in \(S_{12}\). Then we use
the fact that \(A_{12}\) is generated by the 5-cycles.

Anyway, what this indicates is that there is an interesting relation
between the icosahedron and a certain finite group, the Mathieu group
\(M_{12}\). This group has \[12!/7! = 95040\] elements and it is a
``simple'' group, in the technical sense. The simple groups are to
finite groups roughly as the prime numbers are to the counting numbers;
that is, they are the elementary building blocks from which other finite
groups are made (although one has to specify \emph{how} one gloms them
together to get other groups). One of the remarkable achievements of
this century is the classification of these simple groups. In addition
to various infinite families of simple groups, like the alternating
groups \(\mathrm{A}_n\) (consisting of even permutations) there are a
finite number of ``sporadic'' simple groups such as the Mathieu groups,
the Fischer groups, the Suzuki groups, and, biggest of all, the Monster
group, which has
\[\begin{gathered}246\cdot320\cdot59\cdot76\cdot112\cdot133\cdot17\cdot19\cdot23\cdot29\cdot31\cdot41\cdot47\cdot59\cdot71 \\= 808017424794512875886459904961710757005754368000000000\end{gathered}\]
elements!

Here's a fun story about this number.

In the 1970s, the mathematicians Fricke, Ogg and Thompson were studying
the quotient of the hyperbolic plane by various subgroups of
\(\mathrm{SL}(2,\mathbb{R})\) --- the group of \(2\times2\) real
matrices with determinant one --- which acts as isometries of the
hyperbolic plane. Sitting inside \(\mathrm{SL}(2,\mathbb{R})\) is the
group of \(2\times2\) integer matrices with determinant one, called
\(\mathrm{SL}(2,\mathbb{Z})\). Sitting inside that is the group
\(\Gamma_0(p)\) consisting of matrices whose lower left corner is
congruent to zero \(\mod p\) for the prime \(p\). But Fricke, Ogg and
Thompson were actually considering a somewhat larger group
\(\Gamma_0(p)+\), which is the normalizer of \(\Gamma_0(p)\) inside
\(\mathrm{SL}(2,\mathbb{R})\).

If you don't know what this stuff means, don't worry! The point is that
they asked this question: if we take the quotient of the hyperbolic
plane by this group \(\Gamma_0(p)\)+, when does the resulting Riemann
surface have genus zero? And the \emph{real} point is that they found
the answer was: precisely when \(p\) is 2, 3, 5, 7, 11, 13, 17, 19, 23,
29, 31, 41, 47, 59 or 71.

Later, Ogg went to a talk on the Monster and noticed that these primes
were precisely the prime factors of the size of the Monster! He wrote a
paper offering a bottle of Jack Daniels whiskey to anyone who could
explain this fact. This was the beginning of a subject which Conway
dubbed ``Monstrous Moonshine'': the mysterious relation between the
Monster group, the group \(\mathrm{SL}(2,\mathbb{R})\), and Riemann
surfaces.

It turns out that lying behind Monstrous Moonshine is a certain string
theory having the Monster as symmetries, and this was the key to
understanding many strange ``coincidences''.

So, the sporadic groups are telling us something very deep and
mysterious about the universe, since they are very complicated and yet
somehow a basic, intrinsic part of the weave of mathematics. Conway and
Sloane have a lot to say about them and their relations to lattices and
error-correcting codes. For more about Monstrous Moonshine and string
theory, the reader should try this:

\begin{enumerate}
\def\labelenumi{\arabic{enumi})}
\setcounter{enumi}{1}
\tightlist
\item
  Igor Frenkel, James Lepowsky, and Arne Meurman, \emph{Vertex Operator
  Algebras and the Monster}, Academic Press, New York, 1988.
\end{enumerate}

Rather than attempt to describe this work, which I am not really
qualified to do (not that that usually stops me!), I think I will finish
up by describing a charming connection between the beloved icosahedron
and a lattice in 8 dimensions that goes by the name of \(\mathrm{E}_8\).
I'll be a bit more technical here.

The group of rotational symmetries of the icosahedron is, as we have
said, \(A_5\). This is a subgroup of the 3d rotation group
\(\mathrm{SO}(3)\). As all physicists know, whether they know it or not,
the group \(\mathrm{SO}(3)\) has the group \(\mathrm{SU}(2)\) of
\(2\times2\) unitary matrices with determinant \(1\) as its double
cover. So we can find a corresponding double cover of \(A_5\) as a
subgroup of \(\mathrm{SU}(2)\); this has twice as many elements as
\(A_5\), for a total of 120.

Now the group \(\mathrm{SU}(2)\) has a nice description as the group of
\emph{unit quaternions}, that is, things of the form
\[a + bI + cJ + dK\] where \(a,b,c,d\) are real numbers with
\(a^2 + b^2 + c^2 + d^2 = 1\), and \(I\),\(J\), and \(K\) satisfy
\[IJ = -JI = K, \quad JK = -KJ = I, \quad KI = -IK = J, \quad I^2 = J^2 = K^2 = -1\]
(Physicists in the 20th century usually use the Pauli matrices instead,
which are basically the same thing; for the relationship, read
\protect\hyperlink{week5}{``Week 5''}.)

It's natural to ask what the double cover of \(A_5\) looks like
explicitly in terms of the unit quaternions. Conway and Sloane give a
nice description. Let's write \((a,b,c,d)\) for \(a + bI + cJ + dK\),
write \(\Phi\) for the golden ratio as before, and \(\phi\) for the inverse of
the golden ratio:
\[\phi = \Phi^{-1} = \Phi -1 = 0.61803398874989484820458683437\ldots\] Then the
elements of the double cover of \(A_5\) are of the form \[
  \begin{gathered}
    (\pm 1, 0, 0, 0),
  \\(\pm \frac{1}{2}, \pm \frac{1}{2}, \pm \frac{1}{2}, \pm \frac{1}{2}),
  \\(0, \pm\frac{1}{2}, \pm \phi/2, \pm \Phi/2).
  \end{gathered}
\] and everything else that can be gotten by \emph{even} permutations of
the coordinates. (Check that there are 120 and that they are closed
under multiplication!)

Charming, but what does it have to do with \(\mathrm{E}_8\)? Well, note
that if we take all finite sums of elements of the double cover of
\(A_5\) we get a subring of the quaternions that Conway and Sloane calls
the ``icosians.'' Any icosian is of the form \[a + bI + cJ + dK\] where
\(a\),\(b\),\(c\), and \(d\) live in the ``golden field''
\(\mathbb{Q}(\Phi)\) --- this is the field of numbers of the form
\[x + \sqrt{5} y\] where \(x\) and \(y\) are rational. Thus we can think
of an icosian as an 8-tuple of rational numbers. We don't get all
8-tuples, however, but only those lying in a given lattice.

In fact, we can put a norm on the icosians as follows. First of all,
there is usual quaternionic norm
\[\|a + bI + cJ + dK\|^2 = a^2 + b^2 + c^2 + d^2\] But for an icosian
this is always of the form \(x + \sqrt{5} y\) for some rational \(x\)
and \(y\). It turns out we can define a new norm on the icosians by
setting 
\[|a + bI + cJ + dK|^2 = x + y.\] 
With respect to this norm, the
icosians form a lattice that fits isometrically in \(8\)-dimensional
Euclidean space and is the famous one called \(\mathrm{E}_8\)!
\(\mathrm{E}_8\) is known to yield the densest lattice packing of
spheres in 8 dimensions, a fact that is not only useful for
8-dimensional greengrocers, but also is apparently used in
error-correcting codes in a number of commercially available modems! (If
anyone knows \emph{which} modems use \(\mathrm{E}_8\), let me know --- I
might just buy one!) The density with which one can pack spheres in 8
dimensions using \(\mathrm{E}_8\), by the way, is \(\pi^4/384\), or
about \(.2537\).

Group theorists and some physicists, of course, will know that
\(\mathrm{E}_8\) is also the root lattice of the largest exceptional Lie
group, also known as \(\mathrm{E}_8\). This appears as the gauge group
in some string theories. While I find those string theories a bit
baroque for my taste, there is clearly a lot of marvelous mathematics
floating around here, and anyone who wants more might start with Conway
and Sloane.

\begin{center}\rule{0.5\linewidth}{0.5pt}\end{center}

\textbf{Addendum:} I had written:

\begin{quote}
Now here's an interesting question: say we labelled the 12 spheres
touching the central one with numbers 1-12. Is there enough room to roll
these spheres around, always touching the central one, and permute the
spheres 1-12 in an interesting way? Notice that it's trivially easy to
roll them around in a way that amounts to rotating the icosahedron, but
are there more interesting permutations one can get?

In fact, one can achieve all \textrm{even} permutations of the 12 spheres this
way!
\end{quote}
\noindent
I got some email from Conway saying that in fact one get \emph{all}
permutations of the 12 spheres. The point is this. One can start with
the 12 spheres touching the central one arranged so their centers are at
the vertices of an icosahedron. Then one can roll them around so their
centers lie at the vertices of a cuboctahedron:
\[\includegraphics[max width=0.65\linewidth]{wk20_fig5.pdf}\quad\includegraphics[max width=0.65\linewidth]{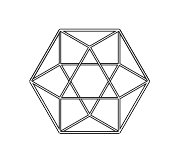}\]
At this point they touch each other, but one can indeed get to this
position. There is a nice picture of \emph{how} to do this in Conway and
Sloane's book, but you might enjoy figuring it out.

Now if we rotate the cuboctahedron around the axis pointing towards the
reader, we get an odd permutation of the 12 spheres, with cycle
structure \((1 2 3)(4 5 6)(7 8 9 10 11 12)\). So one can in fact get all
of \(S_{12}\) if one lets the 12 spheres ``just touch.'' If one thinks
about it a bit, the method described in the previous post gets all of
\(A_{12}\) without the 12 spheres touching at all. Conway says he
doesn't know if one can get all of \(S_{12}\) without the 12 spheres
touching each other at all. So this might be a fun problem to work on. I
bet the answer is ``no''.

Conway says the permutation problem came up before a lecture he gave at
the University of Pennsylvania a while ago, that he solved it during the lecture, and told
Jim Propp about it afterwards.

Someone with much better intuition about these things than I have might
want to consider similar ``rolling spheres permutation problems'' in
higher dimensions. E.g., what permutations can one can achieve in 4
dimensions, where it is possible to get 24 spheres to touch a central
one? (In 4d it is not known whether one can get 25 to touch, but 25 is
an upper bound.) Greg Kuperberg informs me that the known method for
having 24 spheres touch the central sphere is rigid, so that only the
``obvious'' permutations are possible starting from this arrangement.

The only other cases where the rolling spheres permutation problem
appears to be solved is in dimensions 1 and 2 (boring), 8, and 24. I
described the \(\mathrm{E}_8\) lattice in 8 dimensions in ``This Week's
Finds.'' In the corresponding lattice packing, each sphere touches 240
others. This is, up to rotation, the \emph{only} way to get 240 spheres
to touch a central one in 8 dimensions. (Also, one cannot get
\emph{more} than 240 to touch the central one.) So the subgroup of
\(S_{240}\) that one can get by rolling 240 spheres around a central one
is precisely the ``obvious'' subgroup, that is, the subgroup of
\(\mathrm{SO}(8)\) that preserves the \(\mathrm{E}_8\) lattice. This
turns out to be just the Weyl group of \(\mathrm{E}_8\), which has
\[2^{14}\cdot 3^5\cdot 5^2\cdot 7 = 696729600\] elements.

Dimension 24 is probably the most interesting for lattice theory, and
here the densest lattice packing is the Leech lattice. I am somewhat
sorry not to have even mentioned the Leech lattice in my article, since
this is the real star of Conway and Sloan's book. The Leech lattice is
probably related to the appearance of \(26\)-dimensional spacetime in
string theory; to get it, start with the unique even unimodular lattice
in 26-dimensional Minkowski space, and then look at \(M\), the set of
vectors in the lattice perpendicular to the vector
\(w = (70,1,2,3,\ldots,24)\). This is a null vector since
\[1^2 + 2^2 + \cdots + 24^2 = 70^2,\] i.e., it is perpendicular to itself.
Taking the quotient of \(M\) by \(w\) itself, we get the Leech lattice.
In this lattice each sphere touches 196560 others, which is the most one
can attain, and again this is the \emph{only} way to get 196560 spheres
to touch a central one in 24 dimensions. This should be obvious by
visualizing it. :-) So again the answer to the rolling spheres
permutation problem is the subgroup of \(\mathrm{SO}(24)\) that
preserves the Leech lattice. The isomorphism group of the Leech lattice
is an interesting group called \(\mathrm{Co}_0\) or \(.0\) (pronounced
``dotto''). It has
\[2^{22}\cdot 3^9\cdot 5^4\cdot 7^2\cdot 11\cdot 13\cdot 23 = 8315553613086720000\]
elements. The Monster group can also be produced using the Leech
lattice.

Please don't be fooled into thinking I understand this stuff!

Jim Buddenhagen raised another interesting point:

\begin{quote}
Since 13 unit spheres can't quite all touch a unit sphere, one may ask
how much bigger the central sphere must be to allow all 13 to touch.

In 1951 K.\ Schutte and B.\ L.\ van der Waerden found an arrangement of the
13 unit spheres that allows all of them to touch a central sphere of
radius \(r=1.04557\ldots\) This is thought to be optimal but has not
been proved optimal.

This \(r\) is an algebraic number and is a root of the polynomial
\[4096 x^{16} -18432 x^{12} +24576 x^{10} -13952 x^8 +4096 x^6 -608x^4 +32 x^2 +1\]
They computed \(r\) but did not publish the polynomial.
\end{quote}

\hypertarget{week21}{%
\section{October 10, 1993}\label{week21}}

Louis Kauffman is editing a series of volumes called ``Series on Knots
and Everything,'' published by World Scientific. The first volume was
his own book, \emph{Knots and Physics}.  Right now I'd like to talk about
the second volume, by Carter. I got to know Carter and Saito when it
started seeming that a deeper understanding of string theory and the
loop representation of quantum gravity might require understanding how
2-dimensional surfaces can be embedded in \(4\)-dimensional spacetime.
The study of this subject quickly leads into some very fascinating
algebra, such as the ``Zamolodchikov tetrahedron equations'' (which
first appeared in string theory). A nice review of this subject and
their work on it will appear in a while:

\begin{enumerate}
\def\labelenumi{\arabic{enumi})}
\tightlist
\item
 J.\ Scott Carter and Masahico Saito,  ``Knotted surfaces, braid movies, and beyond'', in \emph{Knots and Quantum Gravity}, ed.~John
  Baez, Oxford U.\ Press, Oxford, 1994.
\end{enumerate}
\noindent
but for the non-expert, a great way to get started is:

\begin{enumerate}
\def\labelenumi{\arabic{enumi})}
\setcounter{enumi}{1}
\tightlist
\item
  J.\ Scott Carter, \emph{How Surfaces Intersect in Space: An Introduction to Topology},
   World Scientific Press, Singapore, 1993.
\end{enumerate}

You can tell this isn't a run-of-the-mill introductory topology book as
soon as you read the little blurb about the author on the back
dustjacket. Occasionally there will be tantalizing personal details in
these blurbs that indicate that the author is not just a mathematical
automaton; for example, on the back of Hartshorne's famous text on
algebraic topology it says ``He has travelled widely, speake several
foreign languages, and is an experienced mountain climber. He is also an
accomplished amateur musician; has played the flute for many years, and
during his last visit to Kyoto, he began studying the shakuhachi.'' This
somehow fits with the austere and slightly intimidating quality of the
text itself. The tone of the blurb on the back of Scott Carter's book
could not be more different: ``When he is not drawing pictures, cooking,
or playing with Legos, he is writing songs and playing guitar for his
band The Anteaters who have recorded an eight-song cassette published by
Lobe Current Music.'' This is a book that invites the reader into
topology without taking itself too seriously.

I remember first reading about topology as the study of doughnuts,
M\"obius strips and the like, and then being in a way disappointed as an
undergrad -- although in another way quite excited -- when it seemed
that what topologists \emph{really} did was a lot of
``diagram-chasing,'' the algebraic technique widely used in homology and
homotopy theory. Once, however, as a grad student, I took a course in
``geometric topology'' by Tim Cochran, and was immensely pleased to find
that \emph{some} topologists really did draw wild pictures of
many-handled doughnuts and the like in 4 dimensions, and prove things by
sliding handles around. The nice thing about this book is that it is
readable by any undergraduate --- it doesn't assume or even mention the
definition of a topological space! --- but covers some very nontrivial
geometric topology. It is not a substitute for the usual introductory
course; instead, it concentrates on the study of surfaces embedded or
immersed in 3 and 4 dimensional space, and shows how much there is to
ponder about them. It is \emph{packed} with pictures and is lots of fun
to read.

The intrinsic topology of surfaces is very simple. The simplest one is
the sphere (by which, of course, mathematicians mean the \emph{surface}
of a ball, not the ball itself). The next is the torus, that is, the
surface of the doughnut. One can also think of the torus as what you get
by taking a square and gluing together the edges as below: \[
  \begin{tikzpicture}
    \draw[thick] (0,0) to node{$\blacktriangleright\blacktriangleright$} (2,0);
    \draw[thick] (0,-2) to node{$\blacktriangleright\blacktriangleright$} (2,-2);
    \draw[thick] (0,0) to node{$\blacktriangledown$} (0,-2);
    \draw[thick] (2,0) to node{$\blacktriangledown$} (2,-2);
  \end{tikzpicture}
\] gluing the two horizontal edges together so the double arrows match
up, and gluing the two vertical edges together so the single arrows
match up. There is also a two-handled torus, and so on. The number of
handles is called the ``genus.'' All these surfaces are orientable, that
is, one can define a consistent notion of ``right'' and ``left'' on
them, so that if one writes a little word on them and slides the word
around it'll never come back mirror-imaged. And in fact, all orientable
surfaces are just \(n\)-handled tori, so they are classified by their
genus.

A nice example of a nonorientable surface is the projective plane. One
way to visualize this is to take the surface of the sphere and
``identify'' opposite points, that is, decree them ``the same'' by fiat.
Imagine, for example, a globe in which antipodal points have been
identified. If one writes a word on the north pole and then slides it
down through the Americas to Ecuador, since the southern hemisphere has
been identified with the northern one, we can think of it popping out
over around India somewhere (sorry, my geography is a little rusty when
it comes to antipodes!), but we will see when we slide it back to the
north pole that it has been reversed, and is now written backwards! We
see from this not only that the projective plane is nonorientable, but
that it has another description: simply take a disc and identify
opposite points along the boundary. Since we're doing topology, a square
is just as good as a disc, so we can think of the projective plane as
the result of identifying the points on the boundary of a square as
follows: \[
  \begin{tikzpicture}
    \draw[thick] (0,0) to node{$\blacktriangleleft\blacktriangleleft$} (2,0);
    \draw[thick] (0,-2) to node{$\blacktriangleright\blacktriangleright$} (2,-2);
    \draw[thick] (0,0) to node{$\blacktriangledown$} (0,-2);
    \draw[thick] (2,0) to node{$\blacktriangle$} (2,-2);
  \end{tikzpicture}
\]

Another famous example of a nonorientable surface is the Klein bottle,
which is given by \[
  \begin{tikzpicture}
    \draw[thick] (0,0) to node{$\blacktriangleleft\blacktriangleleft$} (2,0);
    \draw[thick] (0,-2) to node{$\blacktriangleright\blacktriangleright$} (2,-2);
    \draw[thick] (0,0) to node{$\blacktriangledown$} (0,-2);
    \draw[thick] (2,0) to node{$\blacktriangledown$} (2,-2);
  \end{tikzpicture}
\]

We can take either the Klein bottle or the projective plane and get more
nonorientable surfaces by adding handles. Every nonorientable surface is
of one of these forms. I've included a few more basic facts about the
classification of surfaces as puzzles at the end of this article.

Now, the intrinsic topology of surfaces considers them as abstract
spaces in their own right, but the ``extrinsic topology'' of them
considers the ways they may be mapped into other spaces --- for example,
3- or \(4\)-dimensional Euclidean space. And here things get much more
interesting and subtle. For example, while one can embed any orientable
surface in 3d space, one cannot embed any of the nonorientable ones.
Here an embedding is a 1-1 continuous map. However, one can immerse the
non-orientable ones. An immersion is a map that is locally an embedding,
but not necessarily globally; e.g., a figure 8 is an immersion of the
circle in the plane. There's a standard way of immersing the Klein
bottle in 3d space with a circle of ``double points,'' that is, places
where the immersion is 2-1. One can easily turn this immersion into an
embedding of the Klein bottle in 4d space by representing the 4th
coordinate by how \emph{red} the surface is and having the Klein bottle
blush as it passes through itself. In fact, one can embed any surface
into 4d space.

While one can't embed the nonorientable surfaces in 3d space, it is
interesting to see how close one can come. The simplest way an immersion
can fail to be an embedding is by having double points. Another simple
way is to have triple points. Carter discusses a charming immersion of
the projective plane in 3d space that only has curves of double points
and a single triple point. This is known as ``Boy's surface.'' A
somewhat sneakier way immersions fail to be embeddings is by having
``branch points.'' Think, for example, of the function \(\sqrt{z}\) on
the complex plane. This is a two-valued function, so its graph consists
of two ``sheets'' which glom together in a funny way at \(z = 0\), the
branch point. Carter also talks about another neat immersion of the
projective plane in \(\mathbb{R}^3\) that just has double points and a
branch point --- the ``cross cap.'' Another immersion, the ``Roman
surface,'' has both triple points and a branch point.

The general question, then, is what sort of embeddings and immersions
different surfaces admit in 3 and 4 dimensions, and how to classify
these. If we are studying embeddings into 4 dimensions, a nice technique
is that of movies. Calling the 4th coordinate ``time,'' we can draw
slices at different times and get frames of a movie. Most of the frames
of a movie of an embedded surface will show simply a bunch of knots. At
a few times, however, a ``catastrophe'' will occur, e.g.: \[
  \begin{array}{p{8em}p{8em}p{8em}}
    \begin{tikzpicture}
      \begin{knot}
        \strand[thick] (0,0)
          to (0,-2);
        \strand[thick] (1,0)
          to (1,-2);
      \end{knot}
      \node at (0.5,-2.5) {Frame 1};
    \end{tikzpicture}
    &
    \begin{tikzpicture}
      \begin{knot}
        \strand[thick] (0,0)
          to[out=right,in=right,looseness=0.4] (0,-2);
        \strand[thick] (1,0)
          to[out=left,in=left,looseness=0.4] (1,-2);
      \end{knot}
      \node at (0.5,-2.5) {Frame 2};
    \end{tikzpicture}
    &
    \begin{tikzpicture}
      \begin{knot}
        \strand[thick] (0,0)
          to[out=right,in=right,looseness=0.7] (0,-2);
        \strand[thick] (1,0)
          to[out=left,in=left,looseness=0.7] (1,-2);
      \end{knot}
      \node at (0.5,-2.5) {Frame 3};
    \end{tikzpicture}
    \\[1em]
    \hspace*{-1.8em}
    \raisebox{-1.6em}{
    \begin{tikzpicture}
      \begin{knot}[clip width=0]
        \strand[thick] (0,0)
          to (1,-2);
        \strand[thick] (1,0)
          to (0,-2);
      \end{knot}
      \node at (0.5,-2.5) {Frame 4};
      \node at (0.5,-3) {(catastrophe!)};
    \end{tikzpicture}
    }
    &
    \begin{tikzpicture}
      \begin{knot}
        \strand[thick] (0,0)
          to[out=down,in=down,looseness=3] (1,0);
        \strand[thick] (0,-2)
          to[out=up,in=up,looseness=3] (1,-2);
      \end{knot}
      \node at (0.5,-2.5) {Frame 5};
    \end{tikzpicture}
    &
    \begin{tikzpicture}
      \begin{knot}
        \strand[thick] (0,0)
          to[out=down,in=down,looseness=1.5] (1,0);
        \strand[thick] (0,-2)
          to[out=up,in=up,looseness=1.5] (1,-2);
      \end{knot}
      \node at (0.5,-2.5) {Frame 6};
    \end{tikzpicture}
  \end{array}
\] However, there are always many different movies of essentially the
same embedding. We can, however, always relate these by a sequence of
transformations called ``movie moves.'' I wish I could draw these, but
it would take too long, so look at Carter's book!

And while you're at it, check out the index. You will enjoy finding the
excuses he has for such entries as ``hipster jive,'' ``math jail,''
``basket shaped thingy,'' and ``chocolate.'' Heck, I can't resist
one\ldots{} on page 81: ``Mathematicians use the term ``word" to mean any
finite sequence of letters or numbers. This practice can freak out
(disturb) people who are not hip to the lingo (aware of the
terminology)."

I should add that the following book also has a lot of interesting
pictures of surfaces in it:

\begin{enumerate}
\def\labelenumi{\arabic{enumi})}
\setcounter{enumi}{2}
\tightlist
\item
  George Francis, \emph{A Topological Picturebook}, Springer, Berlin,
  1987.
\end{enumerate}

\textbf{Problems:}

\vskip 1em
A. Take a projective plane and cut out a little disc. Show that what's
left is a M\"obius strip. 

B. Take two projective planes, cut out a little
disc from each one and attach them along the resulting circles. This is
called taking the ``connected sum'' of two projective planes. Show that
the result is a Klein bottle. In symbols, \(P + P = K\), or \(2P = K\).

C. Now take the connected sum of a projective plane and a Klein bottle.
Show that this is the same as a projective plane with a handle attached.
A projective plane with a handle attached is just the connected sum of a
projective plane and a torus, so we have: \(3P = P + K = P + T\). 

D. Show: \(4P = K + K = K + T\). 

E. Show: \((2n+1)P = P + nT\). 

F. Show: \((2n+2)P = K + nT\).

\hypertarget{week22}{%
\section{October 16, 1993}\label{week22}}

Lately I've been having fun in this series discussing some things that I
don't really know much about, like lattice packings of spheres. Next
week I'll get back to subjects that I actually know something about, but
today I want to talk about the 4-color theorem, the golden mean, the
silver root, knots and quantum field theory. I know a bit about
\emph{some} of these subjects, but I've only become interested in the
4-color theorem recently, thanks to my friend Bruce Smith, who has a
hobby of trying to prove it, and Louis Kauffman's recent work connecting
it to knot theory. The sources for what follows are:

\begin{enumerate}
\def\labelenumi{\arabic{enumi})}
\tightlist
\item
Thomas L. Saaty
  and Paul C. Kainen, 
  \emph{The Four-Color Problem: Assault and Conquest}, McGraw-Hill, 1977.
\end{enumerate}

and

\begin{enumerate}
\def\labelenumi{\arabic{enumi})}
\setcounter{enumi}{1}
\item
   Louis Kauffman, ``Map coloring and the vector cross product'',
  \emph{J.\ Comb.\ Theory} B \textbf{48} (1990), 45.

  Louis Kauffman, ``Map coloring, 1-deformed spin networks, and Turaev--Viro invariants
  for 3-manifolds'', \emph{Int.\ Jour.\ of Mod.\ Phys.}
  B, \textbf{6} (1992), 1765--1794.

  Louis Kauffman and H. Saleur,
  ``An algebraic approach to the planar colouring problem'', Yale University preprint \texttt{YCTP-P27-91},
  November 8, 1991.
\end{enumerate}
\noindent
(I discussed this work of Kauffman already in
\protect\hyperlink{week8}{``Week 8''}, where I described a way to
reformulate the 4-color theorem as a property of the vector cross
product.)

Where to start? Well, probably back in October, 1852. When Francis
Guthrie was coloring a map of England, he wondered whether it was always
possible to color maps with only 4 colors in such a way that no two
countries (or counties!) touching with a common stretch of boundary were
given the same color. Guthrie's brother passed the question on to De
Morgan, who passed it on to students and other mathematicians, and in
1878 Cayley publicized it in the \emph{Proceedings of the London Mathematical
Society}.

In just one year, Kempe was able to prove it. Whoops! In 1890 Heawood
found an error in Kempe's proof. And then the real fun starts\ldots.

But I don't want to tell the whole story leading up to how Appel and
Haken proved it in 1976 (with the help of a computer calculation
involving \(10^{10}\) operations and taking 1200 hours). I don't even
understand the structure of the Appel--Haken proof --- for that, one
should probably try:

\begin{enumerate}
\def\labelenumi{\arabic{enumi})}
\setcounter{enumi}{2}
\tightlist
\item
 Kenneth Appel and Wolfgang
  Haken,  \emph{Every Planar Map is Four Colorable}, Contemporary Mathematics
  \textbf{98}, American Mathematical Society, Providence, Rhode Island, 1989.
\end{enumerate}
\noindent
Instead, I'd like to talk about some tantalizing hints of relationships
between the 4-color theorem and physics!

First, note that to prove the 4-color theorem, it suffices to consider
the case where only three countries meet at any ``corner,'' since if
more meet, say four: 
\[
  \begin{tikzpicture}
    \draw[thick] (0,0) to (3,0);
    \draw[thick] (1.5,1.5) to (1.5,-1.5);
  \end{tikzpicture}
\] we can stick in a little country at each corner: \[
  \begin{tikzpicture}
    \draw[thick] (0,0) to (1,0);
    \draw[thick] (2,0) to (3,0);
    \draw[thick] (1.5,1.5) to (1.5,0.5);
    \draw[thick] (1.5,-0.5) to (1.5,-1.5);
    \draw[thick] (1,0) to (1.5,0.5) to (2,0) to (1.5,-0.5) to cycle;
  \end{tikzpicture}
\] so that now only three meet at each corner. If we can color the
resulting map, it's easy to check that the same coloring with the little
countries deleted gives a coloring of the original map.

Let us talk in the language of graph theory, calling the map a
``graph,'' the countries ``faces,'' their borders ``edges,'' and the
corners ``vertices.'' What we've basically shown is it suffices to
consider trivalent planar graphs without loops --- that is, graphs on
the plane that have three edges meeting at any vertex, and never have
both ends of the same edge incident to the same vertex.

Now, it's easy to see that 4-coloring the faces of such a graph is
equivalent to 3-coloring the \emph{edges} in such a way that no two
edges incident to the same vertex have the same color. For suppose we
have a 4-coloring of faces with colors \(1\), \(i\), \(j\), and \(k\).
Wait --- you say --- those don't look like colors, they look like the
quaternions. True! Now color each edge either \(i\), \(j\), or \(k\)
according to product of the colors of the two faces it is incident
to, where we define products by:
\[\begin{gathered}1i = i1 = i, \quad 1j = j1 = j, \quad 1k = k1 = k \\ ij = ji = k, \quad jk = kj = i, \quad ki = ik = j.\end{gathered}\]
These are \emph{almost} the rules for multiplying quaternions, but with
some minus signs missing. Since today (October 16th, 1993) is the 150th
birthday of the quaternions, I suppose I should remind the reader what
the right signs are:
\[ij = -ji = k, \quad jk = -kj = i, \quad  ki = -ik = j, \quad i^2 = j^2 = k^2 = -1.\]
Anyway, I leave it to the reader to check that this trick really gives
us a 3-coloring of the edges, and conversely that a 3-coloring of the
edges gives a 4-coloring of the faces.

So, we see that the edge-coloring formulation of the 4-color problem
points to some relation with the quaternions, or, pretty much the same
thing, the group \(\mathrm{SU}(2)\)! (For what \(\mathrm{SU}(2)\) has to
do with quaternions, see \protect\hyperlink{week5}{``Week 5''}.) Those
wrong signs look distressing, but in the following paper Penrose showed
they weren't really so bad:

\begin{enumerate}
\def\labelenumi{\arabic{enumi})}
\setcounter{enumi}{3}
\tightlist
\item
   Roger Penrose, ``Applications of negative dimensional tensors'', in
  \emph{Combinatorial Mathematics and its Applications}, ed.~D. J. A.
  Welsh, Academic Press, 1971.
\end{enumerate}

Namely, he showed one could count the number of ways to 3-color the
edges of a planar graph as follows. Consider all ways of labelling the
edges with the quaternions \(i\), \(j\), and \(k\). For each vertex,
take the product of the quaternions at the three incident edges in
counterclockwise order and then multiply by \(i\), getting either \(i\)
or \(-i\). Take the product of these plus-or-minus-i's over all vertices
of the graph. And THEN sum over all labellings!

This recipe may sound complicated, but only if you haven't ever studied
statistical mechanics of lattice systems. It's exactly the same as how
one computes the ``partition function'' of such a system --- the
partition function being the philosopher's stone of statistical
mechanics, since one can squeeze out so much information from it. (If we
could compute the partition function of water we could derive its
melting point.) To compute a partition one sums over states (labellings
of edges) the product of the exponentials of interaction energies
(corresponding to vertices). The statistical mechanics of
\(2\)-dimensional systems is closely connected to all sorts of nice
subjects like knot theory and quantum groups, so we should suspect
already that something interesting is going on here. It's especially
nice that Penrose's formula makes sense for arbitrary trivalent graphs
(although it does not count their 3-colorings unless they're planar),
and satisfies some juicy ``skein relations'' reminiscent of those
satisfied by the quantum group knot invariants. Namely, we can
recursively calculate Penrose's number for any trivalent graph using the
following three rules:

\begin{enumerate}
\def\labelenumi{\arabic{enumi})}
\item
  Wherever you see \[
    \begin{tikzpicture}
   \draw[thick] (0,0) to (1,-0.5) to (2,0);
   \draw[thick] (1,-0.5) to (1,-1.5);
   \draw[thick] (0,-2) to (1,-1.5) to (2,-2);
    \end{tikzpicture}
  \] you can replace it with \[
    \begin{tikzpicture}
   \draw[thick] (0,0) to (0,-2);
   \draw[thick] (1.5,0) to (1.5,-2);
    \end{tikzpicture}
    \raisebox{2.5em}{$\qquad-\qquad$}
    \begin{tikzpicture}
   \draw[thick] (1.5,0) to node[fill=white]{} (0,-2);
   \draw[thick] (0,0) to (1.5,-2);
    \end{tikzpicture}
  \] In other words, replace the problem of computing Penrose's number
  for the original graph by the problem computing the difference of the
  Penrose numbers for the two graphs with the above changes made. For
  knot theory fans I should emphasize that we are talking about abstract
  graphs here, not graphs in 3d space, so there's no real difference
  between an ``overcrossing'' and an ``undercrossing'' --- i.e., we
  could have said \[
    \begin{tikzpicture}
   \draw[thick] (0,0) to node[fill=white]{} (1.5,-2);
   \draw[thick] (1.5,0) to (0,-2);
    \end{tikzpicture}
  \] instead of \[
    \begin{tikzpicture}
   \draw[thick] (1.5,0) to node[fill=white]{} (0,-2);
   \draw[thick] (0,0) to (1.5,-2);
    \end{tikzpicture}
  \] above, and it wouldn't matter.
\item
  If you do this you will start getting weird loops that have NO
  vertices on them. You are allowed to dispose of such a loop if you
  correct for that by multiplying by 3. (This is not magic, this is just
  because there were 3 ways to color that loop!)
\item
  Finally, when you are down to the empty graph, use the rule that the
  empty graph equals 1.
\end{enumerate}

Greg Kuperberg pointed out to me that this is a case of the quantum
group knot invariant called the Yamada polynomal. This is associated to
the spin-\(1\) representation of the quantum group \(\mathrm{SU}(2)\),
and it is a polynomial in a variable \(q\) that represents \(e^\hbar\),
where \(\hbar\) is Planck's constant. But the ``Penrose number'' is just
the value at \(q = 1\) of the Yamada polynomial --- the ``classical
case'' when \(\hbar = 0\). This makes perfect sense if one knows about
quantum group knot invariants: the factor of 3 in rule B above comes
from the fact that the spin-\(1\) representation of \(\mathrm{SU}(2)\)
is \(3\)-dimensional; this representation is really just another way of
talking about the vector space spanned by the quaternions \(i\), \(j\),
and \(k\). Also, quantum group knot invariants fail to distinguish
between overcrossings and undercrossings when \(\hbar = 0\).

Now let me turn to a different but related issue. Consider the problem
of trying to color the \emph{vertices} of a graph with \(n\) colors in
such a way that no two vertices at opposite ends of any given edge have
the same color. Let \(P(n)\) denote the number of such \(n\)-colorings.
This turns out to be a polynomial in \(n\) --- it's not hard to see
using recursion relations similar to the skein relations above. It also
turns out that the 4-color theorem is equivalent to saying that the
vertices of any planar graph can be 4-colored. (To see this, just use
the idea of the ``dual graph'' of a graph --- the vertices of the one
being in 1-1 correspondence with the edges of the other.) So another way
to state the 4-color theorem is that for no planar graph does the
polynomial \(P(n)\) have a root at \(n = 4\).

\(P(n)\) is called the ``chromatic polynomial'' and has been intensively
investigated. One very curious thing is this. Remember the golden mean
\[\Phi = \frac{\sqrt{5} + 1}{2} = 1.61803398874989484820458683437\ldots?\]
Well, \(\Phi + 1\) is never a root of the chromatic polynomial of a graph!
(Unless the polynomial vanishes identically, which happens just when the
graph has loops.) The proof is not all that hard, and it's in Saaty and
Kainen's book. However --- and here's where things get \emph{really}
interesting --- in 1965, Hall, Siry and Vanderslice figured out the
chromatic polynomial of a truncated icosahedron. (This looks like a
soccer ball or buckyball.) They found that of the four real roots that
weren't integers, one agreed with \(\Phi + 1\) up to 8 decimal places! Of
course, here one might think the 5-fold symmetry of the situation was
secretly playing a role. But in 1966 Barri tabulated a bunch of
chromatic polynomials in her thesis, and in 1969 Berman and Tutte
noticed that most of them had a root that agreed with \(\Phi + 1\) up to at
least 5 decimal places.

This curious situation was at least partially explained by Tutte in
1970. He showed that for a triangular planar graph (that is, one all of
whose faces are triangles) with \(n\) vertices one has
\[|P(\Phi + 1)| \leqslant \Phi^{5-n}.\] This is apparently not a
\emph{complete} explanation, though, because the truncated icosahedron
is not triangular.

This is not an isolated freak curiosity, either! In 1974 Beraha
suggested checking out the behavior of chromatic polynomials at what are
now called the ``Beraha numbers'' \[B(n) = 4 \cos^2(\pi/n).\] These are
\[\begin{aligned}B(1) &= 4 \\ B(2) &= 0 \\ B(3) &= 1 \\ B(4) &= 2 \\ B(5) &= \Phi+1 \\ B(6) &= 3 \\ B(7) &= S \\ \end{aligned}\]
etc. Note by the way that \(B(n)\) approaches 4 as \(n\) approaches
\(\infty\). (What's \(S\), you ask? Well, folks call \(B(7)\) the
``silver root,'' a term I find most poetic and eagerly want to spread!
\[S =  3.246979603717467061050009768008479621265\ldots\] If anyone knows
charming properties of the silver root, I'd be interested.) Anyway, it
turns out that the roots of chromatic polynomials seem to cluster near
Beraha numbers. For example, the four nonintegral real roots of the
chromatic polynomial of the truncated icosahedron are awfully close to
\(B(5)\), \(B(7)\), \(B(8)\) and \(B(9)\). Beraha made the following
conjecture: let \(P_i\) be a sequence of chromatic polynomials of graphs
such whose number of vertices approaches \(\infty\) as \(i\to\infty\).
Suppose \(r_i\) is a real root of \(P_i\) and suppose the \(r_i\)
approach some number \(x\). Then \(x\) is a Beraha number.

In work in the late 60's and early 70's, Tutte proved some results
showing that there really was a deep connection between chromatic
polynomials and the Beraha numbers.

Well, to make a long story short (I'm getting tired), the Beraha numbers
\emph{also} have a lot to do with the quantum group \(\mathrm{SU}(2)\).
This actually goes back to some important work of Jones right before he
discovered the first of the quantum group knot polynomials, the Jones
polynomial. He found that -- pardon the jargon burst -- the Markov trace
on the Temperley--Lieb algebra is only nonnegative when the Markov
parameter is the reciprocal of a Beraha number or less than \(1/4\).
When the relationship of all this stuff to quantum groups became clear,
people realized that this was due to the special nature of quantum
groups when \(q\) is an \(n\)th root of unity (this winds up
corresponding to the Beraha number \(B(n)\)).

This all leads up to a paper that, unfortunately, I have not yet read,
in part because our library doesn't get this journal!

\begin{enumerate}
\def\labelenumi{\arabic{enumi})}
\setcounter{enumi}{4}
\tightlist
\item
  H.\ Saleur, ``Zeroes of chromatic polynomials: a new approach to the Beraha
  conjecture using quantum groups'', \emph{Commun.\ Math.\
  Phys.} \textbf{132} (1990) 657.
\end{enumerate}
\noindent
This apparently gives a ``physicist's proof'' of the Beraha conjecture,
and makes use of conformal field theory, that is, quantum field theory
in 2 dimensions that is invariant under conformal transformations.

I should say more: about what quantum groups have to do with conformal
field theory and knot polynomials, about the Kauffman/Saleur translation
of the 4-color theorem into a statement about the Temperley--Lieb
algebra, etc. But I won't! It's time for dinner. Next week, if all goes
according to plan, I'll move on to another puzzle in \(2\)-dimensional
topology --- the Andrews--Curtis conjecture --- and Frank Quinn's ideas on
tackling \emph{that} using quantum field theory.

\hypertarget{week23}{%
\section{October 24, 1993}\label{week23}}

I will soon revert to my older style, in which I list piles of new
papers as they accumulate on my desk. This time, though, I want to
describe Frank Quinn's work on the Andrews--Curtis conjecture using
topological quantum field theories (TQFTs), as promised. Then, if you'll
pardon me, I'll list the contents of a book I've just finished editing.
It is such a relief to be done that I cannot resist.

So ---

\begin{enumerate}
\def\labelenumi{\arabic{enumi})}
\item Frank Quinn,
  ``Topological quantum invariants and the Andrews--Curtis conjecture
  (Progress report)'', preprint, Sept.~1993.
\item
  Frank Quinn, ``Lectures on axiomatic topological quantum field theory'', to appear in the proceedings of the Park City Geometry
  Institute.
\item
  Wolfgang
  Metzler, ``On the Andrews--Curtis conjecture and related problems'', in \emph{Combinatorial Methods in Topology and Algebraic
  Geometry}, Contemporary Mathematics \textbf{44}, AMS, 1985.
\end{enumerate}

Last week I described --- in a pretty sketchy way --- how the 4-color
theorem and the Beraha conjecture are related to TQFTs. These can be
regarded as two very hard problems in \(2\)-dimensional topology --- one
solved by a mixture of cleverness and extreme brute force, the other
still open. There is another hard problem in \(2\)-dimensional topology
called the Andrews--Curtis conjecture, which Quinn is working on using
TQFT methods, which I'll talk about this time. I don't know too much
about this stuff, so I hope any experts out there will correct my
inevitable mistakes.

Actually, this conjecture is easiest to describe in a purely algebraic
way, so I'll start there. Hopefully most of you know the concept of a
``presentation'' of a group in terms of generators and relations. For
example, the group \(\mathbb{Z}_n\) (integers \(\mod n\)) has the
presentation \(\langle x \mid x^n\rangle\). This means, roughly, that we
form all products of the ``generator'' \(x\) and its inverse, and then
mod out by the ``relation'' \(x^n = 1\). A bit more interesting is the
dihedral group \(\mathrm{D}_n\) of symmetries of a regular \(n\)-gon,
counting rotations and reflections, with presentation
\(\langle x,y \mid x^n, y^2, (xy)^2 \rangle\). Here \(x\) corresponds to
a clockwise rotation by \((1/n)\)-th of a turn, and \(y\) corresponds to
a reflection.

A group always has lots of different presentations, so a natural problem
is to decide whether two different presentations give the same group
(or, strictly speaking, isomorphic groups). It'd be nice to have an
algorithm for deciding this question. But it's a famous result of
mathematical logic that there is no such algorithm!

If two presentations give the same group, one can get from one to the
other by a sequence of the following easy steps, called Tietze moves:

\begin{enumerate}
\def\labelenumi{\arabic{enumi})}
\tightlist
\item
  Throw in an extra new generator \(x\) together with the extra new
  relation \(xg^{-1}\) where \(g\) is a product of the previous
  generators and their inverses.
\item
  The inverse of 1) --- remove a generator \(x\) together with the
  relation \(xg^{-1}\), if possible (the relation \(xg^{-1}\) needs to
  be there!).
\item
  Throw in a new relation that's a consequence of existing relations.
\item
  The inverse of 3) --- remove a relation that's a consequence of other
  relations.
\end{enumerate}

So if one has two presentations and wants to see if they give same
group, you could always set up a program that blindly tries using these
Tietze moves in all possible ways to transform one presentation into the
other. If they are the same it'll eventually catch on! But if they're
not it'll chug on forever. There's no algorithmic way to tell
\emph{when} it should give up and admit the two presentations give
different groups! --- which is why we say there is no ``decision
procedure'' for this problem.

In one form, the Andrews--Curtis conjecture goes as follows. Remember
that the trivial group is the group with just the identity element; it
has a presentation \(\langle x \mid x \rangle\). Suppose we have some
other ``balanced'' presentation of the trivial group, that is a
presentation with just as many generators as relations:
\(\langle x_1,...,x_n \mid r_1,...,r_n \rangle\). Then the conjecture is
that it can be reduced to the presentation \(\langle x \mid x \rangle\)
by a sequence of the following moves that keep the presentation
balanced:

\begin{enumerate}
\def\labelenumi{\arabic{enumi})}
\tightlist
\item
  Throw in an extra new generator \(x\) together with the extra new
  relation \(x\).
\item
  The inverse of 1).
\item
  Permute the relations
\item
  Change \(r_1\) to \(r_1^{-1}\)
\item
  Change \(r_1\) to \(r_1r_2\)
\item
  Change \(r_1\) to \(gr_1g^{-1}\) for any \(g\).
\end{enumerate}

The experts seem to think this conjecture is probably false --- but
nobody has disproved it. Metzler lists a few presentations of the
trivial group that might be counterexamples: nobody has ever found a way
to use moves 1)-6) to boil them down to the presentation
\(\langle x \mid x \rangle\). For example,
\[\langle a, b \mid b^5a^{-4}, aba(bab)^{-1}\rangle.\] Try it!

The Andrews--Curtis conjecture is interesting mainly for its implications
in topology. When they first stated their conjecture they noted a number
of topological consequences, and the referee of the paper noted one
more. For example, it would shed some light on the Poincar\'e conjecture
(although not settle it) as follows. Recall that the Poincar\'e conjecture
says every \(3\)-dimensional manifold homotopic to a 3-sphere is
homeomorphic to a 3-sphere. The Andrews--Curtis conjecture implies that
if the Poincar\'e conjecture is false, any counterexample can in fact be
embedded (topologically) in \(\mathbb{R}^4\)!

It was the referee (does anyone know who that was?) who noted that the
Andrews--Curtis conjecture can be formulated in terms of ``CW
complexes.'' This is how Quinn thinks about it, so I suppose I should
say what those are.

A 0-complex is simply a set of points given the discrete topology. We
call the points ``0-cells.'' To get a 1-complex, we take a set of
``1-cells,'' that is, closed unit intervals, and glue their ends on to
the 0-cells in any way we want. In other words, we get a graph, possibly
with some edges having both ends at the same vertex. To get a 2-complex,
we take a set of ``2-cells,'' that is, \(2\)-dimensional closed disks,
and glue their boundaries onto our 1-complex by any continuous map. And
so on, with the ``\(n\)-cells'' being just copies of the closed unit
ball in \(\mathbb{R}^n\).

CW complexes were invented by J. H. C. Whitehead in 1949 and are a key
tool in algebraic topology. (The word ``CW,'' by the way, seems to come
from ``closure-finite'' and ``weak'' -- as in ``weak topology.'') They
are a nice class of topological spaces since on the one hand, being
built up by gluing simple pieces together, one can really understand
them, and on the other hand, they are actually quite general. In fact,
if one is interested in the usual invariants studied in algebraic
topology (homology and cohomology groups, homotopy groups and the like),
CW complexes are pretty much good enough. More precisely, Whitehead
proved a ``CW approximation theorem'' saying that any halfway decent
topological space (i.e., any ``compactly generated'' space) is ``weakly
homotopy equivalent'' to a CW complex. I won't burden you with the
definitions here; I learned this stuff once upon a time from

\begin{enumerate}
\def\labelenumi{\arabic{enumi})}
\setcounter{enumi}{3}
\tightlist
\item
   George W.\ Whitehead,  \emph{Elements of Homotopy Theory},
  Springer, Berlin, 1978. 
\end{enumerate}

Anyway, the Andrews--Curtis conjecture can be thought of as being about
2-complexes. In fact, a group presentation can be regarded as
instructions for building up a 2-complex --- start with a point, glue on
1-cells, one for each generator (obtaining a ``bouquet of circles'') and
then glue on 2-cells, one for each relation, attaching their boundaries
to the 1-cells in the manner prescribed by the relation. This 2-complex
will have fundamental group equal to the group given by the
presentation. The moves 1)-6) above can be thought of as operations on
these 2-complexes. So one can translate the Andrews--Curtis conjecture
into a statement about 2-complexes. And at this point I guess I'm going
to start getting more technical....

One topological statement of the Andrews--Curtis conjecture is that ``if
two 2-complexes are simply equivalent then one can be 2-deformed to the
other.'' I don't understand this as well as I want, so I won't explain
it; instead, I'll briefly explain the (weaker?) version corresponding
more closely to the algebraic statement above, namely ``if \(X\) is a
contractible 2-complex, it can be 2-deformed to a point.'' Being
``contractible'' means that as far as homotopy theory goes \(X\) is just
like a point. (E.g., the unit disk is contractible, while the circle is
not.) And a ``2-deformation'' roughly means a sequence of moves
consisting of adding or deleting 1-cells or 2-cells in a way that
doesn't affect things as far as homotopy goes, or doing homotopies of
attaching maps of 2-cells. The interesting thing about these
formulations of the Andrews--Curtis conjecture is that their analogs for
\(n > 2\) are true and in fact were shown by J.\ H.\ C.\ Whitehead in 1939!

Quinn's goal is to cook up invariants of 2-complexes that might detect
counterexamples to the Andrews--Curtis conjecture, i.e., invariants under
2-deformation. He wants to do it using 1+1-dimensional TQFTs of a sort
that assign vector spaces to 1-complexes and linear maps to 2-complexes.
Traditionally, TQFTs assign vector spaces to \(n\)-manifolds and linear
maps to \((n+1)\)-manifolds. Quinn calls his TQFTs ``modular'' because
they have a lot of formal similarities to the kind of TQFTs that come up
in string theory (where the modular group reigns supreme). He gives a
thorough axiomatic description of modular TQFTs in his lecture notes,
and this is actually the most fascinating aspect for me, more so than
the Andrews--Curtis conjecture per se, since it bears on physics.

The problem with coming up with an TQFT invariant that can catch
counterexamples to the Andrews--Curtis conjecture is an interesting
``stabilization'' property that 2-complexes have. Namely, if two
2-complexes are simply equivalent, one can can wedge them both with some
large number \(k\) of 2-spheres and get complexes which are 2-deformable
to each other. It turns out that this means we want to find a TQFT such
that \(Z(S^2)^k = 0\). And so Quinn considers TQFTs based, not on the
complex numbers, but on integers \(\mod p\).

A TQFT of his sort amounts to finding a symmetric tensor category of
vector spaces and an object A in this category with some special
properties corresponding to the fact that it is the vector space
corresponding to the unit interval \([0,1]\), which is the basic
1-complex from which one can build up more fancy ones. The kind of
category he uses has been described by:

\begin{enumerate}
\def\labelenumi{\arabic{enumi})}
\setcounter{enumi}{4}
\tightlist
\item
  Sergei Gelfand and David Kazhdan, ``Examples of tensor categories'',
  \emph{Invent.\ Math.} \textbf{109} (1992), 595--617.
\end{enumerate}

It is formed by starting with the category of representations of an
algebraic group in characteristic \(p\), and then making a semisimple
category out of this in a manner strongly reminiscent of what they do in
the theory of quantum groups at roots of unity. (See
\protect\hyperlink{week5}{``Week 5''} for a bit more about this.) The
object \(A\) is taken to be the sum of one copy of each irreducible
representation. (Again, this is strikingly reminiscent, and no doubt
based on, what occurs in the physics of the Wess--Zumino--Witten model,
where quantum groups at roots of unity play the role a finite group is
playing here.)

So, to round off a long story, Quinn and Ivelini Bobtcheva are currently
engaged in some rather massive computer calculations in order to
actually explicitly obtain the data necessary to calculate in the TQFTs
of this form. They have been looking at the groups \(\mathrm{SL}(2)\),
\(\mathrm{SL}(3)\), \(\mathrm{Sp}(4)\) and \(\mathrm{G}_2\) over
\(\mathbb{Z}_p\), where \(p\) is small (up to 19 for the
\(\mathrm{SL}(2)\) case). They are finding some interesting stuff just
by calculating the TQFT invariants of the 2-complexes corresponding to
the presentations \(\langle x \mid x^n \rangle\). (Note that \(n = 0\)
gives a space that's a wedge of a circle and \(S^2\), while \(n = 1\)
gives a disk.) Namely, they are finding periodicity in \(n\).

But they haven't found any counterexamples to the Andrews--Curtis
conjecture yet!

\begin{enumerate}
\def\labelenumi{\arabic{enumi})}
\setcounter{enumi}{5}
\tightlist
\item
  \emph{Knots and Quantum Gravity}, ed.~John Baez, Oxford University
  Press (to appear).
\end{enumerate}

This is the proceedings of a workshop held at U.C. Riverside; a large
percentage of the papers contain new results. Let me simply list them:

\begin{itemize}
\tightlist
\item
   Renate Loll, ``The loop formulation of gauge theory and gravity''.
\item
    Abhay Ashtekar and Jerzy Lewandowski, 
  ``Representation theory of Analytic Holonomy \(C^*\) Algebras''.

\item
  Rodolfo Gambini
  and Jorge Pullin , ``The Gauss linking number in quantum gravity'', 
  available as \href{https://www.arxiv.org/abs/gr-qc/9310025}{\texttt{gr-qc/9310025}}.
\item
  Louis Kauffman,``Vassiliev invariants and the loop states in quantum gravity'', 
  available as \href{https://arxiv.org/abs/gr-qc/9310035}{\texttt{gr-qc/9310035}}.
\item
  Steven Carlip, ``Geometric structures and loop variables in (2+1)-dimensional
  gravity'', available as \href{https://arxiv.org/abs/gr-qc/9309020}{\texttt{gr-qc/9309020}}. 
\item
  Dana S.\ Fine, ``From Chern--Simons to WZW via path integrals''.
\item
  Louis  Crane, ``Topological field theory as the key to quantum gravity'', 
  available as \href{https://arxiv.org/abs/hep-th/9308126}{\texttt{hep-th/9308126}}.
\item
  John Baez, ``Strings, loops, knots and gauge fields'', available as \href{https://arxiv.org/abs/hep-th/9309067}{\texttt{hep-th/9309067}}.
\item
  Paolo Cotta-Ramusino and Maurizio Martellini, ``BF theories and 2-knots''.
\item
   J. Scott Carter and
  Masahico Saito, ``Knotted surfaces, braid movies, and beyond''.
\end{itemize}

\hypertarget{week24}{%
\section{October 31, 1993}\label{week24}}

I will now revert to topics more directly connected to physics and start
catching up on the papers that have been accumulating. First, two very nice
review papers:

\begin{enumerate}
\def\labelenumi{\arabic{enumi})}
\tightlist
\item
  Chris Isham, ``Prima facie questions in quantum gravity'', lecture
  at Bad Honeff, September 1993, available as
  \href{https://arxiv.org/abs/gr-qc/9310031}{\texttt{gr-qc/9310031}}.
\end{enumerate}

If one wants to know why people make such a fuss about quantum gravity,
one could not do better than to start here. There are many approaches to
the project of reconciling quantum mechanics with gravity, all of them
rather technical, but here Isham focuses on the ``prima facie''
questions that present themselves no matter \emph{what} approach one
uses. He even explains why we should study quantum gravity --- a
nontrivial question, given how difficult it has been and how little
practical payoff there has been so far! Let me quote his answers and
urge you to read the rest of this paper:

\begin{quote}

\textbf{We must say something.} The value of the Planck length suggests
that quantum gravity should be quite irrelevant to, for example, atomic
physics. However, the non-renormalisability of the perturbative theory
means it is impossible to actually compute these corrections, even if
physical intuition suggests they will be minute. Furthermore, no
consistent theory is known in which the gravitational field is left
completely classical. Hence we are obliged to say \emph{something} about
quantum gravity, even if the final results will be negligible in all
normal physical domains.

\textbf{Gravitational singularities.} The classical theory of general
relativity is notorious for the existence of unavoidable spacetime
singularities. It has long been suggested that a quantum theory of
gravity might cure this disease by some sort of `quantum smearing'.

\textbf{Quantum cosmology.} A particularly interesting singularity is
that at the beginning of a cosmological model described by, say, a
Robertson--Walker metric. Classical physics breaks down here, but one of
the aims of quantum gravity has always been to describe the `origin' of
the universe as some type of quantum event.

\textbf{The end state of the Hawking radiation process.} One of the most
striking results involving general relativity and quantum theory is
undoubtedly Hawking's famous discovery of the quantum thermal radiation
produced by a black hole. Very little is known of the final fate of such
a system, and this is often taken to be another task for a quantum
theory of gravity.

\textbf{The unification of fundamental forces.} The weak and
electromagnetic forces are neatly unified in the Salam--Weinberg model,
and there has also been a partial unification with the strong force. It
is an attractive idea that a consistent quantum theory of gravity
\emph{must} include a unification of all the fundamental forces.

\textbf{The possibility of a radical change in basic physics.} The deep
incompatibilities between the basic structures of general relativity and
of quantum theory have lead many people to feel that the construction of
a consistent theory of quantum gravity requires a profound revision of
the most fundamental ideas of modern physics. The hope of securing such
a paradigm shift has always been a major reason for studying the
subject.

\end{quote}

\begin{enumerate}
\def\labelenumi{\arabic{enumi})}
\setcounter{enumi}{1}
\tightlist
\item
  Matthias Blau and George Thompson, ``Lectures on 2d gauge theories: topological aspects and     
  path integral techniques'', available as
  \href{https://arxiv.org/abs/hep-th/9310144}{\texttt{hep-th/9310144}}.
\end{enumerate}

Most of the basic laws of physics appear to be gauge theories. Gauge
theories are tricky to deal with because they are inherently nonlinear.
(At least the ``nonabelian'' ones are --- the main example of an abelian
gauge theory is Maxwell's equations.) People have been working hard for
quite some time trying to develop tools to study gauge theories on their
own terms, and \emph{one} reason for the interest in gauge theories in
2-dimensional spacetime is that life is simple enough in this case to
exactly solve the theories and see precisely what's going on. Another
reason is that in string theory one becomes interested in gauge fields
living on the 2-dimesional ``string worldsheet.''

This paper is a thorough review of two kinds of gauge theories in 2
dimensions: topological Yang--Mills theory (also called \(BF\) theory)
and the \(G/G\) gauged Wess--Zumino--Witten model. Both of these are of
great mathematical interest in addition to their physical relevance.
Studying the \(BF\) theory gives a way to do integrals on the moduli
space of flat connections on a bundle over a Riemann surface, while
studying the \(G/G\) model amounts to a geometric construction of the
categories of representations of quantum groups at roots of unity. (Take
my word for it, mathematicians find these important!)

I have found this review a bit rough going so far because the authors
like to use supersymmetry to study these models. But I will continue
digging in, since the authors consider the following topics (and I
quote): solution of Yang--Mills theory on arbitrary surfaces; calculation
of intersection numbers of moduli spaces of flat connections; coupling
of Yang--Mills theory to coadjoint orbits and intersection numbers of
moduli spaces of parabolic bundles; derivation of the Verlinde formula
from the \(G/G\) model; derivation of the shift \(k\) to \(k+h\) in the
\(G/G\) model via the index of the twisted Dolbeault complex.

\begin{enumerate}
\def\labelenumi{\arabic{enumi})}
\setcounter{enumi}{2}
\tightlist
\item
 J.\ W.\ Barrett and T.\ J.\ Foxon, ``Semi-classical limits of simplicial quantum gravity'', 
  available as
  \href{https://arxiv.org/abs/gr-qc/9310016}{\texttt{gr-qc/9310016}}.
\end{enumerate}

This paper looks at quantum gravity in 3 spacetime dimesions formulated
along the lines of Ponzano and Regge, that is, with the spacetime
manifold replaced by a bunch of tetrahedra (a ``simplicial complex''). I
describe some work along these lines in
\protect\hyperlink{week16}{``Week 16''}. Here the Feynman path integral
is replaced by a discrete sum over states, in which the edges of the
tetrahedra are assigned integer or half-integer lengths, which really
correspond to ``spins,'' and the formula for the action is given in
terms of \(6j\)-symbols. The authors look for stationary points of this
action and find that some correspond to Riemannian metrics and some
correspond to Lorentzian metrics. This is strongly reminiscent of Hartle
and Hawking's work on quantum cosmology,

\begin{enumerate}
\def\labelenumi{\arabic{enumi})}
\setcounter{enumi}{3}
\tightlist
\item
     J.\ B.\ Hartle and S.\ W.\ Hawking,
  ``Wave function of the universe'',
  \emph{Phys. Rev.} \textbf{D28} (1983), 2960.
\end{enumerate}

\noindent
in which there is both a Euclidean and a Lorentzian regime (providing a
most fascinating answer to the old question, ``what came before the Big
Bang?''). Here, however, the path integral is oscillatory in the Euclidean
regime and exponential in the Lorentzian one --- the opposite of what
Hartle and Hawking had. This puzzles me.

\begin{enumerate}
\def\labelenumi{\arabic{enumi})}
\setcounter{enumi}{4}
\tightlist
\item
  John Baez, ``Generalized measures in gauge theory'',
  \emph{Lett. Math. Phys.} \textbf{31} (1994), 213--223.  Also available as
  \href{https://arxiv.org/abs/hep-th/9310201}{\texttt{hep-th/9310201}}.
\end{enumerate}

Path integrals in gauge theory typically invoke the concept of Lebesgue
measure on the space of connections. This is roughly an
infinite-dimensional vector space, and there \emph{is} no ``Lebesgue
measure'' on an infinite-dimensional vector space. So what is going on?
Physicists are able to do calculations using this concept and get useful
answers --- mixed in with infinities that have to be carefully
``renormalized.'' Some of the infinities here are supposedly due to the
fact that one should really be working on the space of connections,
but on a quotient space, the connections modulo gauge transformations.
But not all the infinities are removed this way, and mathematically the
whole situation is enormously mysterious.

Recently Ashtekar, Isham, Lewandowski and myself have been looking at a
way to generalize the concept of measure, suggested by earlier work on
the ``loop representation'' of gauge theories. Ashtekar and Lewandowski
managed to rigorously construct a kind of ``generalized measure'' on the
space of connections modulo gauge transformations that acts formally
quite a bit like what one might hope for. In this paper I show how one can
define generalized measures directly on the space of connections. All of
these project down to generalized measures on the space of connections
modulo gauge transformations, but even when one is interested in
gauge-invariant quantities, it is sometimes easier to work ``upstairs.''
In particular, when the gauge group is compact, there is a ``uniform''
generalized measure on the space of connections that projects down to
the measure constructed by Ashtekar and Lewandowski. This generalized
measure is in some respects a rigorous substitute for the ill-defined
``Lebesgue measure,'' but it is actually built using Haar measure on
\(G\). I also define generalized measures on the group of gauge
transformations (which is an infinite-dimensional group), and when \(G\)
is compact I construct a natural example that is a rigorous substitute
for Haar measure on the group of gauge transformations.  As an
application of this ``generalized Haar measure'' I show that any
generalized measure on the space of connections can be averaged against
generalized Haar measure to give a gauge-invariant generalized measure
on the space of connections.

This doesn't, by the way, mean the problems I mentioned at the beginning
are solved!

\hypertarget{week25}{%
\section{November 14, 1993}\label{week25}}

Lately, many things give me the feeling that we're on the brink of some
deeper understanding of the relations between geometry, topology, and
category theory. It is very tantalizing to see the array of clues
pointing towards the fact that many seemingly disparate mathematical
phenomena are aspects of some underlying patterns that we don't really
understand yet. Louis Crane expressed it well when he said that it's as
if we are a bunch of archeologists digging away at different sites, and
are all starting to find different parts of the skeleton of some
gigantic prehistoric creature, the full extent of which is still
unclear.

I want to keep studying the following book until I understand it,
because I think it makes a lot of important connections\ldots{} pardon
the pun:

\begin{enumerate}
\def\labelenumi{\arabic{enumi})}
\tightlist
\item
   Jean-Luc Brylinski, \emph{Loop Spaces, Characteristic Classes and Geometric Quantization},
  Birkh\"auser, Basel, 1993.
\end{enumerate}
\noindent
The title of this book, while accurate, really does not convey the
\emph{novelty} of the ideas it contains. All three subjects listed have
been intensively studied by many people for at least several decades,
but Brylinski's book is not so much a summary of what is understood
about these subjects, as a plan to raise the subjects to a whole new
level.

I can't really describe the full contents of the book, since I haven't
had time to really absorb some of the most interesting parts, but let me
start by listing the contents, and then talk about it a bit.

\begin{enumerate}
\def\labelenumi{\arabic{enumi}.}
\tightlist
\item
  Complexes of Sheaves and Their Hypercohomology
\item
  Line Bundles and Geometric Quantization
\item
  K\"ahler Geometry of the Space of Knots
\item
  Degree 3 Cohomology --- The Dixmier-Douady Theory
\item
  Degree 3 Cohomology --- Sheaves of Groupoids
\item
  Line Bundles over Loop Spaces
\item
  The Dirac Monopole
\end{enumerate}
\noindent
It should be clear that while this is a very mathematical book, it is
informed by ideas from physics. As usual, the physical universe is
serving to goad mathematics to new heights!

The first two chapters are largely, but not entirely, ``standard''
material. I put the word in quotes because while Brylinski's treatment
of it starts with the basics --- the definition of sheaves, sheaf
cohomology, \v Cech cohomology, de Rham theory and the like --- even these
``basics'' are rather demanding, and the slope of the ascent is rather
steep. Really, the reader should already be fairly familiar with these
ideas, since Brylinski is mainly introducing them in order to describe a
remarkable generalization of them in the next chapters.

Let me quickly give a thumbnail sketch of the essential ideas behind
this ``standard'' material. In classical mechanics the main stage is the
phase space of a physical system. Points in this space represent
physical states; smooth functions on it represent observables. Time
evolution acts on this space as a one-parameter group of
diffeomorphisms. The remarkable fact is that time evolution is
determined by an observable, the Hamiltonian, or energy function, by
means of a geometric structure on phase space called a symplectic
structure. This is a nondegenerate closed \(2\)-form. The idea is that
the differential of the Hamiltonian is a \(1\)-form; since the
symplectic structure is nondegenerate it sets up an isomorphism of the
tangent and cotangent bundles of phase space, allowing us to turn the
differential of the Hamiltonian into a vector field; this vector field
generates the 1-parameter group of diffeomorphisms representing time
evolution; and by the magic of symplectic geometry, these
diffeomorphisms automatically preserve the symplectic structure.

This is the starting-point of the beautiful approach to quantum theory
known as geometric quantization, founded by Kostant in the early 1970's.
His first paper is still a good place to start:

\begin{enumerate}
\def\labelenumi{\arabic{enumi})}
\setcounter{enumi}{1}
\tightlist
\item
  Bertram Kostant, ``Quantization and unitary representations'', in
  \emph{Lectures in Modern Analysis and Applications III},
  Springer Lecture Notes in Mathematics \textbf{170} (1970),
  87--208.
\end{enumerate}

Here the idea is to construct a Hilbert space of states of the
\emph{quantum} system corresponding to the classical system, and turn
time evolution into a one-parameter group of unitary operators on this
Hilbert space. Extremely roughly, the idea is to first look at the space
of all \(L^2\) complex functions on phase space, and then use a
``polarization'' to cut down this ``prequantum'' Hilbert space to ``half
the size,'' by which one means something vaguely like how
\(L^2(\mathbb{R}^n)\) is ``half the size'' of \(L^2(\mathbb{R}^{2n})\)
--- this being the classic example. But in fact, it turns out one
doesn't really want to use \emph{functions} on phase space, but instead
sections of a certain complex line bundle. The point is that the
classification of line bundles fits in beautifully with symplectic
geometry. We can equip any line bundle with a hermitian connection; the
curvature of this connection is a closed \(2\)-form; this determines an
element of the 2nd cohomology of phase space called the first Chern
class. An important theorem says this class is necessarily an
\emph{integral} class, that is, it comes from an element of the 2nd
cohomology with integer coefficients; moreover, isomorphism classes of
line bundles over a manifold are in one-to-one correspondence with
elements of its 2nd cohomology with integer coefficients. The trick,
then, is to try to cook up a line bundle over phase space with a
connection whose curvature is the symplectic structure! This will be
possible precisely when the symplectic structure defines an integral
cohomology class. In fact, this integrality condition is nothing but the
old Bohr-Sommerfeld quantization condition dressed up in spiffy new
clothes (and made far more precise).

So: the moral I want to convey here is just that if the symplectic
structure on phase space defines an integral class in the 2nd cohomology
group, then we get a line bundle over phase space which helps us get
going with quantization. It then turns out that the one-parameter group
of diffeomorphisms defined by any Hamiltonian on phase space lifts to a
one-parameter group of transformations of this line bundle, which allows
us to get a unitary operator on the space of \(L^2\) sections of the
line bundle. This is not the end of the quantization story; one still
needs to chop down this ``prequantum'' space to half the size, etc.; but
let me leave off here.

What Brylinski wants to do is to find analogs of all these phenomena
involving the \textbf{third} cohomology groups of manifolds.

At first glance, this might seem to be a very artificial desire. Note
that importance of the \textbf{second} cohomology group in the above
story is twofold: 1) symplectic structures give elements of the second
cohomology, 2) the curvature of a connection gives an element of the
second cohomology, and in fact 2') line bundles are classified by
elements of second cohomology. None of these beautiful things seem to
have analogs in third cohomology! Of course, one can use the curvature
of a connection to get, not just the first Chern class, but higher Chern
classes. But the \(n\)th Chern class is an element of the \(2n\)th
cohomology group, so the odd cohomology groups don't play a major role
here. Of course, experts will immediately reply that there are also
Chern--Simons ``secondary characteristic classes'' that live in odd
cohomology, at least when one has a flat bundle around. And the same
experts will immediately guess that, because Chern--Simons theory has
been near the epicenter of the explosion of new mathematics relating
quantum groups, topological quantum field theories, conformal field
theory and all that stuff, I must be leading up to something along these
lines. Well, there \emph{must be} a relationship here, but
actually it is not emphasized in Brylinski's book! He takes a different
tack, as follows.

The basic point is that given a manifold \(M\), the space of loops in
\(M\), say \(LM\), is a space of great interest in its own right. It is
infinite-dimensional, but that should not deter us. When \(G\) is a Lie
group, \(LG\) is also a group (with pointwise operations); these are the
famous loop groups, which appear as groups of gauge transformations in
conformal field theory. When \(M\) is a \(3\)-dimensional manifold,
\(LM\) contains within it the space of all knots in \(M\); also, we may
think of \(LM\) as the configuration space for the simplest flavor of
string theory in the spacetime \(\mathbb{R} \times M\). Loops also serve to
define observables called ``Wilson loops'' in gauge theories, and these
are the basis of the loop representation of quantum gravity. So there is
a lot of interesting mathematics and physics to be found in the loop
space.

What does this have to do with the 3rd cohomology group of \(M\)? Well,
\(LM\) is a bundle over \(M\), so according to algebraic topology there
is a natural map from the 3rd cohomology of \(M\) to the 2nd cohomology
of \(LM\)! The ramifications of this are multiple.

First, every compact simple Lie group G has 3rd cohomology equal to
\(\mathbb{Z}\). (In fact, Brylinski notes that the cohomology group is
not merely isomorphic to \(\mathbb{Z}\), but canonically so --- and this
extra nuance turns out to be quite significant!) This gives rise to a
special element in the 2nd cohomology of \(LG\). This then gives a line
bundle over \(LG\). Alternatively, it gives a circle bundle over \(LG\),
in fact a central extension of \(LG\), that is, a bigger group
\(\widehat{LG}\) and an exact sequence
\[1 \to S^1 \to \widehat{LG} \to LG \to 1\] This group is called a
Kac--Moody group, and these are well-loved by string theorists since it
turns out that when one wants to quantize a gauge theory on the string
worldsheet (a kind of conformal field theory) one gets, not a
representation of the gauge group \(LG\) on the Hilbert space of quantum
states, but merely a projective representation, or in other words, a
representation of the central extension \(\widehat{LG}\). Brylinski also
notes that in some sense the canonical element in the 3rd cohomology of
\(G\) is responsible for the existence of quantum groups; this is
probably the deep reason for the association between quantum group
representations and Kac--Moody group representations, but, alas, this is
still quite murky to me.

Second, we can do better if we restrict ourself to knots (possibly with
nice self-intersections) rather than loops. Namely, given a 3-manifold
\(M\) equipped with a \(3\)-form, one gets, not just an element of the
2nd cohomology of \(LM\), but a symplectic structure on the space of
knots in \(M\), say \(KM\). It may seem odd to think of the space of
knots as a physical \emph{phase} space, but Brylinski shows that this
idea is related to the work of Marsden and Weinstein on ``vortex
filaments,'' an idealization of fluid dynamics in which all the fluid
motion is concentrated along some curves. Brylinski also notes that if
\(M\) is equipped with a Riemannian structure then \(KM\) inherits a
Riemannian structure (this is easy), and that if \(M\) has a conformal
structure \(KM\) has an almost complex structure. In fact, in the
Riemannian case all these structures on \(KM\) fit together to make it a
sort of K\"ahler manifold (although one must be careful, since the almost
complex structure is only integrable in a certain formal sense).
Brylinski hints that all this geometry may give a nice approach to the
study of knot invariants; I will have to look at the following papers
sometime:

\begin{enumerate}
\def\labelenumi{\arabic{enumi})}
\setcounter{enumi}{2}
\item
  M.\ Rasetti and T.\ Regge,``Vortices in He II, current algebras and quantum knots'', \emph{Physica} \textbf{80A} (1975) 217--233.
\item
   V. Penna and M. Spera, ``A geometric approach to quantum vortices'',
  \emph{Jour. Math. Phys.} \textbf{30} (1989), 2778--2784.
\end{enumerate}

However, Brylinski's real goal is something much more radical! The
beauty of 2nd cohomology is that integer classes in the 2nd cohomology
of \(M\) correspond to line bundles on \(M\); there is, in other words,
a very nice geometrical picture of 2nd cohomology classes. What is the
natural analog for 3rd cohomology? Instead of just working with \(LM\),
it would be nice to have some sort of geometrical objects on \(M\) that
correspond to integer classes in 3rd cohomology. What should they be?

Brylinski gives two answers, one in Chapter 4 and another in Chapter 5.
The first one, due mainly to Dixmier and Douady, is very appealing for a
quantum field theorist such as myself. Just as elements of
\(\mathrm{H}^2(M,\mathbb{Z})\) correspond to line bundles over \(M\), elements of
\(\mathrm{H}^3(M,\mathbb{Z})\) correspond to projective Hilbert space bundles
over \(M\)! Recall that in physics two vectors in a Hilbert space
correspond to the same physical state if one is a scalar multiple of the
other; the space of equivalence classes (starting with a
countable-dimensional Hilbert space) is what I'm calling ``projective
Hilbert space,'' and it is bundles of such rascals that correspond to
elements of \(\mathrm{H}^3(M,\mathbb{Z})\). The reason is roughly this: the
structure group \(G\) for such bundles is the group
\(\mathrm{Aut}(H)/\mathbb{C}^*\), that is, invertible operators on the Hilbert
space \(H\), modulo invertible complex numbers. In other words, we have
an exact sequence 
\[1 \to \mathbb{C}^* \to \mathrm{Aut}(H) \to G \to 1\] 
This gives an exact sequence of sheaves on \(M\), which, combined with the
marvelous fact that \(\mathrm{Aut}(H)\) is contractible, gives an isomorphism
between \(\mathrm{H}^1(M,\mathrm{sh}(G))\) (the cohomology of the sheaf of smooth
\(G\)-valued functions on \(M\)) and \(\mathrm{H}^2(M,\mathrm{sh}(\mathbb{C}^*))\).
But the latter is isomorphic to \(\mathrm{H}^3(M,\mathbb{Z})\).

Brylinski pushes the analogy to the line bundle case further by showing
how to realize the element of \(\mathrm{H}^3(M,\mathbb{Z})\) starting from a
connection on a projective Hilbert space bundle. But in Chapter 5 he
takes a more abstract approach that I want to sketch very vaguely, since
I don't understand it very well yet. This approach is exciting because
it connects to recent work on \(2\)-categories (and higher
\(n\)-categories), which I am convinced will play a role in unifying the
wild profusion of mathematics we are seeing in this tail end of the
twentieth century.

Here the best way to see the analogy to the line bundle case is through
\v Cech cohomology. Recall that we can patch a line bundle together by
covering our manifold \(M\) with charts \(O(i)\) and assigning to each
intersection \(O(i) \cap O(j)\) a \(\mathbb{C}^*\)-valued function
\(g_{ij}\). These ``transition functions'' must satisfy the compatibility
condition 
\[g_{ij} g_{jk} g_{ki} = 1\] 
We say then that the functions
\(g_{ij}\) define a 1-cocycle in \v Cech cohomology --- think of this as
just jargon, if you like. Note that we will get an isomorphic line
bundle if we take some \(\mathbb{C}^*\)-valued functions \(f_i\), one
on each chart \(O(i)\), and multiply \(g_{ij}\) by \(f_ if_j^{-1}\).
This simply amounts to changing the trivialization of the bundle on each
chart. We say that the new \v Cech cocycle differs by a coboundary. So line
bundles are in 1-1 correspondence with the 1st \v Cech cohomology with
values in \(\mathrm{sh}(\mathbb{C}^*)\). This turns out to be the same
thing as \(\mathrm{H}^2(M,\mathbb{Z})\), as noted above.

Now, there is a marvelous thing called a gerbe, which is like a bundle,
but is pieced together using \v Cech 2-cocyles! These will be classified by
the 2nd \v Cech cohomology with values in \(\mathrm{sh}(\mathbb{C}^*)\),
which is nothing but \(\mathrm{H}^3(M,\mathbb{Z})\).

What are these gerbes? Well, I wish I really understood them. Let me
just say what I know. The basic idea is to boost everything up a notch
using category-theoretic thinking. When we were getting ready to define
bundles, we needed to have the concept of a group at our disposal (to
have a structure group).  For gerbes, we need something called a category
of torsors. What is a group? Well, it is a \textbf{set} equipped with
various \textbf{maps} satisfying various properties. What is a category
of torsors? Well, it is a \textbf{category} equipped with various
\textbf{functors} satisfying utterly analogous properties. Note how we
are ``categorifying'' here. We have more structure, since while a set is
just a bunch of naked points, a category is a bunch of points, namely
objects, which are connected by arrows, namely morphisms. Given the
group \(\mathbb{C}^*\) we can get a corresponding category of torsors as
follows: the category of all manifolds with a simply transitive
\(\mathbb{C}^*\)-action (which are called torsors). A nice account of
why this category looks so much like a group appears in

\begin{enumerate}
\def\labelenumi{\arabic{enumi})}
\setcounter{enumi}{4}
\tightlist
\item
  Dan Freed, ``Higher algebraic structures and quantization'',
  preprint, available as
  \href{https://arxiv.org/abs/hep-th/9212115}{\texttt{hep-th/9212115}}.
\end{enumerate}
\noindent
which I already mentioned in \protect\hyperlink{week12}{``Week 12''}.

Just as a group can act on a set, a category of torsors can act on a
category. If we ``sheafify'' this notion, we get the concept of a gerbe.
Clear? Well, part of why I am interested in these ideas is the way they
make me a bit dizzy, so don't feel bad if you are a bit dizzy too now. I
really think that overcoming this dizziness will be necessary for
certain advances in mathematics and physics, though.

Instead of actually coming clean and defining the concept of a gerbe,
let me finish by saying what Brylinski does next. He defines an analog
of connections on bundles, called ``connective structures'' on gerbes.
And he defines an analog of the curvature, the ``curving'' of a
connective structure. This turns out to give an element of
\(\mathrm{H}^3(M,\mathbb{Z})\) in a natural way. He concludes in a blaze of glory
by showing how the Dirac monopole gives a gerbe on \(S^3\) whose curving
is the volume form. The integrality condition turns out to be related to
Dirac's original argument for quantization of electric charge. Whew!

To wrap up, let me note that the following paper, mentioned in
\protect\hyperlink{week23}{``Week 23''}, has shown up on \texttt{gr-qc}:

\begin{enumerate}
\def\labelenumi{\arabic{enumi})}
\setcounter{enumi}{5}
\tightlist
\item
  Abhay Ashtekar and Jerzy Lewandowski, ``Representation theory of analytic holonomy \(C^*\) Algebras'', in \emph{Knots and
  Quantum Gravity}, ed.~J. Baez, available as
  \href{https://arxiv.org/abs/gr-qc/9311010}{\texttt{gr-qc/9311010}}.
\end{enumerate}
\noindent
Ashtekar and Lewandowski are my friendly competitors in the business of making the loop representation of quantum gravity more rigorous by
formalizing the idea of a generalized measure on the space of
connections modulo gauge transformations.

\hypertarget{week26}{%
\section{November 21, 1993}\label{week26}}

I have been struggling to learn the rudiments of Teichm\"uller theory,
and it's almost time for me to face up to my ignorance of it by posting
a ``This Week's Finds'' attempting to explain the stuff, but I am going
to put off the inevitable and instead describe a variety of papers on
different subjects\ldots{}

\begin{enumerate}
\def\labelenumi{\arabic{enumi})}
\tightlist
\item
  Huw Price, ``Cosmology, time's arrow, and that old double standard'', in S.\ Savitt, ed., \emph{Time's Arrows Today},
  Cambridge University Press, 1994, pp.\ 66--94.  Available as
  \href{https://arxiv.org/abs/gr-qc/9310022}{\texttt{gr-qc/9310022}}.
\end{enumerate}

Why is the future different from the past? Because it hasn't happened
yet? Well, sure, but that's not especially enlightening, in fact, it's
downright circular. Unfortunately, a lot of work on the ``arrow of
time'' is just as circular, only so erudite that it is hard to spot it!
That's what this article takes some pains to clarify. I think I will be
lazy and quote the beginning of the paper:

\begin{quote}
A century or so ago, Ludwig Boltzmann and others attempted to explain
the temporal asymmetry of the second law of thermodynamics. The hard-won
lesson of that endeavour --- a lesson still commonly misunderstood ---
was that the real puzzle of thermodynamics lies not in the question why
entropy increases with time, but in that as to why it was ever so low in
the first place. To the extent that Boltzmann himself appreciated that
this was the real issue, the best suggestion he had to offer was that
the world as we know it is simply a product of a chance fluctuation into
a state of very low entropy. (His statistical treatment of
thermodynamics implied that although such states are extremely
improbable, they are bound to occur occasionally, if the universe lasts
a sufficiently long time.) This is a rather desperate solution to the
problem of temporal asymmetry, however, and one of the great
achievements of modern cosmology has been to offer us an alternative. It
now appears that temporal asymmetry is cosmological in origin, a
consequence of the fact that entropy is much lower than its theoretical
maximum in the region of the Big Bang --- i.e., in what we regard as the
\emph{early} stages of the universe.

The task of explaining temporal asymmetry thus becomes the task of
explaining this condition of the early universe. In this paper I want to
discuss some philosophical constraints on the search for such an
explanation. In particular, I want to show that cosmologists who discuss
these issues often make mistakes which are strikingly reminiscent of
those which plagued the nineteenth century discussions of the
statistical foundations of thermodynamics. The most common mistake is to
fail to recognise that certain crucial arguments are blind to temporal
direction, so that any conclusion they yield with respect to one
temporal direction must apply with equal force with respect to the
other. Thus writers on thermodynamics often failed to notice that the
statistical arguments concerned are inherently insensitive to temporal
direction, and hence unable to account for temporal asymmetry. And
writers who did notice this mistake commonly fell for another:
recognising the need to justify the double standard --- the application
of the arguments in question `towards the future' but not `towards the
past' --- they appealed to additional premises, without noticing that
in order to do the job, these additions must effectively embody the very
temporal asymmetry which was problematic in the first place. To assume
the uncorrelated nature of initial particle motions (or incoming
`external influences'), for example, is simply to move the problem from
one place to another. (It may \emph{look} less mysterious as a result,
but this is no real indication of progress. The fundamental lesson of
these endeavours is that much of what needs to be explained about
temporal asymmetry is so commonplace as to go almost unnoticed. In this
area more than most, folk intuition is a very poor guide to explanatory
priority.)

One of the main tasks of this paper is to show that mistakes of these
kinds are widespread in modern cosmology, even in the work of some of
the contemporary physicists who have been most concerned with the
problem of the cosmological basis of temporal asymmetry --- in the
course of the paper we shall encounter illicit applications of a
temporal double standard by Paul Davies, Stephen Hawking and Roger
Penrose, among others. Interdisciplinary point-scoring is not the
primary aim, of course: by drawing attention to these mistakes I hope to
clarify the issue as to what would count as adequate cosmological
explanation of temporal asymmetry.

I want to pay particular attention to the question as to whether it is
possible to explain why entropy is low near the Big Bang without thereby
demonstrating that it must be low near a Big Crunch, in the event that
the universe recollapses. The suggestion that entropy might be low at
both ends of the universe was made by Thomas Gold in the early 1960s.
With a few notable exceptions, cosmologists do not appear to have taken
Gold's hypothesis very seriously. Most appear to believe that it leads
to absurdities or inconsistencies of some kind. However, I want to show
that cosmologists interested in time asymmetry continue to fail to
appreciate how little scope there is for an explanation of the low
entropy Big Bang which does not commit us to the Gold universe. I also
want criticise some of the objections that are raised to the Gold view,
for these too often depend on a temporal double standard. And I want to
discuss, briefly and rather speculatively, some issues that arise if we
take the view seriously. (Could we observe a time-reversing future, for
example?)
\end{quote}

And now let me jump forward to a very interesting issue, Hawking's
attempt to derive the arrow of time from his ``no-boundary boundary
conditions'' choice of the wavefunction of the universe. (See
\protect\hyperlink{week3}{``Week 3''}.) I found this rather unsatisfying
when I read it, and had a sneaking suspicion that he was falling into
the fallacy Gold describes above. Let's hear what Price has to say about
it.

\begin{quote}
Our second example is better known, having been described in Stephen
Hawking's best seller, \emph{A Brief History of Time}. It is Hawking's proposal
to account for temporal asymmetry in terms of what he calls the No
Boundary Condition (NBC) --- a proposal concerning the quantum wave
function of the universe. To see what is puzzling about Hawking's claim,
let us keep in mind the basic dilemma. It seemed that provided we avoid
double standard fallacies, any argument for the smoothness of the
universe would apply at both ends or at neither. So our choices seemed
to be to accept the globally symmetric Gold universe, or to resign
ourselves to the fact that temporal asymmetry is not explicable (without
additional assumptions or boundary conditions) by a time-symmetric
physics. The dilemma is particularly acute for Hawking, because he has a
more reason than most to avoid resorting to additional boundary
conditions. They conflict with the spirit of his NBC, namely that one
restrict possible histories for the universe to those that `are finite
in extent but have no boundaries, edges, or singularities.'

Hawking tells us how initially he thought that this proposal favoured
the former horn of the above dilemma: `I thought at first that the no
boundary condition did indeed imply that disorder would decrease in the
contracting phase.' He changed his mind, however, in response to
objections from two colleagues: `I realized that I had made a mistake:
the no boundary condition implied that disorder would in fact continue
to increase during the contraction. The thermodynamic and psychological
arrows of time would not reverse when the universe begins to contract or
inside black holes.'

This change of mind enables Hawking to avoid the apparent difficulties
associated with reversing the thermodynamic arrow of time. What is not
clear is how he avoids the alternative difficulties associated with
\emph{not} reversing the thermodynamic arrow of time. That is, Hawking
does not explain how his proposal can imply that entropy is low near the
Big Bang, without equally implying that it is low near the Big Crunch.
The problem is to get a temporally asymmetric consequence from a
symmetric physical theory. Hawking suggests that he has done it, but
doesn't explain how. Readers are entitled to feel a little dissatisfied.
As it stands, Hawking's account reads a bit like a suicide verdict on a
man who has been stabbed in the back: not an impossible feat, perhaps,
but we'd like to know how it was done!

It seems to me that there are three possible resolutions of this
mystery. The first, obviously, is that Hawking has found a way round the
difficulty. The easiest way to get an idea of what he would have to have
established is to think of three classes of possible universes: those
which are smooth and ordered at both temporal extremities, those which
are ordered at one extremity but disordered at the other, and those
which are disordered at both extremities. If Hawking is right, then he
has found a way to exclude the last class, without thereby excluding the
second class. In other words, he has found a way to exclude disorder at
one temporal extremity of the universe, without excluding disorder at
both extremities. Why is this combination the important one? Because if
we can't exclude universes with disorder at both extremities, then we
haven't explained why our universe doesn't have disorder at both
extremities --- we know that it has order at least one temporal
extremity, namely the extremity we think of as at the beginning of time.
And if we do exclude disorder at both extremities, we are back to the
answer that Hawking gave up, namely that order will increase when the
universe contracts.

Has Hawking shown that the second class of universal histories, the
order-disorder universes, are overwhelmingly probable? It is important
to appreciate that this would not be incompatible with the underlying
temporal symmetry of the physical theories concerned. A symmetric
physical theory might be such that all or most of its possible
realisations were asymmetric. Thus Hawking might have succeeded in
showing that the NBC implies that any (or almost any) possible history
for the universe is of this globally asymmetric kind. If so, however,
then he hasn't yet explained to his lay readers how he managed it. In a
moment I'll describe my attempts to find a solution in Hawking's
technical papers. What seems clear is that it can't be done by
reflecting on the consequences of the NBC for the state of one temporal
extremity of the universe, considered in isolation. For if that worked
for the `initial' state it would also work for the `final' state; unless
of course the argument had illicitly assumed an objective distinction
between initial state and final state, and hence applied some constraint
to the former that it didn't apply to the latter. What Hawking needs is
a more general argument, to the effect that disorder-disorder universes
are impossible (or at least overwhelmingly improbable). It needs to be
shown that almost all possible universes have at least one ordered
temporal extremity --- or equivalently, at most one disordered
extremity. (As Hawking points out, it will then be quite legitimate to
invoke a weak anthropic argument to explain why we regard the ordered
extremity thus guaranteed as an \emph{initial} extremity. In virtue of
its consequences for temporal asymmetry elsewhere in the universe,
conscious observers are bound to regard this state of order as lying in
their past.)

That's the first possibility: Hawking has such an argument, but hasn't
told us what it is (probably because he doesn't see why it is so
important). As I see it, the other possibilities are that Hawking has
made one of two mistakes (neither of them the mistake he claims to have
made). Either his NBC does exclude disorder at both temporal extremities
of the universe, in which case his mistake was to change his mind about
contraction leading to decreasing entropy; or the proposal doesn't
exclude disorder at either temporal extremity of the universe, in which
case his mistake is to think that the NBC accounts for the low entropy
Big Bang.
\end{quote}

And, eventually, Price concludes that there is indeed something
lacking in Hawking's attempts to derive an arrow of time.

I have said this many times here and there, but I'll say it again. For a
good introduction to these issues, read:

\begin{enumerate}
\def\labelenumi{\arabic{enumi})}
\setcounter{enumi}{1}
\tightlist
\item
 H.\ D.\ Zeh, \emph{The Physical Basis of the Direction of Time},
  Second Edition, Springer, Berlin, 1992.
\end{enumerate}
\noindent
Zeh is one of the most clear-headed writers I know on this vexing
problem. Interestingly, Price acknowledges Zeh in his paper.

\begin{enumerate}
\def\labelenumi{\arabic{enumi})}
\setcounter{enumi}{2}
\tightlist
\item
  Renate Loll, ``Chromodynamics and gravity as theories on loop space'', available as
  \href{https://arxiv.org/abs/hep-th/9309056}{\texttt{hep-th/9309056}}.
\end{enumerate}

This is an especially thorough review of work on the loop representation
of gauge theories, especially the theories of the strong force and
gravity. A lot of work has been done on this subject but there are still
very many basic mathematical problems when it comes to making any of
this work rigorous, and one nice thing about Loll's work is that she is
just as eager to point out the problems as the accomplishments. It can
be dangerous when people become complacent and simply shrug off various
problems just because they are difficult and can be temporarily ignored.

\begin{enumerate}
\def\labelenumi{\arabic{enumi})}
\setcounter{enumi}{3}
\tightlist
\item
 Daniel Armand-Ugon, Rodolfo Gambini and Pablo Mora, ``Intersecting braids and intersecting knot theory'', available as
  \href{https://arxiv.org/abs/hep-th/9309136}{\texttt{hep-th/9309136}}.
\end{enumerate}
\noindent
There are a lot of hints that classical knot theory, which only
considers a circle smoothly embedded in space, is only the tip of a very
interesting iceberg. Namely, if one looks at the space of all loops,
this has the knots as an open dense subset, but then it has loops with a
single ``transverse double point'' like \[
  \begin{tikzpicture}
    \begin{knot}[clip width=0]
      \strand[thick] (0,0)
      to [out=down,in=up] (1,-1);
      \strand[thick] (1,0)
      to [out=down,in=up] (0,-1);
    \end{knot}
  \end{tikzpicture}
\] as a codimension 1 subset (like a hypersurface in the space of all
loops), and more fancy singularities appear as still smaller subsets, or
as the jargon has it, strata of higher codimension. The recent flurry of
work on Vassiliev invariants points out the importance of these other
strata --- or at least a few --- to knot theory. Namely, knot invariants
that extend nicely to knots with arbitrarily many transverse double
points include the famous quantum group knot invariants like the Jones
polynomial, and there is a close relationship between these ``Vassiliev
invariants'' and Lie algebra theory.

Meanwhile, the physicists have been forging ahead into more complicated
strata, motivated mainly by the loop representation of quantum gravity.
Gambini is one of the originators of the loop representation of gauge
theory, so it is not surprising that he is ahead of the game on this
business. Together with Pullin and Br\"ugmann he has been working on
extending the Jones polynomial, for example, to loops with various sorts
of self-intersections, calculating these extensions directly from the
path-integral formula for the Jones polynomial as a Wilson loop
expectation value in Chern--Simons theory. The relationship of this
extension to the theory of Vassiliev invariants was recently clarified
by Kauffman (see \protect\hyperlink{week23}{``Week 23''}), but there is
much more to do. Here Gambini and collaborators look at loops with
transverse triple points. I guess I'll just quote the abstract:

\begin{quote}
An extension of the Artin Braid Group with new operators that generate
double and triple intersections is considered. The extended Alexander
theorem, relating intersecting closed braids and intersecting knots is
proved for double and triple intersections, and a counter example is
given for the case of quadruple intersections. Intersecting knot
invariants are constructed via Markov traces defined on intersecting
braid algebra representations, and the extended Turaev representation is
discussed as an example. Possible applications of the formalism to
quantum gravity are discussed.
\end{quote}

\hypertarget{week27}{%
\section{December 16, 1993}\label{week27}}

This week I would like to describe some of the essays from the following
volume:

\begin{enumerate}
\def\labelenumi{\arabic{enumi})}
\tightlist
\item
  \emph{Conceptual Problems of Quantum Gravity}, eds. 
   Abhay Ashtekar and John Stachel, based on the proceedings of the 1988 Osgood
  Hill Conference, 15--19 May 1988, Birkh\"auser, Basel, 1991.
\end{enumerate}
\noindent
As the title indicates, this conference concentrated not on technical,
mathematical aspects of quantum gravity but on issues with a more
philosophical flavor. The proceedings make it clear how many problems we
still have in understanding how to fit quantum theory and gravity
together. Indeed, the book might be a bit depressing to those who
thought we were close to the ``theory of everything'' which some
optimists once assured us would be ready by the end of the millenium.
But to those like myself who enjoy the fact that there is so much left
to understand about the universe, this volume should be exciting (if
perhaps a bit daunting).

The talks have been divided into a number of groups:

\begin{itemize}
\tightlist
\item
  Quantum mechanics, measurement, and the universe
\item
  The issue of time in quantum gravity
\item
  Strings and gravity
\item
  Approaches to the quantization of gravity
\item
  Role of topology and black holes in quantum gravity
\end{itemize}

Let me describe a few of the talks, or at least their background, in
some detail rather than remaining general and vague.

\begin{enumerate}
\def\labelenumi{\arabic{enumi})}
\setcounter{enumi}{1}
\item
  Wojciech H.\ Zurek, ``Quantum measurements and the environment-induced transition from quantum to classical'', the volume above.

 W.\ G.\ Unruh, ``Loss of quantum coherence for a damped oscillator'', the volume above.
\end{enumerate}

These talks by Zurek and Unruh fit into what one might call the
``post-Everett school'' of research on the foundations of quantum
theory. To understand what Everett did, and what the post-Everett work
is about, you will need to be comfortable with the notions of pure
versus mixed states, and superpositions of states versus mixtures of
states (which are very different things). So, rather than discussing the
talks above, it probably makes more sense for me to talk about these
basic notions. A brief mathematical discussion appears below; one really
needs the clarity of mathematics to get anywhere with this sort of
issue. First, though, let me describe them vaguely in English.

In quantum theory, associated to any physical system there are states
and observables. An observable is a real-valued quantity we might
conceivably measure about the system. A state represents what we might
conceivably know about the system. The previous sentence is quite vague;
all it really means is this: given a state and an observable there is a
mathematical recipe that lets us calculate a probability distribution on
the real number line, which represents the probability of measuring the
observable to have a value lying in any subset of the real line. We call
this the probability distribution of the observable in the state. Using
this we can, if we want, calculate the mean of this probability
distribution (let us assume it exists!), which we call the expectation
value of the observable in the state.

Given two states \(\psi\) and \(\phi\), and a number \(c\) between \(0\)
and \(1\) there is a recipe for getting a new state, called
\(c \psi + (1-c)\phi\). This can be described roughly in words as
follows: ``with probability \(c\), the system is in state \(\psi\); with
probability \(1-c\) it is in state \(\phi\).'' This is called a
\textbf{mixture} of the states \(\psi\) and \(\phi\). If a state is a
mixture of two different states, with \(c\) not equal to \(0\) or \(1\),
we call that state a \textbf{mixed} state. If a state is not mixed it is
\textbf{pure}. Roughly speaking, a pure state is a state with as little
randomness as possible. (More precisely, it has as little entropy as
possible.)

All the remarks so far apply to classical mechanics as well as quantum
mechanics. A simple example from classical mechanics is a 6-sided die.
If we ignore everything about the die except which side is up, we can
say there are six pure states: the state in which the side of the die
showing one dot is up, the state in which the side showing two dots is
up, etc. Call these states 1,2,3,4,5, and 6. If it's a fair die, and we
roll it and don't look at it, the best state we can use to describe what
we know about the die is a mixed state which is a mixture: \(1/6\) of
state 1 plus \(1/6\) of state 2, etc. Note that if you peek at the die
and see that side 4 is actually up, you will be inclined to use a
different state to describe your knowledge: a pure state, state 4. Your
honest friend, who didn't peek, will still want to use a mixed state.
There is no contradiction here; the state simply is a way of keeping
track of what you know about the system, or more precisely, a device for
calculating expectation values of observables; which state you use
reflects your knowledge, and some people may know more than others.

Things get trickier in quantum mechanics. They also get trickier when
the system being described includes the person doing the describing.
They get even trickier when the system being described is the whole
universe -- for example, some people rebel at the thought that the
universe has ``many different states'' -- after all, it is how it is,
isn't it? (Gell-Mann gave a talk at this conference, which unfortunately
does not appear in this volume, entitled ``Quantum mechanics of this
specific universe.'' I have a hunch it deals with this issue, which
falls under the heading of ``quantum cosmology.'')

The first way things get trickier in quantum mechanics is that something
we are used to in classical mechanics fails. In classical mechanics,
pure states are always dispersion-free -- that is, for \emph{every}
observable, the probability measure assigned by the state to that
observable is a Dirac delta measure, that is, the observable has a 100\%
chance of being some specific value and a 0\% chance of having any other
value. (Consider the example of the dice, with the observable being the
number of dots on the face pointing up.) In quantum mechanics, pure
states need \emph{not} be dispersion-free. In fact, they usually aren't.

A second, subtler way things get trickier in quantum mechanics concerns
systems made of parts, or subsystems. Every observable of a subsystem is
automatically an observable for the whole system (but not all
observables of the whole system are of that form; some involve, say,
adding observables of two different subsystems). So every state of the
whole system gives rise to, or as we say, ``restricts to,'' a state of
each of its subsystems. In classical mechanics, pure states restrict to
pure states. For example, if our system consisted of 2 dice, a pure
state of the whole system would be something like ``the first die is in
state 2 and the second one is in state 5;'' this restricts to a pure
state for the first die (state 2) and a pure state for the second die
(state 5). In quantum mechanics, it is \emph{not} true that a pure state
of a system must restrict to a pure state of each subsystem.

It is this latter fact that gave rise to a whole bunch of quantum
puzzles such as the Einstein--Podolsky--Rosen puzzle and Bell's
inequality. And it is this last fact that makes things a bit tricky when
one of the two subsystems happens to be \emph{you}. It is possible, and
indeed very common, for the following thing to happen when two
subsystems interact as time passes. Say the whole system starts out in a
pure state which restricts to a pure state of each subsystem. After a
while, this need no longer be the case! Namely, if we solve
Schr\"odinger's equation to calculate the state of the system a while
later, it will necessarily still be a pure state (pure states of the
whole system evolve to pure states), but it need no longer restrict to
pure states of the two subsystems. If this happens, we say that the two
subsystems have become ``entangled.''

In fact, this is the sort of thing that often happens when one of the
systems is a measuring apparatus and the other is something measured.
Studying this issue, by the way, does \emph{not} require a general
definition of what counts as a ``measuring apparatus'' or a
``measurement'' --- on the contrary, this is exactly what is not needed,
and is probably impossible to attain. What is needed is a description in
quantum theory of a \emph{particular} kind of measuring apparatus,
possibly quite idealized, but hopefully reasonably realistic, so that we
can study what goes on using quantum mechanics and see what it actually
predicts will occur. For example, taking a very idealized case for
simplicity:

Our system consists of two subsystems, the ``detector'' and an
``electron.'' The systems starts out, let's suppose, in a pure state
which restricts to a pure state of each subsystem: the detector is
``ready to measure the electron's spin in the z direction'' and the
electron is in a state with its spin pointing along the x axis. After a
bit of time passes, if we restrict the state of the whole system to the
first subsystem, the detector, we get a mixed state like ``with 50\%
probability it has measured the spin to be up, and with 50\% probability
it has measured the spin to be down.'' Meanwhile, if we restrict the
state to the second subsystem, the electron is in the mixed state
``with 50\% change it has spin up, and with 50\% chance it has spin
down.'' In fact these two mixed states are \emph{correlated} in an
obvious sense. Namely, the observable of the \emph{whole} system that
equals 1 if the reading on the detector agrees with the spin of the
electron, and 0 otherwise, will have expectation value 1 (if the
detector is accurate). The catchy term ``entangled,'' which is a little
silly, really just refers to this correlation. I don't want to delve
into the math of correlations, but it is perhaps not surprising that, in
classical or quantum mechanics, interesting correlations can only occur
between subsystems if both of them are in mixed states. What's sneaky
about quantum mechanics is that the whole system can be in a pure state
which when restricted to each subsystem gives a mixed state, and that
these mixed states are then correlated (necessarily, as it turns out).
That's what ``entanglement'' is all about.

It was through analyses like this, but more detailed, that Everett
realized what was going on in a quantum system composed of two
subsystems, one of which was a measuring apparatus (or person, for that
matter), the other of which was something measured. The post-Everett
work amounts to refining Everett's analysis by looking at more realistic
examples, and more varied examples. In particular, it is interesting to
study situations where nothing very controlled like a scientific
``measurement'' is going on. For example, one subsystem might be an atom
in outer space, and the other subsystem might be its environment (a
bunch of other atoms or radiation). If one started out in a state which
restricted to a pure state of each subsystem, how fast would the
subsystems become entangled? And exactly \emph{how} would they become
entangled? -- this is very interesting. When we are doing a scientific
measurement, it's pretty clear what sort of correlation is involved in
the entanglement. In the above example, say, the detector reading is
becoming correlated to the electron's spin about the z axis. If all we
have is an atom floating about in space, it's not so clear. Can we think
of the environment as doing something analogous ``measuring'' something
about the atom, which establishes correlations of a particular kind?
This is the kind of thing Zurek and Unruh are studying.

In my description above I have tried to be very matter-of-fact, but
probably you all know that this subject is shrouded in mystery, largely
because of the misty and dramatic rhetoric people like to use, which
presumably makes it seem more profound. At least ``entangled'' has a
precise technical meaning. But anyone studying this subject will soon
run into ``collapse of the wavefunction,'' ``branches,'' ``the
many-worlds interpretation,'' the ``observer,'' and so on. These things
mean many things to many people, and nothing in particular to many more,
so one must always be on the alert.

Now for a little math to ground the above discussion. To keep life
simple suppose we have a quantum system described by a \(n\)-dimensional
Hilbert space \(H\) which we just think of as \(\mathbb{C}^n\),
\(n\)-dimensional complex space. The main thing to get straight is the
difference between superpositions and mixtures of quantum states. An
observable in quantum theory is described by a self-adjoint operator
\(A\), which for us is just an \(n\times n\) self-adjoint matrix. A
state is something that assigns to each observable a number called its
expectation value, in a manner that is 1) linear, 2) positive, and 3)
normalized. To explain this let us call our state \(\psi\). Linearity
means \(\psi(A + B) = \psi(A) + \psi(B)\) and \(\psi(cA) = c \psi(A)\)
for all observables \(A\),\(B\) and real numbers \(c\). Positivity means
\(\psi(A) > 0\) when \(A\) is a nonzero matrix that has non-negative
eigenvalues (a so-called non-negative matrix). And the normalization
condition is that \(\psi(1) = 1\).

This may seem unfamiliar, and that is because elementary quantum
mechanics only considers states of the form
\[\psi(A) = \langle v, Av\rangle\] where \(v\) is a unit vector in
\(H\). Not all states are of this form, but they are an extremely
important special class of states. It is also important to consider
states that are represented as ``density matrices,'' which are
non-negative matrices \(D\) with trace 1:
\[\operatorname{tr}(D) = \sum_i D_{ii} = 1\] Such a density matrix
defines a state \(\psi\) by \[\psi(A) = \operatorname{tr}(AD).\] It's
worth checking that this really meets the definition of a ``state''
given above!

The states corresponding to unit vectors in \(H\) are in fact a special
case of the density matrices. Namely, if \(v\) is a unit vector in \(H\)
we can let \(D\) be the self-adjoint matrix corresponding to projection
onto \(v\). I.e., the matrix \(D\) acts on any other vector, say \(w\),
by \[Dw = \langle v,w \rangle v.\] It's not to hard to check that the
matrix \(D\) really is a density matrix (do it!) and that this density
matrix defines the same state as does the vector \(v\), that is,
\[\operatorname{tr}(AD) = \langle v, Av\rangle\] for any observable
\(A\).

The entropy of a state \(\psi\) corresponding to the density matrix
\(D\) is defined to be \[S(\psi) = -\operatorname{tr}(D \ln D)\] where
one calculates \(D \ln D\) by working in a basis where \(D\) is diagonal
and replacing each eigenvalue \(x\) of \(D\) by the number \(x \ln x\),
which we decree to be \(0\) if \(x = 0\). Check that if \(D\)
corresponds to a \emph{pure} state as above then \(D \ln D = 0\) so the
entropy is zero.

Now about superpositions versus mixtures. They teach you how to take
superpositions in basic quantum mechanics. They usually don't tell you
about density matrices; all they teach you about is the states that
correspond to unit vectors in Hilbert space. Given two unit vectors in
H, one can take any linear combination of them and, if it's not zero,
normalize it to be a unit vector again, which we call a superposition.

Mixtures are an utterly different sort of linear combination. Given two
states \(\psi\) and \(\phi\) -- which recall are things that assign
numbers to observables in a linear way -- and given any number \(c\)
between \(0\) and \(1\), we can form a new state by taking
\[c \psi + (1-c) \phi\] This is called a mixture of \(\psi\) and
\(\phi\). Finally, some nontrivial exercises:

\textbf{Exercise:} Recall that a pure state is defined to be a state
which is not a mixture of two different states with \(0 < c < 1\). Show
that the states corresponding to unit vectors in Hilbert space are pure.

\textbf{Exercise:} Conversely, show (in the finite-dimensional case we
are considering) that all the pure states correspond to unit vectors in
Hilbert space.

\textbf{Exercise:} Show that every density matrix is a mixture of states
corresponding to unit vectors in Hilbert space.

\textbf{Exercise:} Show (in the finite-dimensional case we are considering) that
all states correspond to density matrices. Show that such a state is
pure if and only if its entropy is zero.

Well, this took longer than expected, so let me quickly say a bit more
about a few other papers in the conference proceedings\ldots.

\begin{enumerate}
\def\labelenumi{\arabic{enumi})}
\setcounter{enumi}{2}
\item
  Carlo Rovelli, ``Is there incompatibility between the ways time is treated in general
  relativity and in standard quantum mechanics?'', the
  volume above.

 Karel V.\ Kuchar, ``The problem of time in canonical quantization of relativistic systems'', the volume above.

  James B.\ Hartle, ``Time and prediction in quantum cosmology'', the
  volume above.

  Lee Smolin, ``Space and time in the quantum universe'', the volume
  above.
\end{enumerate}
\noindent
In the section on the problem of time in quantum gravity, these papers
in particular show a lively contrast between points of view. One nice
thing is that discussions after the papers were presented have been
transcribed; these make the disagreements even more clear. Let me simply
give some quotes that highlight the issues:

\begin{quote}
\textbf{Rovelli:} A \textbf{partial observable} is an operation on the
system that produces a number. But this number may be totally
unpredictable even if the state is perfectly known. Equivalently, this
number by itself may give no information on the state of the system
{[}in the Heisenberg picture --- jb{]}. For example, the reading of a
clock, or the value of a field, not knowing where and when it has been
measured, are partial observables.

A \textbf{true observable} or simply an \textbf{observable} is an
operation on the system that produces a number than can be predicted (or
whose probability distribution may be predicted) if the (Heisenberg)
state is known. Equivalently, it is an observable that gives information
about the state of the system.

\ldots.

Time is an experimental fact of nature, a very basic and general
experimental fact, but just an experimental fact. The formal development
of mechanics, and in particular Heisenberg quantum mechanics and the
presymplectic formulation of classical mechanics, suggests that it is
possible to give a coherent description of the world that is independent
of the presence of time.

\ldots.

From the mathematical point of view, \textbf{time} is a structure on the
set of observables (the foliation that I called a time structure).

From the physical point of view, time is the \emph{experimental fact}
that, in the nature as we see it, meaningful observables are always
constructed out of two partial observables. That is, it is the
experimental fact (not a priori required), that knowing the position of
a particle is meaningless unless we also know ``at what time'' a
particle was at that position.

In the formulation of the theory, this experimental fact is coded in the
time structure of the set of observables. If true observables are
composed of correspondences of partial observables, one of which is the
reading of a clock, then the set of true observables can be foliated
into one-parameter families that are given by the same partial
observables at different clock readings.

From an operational point of view, mechanics is perfectly well defined
in the absence of this time structure. It will describe a world (maybe
one slightly unfamiliar to us) in which observables are not arranged
along one-parameter lines, in which they have no such time structure (a
kind of fixed-time world), or have more complicated structures. We must
not confuse the psychological difficulties in visualizing such worlds
with their logical impossibility.

\ldots{} Heisenberg states, observables, measurement theory -- none of
these require time.

The notion of probability does not require time\ldots.

What I am proposing is that there may exist a coherent description of a
system in the framework of standard quantum mechanics even if it does
not have a standard ``time evolution.''

Why should we be interested in mechanics with no time structure? Because
general relativity \emph{is} a system (a classical system) with no time
structure. At least, it has no clearly defined time structure.

\ldots{} What we have to do is simple: ``forget time.''
\end{quote}

\begin{quote}
\textbf{Kuchar:} For myself, I want to see observables changing along my
world line and therefore associated with individual leaves of a
foliation. In that sense, the problem of time is shifted to the problem
of constructing an appropriate class of quantities one would like to call
observables. Now, what I would like to call observables probably differs
from what Carlo Rovelli would like to call observables. Carlo may like
to restrict that term to constants of motion, while I would like to use
variables that depend on a time hypersurface. Of course, both of us know
that there is a technical way of translating my observables into his
observables. However, it is difficult to subject such a translated
observable to an actual observation. In principle, of course, it does
not matter at what instant of time one measures a constant of motion.
But the constants of motion that are translations of my observables are
much too complicated when expressed in terms of the coordinates and the
momenta at the time of measurement. You thus have a hard time to design
an apparatus that would measure such a constant of motion at a time
different from the moment for which it was originally designed.
\end{quote}

\begin{quote}
\textbf{Smolin:} Now, as I discussed above, and as Jim Hartle argues at
length, there can be no strict implementation of the principle of
conservation of probability for a time that is the value of a
dynamical variable of a quantum system. Therefore, a sensible
measurement theory for quantum cosmology can only be constructed if
there is a time variable that is not a dynamical variable of the quantum
system that describes the universe.

Does this mean that quantum cosmology is impossible, since there is no
possibility of a clock outside of the system?

There is, as far as I know, exactly one loophole in this argument, which
is the one exploited by the program of intrinsic time. This is that one
coordinate on the phase space of general relativity might be singled out
and called time in such a way that the states, represented by functions
on the configuration space, could be read as time-dependent functions
over a reduced configuration space from which the privileged time
coordinate is excluded.
\end{quote}

\noindent Jokes:

\begin{quote}
Kuchar: Because Leibniz didn't believe in the ontological significance
of time, he dropped the letter ``t'' from his name.

Smolin: Is that true?

Kuchar: Yes! He spelled his name with a ``z''.

DeWitt: It's a good thing that he did believe in space because the ``z''
would've gone too.
\end{quote}
Of course, time is ``Zeit'' in German, which complicates things. 

\begin{enumerate}
\def\labelenumi{\arabic{enumi})}
\setcounter{enumi}{3}
\item
  Abhay Ashtekar, ``Old problems in the light of new variables'', the
  volume above.

  Carlo Rovelli, ``Loop representation in quantum gravity'', the
  volume above.

  Lee Smolin, ``Nonperturbative quantum gravity via the loop representation'',  the volume above.
\end{enumerate}

\noindent
These are a bit more technical papers that give nice introductions to
various aspects of the loop representation.

\hypertarget{week28}{%
\section{January 4, 1994}\label{week28}}

I think I'll finally break my New Year's resolution to stop making a
fool of myself on the net, and attempt an explanation of some things I
am just learning a bit about, namely, Teichm\"uller space and moduli
space. These are concepts that string theorists often throw around, and
when I first heard of them in that context I immediately dimissed them
as just another example of how physicists were learning far too much
mathematics for their own good. I take it all back! They are, in fact,
beautifully simple pieces of mathematics suited to physics in
2-dimensional spacetime. Two dimensions is low enough that one can often
actually understand exactly what's going on in problems that become
infinitely more tricky in higher dimensions. So even if one doesn't
``believe'' in \(2\)-dimensional physics the way the string theorists do
-- for them, physics happens on the worldsheet of the string, which is
2-dimensional -- it's worth learning as a kind of textbook case.

Everything I know so far is culled from the following sources, which by
no means form an exhaustive or even optimal set of references:

\begin{enumerate}
\def\labelenumi{\arabic{enumi})}
\item
  Y.\ Imayoshi and M.\  Taniguchi, \emph{An Introduction to Teichm\"uller spaces}, Springer, Berlin, 1991.
\item
  Joe Harris, ``An introduction to the moduli space of curves'', in
  \emph{Mathematical Aspects of String Theory} (proceedings of a
  conference at UC San Diego in 1986), ed.~S. T. Yau, World Scientific
  Press, 1987.
\item 
R.\ C.\ Penner, ``The moduli space of punctured surfaces'', same
volume.
\item
   John L.\ Harer, ``The cohomology of the moduli space of curves'', in
  \emph{Theory of Moduli (lectures given at the 3rd 1985 session of
  C.I.M.E. at Mondecatini Terme, Italy)}, ed.~E. Sernesi,
  Springer Lecture Notes in Mathematics \textbf{1337}, 1988.
\end{enumerate}

There are a number of ways of describing Teichm\"uller space and moduli
space. Maybe the easiest is this. Start with the surface of a doughnut
with \(g\) handles, or as the experts say, a ``surface of genus \(g\)''.
We can make this into a ``Riemann surface'' if we cover it with lots of
patches, or ``charts,'' each of which looks just like part the complex
plane (imagine a little piece of graph paper), and such that the change
of coordinates function relating overlapping patches is analytic, in the
usual sense of complex variables. The simplest Riemann surface is the
Riemann sphere, which is of genus zero; one gets this by taking the
complex plane and sticking on one more point, ``\(\infty\)'' -- think of
a sphere with ``\(\infty\)'' as the north pole. If we have a Riemann
surface, we can tell whether a complex-valued function on it is
analytic, simply by working locally in charts, so we can do complex
analysis as usual on a doughnut with lots of handles as long as we make
it into a Riemann surface! Since the main point of my article will be to
provide the reader with lots of buzzwords, I should add that making a
surface of genus \(g\) into a Riemann surface is called ``giving it a
complex structure.''

Now, how many ways are there to give a surface of genus \(g\) a complex
structure? First of all we need a good notion of when two Riemann
surfaces are ``the same''. They must, of course, have the same genus,
but there must also be a 1-1 and onto function from one to the other
that is everywhere analytic, with an analytic inverse. (Again, it makes
sense to say such a function is analytic, since we can cover each
Riemann surface with charts that we can think of as bits of the complex
plane.) Such a mapping is called a ``biholomorphic mapping'' --
holomorphic just being another word for analytic -- and if we want to
sound fancy, we say that the two Riemann surfaces are
``biholomorphically equivalent.''

Well, there's a famous old theorem of Riemann that for genus 0, there is
only ONE way to do it; any Riemann surface of genus 0 is
biholomorphically equivalent with the Riemann sphere. But for higher
genus there are infinitely many essentially different ways to give a
surface of genus \(g\) a complex structure. In fact, we can imagine a
big fat ``space'' of all ways. This is the moduli space of genus \(g\)!

The first key problem in the theory is to ``get our hands on'' moduli
space; that is, to describe it quite concretely. String theory provided
a lot of motivation for doing this very well, since the worldsheet of a
string -- that is, the string viewed in spacetime -- is just a surface,
and Feynman path integrals in string theory involve integrating over all
complex structures for this surface. To do integrals over moduli space
we need to bring it down to earth!

To do so, it's awfully handy to get involved with Teichm\"uller space.
Note that moduli space can be thought of as the space of equivalence
classes of complex structures on a fixed surface of genus \(g\), where
two complex structures are deemed ``the same'' if they are
biholomorphically equivalent. Teichm\"uller space is defined using a more
fine-grained notion of ``the same''. Note that any biholomorphic mapping
is a diffeomorphism, that is, a smooth mapping with a smooth inverse. In
fact, it must also be orientation-preserving, since an
orientation-reversing map like complex conjugation can never be
holomorphic! Henceforth I will always mean ``orientation-preserving
diffeomorophism'' when I speak of diffeomorphisms of a surface.

Now, some diffeomorphisms are ``connected to the identity'' and some
aren't. We say a diffeomorphism \(f\) is connected to the identity if
there is a smooth 1-parameter family of diffeomorphisms starting at
\(f\) and ending at the identity diffeomorphism. In other words, a
diffeomorphism is connected to the identity if you can do it
``gradually'' without ever having to cut the surface. To really
understand this you need to know some diffeomorphisms that \emph{aren't}
connected to the identity. Here's how to get one: start with your
surface of genus \(g > 0\), cut apart one of the handles along a circle,
give one handle a 360-degree twist, and glue the handles back together!
This is called a Dehn twist.

Anyway, Teichm\"uller space may be defined as the space of equivalence
classes of complex structures on a fixed surface of genus \(g\), where
two complex structures are counted as the same if they are
biholomorphically equivalent \emph{by a diffeomorphism connected to the
identity}.

A good way of understanding the relation between Teichm\"uller space and
moduli space is this. Define the mapping class group (of genus \(g\)) to
be the group of diffeomorphisms of a surface of genus \(g\) modulo the
subgroup of those connected to the identity. A beautiful fact is that
this group is generated by Dehn twists! In other words, given any
diffeomorphism of a surface, you can get it by first doing a bunch of
Dehn twists and then doing a diffeomorphism connected to the identity.
Since the mapping class group is finitely generated one should think of
it as a discrete group. In fact, folks know what the relations between
the generators are, too, and these are also very beautiful. I guess good
places to read about this stuff are the first paper that gave an
explicit presentation of mapping class groups:

\begin{enumerate}
\def\labelenumi{\arabic{enumi})}
\setcounter{enumi}{4}
\tightlist
\item
  A.\ Hatcher and W.\ Thurston, ``A presentation of the mapping class group of a closed, orientable
  surface'', \emph{Topology} \textbf{19}
  (1980), 221--237.
\end{enumerate}
\noindent
and the simplified treatment in

\begin{enumerate}
\def\labelenumi{\arabic{enumi})}
\setcounter{enumi}{5}
\tightlist
\item
  Bronislaw Wajnryb, ``A simple presentation for the mapping class group of an orientable
  surface'', \emph{Israel J.\ Math.} \textbf{45} (1983), 157--174.
\end{enumerate}
\noindent
Actually, though, I must admit my only acquaintance with mapping class
groups comes from leafing through

\begin{enumerate}
\def\labelenumi{\arabic{enumi})}
\setcounter{enumi}{6}
\tightlist
\item
  Joan S.\ Birman, \emph{Braids, Links, and Mapping Class Groups},
  Annals of Mathematics Studies \textbf{82}, Princeton
  University Press, 1974.
\end{enumerate}

As one can gather from this title there is a close connection between
mapping class groups and the braid group and knot theory, which is the
main reason why string theory allowed Witten to get new insights into
knots. (The more mundane connection, namely that one ties knots out of
string, seems largely unexplored, but see
\protect\hyperlink{week18}{``Week 18''}.) Let me not digress into this
fascinating realm, however! The point I want to make here is just that:

\begin{quote}
The mapping class group acts on Teichm\"uller space, and the quotient by
this group action is moduli space.
\end{quote}

Anyone used to how math goes should find this pretty believable, but let
me explain: given a diffeomorphism of our surface of genus \(g\), we can use
it as a kind of ``coordinate transformation'' to turn one complex
structure into another. So the group of diffeomorphisms acts on
Teichm\"uller space, but, given how Teichmuller space is defined, the
subgroup of diffeomorphisms connected to the identity acts trivially.
Thus the mapping class group acts on Teichm\"uller space. By how moduli
space is defined, two points in Teichm\"uller space define the same point
in moduli space iff one is obtained from another by an element of the
mapping class group.

Now the good thing about Teichm\"uller space is that it has very nice
coordinates on it, called Fenchel-Nielsen coordinates, which reveal it
to be diffeomorphic to \(\mathbb{R}^{6g-6}\) when \(g > 1\). (The case
\(g = 0\) is utterly dull, since Teichm\"uller space is a point, and the
case \(g = 1\) is beautifully treated using the fact that any Riemann
surface of genus 1 is biholomorphically equivalent to the quotient of
the complex plane by a lattice, relating this case to the subject of
elliptic functions, as touched upon in \protect\hyperlink{week13}{``Week
13''}. I should also add that \protect\hyperlink{week13}{``Week 13''}
indicates, at least in the \(g = 1\) case, why moduli space is often
called the ``moduli space of curves.'')

Let me say how these coordinates go, rather sketchily, just so the
mysterious number \(6g-6\) becomes not so mysterious! Take your surface
of genus \(g\) -- just think of it as a doughnut with \(n\) holes -- and
cut it up into ``pairs of pants,'' that is, pieces that look like \[
  \begin{tikzpicture}
    \draw[thick] (0,0) to (0,-1) to (-0.5,-2);
    \draw[thick] (0.75,0) to (0.75,-1) to (1.25,-2);
    \draw[thick] (0.125,-2) to (0.375,-1.5) to (0.625,-2);
  \end{tikzpicture}
\] from above. Topologically, a pair of pants is just a sphere with
three discs cut out of it! A more dignified term for a pair of pants is
a ``trinion,'' by the way.

The idea now is to describe the complex structure on each pair of pants
separately, and then describe how the pairs of pants are glued together.
Now, it turns out that the complex structures on each pair of pants are
very easily described (up to biholomorphic equivalences connected to the
identity). It takes 3 positive real numbers. There's a unique metric on
the original surface that is compatible with the complex structure
(i.e.~is a ``K\"ahler metric") and has curvature equal to \(-1\). This is
called a hyperbolic metric, as in hyperbolic geometry. Then, we can cut
the surface into pairs of pants along circles that are geodesics
relative to this metric. To describe the complex structure on each pair
of pants we simply need to measure the lengths of the 3 bounding
circles; these are called the ``geodesic length functions". In other
words, if your pair of pants was hyperbolic, a tailor would only need to
measure you waistlength and the lengths around the two cuffs, not the
inseams!

Now it's a fun exercise to show that we can chop up a surface of genus
\(g > 1\) into exactly \(2g-2\) pairs of pants. Doing so, moreover,
requires that we cut the surface along \(3g-3\) circles. (Draw some
pictures!) Thus, we have a total of \(3g-3\) geodesic length functions.
However, we also need to describe how the pairs of pants are attached to
each other. We can glue them together with any amount of twisting, and
this twisting is a real-valued parameter. So there are \(3g-3\)
``twisting parameters'' required to describe how the pairs of pants are
attached. We thus have a total of \(6g-6\) parameters, the
Fenchel-Nielsen coordinates, and Teichm\"uller space is diffeomorphic to
\(\mathbb{R}^{6g-6}\) (since the positive real numbers are diffeomorphic
to \(\mathbb{R}\) itself).

I think I want to quit here but not before making a few random remarks.

First, there's another description of Teichm\"uller space which gives it
a triangulation, i.e., describes it as a bunch of high-dimensional
tetrahedra (simplices) glued together. Harer's paper gives a nice quick
sketch of this construction; the buzzword to look for is ``quadratic
differentials.'' The nice thing about this is that the mapping class
group action respects this triangulation so we get a triangulation of
moduli space.

Second, quite recent work by Penner:

\begin{enumerate}
\def\labelenumi{\arabic{enumi})}
\setcounter{enumi}{7}
\tightlist
\item
  R.\ C.\ Penner, ``Universal constructions in Teichm\"uller theory'', 
  \emph{Adv. Math.} \textbf{98} (1993), 143--215.
\end{enumerate}
\noindent
shows how to fit the Teichm\"uller spaces for different genus \(g\) into
a single universal object. This was directly motivated by string theory,
but the basic idea is (I think) simply that there should be some sense
in which one can go ``continously'' from genus \(g+1\) to genus \(g\) by
making one handle smaller and smaller (with respect to the hyperbolic
metric) until it goes away.

Third, when studying the triangulation of Teichm\"uller space one
repeatedly runs across a certain ``Pachner move'' which goes from one
triangulation of a surface to another: \[
  \begin{tikzpicture}
    \draw[thick] (0,0) node{$\bullet$} to (1,1.5) node{$\bullet$} to (2,0) node{$\bullet$} to (1,-1.5) node{$\bullet$} to cycle;
    \draw[thick] (1,1.5) to (1,-1.5);
  \end{tikzpicture}
  \raisebox{4.5em}{$\qquad\longleftrightarrow\qquad$}
  \begin{tikzpicture}
    \draw[thick] (0,0) node{$\bullet$} to (1,1.5) node{$\bullet$} to (2,0) node{$\bullet$} to (1,-1.5) node{$\bullet$} to cycle;
    \draw[thick] (0,0) to (2,0);
  \end{tikzpicture}
\] which reminds me of lattice field theories in 2 dimensions (see
\protect\hyperlink{week16}{``Week 16''} for an explanation) and the
``pentagon diagram'' showing how to get between the 5 simplest ways to
triangulate a pentagon using this move (see, for example, Figure 3 in
Penner's paper above, or Figure 3.2 of Harer's paper). The pentagon
diagram appears both in Moore and Seiberg's famous paper on string
theory and chopping up surfaces into pairs of pants:

\begin{enumerate}
\def\labelenumi{\arabic{enumi})}
\setcounter{enumi}{8}
\tightlist
\item
  G.\ Moore and S.\ Seiberg, ``Classical and quantum conformal field theory'', \emph{Commun.\ Math.\ Phys.} \textbf{123} (1989), 177--254
\end{enumerate}
\noindent
and in category theory, with the relationship \emph{there} now pretty
well understood in terms of ``modular tensor categories'' --- see e.g.

\begin{enumerate}
\def\labelenumi{\arabic{enumi})}
\setcounter{enumi}{9}
\tightlist
\item
  Louis Crane, ``2-d physics and 3-d topology'', \emph{Commun.\ Math.\
  Phys.} \textbf{135} (1991) 615-640.
\end{enumerate}

Unfortunately, I haven't been able to get any string guru to sit down
and really clarify how the triangulations of Teichm\"uller space fit in.
Lest the reader wonder what the heck I'm going on about, the idea is
that category theory provides a marvelous way to unify the profusion of
mathematical structures that are coming up these days, and if we ever
understood everything in those terms, it would all seem much less
confusing.

\hypertarget{week29}{%
\section{January 14, 1994}\label{week29}}

I'm awfully busy this week, but feel like attempting to keep up with the
pile of literature that is accumulating on my desk, so this will be a
rather terse description of papers. All of these papers are related to
my current obsession with ``higher-dimensional algebra'' and its
applications to physics.

\begin{enumerate}
\def\labelenumi{\arabic{enumi})}
\item
   Ruth J.\ Lawrence, ``On algebras and triangle relations'', to appear
  in \emph{Proc. Top. \& Geom. Methods in Field Theory} (1992), eds.~J.\
  Mickelsson and O.\ Pekonen, World Scientific, Singapore.

  Ruth J.\ Lawrence, ``A presentation for Manin and Schechtman's higher braid groups'',
  available as MSRI preprint \texttt{04129-91}.

   Ruth J.\ Lawrence, ``Triangulations, categories and extended topological field
  theories'', to appear in \emph{Quantum Topology},
  eds.\ L.\ Kauffman and R.\ A.\ Baadhio, World Scientific, Singapore, 1993.

  Ruth J.\ Lawrence, ``Algebras and triangle relations'', Harvard U.\
  preprint.
\end{enumerate}

Many people are busily trying to extend the remarkable relationship
between knot theory and physics, which is essentially a feature of 3
dimensions, to higher dimensions. Since the \(3\)-dimensional case
required the development of new branches of algebra (namely, quantum
groups and braided tensor categories), it seems that the
higher-dimensional cases will require still further ``higher-dimensional
algebra.'' One approach, which is still being born, involves the use of
``\(n\)-categories,'' which are generalizations of braided tensor
categories suited for higher-dimensional physics. (See for example the
papers by Crane in \protect\hyperlink{week2}{``Week 2''}, by Kapranov
and Voevodsky in \protect\hyperlink{week4}{``Week 4''}, by Fischer and
Freed (separately) in \protect\hyperlink{week12}{``Week 12''}, and the
one by Gordon, Power, and Street below.) Lawrence has instead chosen to
invent ``\(n\)-algebras,'' which are vector spaces equipped with
operations corresponding to the ways one can subdivide
\((n-1)\)-dimensional simplices into more such simplices. (See the paper
by Chung, Fukuma and Shapere in \protect\hyperlink{week16}{``Week 16''}
for some of the physics motivation here.)

These alternative approaches should someday be seen as different aspects
of the same thing, but there as yet I know of no theorems to this
effect, so there is a lot of work to be done. Even more importantly,
there is a lot of work left to be done about inventing \emph{examples}
of these higher-dimensional structures. For example, there may
eventually be general results on ``boosting'' \(n\)-algebras to
\((n+1)\)-algebras, or \(n\)-categories to \((n+1)\)-categories, which
will explain how generally covariant physics in \(n\)-dimensional
spacetime relates to the same thing in one higher dimension. So far,
however, all we have is a few examples, which are not even clearly
related to each other. For example, Crane calls this boosting process
``categorification'' and has done it starting with the braided tensor
category of representations of a quantum group. Lawrence, on the other
hand, shows how to construct some 3-algebras from quantum groups. And
Freed has given a general procedure for ``boosting'' using path integral
methods that are not yet rigorous in the most interesting cases.

\begin{enumerate}
\def\labelenumi{\arabic{enumi})}
\setcounter{enumi}{1}
\tightlist
\item
   R.\ Gordon, A.\ J.\ Power, and Ross
  Street, \emph{Coherence for Tricategories}, Memoirs of the American
  Mathematical Society \textbf{558}, Providence, Rhode Island, 1995.
\end{enumerate}

An ``\(n\)-category'' is a kind of algebraic structure that has
``objects,'' ``morphisms'' between objects, ``2-morphisms'' between
morphisms, and so on up to ``\(n\)-morphisms.'' However, the
\emph{correct} definition of an \(n\)-category for the purposes of
physics is still unclear! I gave a rough explanation of the importance
of \(2\)-categories in physics in \protect\hyperlink{week4}{``Week 4''},
where I discussed Kapranov and Voevodsky's nice definition of braided
tensor \(2\)-categories. However, it seems likely that we will need to
understand the situation for larger \(n\) as well. This paper makes a big
step in this direction, by defining ``tricategories'' (what I might call
``weak \(2\)-categories'') and proving a ``strictification'' or
``coherence'' result analogous to the result that every braided tensor
category is equivalent to a ``strict'' one. The result is, however,
considerably more subtle, as it involves a special way of defining the
tensor product of \(2\)-categories due to Gray:

\begin{enumerate}
\def\labelenumi{\arabic{enumi})}
\setcounter{enumi}{2}
\item
   John W.\ Gray, \emph{Formal Category Theory: Adjointness for \(2\)-Categories}, 
  Springer Lecture Notes in Mathematics \textbf{391},
  Springer, Berlin, 1974.

    John W.\ Gray, ``Coherence for the tensor product of \(2\)-categories, and braid
  groups'', in \emph{Algebras, Topology, and Category Theory}, eds.~A.
  Heller and M. Tierney, Academic Press, New York, 1976, pp.~63--76.
\end{enumerate}

Briefly speaking, Gordon--Power--Street use a category they call ``Gray,''
the category of all \(2\)-categories, made into a symmetric monoidal
closed category using a modified version of Gray's tensor product. Then
they show that every tricategory (as defined by them) is
``triequivalent'' to a category enriched over Gray.

\begin{enumerate}
\def\labelenumi{\arabic{enumi})}
\setcounter{enumi}{3}
\tightlist
\item
  J.\ M.\ Maillet, ``On pentagon and tetrahedron equations'',
  available as
  \href{https://arxiv.org/abs/hep-th/9312037}{\texttt{hep-th/9312037}}.
\end{enumerate}

Maillet shows how to obtain solutions of the tetrahedron equations from
solutions of pentagon equations; all these geometrical equations are
part of the theory of \(2\)-categories, and this is yet another example
of a ``boosting'' construction as alluded to above.

\begin{enumerate}
\def\labelenumi{\arabic{enumi})}
\setcounter{enumi}{4}
\tightlist
\item
    David Yetter, ``Homologically twisted invariants related to (2+1)- and
  (3+1)-dimensional state-sum topological quantum field theories'', 
 available as
  \href{https://arxiv.org/abs/hep-th/9311082}{\texttt{hep-th/9311082}}.
\end{enumerate}

Let me simply quote the abstract: ``Motivated by suggestions of Paolo
Cotta-Ramusino's work at the physical level of rigor relating \(BF\)
theory to the Donaldson polynomials, we provide a construction
applicable to the Turaev--Viro and Crane--Yetter invariants of \emph{a
priori} finer invariants dependent on a choice of (co)homology class on
the manifold.'' The dream is that this would give a state-sum formula
for the Donaldson polynomials, but Yetter is careful to avoid claiming
this! A while back, Crane and Yetter showed how to get \(4\)-dimensional
TQFTs from certain 3d TQFTs by another kind of ``boosting'' procedure
related to those mentioned above, but the resulting TQFT in
\(4\)-dimensions did not by itself yield interesting new invariants of
4-manifolds. The procedure Yetter describes here generalizes the earlier
work by allowing the inclusion of an embedded 2-manifold.

\hypertarget{week30}{%
\section{January 14, 1994}\label{week30}}

For the most part, this is a terse description of some papers dealing
with quantum gravity. Some look to be quite important, but as I have not
had time to read them as thoroughly as I would like, I won't say much.

First, however, let me note some books:

\begin{enumerate}
\def\labelenumi{\arabic{enumi})}
\tightlist
\item
   Silvan S. Schweber, \emph{QED and the Men Who Made It: Dyson, Feynman, Schwinger and
  Tomonaga}, Princeton U.\ Press, Princeton, 1994.
\end{enumerate}

Back in the 1930s there was a crisis in physics: nobody knew how to
reconcile quantum theory with special relativity. This book describes
the history of how people struggled with this problem and achieved a
marvelous, but flawed, solution: quantum electrodynamics (QED).
Marvelous, because it made verified predictions of unparalleled
accuracy, involves striking new concepts, and gave birth to beautiful
new mathematics. Flawed, only because nobody yet knows for sure whether
the theory is mathematically well-defined --- for reasons profoundly
related to physics at ultra-short distance scales. This story should
give some inspiration to those currently attempting to reconcile quantum
theory with general relativity! Feynman, Schwinger, and Tomonaga won
Nobel prizes for QED, but Dyson was also instrumental in inventing the
theory, and the book is mainly a story of these 4 men.

\begin{enumerate}
\def\labelenumi{\arabic{enumi})}
\setcounter{enumi}{1}
\item
  Bruce Stephenson, \emph{The Music of the Heavens: Kepler's Harmonic Astronomy},  
  Princeton U.\ Press, Princeton, 1994.

   Bruce Stephenson, \emph{Kepler's Physical Astronomy}, Princeton U.\
  Press, Princeton, 1994.
\end{enumerate}

Kepler's work on astronomy was in part based on the notion of the
``music of the spheres,'' and in his Harmonice Mundi (1619) he sought to
relate planetary velocities to the notes of a chord. He was also
fascinated with geometry, and sought to relate the radii of the
planetary orbits to the Platonic solids. While this may seem a bit silly
nowadays, it's clear that this faith that mathematical patterns pervade
the heavens was a crucial part of how Kepler found his famous laws of
planetary motion. Also important, of course, was his use of what we
would now call ``physical'' reasoning to understand the heavens --- that
is, the use of analogies between the motions of heavenly bodies and that
of ordinary terrestrial matter. But even this is not as straightforward
as one might hope, since (Stephenson argues in the second book) this
physical reasoning was what we would now consider incorrect, even though
it led to valid laws. More inspiration for those now struggling amid
error to understand what the universe is really like!

\begin{enumerate}
\def\labelenumi{\arabic{enumi})}
\setcounter{enumi}{2}
\tightlist
\item
  Louis Kauffman and Sostenes Lins, 
  \emph{Temperley--Lieb Recoupling Theory and Invariants of 3-Manifolds},
  Annals of Mathematics Studies
  \textbf{133}, Princeton U.\ Press, Princeton, 1994.
\end{enumerate}

I described this briefly in \protect\hyperlink{week17}{``Week 17''},
before I had spent much time on it. Let me recall the main point: in the
late 80's Jones invented a new invariant of knots and links in ordinary
3d space, but then Witten recognized that this invariant came from a
quantum field theory, and thus could be extended to obtain an invariant
of links in arbitrary 3d manifolds. (In particular, taking the link to
be empty, one obtains a 3-manifold invariant.) In fact, there is a whole
family of such invariants, essentially one for each semisimple Lie
algebra, and Jones original example corresponded to the case
\(\mathfrak{su}(2)\). In this case the combinatorics of the invariants
are so simple that one can write a nice exposition in which one forgets
the underlying, fairly sophisticated, mathematical physics (quantum
groups, conformal field theory and the like) and simply presents the
``how-to'' using a kind of diagrammatic calculus known as ``skein
relations,'' or what Kauffman calls ``Temperley--Lieb recoupling
theory.'' That is the approach the authors take here. The curious reader
will naturally want to know more! For example, anyone familiar with
quantum theory and ``\(6j\) symbols'' will sense that this kind of thing
is lurking in the background, and indeed, it is.

Now for the papers:

\begin{enumerate}
\def\labelenumi{\arabic{enumi})}
\setcounter{enumi}{3}
\item
  C.\ Rovelli and L.\
  Smolin, ``The physical hamiltonian in quantum gravity'', available as
  \href{https://arxiv.org/abs/gr-qc/9308002}{\texttt{gr-qc/9308002}}.

  H.\ A.\ Morales-Tecotl and C.\
  Rovelli, ``Fermions in quantum gravity'', available as
  \href{https://arxiv.org/abs/gr-qc/9401011}{\texttt{gr-qc/9401011}}.
\end{enumerate}

The Rovelli--Smolin loop variables program proceeds apace! In the former
paper, Rovelli and Smolin consider quantum gravity coupled to a scalar
matter field which plays the role of a clock. (Using part of the system
described to play the role of a clock is a standard idea for dealing
with the ``problem of time,'' which arises in quantum theories on
spacetimes having no preferred coordinates, like quantum gravity.
However, getting this idea to actually work is not at all easy. For a
bit on this issue see \protect\hyperlink{week27}{``Week 27''}.) Only
after choosing this ``clock field'' can one work out a Hamiltonian for
the theory, write down the analog of Schr\"odinger's equation, and examine
the dynamics. Before, there is only a ``Hamiltonian constraint''
equation, also known as the Wheeler--DeWitt equation.

In the latter paper, Rovelli and Marales-Tecotl discuss how to include
fermions. The beautiful thing here is that fermions are described in the
loop representation by the \emph{ends} of arcs, while pure gravity is
described by loops. This is completely analogous to the old string
theory of mesons, in which mesons were represented as arcs of ``string''
-- the gluon field -- connecting two fermionic ``ends'' -- the quarks.

\begin{enumerate}
\def\labelenumi{\arabic{enumi})}
\setcounter{enumi}{4}
\tightlist
\item
  C.\ Di Bartolo, R.\ Gambini, J.\ Griego and J.\ Pullin, ``Extended loops: a new arena for           
  nonperturbative quantum gravity'',
  available as
  \href{https://arxiv.org/abs/gr-qc/9312029}{\texttt{gr-qc/9312029}}.
\end{enumerate}

For a while now, Gambini and collaborators have been developing a
modified version of the loop representation that appears to be
especially handy for doing perturbative calculations (perturbing in the
coupling constant, that is, Newton's gravitational constant ---
\emph{not} perturbing about a fixed flat ``background'' spacetime, which
is regarded as a ``no-no'' in this philosophy). The mathematical basis
for this ``extended'' loop representation is something quite charming in
itself: it amounts to embedding the loop group into an
(infinite-dimensional) Lie group. The ``perturbative'' calculations
described above are thus analogous to how one uses Lie algebras to study
Lie groups. In fact, this analogy is a deep one, since the extended loop
representation also permits perturbative calculations in Chern--Simons
theory, allowing one to calculate ``Vassiliev invariants'' starting just
from Lie-algebraic data. In fact this was done by Bar-Natan (cf
\protect\hyperlink{week3}{``Week 3''}), who was using the extended loop
representation without particularly knowing about that fact!

This paper puts the extended loop representation to practical use by
finding some new solutions to the quantum version of Einstein's
equations. These solutions are essentially Vassiliev invariants! (See
also the paper by Gambini and Pullin listed in
\protect\hyperlink{week23}{``Week 23''}.)

\begin{enumerate}
\def\labelenumi{\arabic{enumi})}
\setcounter{enumi}{5}
\tightlist
\item
  Domenico
  Giulini, ``Ashtekar variables in classical general relativity'', available as
  \href{https://arxiv.org/abs/gr-qc/9312032}{\texttt{gr-qc/9312032}}.
\end{enumerate}
\noindent
This was a lecture given at the 117th WE-Heraeus Seminar: ``The
Canonical Formalism in Classical and Quantum General Relativity'', 13-17
September 1993, Bad Honnef, Germany, the goal of which was to give an
introduction to Ashtekar's ``new variables'' for general relativity.

\hypertarget{week31}{%
\section{February 18, 1994}\label{week31}}

Well, I'm really busy these days trying to finish up a big project,
hence the low number of ``Weeks'' per week,
but papers are piling up, and there are some pretty interesting ones, so
I thought I'd quickly mention a few. This bunch will mainly concern
quantum gravity.

\begin{enumerate}
\def\labelenumi{\arabic{enumi})}
\tightlist
\item
  Louis
  Crane, ``Possible implications of the quantum theory of gravity'', available as
  \href{https://arxiv.org/abs/hep-th/9402104}{\texttt{hep-th/9402104}}.
\end{enumerate}
\noindent
This is one paper that everyone can read and enjoy, for although one may
find it too close to science fiction for comfort, it is far more
interesting than most science fiction. Louis Crane has been doing a lot
of excellent work on topological quantum field theory for the last few
years, strongly advocating the use of category theory as a unifying
principle in physics (essentially as an extension of the concept of
symmetry embodied in \emph{group} theory), but this is quite different
in flavor.

To begin with, Lee Smolin, one of the originators of the loop
representation of quantum gravity, has been spending the last year or so
writing a book in a popular style, to be entitled ``Life and Light,''
which tours the cosmos and makes some interesting speculations on
``evolutionary cosmology.'' These speculations are based on 2
hypotheses.

\begin{quote}
A. The formation of a black hole creates ``baby universe,'' the final
singularity of the black hole tunnelling right on through to the initial
``Big Bang'' singularity of the new universe thanks to quantum effects.
\end{quote}

While this must undoubtedly seem outre to anyone unfamiliar with the
sort of thing theoretical physicists amuse themselves with these days,
in a recent review article by John Preskill on the information loss
paradox for black holes, he reluctantly concluded that this was the
\emph{most conservative} solution of that famous problem! Recall the
problem: if a black hole evaporates its mass away via Hawking radiation,
and that radiation is pure blackbody radiation, hence carries none of
the information about the matter that originally formed the black hole,
one does not have conservation of information, or more technically
speaking, the time evolution is not unitary, since a pure state is
evolving into a mixed state. Hawking's original solution to this problem
was to bite the bullet and accept the nonunitarity, even though it goes
against the basic principles of quantum theory. This appears in:

\begin{enumerate}
\def\labelenumi{\arabic{enumi})}
\setcounter{enumi}{1}
\tightlist
\item
  S. W. Hawking, ``Black holes and thermodynamics'', \emph{Phys. Rev. D} \textbf{13} 
  (1976), 191--197.
\end{enumerate}

The ``baby universe'' solution simply says that the matter seeds a baby
universe and the information goes \emph{there}. Many other solutions
have been proposed; two recent review articles are

\begin{enumerate}
\def\labelenumi{\arabic{enumi})}
\setcounter{enumi}{2}
\item
  J. Preskill, ``Do black holes destroy information?'', available as
  \href{https://arxiv.org/abs/hep-th/9209058}{\texttt{hep-th/9209058}}.
\item
   Don Page, ``Black hole information'', in \emph{Proceedings of the 5th Canadian Conference on
  General Relativity and Relativistic Astrophysics, University of
  Waterloo, 13--15 May, 1993}, eds.\ R. B. Mann and R. G. McLenaghan,
  World Scientific, Singapore, 1994.  Available as
  \href{https://arxiv.org/abs/hep-th/9305040}{\texttt{hep-th/9305040}}.
\end{enumerate}
\noindent
Personally, I am a complete agnostic about this problem, since it rests
upon so many phenomena that are hypothesized but not yet observed, and
since any solution would require a theory of quantum gravity. I am
merely reporting the ideas of respected physicists! In any event, the
second hypothesis is:

\begin{quote}
B. Certain parameters of the baby universe are close to but different
than those of the parent universe.
\end{quote}

 The notion that certain physical
facts that appear as ``laws'' are actually part of the state of the
univese has in fact been rather respectable since the application of
spontaneous symmetry breaking to the Weinberg--Salaam model of electroweak
interactions, part of the standard model. (Again, being my usual
cautious self, I must note that a crucial piece of evidence for this
model, the Higgs boson, has not yet been seen.) The notion of
spontaneous symmetry breaking has become quite popular in particle
physics and is a key component of all current theories, such as GUTs or
string theory, that attempt to model the messy heap of observed
particles and interactions by some pristinely symmetrical Lagrangian.
The spontaneous symmetry breaking would be expected to have occured
shortly after the Big Bang, when it got cool enough, much as a hot piece
of iron will randomly settle upon some direction of magnetization as its
temperature fall below the Curie temperature. One application of this
notion to cosmology is already widely popular, namely, inflation. In
fact, pursuing the analogy with magnetic domains, i.e.~small regions
with different directions of magnetization, cosmologists have spend a
fair amount of energy thinking about ``domain walls,'' ``cosmic
strings,'' monopoles and other defects that might occur as residues of
this cooling-down process.

So again, while the idea must seem wild to anyone who has not
encountered it before, physicists these days are fairly comfortable with
the idea that certain ``fundamental constants'' could have been other
than they were. As for the constants of a baby universe being close to,
but different than, those of the parent universe, there is as far as I
know no suggested mechanism for this. This is perhaps the weakest link
in Smolin's argument (though I haven't seen his book yet). But it is at
least conceivable.

Now, given these hypotheses a marvellous consequence ensues: Darwinian
evolution! Those universes whose parameters are such that many black
holes are formed will have many progeny, so the constants of physics can
be expected to be ``tuned'' for the formation of many black holes. As
Smolin emphasizes, while the hypotheses A and B may seem impossible to
test directly at present, we do at least have a hope of testing this
consequence. He has studied the marvellously intricate process of star
formation in the galaxy and attempted to see whether altering the
constants of physics appear ``tuned'' for maximizing black hole
production, and he argues in his book that they do appear so tuned. Of
course, this is an extremely delicate business, since our understanding
of galaxy formation, star formation and black hole formation even in
\emph{this} universe is still rather weak --- much less for other
conceivable universes in which the fundamental constants take different
values.

Crane enters the fray at this point, and proposes an additional
conjecture:

\begin{quote}
SUCCESSFUL INDUSTRIAL CIVILIZATIONS WILL EVENTUALLY CREATE BLACK HOLES.
\end{quote}
\noindent
(The capital letters are his.) He breaks it up into two parts for us:

\begin{quote}
SUBCONJECTURE 1: SUCCESSFUL INDUSTRIAL CIVILIZATIONS WILL EVENTUALLY
WANT TO MAKE BLACK HOLES
\end{quote}
\noindent
and

\begin{quote}
SUBCONJECTURE 2: SUCCESSFUL INDUSTRIAL CIVILIZATIONS WILL EVENTUALLY BE
ABLE TO PRODUCE BLACK HOLES.
\end{quote}
\noindent
and argues for each. The result, as any good evolutionist will
recognize, is a kind of feedback loop whereby intelligence and baby
universe formation both affect each other. Indeed, Crane calls his
hypothesis the ``meduso-anthropic hypothesis,'' after certain jellyfish
with a two-stage life cycle in which medusids produce polyps and vice
versa. This has the charm of completely destroying the usual approach
(dare I say ``paradigm''?) of physics in which the parameters of the
universe are regarded as indifferent to the existence of intelligence.
Of course, the anthropic hypothesis is a previous attempt to breach this
firewall, but a much less dramatic one, since the only role intelligence
plays in that is \emph{noticing} the laws of the universe.

At this point let me leave off with a quote from Crane's paper:

\begin{quote}
``It is not hard to see that if these ideas are true, they will be the
victims of abuse to dwarf quantum healing and even quantum golf. That is
not sufficient reason to ignore them.''
\end{quote}
\noindent
and let me \emph{gradually} turn towards slightly less speculative
realms, eventually finishing with some papers containing rigorous
mathematics! To begin with, some more on black hole entropy:

\begin{enumerate}
\def\labelenumi{\arabic{enumi})}
\setcounter{enumi}{4}
\item
  Leonard Susskind, ``Some speculations about black hole entropy in string theory'',
   available as
  \href{https://arxiv.org/abs/hep-th/9309145}{\texttt{hep-th/9309145}}.

  Leonard Susskind and John Uglum, ``Black hole entropy in canonical quantum gravity and superstring
  theory'', available as
  \href{https://arxiv.org/abs/hep-th/9401070}{\texttt{hep-th/9401070}}.
\end{enumerate}
\noindent
The fact that the entropy of a black hole is (at least under certain
circumstances) proportional to the area of its event horizon is a
curious relationship between general relativity, quantum field theory
and statistical mechanics that many people believe to pointing
somewhere, but unfortunately nobody is sure where. Part of the reason is
that the standard derivations are somewhat indirect, and the event
horizon is not a physical object, so the sense in which it is the locus
of entropy is difficult to understand. These authors suggest that in
string theory it can be explained in terms of open strings having both
ends attached to the horizon.

\begin{enumerate}
\def\labelenumi{\arabic{enumi})}
\setcounter{enumi}{5}
\tightlist
\item
  C.\ R.\ Stephens,
  G.\ 't Hooft and B.\ F.\ Whiting, ``Black hole evaporation without information loss'', available as
  \href{https://arxiv.org/abs/gr-qc/9310006}{\texttt{gr-qc/9310006}}.
\end{enumerate}
\noindent
This is an attempt to make black holes radiate away and disappear in a
manner that preserves unitarity. I've been too busy to read it. And now
for some wormholes:

\begin{enumerate}
\def\labelenumi{\arabic{enumi})}
\setcounter{enumi}{6}
\tightlist
\item
  Hoi-Kwong Lo,
  Kai-Ming Lee, and John Preskill, ``Complementarity in Wormhole Chromodynamics'', available as
  \href{https://arxiv.org/abs/hep-th/9308044}{\texttt{hep-th/9308044}}.
\end{enumerate}
\noindent
Let me just quote the abstract and note that there is probably some
quite interesting topology to be obtained by applying this sort of idea
to mathematics:

\begin{quote}
The electric charge of a wormhole mouth and the magnetic flux ``linked''
by the wormhole are non-commuting observables, and so cannot be
simultaneously diagonalized. We use this observation to resolve some
puzzles in wormhole electrodynamics and chromodynamics. Specifically, we
analyze the color electric field that results when a colored object
traverses a wormhole, and we discuss the measurement of the wormhole
charge and flux using Aharonov-Bohm interference effects. We suggest
that wormhole mouths may obey conventional quantum statistics, contrary
to a recent proposal by Strominger.
\end{quote}

Finally, lest the mathematicians think I have abandoned ship, some
rigorous results:

\begin{enumerate}
\def\labelenumi{\arabic{enumi})}
\setcounter{enumi}{7}
\tightlist
\item
  Piotr T.\
  Chrusciel, ````No hair'' theorems --- folklore, conjectures, results", available 
  as \href{https://arxiv.org/abs/gr-qc/9402032}{\texttt{gr-qc/9402032}}.
\end{enumerate}
\noindent
The famous ``no hair'' theorem says that in general relativity static
black hole solutions are determined by very few parameters --- typically
listed as mass, angular momentum and charge in ``rest frame'' of the
black hole. There have been many attempts to extend this result,
especially because no \emph{actual} black hole is likely to be utterly
static, since it presumably formed at some time. I have not read this
but Chrusciel is a very careful person so I expect it will be up to the
standards of his nice review of work on the cosmic censorship
hypothesis,

\begin{enumerate}
\def\labelenumi{\arabic{enumi})}
\setcounter{enumi}{8}
\tightlist
\item
  Piotr T. Chrusciel, ``On uniqueness in the large of solutions of Einstein's equations
  (``Strong cosmic censorship'')'', in
  \emph{Mathematical Aspects of Classical Field Theory}, Contemp. Math.
  \textbf{132}, eds.~Gotay, Marsden and Moncrief, American Mathematical
  Society, Providence, Rhode Island,
  1992, pp.~235--274.
\end{enumerate}

\begin{center}\rule{0.5\linewidth}{0.5pt}\end{center}

\textbf{Addendum}: See \protect\hyperlink{week33}{``Week 33''} for a
paper by Smolin on his evolutionary cosmology theory. His book came out
in 1997 under the title \emph{The Life of the Cosmos'} --- see
\protect\hyperlink{101}{``Week 101''} for details.

\begin{center}\rule{0.5\linewidth}{0.5pt}\end{center}

\begin{quote}
\emph{The thing that makes things and the thing that makes things fall
apart --- they're the same thing. Entropy maximization!}

--- Chris Lee
\end{quote}

\hypertarget{week32}{%
\section{March 10, 1994}\label{week32}}

Well, I visited Georgia Tech last week to spread the gospel of ``knots
and quantum gravity,'' and came across a most fascinating development.
I'm sure readers of \texttt{sci.math} and \texttt{sci.math.research}
have taken note of the \emph{New York Journal of Mathematics}. This is
one of the first refereed electronic journals of mathematics. Neil
Calkin at Georgia Tech is helping to start up another one --- the
\emph{Electronic Journal of Combinatorics}. Though it's unlikely,
perhaps some among you are still unaware (or unconvinced) of how
essential it is that we develop fully refereed free-of-charge electronic
journals of mathematics and physics. The first and most obvious reason
is that computer-based media offer all sorts of flexibility that print
media lack --- more on this later. But the other reason is that the
monopoly of print journals \emph{must} be broken.

For example, U.\ C.\ Riverside does not subscribe to \emph{Communications
in Mathematical Physics}, despite the fact that this is \emph{the}
crucial journal in that subject, because this journal costs \$3,505 a
year! The ridiculous price is, of course, in part precisely because this
is the crucial journal in that subject, in part because the journal uses
antiquated and expensive production methods involving paper, and in part
because, being a big operation, it is basically run by a publishing
house rather than mathematical physicists. Luckily, with the advent of
the preprint mailing lists \texttt{hep-th} and \texttt{gr-qc}, I don't
\emph{need} to read \emph{Communications in Mathematical Physics} very often! I
simply get my list of abstracts each day by email from Los Alamos, and
send email to get the papers I want, in LaTeX or TeX form. The middleman
has been cut out --- at least for the moment.

One problem with preprint mailing lists, though, is that the preprints
have not gone through the scrutiny of the referee process. This is,
frankly, much less of a problem for the \emph{readers} than is commonly
imagined, because this scrutiny is less intense than people who have
never refereed papers think! Many refereed papers have errors, and I
would personally feel very uncomfortable using a result unless I either
understood the proof or knew that most experts believed it. The real
need for refereed journals, in my slightly cynical opinion, is that
academics need \emph{refereed publications} to advance in their jobs:
the people who give tenure, promotions etc.\ cannot be expected to read
and understand one's papers. This is, of course, also the reason for
other strange phenomena, such as the idea of \emph{counting} somebody's
publications to see how good they are. We need only count Alexander
Abian's publications to see the limitations of this approach.

Eventually, a few birds may be killed with one stone by means of ``seals
of approval'' or SOAPs, which are being widely discussed by people
interested in the ``information superhighway,'' or --- let's call a
spade a spade --- the Internet. For more on these, check out the
\texttt{newsgroup\ comp.interpedia}, or read material about the Xanadu
project. The idea here is that eventually we will have a good system
whereby people can append comments to documents, such as ``there is an
error in the proof of Lemma 1.5, which can be fixed as follows\ldots{}''
or simply various seals of approval, functioning similarly to the seal
of approval ones paper obtains by being published in a journal. E.g.,
one could make ones paper available by ftp or some other protocol, and
``submitting it to a journal'' might amount to asking for a particular
SOAP, with various SOAPs carrying various amounts of prestige, and so
on.

Of course, journals also function as a kind of information ``hub'' or
central access point. We all know that to find out what's the latest
trend in particle physics, it suffices to glance at \emph{Nucl.\ Phys.\ B}
and certain other journals, and so on. It is not clear that the function
of ``hub'' and the function of SOAP need be combined into a single
institution, once the onerous task of transcribing ideas onto dead trees
and shipping them all around the world becomes (at least partially)
obsolete.

It is also not at all certain whether, in the long run, the monopolistic
power of journals to charge large fees for accessing information will be
broken by the new revolutions in technology. This is, of course, just
one small facet of the political/economic struggle for control over
information flow that is heating up these days, at least in the U.S.,
among telephone companies, cable TV stations, computer networks such as
Compuserve, etc. etc. f mathematicians and physicists don't think
about these issues, someone else who has will wind up defining the
future for us.

Anyway, for now it seems to make good sense to start refereed journals
of mathematics and physics that are accessible electronically, free of
charge, over the Internet. Not too long ago one would commonly hear the
remark ``\ldots but of course nobody would ever want to do that,
because\ldots{}'' followed by some reason or other, reminiscent of how
\emph{clearly} nobody would want to switch from horses to automobiles
because then one would have to build \emph{gas stations all over the
place} --- obviously too much bother to be worthwhile. Now, however,
things are changing and the new electronic journals are getting quite
respectable lists of editors, and they seem to have a good chance of
doing well. I urge everyone to support free-of-charge electronic
journals by submitting good papers!

Let me briefly describe the electronic journals I mentioned above. The
New York Journal's chief editor is Mark Steinberger, at SUNY Albany. 
The journal covers algebra, modern
analysis, and geometry/topology. Access is through anonymous ftp, gopher
and listserv, the latter being (I believe) a mailing list protocol. One
can subscribe by sending email to \texttt{listserv@albany.edu} or
\texttt{listserv@albany.bitnet}; if you want abstracts for all the
papers, the body of your email should read

\begin{center}
\begin{verbatim}
subscribe NYJMTH-A <your full name>
\end{verbatim}
\end{center}
\noindent
but you can also subscribe to only certain topics (one of the great
things about electronic journals --- one can only begin to imagine the
possibilities inherent in this concept!), as follows:

\begin{quote}
{\rm 
Algebra:

\texttt{subscribe\ NYJM-ALG\ \textless{}your\ full\ name\textgreater{}}

Analysis:

\texttt{subscribe\ NYJM-AN\ \ \textless{}your\ full\ name\textgreater{}}

Geometry/Topology:

\texttt{subscribe\ NYJM-TOP\ \textless{}your\ full\ name\textgreater{}}
}
\end{quote}
\noindent
Papers are accepted in amstex and amslatex, and when you get papers you
get a \texttt{.dvi} file.

The \emph{Electronic Journal of Combinatorics} is taking a somewhat more
ambitious approach that has me very excited. Namely, they are using
Mosaic, a hypertext interface to the WWW (World-Wide Web). This means,
to technical illiterates such as myself, that if you can ever get your
system manager to get the software running, you can see a ``front page''
of the journal, with the names of the articles and other things
underlined (or in color if you're lucky). To go to any underlined item,
you simply click your mouse on it. In fact, you can use this method to
navigate throughout the whole WWW, which is a vast, sprawling network of
linked files, including --- so I hear --- ``This Week's Finds''! In the
\emph{Electronic Journal of Combinatorics}, when you click on an article
you will see it in postscript form, pretty equations and all. You can
also get yourself a copy and print it out. Neil showed me all this stuff
and my mouth watered! The danger of this ambitious approach is of course
that folks who haven't kept up with things like the WWW may find it
intimidating\ldots{} for a while. It's actually not too complicated.

This journal will be widely announced pretty soon. The editor in chief
is Herbert S. Wilf, \texttt{wilf@central.cis.upenn.edu}, and the
managing editor is Neil Calkin, \texttt{calkin@math.gatech.edu}. It
boasts an impressive slate of editors (even to me, who knows little
about combinatorics), including Graham, Knuth, Rota and Sloane. To get
browse the journal, which is presently under construction, you just do
the following if you can use Mosaic: ``Click on the button marked `Open'
and then type in
\texttt{http://math34.gatech.edu:8080/Journal/journalhome.html}''. To
\emph{get} Mosaic, do anonymous ftp to \texttt{ftp.ncsa.uiuc.edu} and
\texttt{cd} to \texttt{Web/Mosaic\_binaries} --- and then you're on your
own, I just tried it and there were too many people on! --- but Neil
says it's not too hard to get going. I will try as soon as I have a free
day.

``Ahem!'' the reader comments. ``What does this have to do with
mathematical physics?'' Well, seeing how little I'm being paid, I see
nothing wrong with interpreting my mandate rather broadly, but I should
add the following. 

1) There are periodic posts on \texttt{sci.physics}
about physics on the WWW; there's a lot out there, and to get started
one can always try the following. The information below is taken from Scott
Chase's physics FAQ:

\begin{verbatim}
* How to get to the Web

    If you have no clue what WWW is, you can go over the Internet with
telnet to info.cern.ch (no login required) which brings you to the WWW 
Home Page at CERN. You are now using the simple line mode browser. To move 
around the Web, enter the number given after an item. 

* Browsing the Web

    If you have a WWW browser up and running, you can move around
more easily. The by far nicest way of "browsing" through WWW uses the
X-Terminal based tool "XMosaic". Binaries for many platforms (ready for
use) and sources are available via anonymous FTP from ftp.ncsa.uiuc.edu
in directory Web/xmosaic.  The general FTP repository for browser
software is info.cern.ch (including a hypertext browser/editor for
NeXTStep 3.0)

* For Further Information

    For questions related to WWW, try consulting the WWW-FAQ: Its most 
recent version is available via anonymous FTP on rtfm.mit.edu in 
/pub/usenet/news.answers/www-faq , or on WWW at 
http://www.vuw.ac.nz:80/overseas/www-faq.html

    The official contact (in fact the midwife of the World Wide Web) 
is Tim Berners-Lee, timbl@info.cern.ch. For general matters on WWW, try 
www-request@info.cern.ch or Robert Cailliau (responsible for the "physics" 
content of the Web, cailliau@cernnext.cern.ch).
\end{verbatim}

And: 2) there are rumors, which I had better not elaborate on yet, about
an impending electronic journal of mathematical physics! I eagerly await
it!

Okay, just a bit about actual mathematical physics per se this time.

\begin{enumerate}
\def\labelenumi{\arabic{enumi})}
\tightlist
\item
  Carlo Rovelli, ``On quantum mechanics'', available as
  \href{https://arxiv.org/abs/hep-th/9403015}{\texttt{hep-th/9403015}}.
\end{enumerate}
\noindent
This interesting paper suggests that reason why we are constantly
arguing about the meaning of quantum mechanics, despite the fact that it
works perfectly well and is obviously correct, is that we are making a
crucial conceptual error. Rovelli very nicely compares the problem to
special relativity before Einstein did his thing: we had Lorentz
transformations, but they seemed very odd and paradoxical, because the
key notion that the space/time split was only defined \emph{relative to
a frame} (or ``observer'' if we wish to anthropomorphize) was lacking.
Rovelli proposes that in quantum mechanics the problem is that we are
lacking the notion that the \emph{state} of a system is only defined
relative to an observer. (The ``Wigner's friend'' puzzle is perhaps the
most obvious illustration here.) What, though, is an observer? Any
subsystem of a quantum system, says Rovelli; there is no fundamental
``observer-observed distinction.'' This fits in nicely with some recent
work by Crane and myself on quantum gravity, so I like it quite a bit,
though it is clearly not the last word on this issue (nor does Rovelli
claim it to be).

\begin{enumerate}
\def\labelenumi{\arabic{enumi})}
\setcounter{enumi}{1}
\tightlist
\item
  Alan D.\ Rendall, ``Adjointness relations as a criterion for choosing an inner
  product'', available as
  \href{https://arxiv.org/abs/gr-qc/9403001}{\texttt{gr-qc/9403001}}.
\end{enumerate}
\noindent
The inner product problem in quantum gravity is an instance of a
general, very interesting mathematics problem, namely, of determining an
inner product on a representation of a star-algebra, by demanding that
the representation be a star-representation. Rendall has proved some
very nice results on this issue.

\begin{enumerate}
\def\labelenumi{\arabic{enumi})}
\setcounter{enumi}{2}
\tightlist
\item
  Maxim Kontsevich and Yuri Manin, ``Gromov-Witten classes, quantum cohomology, and  
  enumerative  geometry'', available as 
  \href{https://arxiv.org/abs/hep-th/9402147}{\texttt{hep-th/9402147}}.
\end{enumerate}

I will probably never understand this paper so I might as well mention
it right away. Kontsevich's work on knot theory, and Manin's work on
quantum groups and (earlier) instantons is extremely impressive, so I
guess they can be forgiven for their interest in algebraic geometry. (A
joke.) Let me simply quote:

\begin{quote}
``The paper is devoted to the mathematical aspects of topological
quantum field theory and its applications to enumerative problems of
algebraic geometry. In particular, it contains an axiomatic treatment of
Gromov-Witten classes, and a discussion of their properties for Fano
varieties. Cohomological Field Theories are defined, and it is proved
that tree level theories are determined by their correlation functions.
Applications to counting rational curves on del Pezzo surfaces and
projective spaces are given.''
\end{quote}

\hypertarget{week33}{%
\section{May 10, 1994}\label{week33}}

With tremendous relief, I have finished writing a book, and will return
to putting out This Week's Finds on a roughly weekly basis. Let me
briefly describe my book, which took so much more work than I had
expected\ldots{} and then let me start catching up on listing some of
the stuff that's cluttering my desk!

\begin{enumerate}
\def\labelenumi{\arabic{enumi})}
\tightlist
\item
  John Baez and Javier de
  Muniain, \emph{Gauge Fields, Knots and Gravity}, World Scientific Press, Singapore, 
  1994.
\end{enumerate}

This book is based on a seminar I taught in 1992--93. We start out
assuming the reader is familiar with basic stuff --- Maxwell's
equations, special relativity, linear algebra and calculus of several
variables --- and try to prepare the reader to understand recent work on
quantum gravity and its relation to knot theory. It proved difficult to
do this well in a mere 460 pages. Lots of tantalizing loose ends are
left dangling. However, there are copious references so that the reader
can pursue various subjects further.

\begin{quote}
{\rm
Part 1. Electromagnetism

Chapter 1. Maxwell's Equations

Chapter 2. Manifolds

Chapter 3. Vector Fields

Chapter 4. Differential Forms

Chapter 5. Rewriting Maxwell's Equations

Chapter 6. DeRham Theory in Electromagnetism

Part 2. Gauge Fields

Chapter 1. Symmetry

Chapter 2. Bundles and Connections

Chapter 3. Curvature and the Yang--Mills Equations

Chapter 4. Chern--Simons Theory

Chapter 5. Link Invariants from Gauge Theory

Part 3. Gravity

Chapter 1. Semi-Riemannian Geometry

Chapter 2. Einstein's Equations

Chapter 3. Lagrangians for General Relativity

Chapter 4. The ADM Formalism

Chapter 5. The New Variables
}
\end{quote}

\begin{enumerate}
\def\labelenumi{\arabic{enumi})}
\setcounter{enumi}{1}
\tightlist
\item
  Asher Peres, \emph{Quantum Theory: Concepts and Methods}, Kluwer
  Academic Publishers, Amsterdam, 1994.
\end{enumerate}
\noindent
As Peres notes, there are many books that teach students how to solve
quantum mechanics problems, but not many that tackle the conceptual
puzzles that fascinate those interested in the foundations of the
subject. His book aims to fill this gap. Of course, it's impossible not
to annoy people when writing about something so controversial; for
example, fans of Everett will be distressed that Peres' book contains
only a brief section on ``Everett's interpretation and other bizarre
interpretations''. However, the book is clear-headed and discusses a lot
of interesting topics, so everyone should take a look at it.

Schr\"odinger's cat, Bell's inequality and Wigner's friend are old
chestnuts that everyone puzzling over quantum theory has seen, but there
are plenty of popular new chestnuts in this book too, like ``quantum
cryptography'', ``quantum teleportation'', and the ``quantum Zeno
effect'', all of which would send shivers up and down Einstein's spine.
There are also a lot of gems that I hadn't seen, like the
Wigner--Araki--Yanase theorem. Let me discuss this theorem a bit.

Roughly, the WAY theorem states that it is impossible to measure an
operator that fails to commute with an additive conserved quantity. Let
me give an example to clarify this and then give the proof. Say we have
a particle with position \(q\) and momentum \(p\), and a measuring
apparatus with position \(Q\) and momentum \(P\). Let's suppose that the
total momentum \(p + P\) is conserved --- which will typically be the
case if we count as part of the ``apparatus'' everything that exerts a
force on the particle. Then as a consequence of the WAY theorem we can
see that (in a certain sense) it is impossible to measure the particle's
position \(q\); all we can measure is its position \emph{relative} to
the apparatus, \(q-Q\).

Of course, whenever a ``physics theorem'' states that something is
impossible one must peer into it and determine the exact assumptions and
the exact result! Lots of people have gotten in trouble by citing
theorems that seem to show something is impossible without reading the
fine print. So let's see what the WAY theorem \emph{really} says!

It assumes that the Hilbert space for the system is the tensor product
of the Hilbert space for the thing being observed --- for short, let's
call it the ``particle'' --- and the Hilbert space for the measuring
apparatus. Assume also that \(A\) and \(B\) are two observables
belonging to the observed system, while \(C\) is an observable belonging
to the measuring apparatus; suppose that \(B + C\) is conserved, and
let's try to show that we can only measure \(A\) if it commutes with
\(B\). (Our assumptions automatically imply that \(A\) commutes with
\(C\), by the way.)

So, what do we mean when we speak of ``measuring \(A\)''? Well, there
are various things one might mean. The simplest is that if we start the
combined system in some tensor product state \(u(i) \otimes v\), where
\(v\) is the ``waiting and ready'' state of the apparatus and \(u(i)\)
is a state of the observed system that's an eigenvector of \(A\):
\[Au(i) = a(i)u(i),\] then the unitary operator \(U\) corresponding to
time evolution does the following:
\[\mathrm{U}(u(i) \otimes v) = u(i) \otimes v(i)\] where the state
\(v(i)\) of the apparatus is one in which it can be said to have
measured the observable \(A\) to have value \(a(i)\). E.g., the
apparatus might have a dial on it, and in the state \(v(i)\) the dial
reads ``\(a(i)\)''. Of course, we are really only justified in saying a
measurement has occurred if the states \(v(i)\) are \emph{distinct} for
different values of \(i\).

Note: here the WAY theorem seems to be restricting itself to
nondestructive measurements, since the observed system is remaining in
the state \(u(i)\). If you go through the proof you can see to what
extent this is crucial, and how one might modify the theorem if this is
not the case.

Okay, we have to show that we can only ``measure \(A\)'' in this sense
if \(A\) commutes with \(B\). We are assuming that \(B + C\) is
conserved, i.e., \[U^*(B + C)U = B + C.\] First note that
\[\langle u(i), [A,B] u(j)\rangle = (a(i)-a(j)) \langle u(i), Bu(j)\rangle.\]
On the other hand, since \(A\) and \(B\) only act on the Hilbert space
for the particle, we also have
\[\begin{aligned}\langle u(i),[A,B]u(j) \rangle &= \langle u(i)\otimes v,[A,B]u(j) \rangle \\ &= \langle u(i)\otimes v,[A,B+C]u(j) \rangle \\ &= (a(i)-a(j))\langle u(i)\otimes v,(B+C)u(j)\otimes v \rangle\end{aligned}\]
It follows that if \(a(i)-a(j)\) isn't zero,
\[\begin{aligned}\langle u(i),Bu(j) \rangle &= \langle u(i)\otimes v,(B+C)u(j)\otimes v \rangle \\ &= \langle u(i)\otimes v, U^*(B+C)Uu(j)\otimes v \rangle \\ &= \langle u(i)\otimes v(i),(B+C)u(j)\otimes v(j) \rangle \\ &= \langle u(i),Bu(j) \rangle\langle v(i),v(j) \rangle + \langle u(i),u(j) \rangle\langle v(i),Cv(j) \rangle\end{aligned}\]
but the second term vanishes since \(u(i)\) are a basis of eigenvectors
and \(u(i)\) and \(u(j)\) correspond to different eigenvalues, so
\[\langle u(i), Bu(j)\rangle = \langle u(i), Bu(j)\rangle \langle v(i), v(j)\rangle\]
which means that either \(\langle v(i), v(j)\rangle = 1\), hence
\(v(i) = v(j)\) (since they are unit vectors), so that no measurement
has really been done, OR that \(\langle u(i), B u(j)\rangle = 0\), which
means (if true for all \(i,j\)) that \(A\) commutes with \(B\).

So, we have proved the result, using one extra assumption that I didn't
mention at the start, namely that the eigenvalues \(a(i)\) are distinct.

I can't say that I really understand the argument, although it's easy
enough to follow the math. I will have to ponder it more, but it is
rather interesting, because it makes more precise (and general) the
familiar notion that one can't measure \emph{absolute} positions, due to
the translation-invariance of the laws of physics; this translation
invariance is of course what makes momentum be conserved. (What I just
wrote makes me wonder if someone has shown a classical analog of the WAY
theorem.)

Anyway, here's the table of contents of the book:

\begin{quote}
{\rm
Chapter 1: Introduction to Quantum Physics

1-1. The downfall of classical concepts 

1-2. The rise of randomness 

1-3. Polarized photons 

1-4. Introducing the quantum language 

1-5. What is a measurement? 

1-6. Historical remarks 

1-7. Bibliography 

Chapter 2: Quantum Tests

2-1. What is a quantum system? 

2-2. Repeatable tests 

2-3. Maximal quantum tests 

2-4. Consecutive tests 

2-5. The principle of interference 

2-6. Transition amplitudes 

2-7. Appendix: Bayes's rule of statistical inference 

2-8. Bibliography 

Chapter 3: Complex Vector Space

3-1. The superposition principle 

3-2. Metric properties 

3-3. Quantum expectation rule 

3-4. Physical implementation 

3-5. Determination of a quantum state 

3-6. Measurements and observables 

3-7. Further algebraic properties 

3-8. Quantum mixtures 

3-9. Appendix: Dirac's notation 

3-10. Bibliography 

Chapter 4: Continuous Variables

4-1. Hilbert space 

4-2. Linear operators 

4-3. Commutators and uncertainty relations 

4-4. Truncated Hilbert space 

4-5. Spectral theory 

4-6. Classification of spectra 

4-7. Appendix: Generalized functions 

4-8. Bibliography 

Chapter 5: Composite Systems

5-1. Quantum correlations 

5-2. Incomplete tests and partial traces 

5-3. The Schmidt decomposition 

5-4. Indistinguishable particles 

5-5. Parastatistics 

5-6. Fock space 

5-7. Second quantization 

5-8. Bibliography 

Chapter 6: Bell's Theorem

6-1. The dilemma of Einstein, Podolsky, and Rosen 

6-2. Cryptodeterminism 

6-3. Bell's inequalities 

6-4. Some fundamental issues 

6-5. Other quantum inequalities 

6-6. Higher spins 

6-7. Bibliography 

Chapter 7: Contextuality

7-1. Nonlocality versus contextuality 

7-2. Gleason's theorem 

7-3. The Kochen--Specker theorem 

7-4. Experimental and logical aspects of inseparability 

7-5. Appendix: Computer test for Kochen-Specker contradiction 

7-6. Bibliography 

Chapter 8: Spacetime Symmetries

8-1. What is a symmetry? 

8-2. Wigner's theorem 

8-3. Continuous transformations 

8-4. The momentum operator 

8-5. The Euclidean group 

8-6. Quantum dynamics 

8-7. Heisenberg and Dirac pictures 

8-8. Galilean invariance 

8-9. Relativistic invariance 

8-10. Forms of relativistic dynamics 

8-11. Space reflection and time reversal 

8-12. Bibliography 

Chapter 9: Information and Thermodynamics

9-1. Entropy 

9-2. Thermodynamic equilibrium 

9-3. Ideal quantum gas 

9-4. Some impossible processes 

9-5. Generalized quantum tests 

9-6. Neumark's theorem 

9-7. The limits of objectivity 

9-8. Quantum cryptography and teleportation 

9-9. Bibliography 

Chapter 10: Semiclassical Methods

10-1. The correspondence principle 

10-2. Motion and distortion of wave packets 

10-3. Classical action 

10-4. Quantum mechanics in phase space 

10-5. Koopman's theorem 

10-6. Compact spaces 

10-7. Coherent states 

10-8. Bibliography 

Chapter 11: Chaos and Irreversibility

11-1. Discrete maps 

11-2. Irreversibility in classical physics 

11-3. Quantum aspects of classical chaos 

11-4. Quantum maps 

11-5. Chaotic quantum motion 

11-6. Evolution of pure states into mixtures 

11-7. Appendix: PostScript code for a map 

11-8. Bibliography 

Chapter 12: The Measuring Process

12-1. The ambivalent observer 

12-2. Classical measurement theory 

12-3. Estimation of a static parameter 

12-4. Time-dependent signals 

12-5. Quantum Zeno effect 

12-6. Measurements of finite duration 

12-7. The measurement of time 

12-8. Time and energy complementarity 

12-9. Incompatible observables 

12-10. Approximate reality 

12-11. Bibliography
}
\end{quote}

\begin{enumerate}
\def\labelenumi{\arabic{enumi})}
\setcounter{enumi}{2}
\tightlist
\item
  Bernd Br\"ugmann, ``Loop representations'', available as \href{https://arxiv.org/abs/gr-qc/9312001}{\texttt{gr-qc/9312001}}.
\end{enumerate}
\noindent
This is a nice review article on loop representations of gauge theories.
Anyone wanting to jump into the loop representation game would be well
advised to start here.

\begin{enumerate}
\def\labelenumi{\arabic{enumi})}
\setcounter{enumi}{3}
\tightlist
\item
  Lee Smolin, ``The fate of black hole singularities and the parameters of the
  standard models of particle physics and cosmology'', 
  available as
  \href{https://arxiv.org/abs/gr-qc/9404011}{\texttt{gr-qc/9404011}}.
\end{enumerate}
\noindent
This is about Smolin's ``evolutionary cosmology'' scenario, which I
already discussed in \protect\hyperlink{week31}{``Week 31''}. Let me
just quote the abstract:

\begin{quote}
A cosmological scenario which explains the values of the parameters of
the standard models of elementary particle physics and cosmology is
discussed. In this scenario these parameters are set by a process
analogous to natural selection which follows naturally from the
assumption that the singularities in black holes are removed by quantum
effects leading to the creation of new expanding regions of the
universe. The suggestion of J.\ A.\ Wheeler that the parameters change
randomly at such events leads naturally to the conjecture that the
parameters have been selected for values that extremize the production
of black holes. This leads directly to a prediction, which is that small
changes in any of the parameters should lead to a decrease in the number
of black holes produced by the universe. On plausible astrophysical
assumptions it is found that changes in many of the parameters do lead
to a decrease in the number of black holes produced by spiral galaxies.
These include the masses of the proton, neutron, electron and neutrino
and the weak, strong and electromagnetic coupling constants.
Finally, this scenario predicts a natural time scale for cosmology equal
to the time over which spiral galaxies maintain appreciable rates of
star formation, which is compatible with current observations that
\(\Omega = .1-.2\).
\end{quote}

\hypertarget{week34}{%
\section{May 24, 1994}\label{week34}}

A bit of a miscellany this week\ldots.

\begin{enumerate}
\def\labelenumi{\arabic{enumi})}
\tightlist
\item
  ``Algorithms for quantum computation: discrete log and factoring'',
  extended abstract by Peter Shor.
\end{enumerate}
\noindent
There has been a bit of a stir about this paper; since I know Peter
Shor's sister I was able to get a copy and see what it was really all
about.

Quantum computers are so far just a theoretical possibility. It's
easiest to see why machines that take advantage of quantum theory might
be efficient at computation if we think in terms of path integrals. In
Feynman's path-integral approach to quantum theory, the probability of
getting from state \(A\) at time zero to state \(B\) some later time is
obtained by integrating the exponential of the action 
\[\exp(iS/\hbar)\]
over \emph{all} paths from \(A\) to \(B\), and then taking the absolute
value squared. (Here we are thinking of states \(A\) and \(B\) that
correspond to points in the classical configuration space.) We can think
of the quantum system as proceeding along all paths simultaneously; it
is the constructive or destructive interference between paths due to the
phases \(exp(iS/\hbar)\) that makes certain final outcomes \(B\) more
likely than others. In many situations, there is massive destructive
interference except among paths very close to those which are critical
points of the action \(S\); the latter are the \emph{classical} paths.
So in a sense, a classical device functions as it does by executing all
possible motions; motions far from those satisfying Newton's laws simply
cancel out by destructive interference! (There are many other ways of
thinking about quantum theory; this one can be difficult to make
mathematically rigorous, but it's often very handy.)

This raises the idea of building a computer that would take advantage of
quantum theory by trying out all sorts of paths, but making sure that
paths that give the wrong answer cancel out! It seems that Feynman was
the first to seriously consider quantum computation:

\begin{enumerate}
\def\labelenumi{\arabic{enumi})}
\setcounter{enumi}{1}
\tightlist
\item
   Richard Feynman, ``Simulating physics with computers'',
  \emph{International Journal of Theoretical Physics},  \textbf{21},
  (1982), 467--488.
\end{enumerate}
\noindent
but by now quite a bit of work has been done on the subject, e.g.:

\begin{enumerate}
\def\labelenumi{\arabic{enumi})}
\setcounter{enumi}{2}
\item
    P.\  Benioff,
  ``Quantum mechanical Hamiltonian models of Turing machines'', \emph{J.\ Stat.\ Phys.} 
   \textbf{29} (1982), 515--546.

  D.\ Deutsch, ``Quantum theory, the Church--Turing principle and the universal
  quantum computer'',  \emph{Proc. R. Soc. Lond. A}
  \textbf{400} (1985), 96--117.

  D.\ Deutsch, Quantum computational networks, \emph{Proc. R. Soc.
  Lond. A} \textbf{425} (1989), 73--90.

  D.\ Deutsch and R.\ Jozsa, Rapid solution of problems by quantum computation, 
  \emph{Proc. R. Soc. Lond. A } \textbf{439} (1992), 553--558.

  E.\ Bernstein and U.\ Vazirani, Quantum complexity theory, \emph{Proc.
  25th ACM Symp. on Theory of Computation} (1993), 11--20.

   A.\ Berthiaume and G.\ Brassard, editors, ``The quantum challenge to structural complexity theory'', 
   in  \emph{Proc.\ 7th IEEE Conference on
  Structure in Complexity Theory}, 1992.

   A.\ Yao, ``Quantum circuit complexity'', in \emph{Proc.\ 34th IEEE Symp.\
  on Foundations of Computer Science}, 1993.
\end{enumerate}
\noindent
Thanks to this work, there are now mathematical definitions of quantum
Turing machines and the class ``BQP'' of problems that can be solved in
polynomial time with a bounded probability of error. This allows a
mathematical investigation of whether quantum computers can, in
principle, do things more efficiently than classical ones. Shor shows
that factoring integers is in BQP. I won't try to describe how, as it's
a bit technical and I haven't really comprehended it. Instead, I'd like
say a couple things about the \emph{practicality} of building quantum
computers, since people seem quite puzzled about this issue.

There are, as I see it, two basic problems with building quantum
computers. First, it seems that the components must be carefully
shielded from unwanted interactions with the outside world, since such
interactions can cause ``decoherence'', that is, superpositions of the
computer states will evolve into superpositions of the system consisting
of the computer together with what it's interacting with, which from the
point of view of the computer alone are the same as mixed states. This
tends to ruin the interference effects upon which the advantages of
quantum computation are based.

Second, it seems difficult to incorporate error-correction mechanisms in
a quantum computer. Without such mechanisms, slight deviations of the
Hamiltonian of the computer from the design specifications will cause
the computation to drift away from what was intended to occur. Luckily,
it appears that this drift is only \emph{linear} rather than
\emph{exponential} as a function of time. (This impression is based on
some simplifications that might be oversimplifications, so anyone who
wants to build a quantum computer had better ponder this issue
carefully.) Linear increase of error with time sets an upper bound on
how complicated a computation one could do before the answer is junk,
but if the rate of error increase was made low enough, this might be
acceptable.

Certainly as time goes by and computer technology becomes ever more
miniaturized, hardware designers will have to pay ever more attention to
quantum effects, for good or for ill! (Vaughn Pratt estimates that
quantum effects will be a serious concern by 2020.) The question is just
whether they are only a nuisance, or whether they can possibly be
harnessed. Some designs for quantum computers have already been proposed
(sorry, I have no reference for these), and seeing whether they are
workable should be a very interesting engineering problem, even if they
are not good enough to outdo ordinary computers.

\begin{enumerate}
\def\labelenumi{\arabic{enumi})}
\setcounter{enumi}{3}
\tightlist
\item
  Lee Smolin and Chopin Soo,  ``The Chern--Simons invariant as the natural time variable for
  classical and quantum cosmology'', available as
  \href{https://arxiv.org/abs/gr-qc/9405015}{\texttt{gr-qc/9405015}}.
\end{enumerate}
\noindent
Let me just quote the abstract on this one:

\begin{quote}
We propose that the Chern--Simons invariant of the Ashtekar--Sen
connection (or its imaginary part in the Lorentzian case) is the natural
internal time coordinate for classical and quantum cosmology. The
reasons for this are: 1) It is a function on the gauge and
diffeomorphism invariant configuration space, whose gradient is
orthogonal to the two physical degrees of freedom, in the metric defined
by the Ashtekar formulation of general relativity. 2) The imaginary part
of the Chern--Simons form reduces in the limit of small cosmological
constant, \(\Lambda\), and solutions close to DeSitter spacetime, to the
York extrinsic time coordinate. 3) Small matter-field excitations of the
Chern--Simons state satisfy, by virtue of the quantum constraints, a
functional Schr\"odinger equation in which the matter fields evolve on a
DeSitter background in the Chern--Simons time. We then propose this is
the natural vacuum state of the theory for nonzero \(\Lambda\). 4) This
time coordinate is periodic on the Euclidean configuration space, due to
the large gauge transformations, which means that physical expectation
values for all states in non-perturbative quantum gravity will satisfy
the KMS condition, and may then be interpreted as thermal states.
Finally, forms for the physical hamil- tonian and inner product are
suggested and a new action principle for general relativity, as a
geodesic principle on the connection superspace, is found.
\end{quote}

\begin{enumerate}
\def\labelenumi{\arabic{enumi})}
\setcounter{enumi}{4}
\tightlist
\item
  ``Symplectic geometry'', a series of lectures by Mikhail Gromov,
  compiled by Richard Brown, edited by Robert Miner.
 
\end{enumerate}
\noindent
Symplectic geometry is the geometry of classical phase spaces. That is,
it's the geometry of spaces on which one can take Poisson brackets of
functions in a manner given locally by the usual formulas. Gromov has
really revolutionized the subject, and these lectures look like a good
place to begin learning what is going on. There is also an appendix on
contact geometry (another aspect of classical physics) based on a
lecture by Eliashberg.

\hypertarget{week35}{%
\section{June 5, 1994}\label{week35}}

\begin{enumerate}
\def\labelenumi{\arabic{enumi})}
\tightlist
\item
  Alexander Grothendieck, \emph{Pursuing Stacks (A la Poursuite des Champs)}, 1983, 
  593 pages.  Available as
  \href{https://thescrivener.github.io/PursuingStacks/ps-online.pdf}{\texttt{https://thescrivener.github.io/PursuingStacks/ps-online.pdf}}.
\end{enumerate}
\noindent
I owe Ronnie Brown enormous thanks for sending this to me (before it was
available online).   Grothendieck
is mainly famous for his work on algebraic geometry, in which he
introduced the concept of ``schemes'' to provide a modern framework for
the subject. He was also interested in reformulating the foundations of
topology, which is reflected in \emph{Pursuing Stacks}. This is a long
letter to Quillen, inspired by Quillen's 1967 book \emph{Homotopical
Algebra}. It's a fascinating mixture of visionary mathematics, general
philosophy and a bit of personal chat. Let me quote a bit:

\begin{quote}
I write you under the assumption that you have not entirely lost
interest for those foundational questions you were looking at more than
fifteen years ago. One thing which strikes me, is that (as far as I
know) there has not been any substantial progress since --- it looks to
me that understanding of the basic structures underlying homotopy
theory, or even homological algebra only, is still lacking --- probably
because the few people who have a wide enough background and perspective
enabling them to feel the main questions, are devoting their energies to
things which seem more directly rewarding. Maybe even a wind of
disrepute for any foundational matters whatever is blowing nowadays! In
this respect, what seems to me even more striking than the lack of
proper foundations for homological and homotopical algebra, is the
absence I daresay of proper foundations for topology itself! I am
thinking here mainly of the development of a context of ``tame''
topology, which (I am convinced) would have on the everyday technique of
geometric topology (I use this expression in contrast to the topology of
use for analysts) a comparable impact or even a greater one, than the
introduction of the point of view of schemes had on algebraic geometry.
The psychological drawback here I believe is not anything like
messyness, as for homological and homotopical algebra (or for schemes),
but merely the inrooted inertia which prevents us so stubbornly from
looking innocently, with fresh eyes, upon things, without being dulled
and emprisoned by standing habits of thought, going with a familiar
context --- \emph{too} familiar a context!
\end{quote}

One reason why I'm interested in this letter is that Grothendieck seems
to have understood the importance of ``higher algebraic structures''
before most people. Recently, interest in these has been heating up,
largely because of the recent work on ``extended topological quantum
field theories.'' The basic idea is that, just as a traditional quantum
field theory is (among other things) a representation of the symmetry
group of spacetime, a topological quantum field theory is a
representation of a more sophisticated algebraic structure, a
``cobordism \(n\)-category.'' An \(n\)-category is a wonderfully
recursive sort of thing in which there are objects, \(1\)-morphisms
between objects, 2-morphisms between morphisms, and so on up to
\(n\)-morphisms. In a ``cobordism \(n\)-category'' the objects are
0-manifolds, the \(1\)-morphisms are \(1\)-dimensional manifolds that go
between 0-manifolds (as the unit interval goes from one endpoint to
another), the \(2\)-morphisms are 2-dimensional manifolds that go
between 1-manifolds (as a cylinder goes from on circle to another), etc.
In practice one must work with manifolds admitting certain types of
``corners'', and equipped with extra structures that topologists and
physicists like, such as orientations, framings, or spin structures. The
idea is that all the cutting-and-pasting constructions in
\(n\)-dimensional topology can be described algebraically in the
cobordism \(n\)-category. To wax rhapsodic for a moment, we can think of
an \(n\)-category as exemplifying the notion of ``ways to go between
ways to go between ways to go between... ways to go between
things,'' and cobordism \(n\)-categories are the particular
\(n\)-categories that algebraically encode the possibilities along these
lines that are implicit in the notion of \(n\)-dimensional spacetime.

Now, the problem is that the correct \emph{definition} of an
\(n\)-category is a highly nontrivial affair! And it gets more
complicated as \(n\) increases! A 0-category is nothing but a bunch of
objects. In other words, it's basically just a \emph{set}, if we allow
ourselves to ignore certain problems about classes that are too big to
qualify as sets. A \(1\)-category is nothing but a category. Recall the
definition of a
\href{http://math.ucr.edu/home/baez/categories.html}{category}:

A category consists of a set of \textbf{objects} and a set of
\textbf{morphisms}. Every morphism has a \textbf{source} object and a
\textbf{target} object. (A good example to think of is the category in
which the objects are sets and the morphisms are functions. If
\(f\colon X\to Y\), we call \(X\) the source and \(Y\) the target.)
Given objects \(X\) and \(Y\), we write \(\mathrm{Hom}(X,Y)\) for the
set of morphisms from \(X\) to \(Y\) (i.e., having \(X\) as source and
\(Y\) as target).

The axioms for a category are that it consist of a set of objects and
for any 2 objects \(X\) and \(Y\) a set \(\mathrm{Hom}(X,Y)\) of
morphisms from \(X\) to \(Y\), and

\begin{enumerate}
\def\labelenumi{\arabic{enumi}.}
\tightlist
\item
  Given a morphism \(g\) in \(\mathrm{Hom}(X,Y)\) and a morphism \(f\)
  in \(\mathrm{Hom}(Y,Z)\), there is morphism which we call
  \(f \circ g\) in \(\mathrm{Hom}(X,Z)\). (This binary operation
  \(\circ\) is called \textbf{composition}.)
\item
  Composition is associative:
  \((f \circ g) \circ h = f \circ (g \circ h)\).
\item
  For each object \(X\) there is a morphism \(\mathrm{id}_X\) from \(X\)
  to \(X\), called the \textbf{identity on} \(X\).
\item
  Given any \(f\) in \(\mathrm{Hom}(X,Y)\),
  \(f \circ \mathrm{id}_X = f\) and \(\mathrm{id}_Y \circ f = f\).
\end{enumerate}

Now, a \(2\)-category is more complicated. There are objects,
\(1\)-morphisms, and \(2\)-morphisms, and one can compose morphisms and
also compose 2-morphisms. There is, however, a choice: one can make ones
\(2\)-category ``strict'' and require that the rules 2) and 4) above
hold for the 1-morphisms and \(2\)-morphisms, or one can require them
``literally'' only for the \(2\)-morphisms, and allow the
\(1\)-morphisms some slack. Technically, one can choose between
``strict'' \(2\)-categories, usually just called \(2\)-categories, or
``weak'' ones, which are usually called ``bicategories.''

What do I mean by giving the \(1\)-morphisms some ``slack''? This is a
very important aspect of the \(n\)-categorical philosophy\ldots{} I mean
that in a \(2\)-category one has the option of replacing
\emph{equations} between 1-morphisms by \emph{isomorphisms} --- that is,
by \(2\)-morphisms that have inverses! The basic idea here is that in
many situations when we like to pretend things are equal, they are
really just \emph{isomorphic}, and we should openly admit this when it
occurs. So, for example, in a ``weak'' \(2\)-category one doesn't have
associativity of \(1\)-morphisms. Instead, one has ``associators'',
which are \(2\)-morphisms like this:
\[a_{f,g,h}: (f \circ g) \circ h \to f \circ (g \circ h)\] In other
words, the associator is the \emph{process of rebracketing} made
concrete. Now, when one replaces equations between \(1\)-morphisms by
isomorphisms, one needs these isomorphisms to satisfy ``coherence
relations'' if we're going to expect to be able to manipulate them more
or less as if they \emph{were} equations. For example, in the case of
the associators above, one can use associators to go from
\[f \circ (g \circ (h \circ k))\] to \[((f \circ g) \circ h) \circ k\]
in two different ways: either
\[f \circ (g \circ (h \circ k)) \to (f \circ g) \circ (h \circ k) \to ((f \circ g) \circ h) \circ k\]
or
\[f \circ (g \circ (h \circ k)) \to f \circ ((g \circ h) \circ k) \to (f \circ (g \circ h)) \circ k \to ((f \circ g) \circ h) \circ k\]
Actually there are other ways, but in an important sense these are the
basic two. In a ``weak'' \(2\)-category one requires that these two ways
are equal\ldots{} i.e., this is an identity that the associator must
satisfy, known as the pentagon identity. This is one of the first
examples of a coherence relation. It turns out that if this holds,
\emph{all} ways of rebracketing that get from one expression to another
are equal. (Here I'm being rather sloppy, but the precise result is
known as Mac Lane's theorem.)

To learn about weak \(2\)-categories, which as I said people usually
call bicategories, try:

\begin{enumerate}
\def\labelenumi{\arabic{enumi})}
\setcounter{enumi}{1}
\tightlist
\item
  J. Benabou, \emph{Introduction to Bicategories}, Springer Lecture Notes in Mathematics,
  \textbf{47},  Springer, Berlin, 1968.
\end{enumerate}
\noindent
Now, one can continue this game, but it gets increasingly complex if one
goes the ``weak'' route. In a ``weak \(n\)-category'' the idea is to
replace all basic identities that one might expect between
\(j\)-morphisms, such as the associative law, by \((j+1)\)-isomorphisms.
These, in turn, satisfy certain ``coherence relations'' that are really
not equations, but \((j+2)\)-morphisms, and so on\ldots{} up to level
\(n\). This becomes so complicated that only recently have ``weak
\(3\)-categories'' been properly defined, by Gordon, Power and Street,
who call them tricategories (see \protect\hyperlink{week29}{``Week
29''}).

A bit earlier, Kapranov and Voevodsky succeeded in defining a certain
class of weak \(4\)-categories, which happen to be called ``braided
monoidal \(2\)-categories'' (see \protect\hyperlink{week4}{``Week 4''}).
The interesting thing, you see, which justifies getting involved in this
business, is that a lot of topology \emph{automatically pops out} of the
definition of an \(n\)-category. In particular, \(n\)-categories have a
lot to do with \(n\)-dimensional space. A weak \(3\)-category with only
one object and one \(1\)-morphism is usually known as a ``braided
monoidal category,'' and the theory of these turns out to be roughly the
same as the study of knots, links and tangles!  The
``braided monoidal \(2\)-categories'' of Kapranov and Voevodsky are
really just weak \(4\)-categories with only one object and one
1-morphism. (The reason for the term ``\(2\)-category'' here is that
since all one has is \(2\)-morphisms, \(3\)-morphisms, and
\(4\)-morphisms, one can pretend one is in a \(2\)-category in which
those are the objects, morphisms, and \(2\)-morphisms.)

In any event, these marvelous algebraic structures have been cropping up
more and more in physics (see especially Crane's stuff listed in
\protect\hyperlink{week2}{``Week 2''} and Freed's paper listed in
\protect\hyperlink{week12}{``Week 12''}), so I got ahold of a copy of
Grothendieck's letter and have begun trying to understand it.

Actually, it's worth noting that these \(n\)-categorical ideas have been
lurking around homotopy theory for quite some time now. As Grothendieck
wrote:

\begin{quote}
At first sight it had seemed to me that the Bangor group had indeed come
to work out (quite independently) one basic intuition of the program I
had envisioned in those letters to Larry Breen --- namely, that the
study of \(n\)-truncated homotopy types (of semisimplicial sets, or of
topological spaces) was essentially equivalent to the study of so-called
\(n\)-groupoids (where \(n\) is any natural integer). This is expected
to be achieved by associating to any space (say) \(X\) its ``fundamental
\(n\)-groupoid'' \(\Pi_n(X)\), generalizing the familiar Poincar\'e
fundamental groupoid for \(n = 1\). The obvious idea is that 0-objects
of \(\Pi_n(X)\) should be the points of \(X\), 1-objects should be
``homotopies'' or paths between points, 2-objects should be homotopies
between 1-objects, etc. This \(\Pi_n(X)\) should embody the
\(n\)-truncated homotopy type of \(X\), in much the same way as for
\(n = 1\) the usual fundamental groupoid embodies the 1-truncated
homotopy type. For two spaces \(X\), \(Y\), the set of homotopy-classes
of maps \(X \to Y\) (more correctly, for general \(X\), \(Y\), the maps
of \(X\) into \(Y\) in the homotopy category) should correspond to
\(n\)-equivalence classes of \(n\)-functors from \(\Pi_n(X)\) to
\(\Pi_n(Y)\) --- etc. There are some very strong suggestions for a nice
formalism including a notion of geometric realization of an
\(n\)-groupoid, which should imply that any \(n\)-groupoid is
\(n\)-equivalent to a \(\Pi_n(X)\). Moreover when the notion of an
\(n\)-groupoid (or more generally of an \(n\)-category) is relativized
over an arbitrary topos to the notion of an \(n\)-gerbe (or more
generally, an \(n\)-stack), these become the natural ``coefficients''
for a formalism of non commutative cohomological algebra, in the spirit
of Giraud's thesis.
\end{quote}
\noindent
The ``Bangor group'' referred to includes Ronnie Brown, who has done a
lot of work on ``\(\omega\)-groupoids''. A while back he sent me a nice
long list of references on this subject; here are some that seemed
particularly relevant to me (though I haven't looked at all of them).

\begin{enumerate}
\def\labelenumi{\arabic{enumi})}
\setcounter{enumi}{2}
\item
  G.\ Abramson, J.-P.\ Meyer and J.\ Smith, ``A higher dimensional analogue of
  the fundamental groupoid'', in \emph{Recent Developments of Algebraic
  Topology}, RIMS Kokyuroku \textbf{781}, Kyoto, 1992, pp.\ 38--45, 

  F.\ Al-Agl, ``Aspects of multiple categories'', University of Wales Ph.D.
  Thesis, 1989.

  F.\ Al-Agl and R.\ J.\ Steiner, ``Nerves of multiple categories'',
  \emph{Proc. London Math. Soc.} \textbf{66} (1992), 92--128.

  N.\ Ashley, ``Simplicial T-complexes'', University of Wales Ph.D. Thesis,
  1976, published as \textsl{Simplicial T-complexes and Crossed Complexes: a 
  Non-Abelian Version of
  a Theorem of Dold and Kan}, Dissertationes Mathematicae \textbf{265}, 1988.

  H.\ J.\ Baues, \emph{Algebraic Homotopy}, Cambridge University Press, Cambridge,
  1989.

  H.\ J.\ Baues, \emph{Combinatorial Homotopy and \(4\)-Dimensional
  Complexes}, De Gruyter, Berlin, 1991.

  L.\ Breen, ``Bitorseurs et cohomologie non-Ab\'elienne'', in \emph{The
  Grothendieck Festschrift: a collection of articles written in honour
  of the 60th birthday of Alexander Grothendieck, Vol. I}, ed.\
  P.\ Cartier \textit{et al.}, Birkh\"auser, Basel, 1990, pp.\ 401--476.

  R.\ Brown, ``Higher dimensional group theory'', in \emph{Low-dimensional
  topology}, ed.~R.\ Brown and T.L.Thickstun, \emph{London Math. Soc.
  Lect. Notes} \textbf{46}, Cambridge University Press, 215--238, 1982.

  R.\ Brown, ``From groups to groupoids: a brief survey'', \emph{Bull.
  London Math. Soc.}, \textbf{19}, 113--134, 1987.

  R.\ Brown, \emph{Elements of Modern Topology}, McGraw Hill, Maidenhead,
  1968.
  
   R.\ Brown, \emph{Topology: a Geometric Account of General Topology,
  Homotopy Types and the Fundamental Groupoid}, Ellis Horwood,
  Chichester, 1988.

  R.\ Brown, ``Some problems in non-Abelian homological and homotopical
  algebra'', in \emph{Homotopy Theory and Related Topics: Proceedings
  Kinosaki}, ed.\ M.\ Mimura, Springer Lecture Notes in
  Mathematics \textbf{1418}, 1990, Springer, Berlin pp.\ 105--129.

  R.\ Brown and P.\ J.\ Higgins, ``The equivalence of \(\omega\)-groupoids and
  cubical T-complexes'', \emph{Cah. Top. Geom. Diff.} \textbf{22} (1981),
  349--370.

  R.\ Brown and P.\ J.\ Higgins, ``The equivalence of \(\infty\)-groupoids and
  crossed complexes'', \emph{Cah. Top. Geom. Diff.} \textbf{22} (1981),
  371--386.

  R.\ Brown and P.\ J.\ Higgins, ``The algebra of cubes'', \emph{J. Pure Appl.
  Algebra} \textbf{21} (1981), 233--260.

  R.\ Brown and P.\ J.\ Higgins, ``Tensor products and homotopies for
  \(\omega\)-groupoids and crossed complexes'', \emph{J. Pure Appl.
  Algebra}, \textbf{47} (1987), 1--33.

  R.\ Brown and J.\ Huebschmann, ``Identities among relations'', in
  \emph{Low-dimensional topology}, eds.~R.\ Brown and T.\ L.\ Thickstun,
  London Math. Soc. Lect. Notes \textbf{46}, Cambridge University
  Press, Cambridge, 1982, pp.\ 153--202.

  R.\ Brown, ``Generalised group presentations'', \emph{Trans. Amer.
  Math. Soc.}, \textbf{334} (1992), 519--549.

  M.\ Bullejos, A.\ M.\ Cegarra and J.\ Duskin, ``On \(\mathrm{cat}^n\)-groups and
  homotopy types'', \emph{J. Pure Appl. Algebra} \textbf{86} (1993),
  135--154.

  M.\ Bullejos, P.\ Carrasco and A.\ M.\ Cegarra, ``Cohomology with coefficients in
  symmetric \(\mathrm{cat}^n\)-groups. An extension of Eilenberg-Mac
  Lane's classification theorem.'' Granada Preprint, 1992.

  P.\ J.\ Ehlers and T.\ Porter, ``From simplicial groupoids to crossed
  complexes'', UCNW Maths Preprint \texttt{92.19}, 1992.

  D.\ W.\ Jones, ``Polyhedral T-complexes'', University of Wales Ph.D. Thesis,
  1984; published as \emph{A General Theory of Polyhedral Sets and Their
  Corresponding T-complexes}, Dissertationes Mathematicae \textbf{266}, 1988.

  M.\ M.\ Kapranov and V.\ Voevodsky, ``Combinatorial-geometric aspects of
  polycategory theory: pasting schemes and higher Bruhat orders (list of
  results)'', \emph{Cah. Top. Geom. Diff. Cat.} \textbf{32} (1991), 11--27.

  M.\ M.\ Kapranov and V.\ Voevodsky, ``\(\infty\)-groupoids and homotopy
  types'', \emph{Cah. Top. Geom. Diff. Cat.} \textbf{32}, 29--46, 1991.

  M.\ M.\ Kapranov and V.\ Voevodsky, ``\(2\)-categories and Zamolodchikov
  tetrahedra equations'', preprint, 1992.

  J.-L.\ Loday, ``Spaces with finitely many non-trivial homotopy groups'',
  \emph{J. Pure Appl. Algebra} \textbf{24} (1982), 179--202.

  G.\ Nan Tie, ``Iterated W and T-groupoids'', \emph{J. Pure Appl.
  Algebra} \textbf{56} (1989), 195--209.

  T.\ Porter, ``A combinatorial definition of \(\infty\)-types'',
  \emph{Topology} \textbf{22} (1993), 5--24.

  S.\ J.\ Pride, ``Identities among relations of group presentations'', in
  \emph{Proc. Workshop on Group
  Theory from a Geometrical Viewpoint}, eds.\ E.\ Ghys, A.\ Haefliger, A.\ Verjodsky, 
    World Scientific (1991), 687--716.

  R.\ Steiner, ``The algebra of directed complexes'', University of
  Glasgow Math Preprint, 1992.

  A.\ Tonks, ``Cubical groups which are Kan'', \emph{J. Pure Appl.
  Algebra} \textbf{81} (1992), 83--87.

  A.\ Tonks and R.\ Brown, ``Calculation with simplicial and cubical groups
  in Axiom'', UCNW Math Preprint \texttt{93.04}.

  A.\ R.\ Wolf, ``Inherited asphericity, links and identities among
  relations'', \emph{J. Pure Appl. Algebra} \textbf{71} (1991), 99--107.
\end{enumerate}

\begin{center}\rule{0.5\linewidth}{0.5pt}\end{center}

\begin{quote}
\emph{Since the month of March last year, so nearly a year ago, the
greater part of my energy has been devoted to a work of reflection on
the \textbf{foundations of non-commutative (co)homological algebra}, or
what is the same, after all, of homotopic algebra. These reflections
have taken the concrete form of a voluminous stack of typed notes,
destined to for then first volume (now being finished) of a work in two
volumes to be published by Hermann, under the overall title}
\textbf{Pursuing Stacks}\emph{. I now foresee (after successive
extensions of the initial project) that the manuscript of the whole of
the two volumes, which I hope to finish definitively in the course of
this year, will be about 1500 typed pages in length. These two volumes
are moreover for me the first in a vaster series, under the overall
title} \textbf{Mathematical Reflections}\emph{, in which I intend to
develop some of the themes sketched in the present report}

--- Alexander Grothendieck,
\textsl{Sketch of a Program} (1983)
\end{quote}

\hypertarget{week36}{%
\section{July 15, 1994}\label{week36}}

I am attempting to keep my nose to the grindstone these days, in part
since I'm getting ready for the Knots and Quantum Gravity session of the
Marcel Grossman meeting on general relativity, which will take place at
Stanford the week after next (I will report any interesting news I hear
out there), and in part to make up for earlier stretches of laziness on
my part\ldots. Nonetheless, I feel I should describe a few new papers on
topological quantum field theories.

The real reason for physicists' interest in topological quantum field
theories (TQFTs) is that we sorely need a mathematical framework that
quantum gravity might fit into. It's likely, however, that quantum
gravity won't be much like the TQFTs people have studied so far. The
existing TQFTs tend to be ``exactly soluble'' and have
finite-dimensional state spaces; quantum gravity is likely to be
different. Still, any experience in studying quantum field theories that
don't rely on ``fixed background structures'' like a fixed spacetime
metric is likely to be worth having. Also, quantum gravity appears to be
tied mathematically to simpler TQFTs in a variety of ways. In
particular, the Ashtekar formulation of quantum gravity is closely
related to a \(4\)-dimensional TQFT variously known as ``\(B \wedge F\)
theory,'' ``\(BF\) theory,'' ``topological \(2\)-form gravity'' or
``topological quantum gravity''. This in turn is closely related to
Chern--Simons theory in 3 dimensions.

Let me just say what the heck BF theories \emph{are}, and then list a
few references on them. The easiest way to describe them is by giving
the Lagrangian. Say spacetime is an \(n\)-dimensional orientable
manifold \(M\) and we have a principal \(G\)-bundle \(E\) over \(M\),
where \(G\) is a Lie group whose Lie algebra is equipped with an
invariant trace on it. The two fields in BF theory are a connection
\(A\) on \(E\) --- which we can think of locally as a
\(\mathrm{Lie}(G)\)-valued one-form --- and a \(\mathrm{Lie}(G)\)-valued
\((n-2)\)-form called \(B\). If \(F\) denotes the curvature of \(A\),
which is a \(\mathrm{Lie}(G)\)-valued \(2\)-form, we can take the wedge
product \(B\wedge F\) and get a \(\mathrm{Lie}(G)\)-valued \(n\)-form,
which gives the Lagrangian \[\operatorname{tr}(B \wedge F),\] an
\(n\)-form we can integrate over \(M\) to get the action. Since we don't
need any metric to integrate an \(n\)-form over our spacetime \(M\),
this action is ``generally covariant''. (Note also that it's
gauge-invariant.) If we vary \(B\) and \(F\) to get the classical
equations of motion, varying \(B\) gives us \[F = 0,\] that is, the
connection \(A\) is flat, and varying \(A\) gives us \[d_A B = 0,\] that
is, the exterior covariant derivative of \(B\) vanishes.

In 4 dimensions we can soup up our \(BF\) theory a bit by adding terms
proportional to \[\operatorname{tr}(B \wedge B)\] and
\[\operatorname{tr}(F \wedge F)\] to the Lagrangian. If we take as our
Lagrangian
\[\operatorname{tr}(B \wedge F) + a \operatorname{tr}(B \wedge B) + b \operatorname{tr}(F \wedge F),\]
the third term, when integrated over \(M\), is proportional to an
invariant called the second Chern class of \(E\), that is, it's
independent of the connection \(A\), so it really doesn't affect the
equations of motion at all! In a sense it's utterly useless. The second
term does something, though; our equations of motion become
\[F = -2aB, \quad d_A B = 0.\] Note that if we plug the first equation
into the Lagrangian, we get that for solutions of the equations of
motion, the action is a constant times the second Chern class (if \(a\)
is nonzero).

Horowitz showed, in this four-dimensional case, that if \(a\) is
nonzero, there is a single state of the \emph{quantum} version of \(BF\)
theory when spacetime has the form \(\mathbb{R} \times S\) (\(S\) being
some oriented 3-manifold), and that this state, thought of as a
wavefunction on the space of connections on \(S\), is just the
exponential of the Chern--Simons action. This is how Chern--Simons theory
gets into the game.

Moreover, Baulieu and Singer showed that if you take the boring-looking
``FF theory'' with Lagrangian \(\operatorname{tr}(F \wedge F)\), and
quantize it using the BRST approach, you get something that Witten
proved was closely related to Donaldson theory --- an invariant of
4-manifolds. So there seems to be a relation between this stuff and
Donaldson theory. It is a rather mysterious one as far as I'm concerned,
though, because it seems you could just as well have used \emph{zero} as
a Lagrangian, rather than \(\operatorname{tr}(F \wedge F)\), and done
the same things Baulieu and Singer did. (At least, that's how it seems
to me, and I got Scott Axelrod to agree with me on that yesterday.) In
other words, Donaldson theory seems to have more to do with the geometry
of the space of connections on \(M\) than it does with the ``FF''
Lagrangian per se. But still, there are other relationships between
Donaldson theory and Chern--Simons theory (which I don't understand well
enough to want to discuss), so perhaps it is not too silly to say that
\(BF\) theory is related to Donaldson theory in some poorly understood
manner.

Now for some references: the reference that got me started on these was

\begin{enumerate}
\def\labelenumi{\arabic{enumi})}
\tightlist
\item
Gary Horowitz, ``Exactly soluble diffeomorphism-invariant theories'',
\emph{Commun. Math. Phys.} \textbf{125} (1989), 417--437. (Listed in
\protect\hyperlink{week19}{``Week 19''})
\end{enumerate}
\noindent
I got more interested in them when I read

\begin{enumerate}
\def\labelenumi{\arabic{enumi})}
\setcounter{enumi}{1}
\tightlist
\item
Paolo Cotta-Ramusino and Maurizio
Martellini, ``BF Theories and 2-knots'' in \emph{Knots and Quantum Gravity}, ed.~J. Baez,
Oxford U. Press, Oxford, 1994.  (Listed in \protect\hyperlink{week23}{``Week 23''})
\end{enumerate}
\noindent
which indicated that BF theories may give invariants of surfaces
embedded in \(4\)-dimensional manifolds, much as Chern--Simons theory
gives invariants of knots in \(3\)-dimensional manifolds. Actually, BF
theories make sense in any dimension, and in dimension 3 they appear to
give knot invariants, including the Alexander-Conway polynomial:

\begin{enumerate}
\def\labelenumi{\arabic{enumi})}
\setcounter{enumi}{2}
\tightlist
\item
  A.\ S.\ Cattaneo, P.\ Cotta-Ramusino, and M.\ Martellini, ``Three-dimensional 
  BF theories and the Alexander--Conway invariant of knots'', available as
  \href{https://arxiv.org/ps/hep-th/9407070}{\texttt{hep-th/9407070}}.
\end{enumerate}
\noindent
Another nice-looking paper on BF theories is the following:

\begin{enumerate}
\def\labelenumi{\arabic{enumi})}
\setcounter{enumi}{3}
\tightlist
\item
  Henri Waelbroeck, ``\(B \wedge F\) theory and flat spacetimes'', available as
  \href{https://arxiv.org/ps/gr-qc/9311033}{\texttt{gr-qc/9311033}}.
\end{enumerate}
\noindent
Waelbrock also has a recent paper with Zapata on a Hamiltonian
formulation of the theory on a lattice:

\begin{enumerate}
\def\labelenumi{\arabic{enumi})}
\setcounter{enumi}{2}
\tightlist
\item
  Henri  Waelbrock and J. A. Zapata, ``A Hamiltonian formulation of topological gravity'', 
  available as
  \href{https://arxiv.org/ps/gr-qc/9311035}{\texttt{gr-qc/9311035}}.
\end{enumerate}
\noindent
The paper by Baulieu and Singer relating FF theory to Donaldson theory
is:

\begin{enumerate}
\def\labelenumi{\arabic{enumi})}
\setcounter{enumi}{3}
\tightlist
\item
  L. Baulieu and I. M. Singer, ``Topological Yang--Mills symmetry'',
  \emph{Nucl. Phys. (Proc. Suppl.)} \textbf{B5} (1988), 12--19.
\end{enumerate}
\noindent
\(BF\) theory in 2 dimensions is also called ``topological Yang--Mills
theory'', and it's discussed in various places, including

\begin{enumerate}
\def\labelenumi{\arabic{enumi})}
\setcounter{enumi}{4}
\tightlist
\item
   Edward Witten, ``On quantum gauge theories in two dimensions'',
  \emph{Commun.\ Math.\ Phys.} \textbf{141} (1991), 153--209.
\end{enumerate}
\noindent
and

\begin{enumerate}
\def\labelenumi{\arabic{enumi})}
\setcounter{enumi}{5}
\tightlist
\item
 Matthias Blau and George Thompson,  ``Topological gauge theories of antisymmetric tensor fields'', \emph{Ann. Phys.} \textbf{205} (1991), 130--172.
\end{enumerate}

\hypertarget{week37}{%
\section{August 10, 1994}\label{week37}}

Mainly this week I have various bits of news to report from the 7th
Marcel Grossman Meeting on general relativity. It was big and had lots
of talks. Bekenstein gave a nice review talk on entropy/area relations
for black holes, and Strominger gave a talk in which he proposed a
solution to the information loss puzzle for black holes. (Recall that if
one believes, as most people seem to believe, that black holes radiate
away all their mass in the form of completely random Hawking radiation,
then there's a question about where the information has gone that you
threw into the black hole in the form of, say, old issues of \emph{Phys.\
Rev.~Lett.}  Some people think the information goes into a new ``baby
universe'' formed at the heart of the black hole --- see
\protect\hyperlink{week31}{``Week 31''} for more. The information would
still, of course, be gone from \emph{our} point of view in this picture.
Strominger proposed a set up in which one had a quantum theory of
gravity with annihilation and creation operators for baby universes, and
proposed that the universe (the ``metauniverse''?) was in a coherent
state, that is, an eigenstate of the annihilation operator for baby
universes. This would apparently allow handle the problem, though right
now I can't remember the details.) There were also lots of talks on the
interferometric detection of gravitational radiation, other general
relativity experiments, cosmology, etc. But I'll just try to describe
two talks in some detail here.

\begin{enumerate}
\def\labelenumi{\arabic{enumi})}
\tightlist
\item
  L. Lindblom, ``Superfluid hydrodynamics and the stability of rotating
  neutron stars'', talk at MG7 meeting, Monday July 5, 1994, Stanford
  University.
\end{enumerate}

Being fond of knots, tangles, and such, I have always liked knowing that
in superfluids, vorticity (the curl of the velocity vector field) tends
to be confined in ``flux tubes'', each containing an angular momentum
that's an integral multiple of Planck's constant, and that similarly, in
type II superconductors, magnetic fields are confined to magnetic flux
tubes. And I was even more happy to find out that the cores of neutron
stars are expected to be made of neutronium that is \emph{both}
superfluid \emph{and} superconductive, and contain lots of flux tubes of
both types. In this talk, which was really about a derivation of
detailed equations of state for neutron stars, Lindblom began by saying
that the maximum rotation rate of a rotating neutron star is due to some
sort of ``gravitational radiation instability due to internal fluid
dissipation''. I didn't quite understand the details of that, which
weren't explained, but it motivated him to study the viscosity in
neutron star cores, which are superfluid if they are cool enough (less
than a billion degrees Kelvin). There are some protons and electrons
mixed in with the neutrons in the core, and both the protons and
neutrons go superfluid, but the electrons form a normal fluid. That
means that there are actually \emph{two} kinds of superfluid vortices
--- proton and neutron --- in addition to the magnetic vortices. These
vortices mainly line up along the axis of rotation, and their density is
about \(10^6\) per square centimeter. Rather curiously, since the
proton, neutron, and electron fluids are coupled due to \(\beta\) decay
(and the reverse process), even the neutron vortices have electric
currents associated to them and generate magnetic fields. This means
that the electrons scatter off the neutron vortex cores as well as the
proton vortex cores, which is one of the mechanisms that yields
viscosity.

\begin{enumerate}
\def\labelenumi{\arabic{enumi})}
\setcounter{enumi}{1}
\tightlist
\item
  Abhay Ashtekar, ``Mathematical developments in quantum general
  relativity, a sampler'', talk at MG7 meeting, Tuesday July 6, 1994,
  Stanford University.
\end{enumerate}
\noindent
This talk, in addition to reviewing what's been done so far on the
``loop representation'' of quantum gravity, presented two new
developments that I found quite exciting, so I'd like to sketch what
they are. The details will appear in future papers by Ashtekar and
collaborators.

The two developments Ashtekar presented concerned mathematically
rigorous treatments of the ``reality conditions'' in his approach to
quantum gravity, and the ``loop states'' used by Rovelli and Smolin.
First let me try to describe the issue of ``reality conditions''. As I
described in \protect\hyperlink{week7}{``Week 7''}, one trick that's
important in the loop representation is to use the ``new variables'' for
general relativity introduced by Ashtekar (though Sen and Plebanski
already had worked with similar ideas). In the older Palatini approach
to general relativity, the idea was to view general relativity as
something like a gauge theory with gauge group given by the Lorentz
group, \(\mathrm{SO}(3,1)\). But to do this one actually uses two
different fields: a ``frame field'', also called a ``tetrad'',
``vierbein'' or ``soldering form'' depending on who you're talking to,
and the gauge field itself, usually called a ``Lorentz connection'' or
``\(\mathrm{SO}(3,1)\) connection''. Technically, the frame field is an
isomorphism between the tangent bundle of spacetime and some other
bundle having a fixed metric of signature \(+--\,-\), usually called the
``internal space'', and the Lorentz connection is a metric-preserving
connection on the internal space.

The ``new variables'' trick is to use the fact that \(\mathrm{SO}(3,1)\)
has as a double cover the group \(\mathrm{SL}(2,\mathbb{C})\) of
two-by-two complex matrices with determinant one. (For people who've
read previous posts of mine, I should add that the Lie algebra of
\(\mathrm{SL}(2,\mathbb{C})\) is called \(\mathfrak{sl}(2,\mathbb{C})\)
and is the same as the complexification of the Lie algebra
\(\mathfrak{so}(3)\), which allows one to introduce the new variables in
a different but equivalent way, as I did in
\protect\hyperlink{week7}{``Week 7''}.) Ignoring topological niceties
for now, this lets one reformulate \emph{complex} general relativity
(that is, general relativity where the metric can be complex-valued) in
terms of a \emph{complex-valued} frame field and an
\(\mathrm{SL}(2,\mathbb{C})\) connection that is just the Lorentz
connection in disguise. The latter is called either the ``Sen
connection'', the ``Ashtekar connection'', or the ``chiral spin
connection'' depending on who you're talking to. The advantage of this
shows up when one tries to canonically quantize the theory in terms of
initial data. (For a bit on this, try \protect\hyperlink{week11}{``Week
11''}.) Here we assume our \(4\)-dimensional spacetime can be split up
into ``space'' and ``time'', so that space is a \(3\)-dimensional
manifold, and we take as our canonically conjugate fields the
restriction of the chiral spin connection to space, call it \(A\), and
something like the restriction of the complex frame field to a complex
frame field \(E\) on space. (Restricting the complex frame field to one
on space is a wee bit subtle, especially because one doesn't really want
a frame field or ``triad field'', but really a ``densitized cotriad
field'' --- but let's not worry about this here. I explain this in terms
even a mathematician can understand in my paper ``Strings, loops, knots
and gauge fields'' (\href{https://arxiv.org/abs/hep-th/9309067}{\texttt{hep-th/9309067}}). 
The point is, first, that the \(A\) and \(E\)
fields are mathematically very analogous to the vector potential and
electric field in electromagnetism --- or really in
\(\mathrm{SL}(2,\mathbb{C})\) Yang--Mills theory --- and secondly, that
if you compute their Poisson brackets, you really do see that they're
canonically conjugate. Third and best of all, the constraint equations
in general relativity can be written down very simply in terms of \(A\)
and \(E\). Recall that in general relativity, 6 of Einstein's 10
equations act as \emph{constraints} that the metric and its time
derivative must satisfy at \(t = 0\) in order to get a solution at later
times. In quantum gravity, these constraints are a big technical problem
one has to deal with, and the point of Ashtekar's new variables is
precisely that the constraints simplify in terms of these variables.
(There's more on these constraints in \protect\hyperlink{week11}{``Week
11''}.)

The price one has paid, however, is that one now seems to be talking
about \emph{complex-valued} general relativity, which isn't what one had
started out being interested in. One needs to get back to reality, as it
were --- and this is the problem of the so-called ``reality
conditions''. One approach is to write down extra constraints on the
\(E\) field that say that it comes from a \emph{real} frame field. These
are a little messy. Ashtekar, however, has proposed another approach
especially suited to the quantum version of the theory, and in his talk
he filled in some of the crucial details.

Here, to save time, I will allow myself to become a bit more technical.
In the quantum version of the theory one expects the space of
wavefunctions to be something like \(L^2\) functions on the space of
connections modulo gauge transformations --- actually this is the
``kinematical state space'' one gets before writing the constraints as
operators and looking for wavefunctions annihilated by these
constraints. The problem had always been that this space of \(L^2\)
functions is ill-defined, since there is no ``Lebesgue measure'' on the
space of connections. This problem is addressed (it's premature to say
``solved'') by developing a theory of generalized measures on the space
of connections and proving the existence of a canonical generalized
measure that deserves the name ``Lebesgue measure'' if anything does.
One can then define \(L^2\) functions and work with them. For compact
gauge groups, like \(\mathrm{SU}(2)\), this was done by Ashtekar,
Lewandowski and myself; see e.g.~the papers ``Spin network states in gauge
theory'' (\href{https://arxiv.org/abs/gr-qc/9411007}{\texttt{gr-qc/9411007}}) and
``Generalized measures in gauge theory'' (\href{https://arxiv.org/abs/hep-th/9310201}
{\texttt{hep-th/9310201}}).  In the case of
\(\mathrm{SU}(2)\), Wilson loops act as self-adjoint multiplication
operators on the resulting \(L^2\) space. But in quantum gravity we
really want to use gauge group \(\mathrm{SL}(2,\mathbb{C})\), which is
not compact, and we want the adjoints of Wilson loop operators to
reflect that fact that the \(\mathrm{SL}(2,\mathbb{C})\) connection
\(A\) in quantum gravity is really equal to \(\Gamma + iK\), where
\(\Gamma\) is the Levi--Civita connection on space, and \(K\) is the
extrinsic curvature. Both \(\Gamma\) and \(K\) are real in the classical
theory, so the adjoint of the quantum version of \(A\) should be
\(\Gamma - iK\), and this should reflect itself in the adjoints of
Wilson loop operators.

The trick, it turns out, is to use some work of Hall which appeared in
the \emph{Journal of Functional Analysis} in 1994 (I don't have a precise
reference on me). The point is that \(\mathrm{SL}(2,\mathbb{C})\) is the
complexification of \(\mathrm{SU}(2)\), and can also be viewed as the
cotangent bundle of \(\mathrm{SU}(2)\). This allows one to copy a trick
people use for the quantum mechanics of a point particle on
\(\mathbb{R}^n\) --- a trick called the Bargmann--Segal-Fock
representation. Recall that in the ordinary Schr\"odinger representation
of a quantum particle on \(\mathbb{R}^n\), one takes as the space of
states \(L^2(\mathbb{R}^n)\). However, the phase space for a particle in
\(\mathbb{R}^n\), which is the cotangent bundle of \(\mathbb{R}^n\), can
be identified with \(\mathbb{C}^n\), and in the Bargmann representation
one takes as the space of states \(HL^2(\mathbb{C}^n)\), by which I mean
the \emph{holomorphic} functions on \(\mathbb{C}^n\) that are in \(L^2\)
with respect to a Gaussian measure on \(\mathbb{C}^n\). In the Bargmann
representation for a particle on the line, for example, the creation
operator is represented simply as multiplication by the complex
coordinate \(z\), while the annihilation operator is \(d/dz\).
Similarly, there is an isomorphism between \(L^2(\mathrm{SU}(2))\) and a
certain space \(HL^2(\mathrm{SL}(2,\mathbb{C}))\). Using this, one can
obtain an isomorphism between the space of \(L^2\) functions on the
space of \(\mathrm{SU}(2)\) connections modulo gauge transformations,
and the space of holomorphic \(L^2\) functions on the space of
\(\mathrm{SL}(2,\mathbb{C})\) connections modulo gauge transformations.
Applying this to the loop representation, Ashtekar has found a very
natural way to take into account the fact that the chiral spin
connection \(A\) is really \(\Gamma + iK\), basically analogous to the
fact that in the Bargmann representation multiplication by \(z\) is really \(q + ip\)
(well, up to various factors of \(\sqrt{2}\), signs and the like).

Well, that was pretty sketchy and probably not especially comprehensible
to anyone who hasn't already worried about this issue a lot! In any
event, let me turn to the other good news Ashtekar reported: the
construction of ``loop states''. Briefly put (I'm getting worn out), he
and some collaborators have figured out how to \emph{rigorously}
construct generalized measures on the space of connections modulo gauge
transformations, starting from invariants of links. This begins to
provide an inverse to the ``loop transform'' (which is a construction
going the other way).

\hypertarget{week38}{%
\section{August 19, 1994}\label{week38}}

I've been busy, and papers have been piling up; there are lots of
interesting ones that I really should describe in detail, but I had
better be terse and list them now, rather than waiting for the mythical
day when I will have time to do them justice.

So:

\begin{enumerate}
\def\labelenumi{\arabic{enumi})}
\tightlist
\item
  B. Durhuus, H. P. Jakobsen and R. Nest, ``Topological quantum field
  theories from generalized 6j-symbols'', \emph{Rev.\ Math.
  Physics} \textbf{5} (1993), 1--67.
\end{enumerate}
\noindent
In \protect\hyperlink{week16}{``Week 16''} I explained a paper by
Fukuma, Hosono and Kawai in which they obtained topological quantum
field theories in 2 dimensions starting with a triangulation of a 2d
surface. The theories were ``topological'' in the sense that the final
answers one computed didn't depend on the triangulation. One can get
between any two triangulations of a surface by using a sequence of the
following two moves (and their inverses), called the (2,2) move: \[
  \begin{tikzpicture}
    \draw[thick] (0,0) node{$\bullet$} to (1,1.5) node{$\bullet$} to (2,0) node{$\bullet$} to (1,-1.5) node{$\bullet$} to cycle;
    \draw[thick] (1,1.5) to (1,-1.5);
  \end{tikzpicture}
  \raisebox{4.5em}{$\qquad\longleftrightarrow\qquad$}
  \begin{tikzpicture}
    \draw[thick] (0,0) node{$\bullet$} to (1,1.5) node{$\bullet$} to (2,0) node{$\bullet$} to (1,-1.5) node{$\bullet$} to cycle;
    \draw[thick] (0,0) to (2,0);
  \end{tikzpicture}
\] and the (3,1) move: \[
  \begin{tikzpicture}
    \draw[thick] (0,0) node{$\bullet$} to (1.8,3) node{$\bullet$} to (3.6,0) node{$\bullet$} to cycle;
    \node at (1.8,1.1) {$\bullet$};
    \draw[thick] (0,0) to (1.8,1.1);
    \draw[thick] (1.8,3) to (1.8,1.1);
    \draw[thick] (3.6,0) to (1.8,1.1);
  \end{tikzpicture}
  \raisebox{4.5em}{$\qquad\longleftrightarrow\qquad$}
  \begin{tikzpicture}
    \draw[thick] (0,0) node{$\bullet$} to (1.8,3) node{$\bullet$} to (3.6,0) node{$\bullet$} to cycle;
  \end{tikzpicture}
\] Note that in either case these moves amount to replacing one part of
the surface of a tetrahedron with the other part! In fact, similar moves
work in any dimension, and they are often called the Pachner moves.

The really \emph{wonderful} thing is that these moves are also very
significant from the point of view of algebra\ldots{} and especially
what I call ``higher-dimensional algebra'' (following Ronnie Brown), in
which the distinction between algebra and topology is largely erased,
or, one might say, revealed for the sham it always was.

For example, as explained more carefully in
\protect\hyperlink{week16}{``Week 16''}, the (2,2) move is really just
the same as the \emph{associative} law for multiplication. The idea is
that we are in a \(2\)-dimensional spacetime, and a triangle represents
multiplication: two ``incoming states'' go in two sides and their
product, the ``outgoing state'', pops out the third side: \[
  \begin{tikzpicture}
    \draw[thick] (0,0) node{$\bullet$} to node[fill=white]{$A$} (2,3) node{$\bullet$} to node[fill=white]{$B$} (4,0) node{$\bullet$} to node[fill=white]{$AB$} cycle;
  \end{tikzpicture}
\] Then the (2,2) move represents associativity: \[
  \begin{tikzpicture}[scale=1.5]
    \draw[thick] (0,0) node{$\bullet$} to node[fill=white]{$A$} (1,1.5) node{$\bullet$} to node[fill=white]{$(AB)C$} (2,0) node{$\bullet$} to node[fill=white]{$C$} (1,-1.5) node{$\bullet$} to node[fill=white]{$B$} cycle;
    \draw[thick] (1,1.5) to node[fill=white]{$AB$} (1,-1.5);
  \end{tikzpicture}
  \raisebox{6.5em}{$\qquad\longleftrightarrow\qquad$}
  \begin{tikzpicture}[scale=1.5]
    \draw[thick] (0,0) node{$\bullet$} to node[fill=white]{$A$} (1,1.5) node{$\bullet$} to node[fill=white]{$A(BC)$} (2,0) node{$\bullet$} to node[fill=white]{$C$} (1,-1.5) node{$\bullet$} to node[fill=white]{$B$} cycle;
    \draw[thick] (0,0) to node[fill=white]{$BC$} (2,0);
  \end{tikzpicture}
\] Of course, the distinction between ``incoming'' and ``outgoing''
sides of the triangle is conventional, and the more detailed explanation
in \protect\hyperlink{week16}{``Week 16''} shows how that fits into the
formalism. Roughly speaking, what we have is not just any old algebra,
but an algebra that, thought of as a vector space, is equipped with an
isomorphism between it and its dual. This isomorphism allows us to
forget whether we are coming or going, so to speak.

Hmm, and here I was planning on being terse! Anyway, the still
\emph{more} interesting point is that when we think about
\(3\)-dimensional topology and ``3-dimensional algebra,'' we should no
longer think of \[
  \begin{tikzpicture}
    \draw[thick] (0,0) node{$\bullet$} to (1,1.5) node{$\bullet$} to (2,0) node{$\bullet$} to (1,-1.5) node{$\bullet$} to cycle;
    \draw[thick] (1,1.5) to (1,-1.5);
  \end{tikzpicture}
  \raisebox{4.5em}{$\qquad\text{and}\qquad$}
  \begin{tikzpicture}
    \draw[thick] (0,0) node{$\bullet$} to (1,1.5) node{$\bullet$} to (2,0) node{$\bullet$} to (1,-1.5) node{$\bullet$} to cycle;
    \draw[thick] (0,0) to (2,0);
  \end{tikzpicture}
\] as representing \emph{equal} operations (the 3-fold multiplication of
\(A\), \(B\), and \(C\)); instead, we should think of them as merely
\emph{isomorphic}, with the tetrahedron of which they are the front and
back being the isomorphism. The basic philosophy is that in
higher-dimensional algebra, as one ascends the ladder of dimensions,
certain things which had been regarded as \emph{equal} are revealed to
be merely isomorphic. This gets tricky, since certain
\emph{isomorphisms} that were regarded as equal at one level are
revealed to be merely isomorphic at the next level\ldots{} leading us
into a subtle world of isomorphisms between isomorphisms between
isomorphisms\ldots{} which the theory of \(n\)-categories attempts to
systematize. (I should note, however, that in the particular case of
associativity this business was worked out by Jim Stasheff quite a while
back: it's the homotopy theorists who were the ones with the guts to
deal with such issues first.)

Now, it turns out that in \(3\)-dimensional algebra, the isomorphism
corresponding to the (2,2) move is not something marvelously obscure. It
is in fact precisely what physicists call the ``\(6j\) symbol'', a
gadget they've been using to study angular momentum in quantum mechanics
for a long time! In quantum mechanics, the study of angular momentum is
just the study of representations of the group \(\mathrm{SU}(2)\), and
if one has representations \(A\), \(B\), and \(C\) of this group (or any
other), the tensor products \((A \otimes B) \otimes C\) and
\(A \otimes (B \otimes C)\) are not \emph{equal}, but merely
\emph{isomorphic}. It should come as no surprise that this isomorphism
is represented by physicists as a big gadget with 6 indices dangling on
it, the ``\(6j\) symbol''.

Quite a while back, Regge and Ponzano tried to cook up a theory of
quantum gravity in 3 dimensions using the \(6j\) symbols for
\(\mathrm{SU}(2)\). More recently, Turaev and Viro built a
\(3\)-dimensional topological quantum field theory using the
\(6j\)-symbols of the \emph{quantum group} \(SU_q(2)\), and this led to
lots of work, which the above article explains in a distilled sort of
way.

The original Ponzano--Regge and Turaev--Viro papers, and various other
ones clarifying the relation of the Turaev--Viro theory to quantum
gravity in spacetimes of dimension 3, are listed in
\protect\hyperlink{week16}{``Week 16''}. It's also worth checking out
the paper by Barrett and Foxon listed in
\protect\hyperlink{week24}{``Week 24''}, as well as the following paper,
for which I'll just quote the abstract:

\begin{enumerate}
\def\labelenumi{\arabic{enumi})}
\setcounter{enumi}{1}
\tightlist
\item
  Timothy J. Foxon, ``Spin networks, Turaev--Viro theory and the loop
  representation'', available as
  \href{https://arxiv.org/abs/gr-qc/9408013}{\texttt{gr-qc/9408013}}.
\end{enumerate}

\begin{quote}
We investigate the Ponzano--Regge and Turaev--Viro topological field
theories using spin networks and their \(q\)-deformed analogues. I
propose a new description of the state space for the Turaev--Viro theory
in terms of skein space, to which \(q\)-spin networks belong, and give a
similar description of the Ponzano--Regge state space using spin
networks. I give a definition of the inner product on the skein space
and show that this corresponds to the topological inner product, defined
as the manifold invariant for the union of two 3-manifolds. Finally, we
look at the relation with the loop representation of quantum general
relativity, due to Rovelli and Smolin, and suggest that the above inner
product may define an inner product on the loop state space.
\end{quote}

(Concerning the last point I cannot resist mentioning my own paper on
knot theory and the inner product in quantum gravity, ``Quantum gravity
and the algebra of tangles''.)

In addition to the papers by Turaev--Viro and Fukuma--Shapere listed in
\protect\hyperlink{week16}{``Week 16''}, there are some other papers on
Hopf algebras and 3d topological quantum field theories that I should
list:

\begin{enumerate}
\def\labelenumi{\arabic{enumi})}
\setcounter{enumi}{2}
\item
  Greg Kuperberg, ``Involutory Hopf algebras and three-manifold
  invariants'', \emph{Internat. Jour. Math} \textbf{2} (1991), 41--66.

  ``A definition of \(\#(M,H)\) in the non-involutory case'', by Greg
  Kuperberg, unpublished.
\end{enumerate}
\noindent
Greg Kuperberg is one of the few experts on this subject who is often
found on the net; he is frequently known to counteract my rhetorical
excesses with a dose of precise information. The above papers, one of
which is sadly still unpublished, make it beautifully clear how
``algebra knows more about topology than we do'', since various basic
structures on Hopf algebras have a pleasant tendency to interact just as
needed to give 3d topological quantum field theories.

\begin{enumerate}
\def\labelenumi{\arabic{enumi})}
\setcounter{enumi}{3}
\item
  John W. Barrett and Bruce W. Westbury, ``Spherical categories'',
  available as
  \href{https://arxiv.org/abs/hep-th/9310164}{\texttt{hep-th/9310164}}.

  John W. Barrett and Bruce W. Westbury, ``Invariants of
  piecewise-linear 3-manifolds'', \emph{Trans. Amer. Math. Soc.}
  \textbf{348} (1996), 3997--4022. Also available as
  \href{https://arxiv.org/abs/hep-th/9311155}{\texttt{hep-th/9311155}}.

  John W. Barrett and Bruce W. Westbury, ``The equality of 3-manifold
  invariants'', available as
  \href{https://arxiv.org/abs/hep-th/9406019}{\texttt{hep-th/9406019}}.
\end{enumerate}
\noindent
Let me quote the abstract for the first one; the second one gives a
construction of 3-manifold invariants, and the third shows that the
authors' 3-manifold invariants agree with Kuperberg's when both are
defined.

\begin{quote}
This paper is a study of monoidal categories with duals where the tensor
product need not be commutative. The motivating examples are categories
of representations of Hopf algebras and the motivating application is
the definition of \(6j\)-symbols as used in topological field theories.

We introduce the new notion of a spherical category. In the first
section we prove a coherence theorem for a monoidal category with duals
following MacLane (1963). In the second section we give the definition
of a spherical category, and construct a natural quotient which is also
spherical.

In the third section we define spherical Hopf algebras so that the
category of representations is spherical. Examples of spherical Hopf
algebras are involutory Hopf algebras and ribbon Hopf algebras. Finally
we study the natural quotient in these cases and show it is semisimple.
\end{quote}

\begin{enumerate}
\def\labelenumi{\arabic{enumi})}
\setcounter{enumi}{4}
\tightlist
\item
  Louis H. Kauffman and David E. Radford, ``Invariants of 3-Manifolds
  derived from finite dimensional Hopf algebras'', available as
  \href{https://arxiv.org/abs/hep-th/9406065}{\texttt{hep-th/9406065}}.
\end{enumerate}
\noindent
This paper also relates 3d topology and certain finite-dimensional
Hopf algebras, and it shows they give 3-manifold invariants distinct
from the more famous ones due to Witten (and a horde of mathematicians).
I have not had time to think about how they relate to the above ones,
but I have a hunch that they are the same, since all of them make heavy
use of special grouplike elements associated to the antipode.

\begin{enumerate}
\def\labelenumi{\arabic{enumi})}
\setcounter{enumi}{5}
\tightlist
\item
  Louis Crane and Igor Frenkel, ``Four dimensional topological quantum
  field theory, Hopf categories, and the canonical bases'', available as
  \href{https://arxiv.org/abs/hep-th/9405183}{\texttt{hep-th/9405183}}.
\end{enumerate}
\noindent
Work in 4 dimensions is, as one expects, still more subtle than in 3,
since again various things that were equalities become isomorphisms. In
particular, this means that various things one thought were vector
spaces --- which are \emph{sets} that have \emph{elements} that you can
\emph{add} and \emph{multiply by numbers}, and which satisfy
\emph{equations} like \[A + B = B + A\] are now reinterpreted as
``2-vector spaces'', which are \emph{categories} that have
\emph{objects} that you can \emph{direct sum} and \emph{tensor with
vector spaces}, and which have certain \emph{natural isomorphisms} like
the isomorphism \[A \oplus B \cong B \oplus A.\] In particular, using
Lusztig's canonical basis, Crane and Frenkel start with quantum groups
(which are Hopf algebras of a certain sort) and build marvelous ``Hopf
categories'' out of them. While they do not construct a 4d TQFT in this
paper, they indicate the game plan in terms clear enough that they will
probably now have to race other workers in the field to see who can get
the first interesting 4d TQFT\ldots{} or perhaps something a bit subtler
than a 4d TQFT (e.g.~Donaldson theory).

Finally, let me turn to a subject that is closely related (though
unfortunately this has not yet been made sufficiently clear), namely,
holonomy algebras and the loop representation of quantum gravity. Let me
simply list the references now; many of these papers were discussed at
my session on knots and quantum gravity at the Marcel Grossman
conference, so I promise to explain at some later time (and in some
papers I'm writing) a bit more about how the loop representation of a
gauge theory is interesting from the viewpoint of higher-dimensional
algebra!

\begin{enumerate}
\def\labelenumi{\arabic{enumi})}
\setcounter{enumi}{6}
\item
  A. Ashtekar, J. Lewandowski, D. Marolf, J. Mourao and T. Thiemann, ``A
  manifestly gauge-invariant approach to quantum theories of gauge
  fields'', contribution to the Cambridge meeting proceedings, available
  as
  \href{https://arxiv.org/abs/hep-th/9408108}{\texttt{hep-th/9408108}}.

  Jerzy Lewandowski, ``Topological measure and graph-differential
  geometry on the quotient space of connections'', \emph{Proceedings of
  ``Journees Relativistes 1993''}, available as
  \href{https://arxiv.org/abs/gr-qc/9406025}{\texttt{gr-qc/9406025}}.

  Abhay Ashtekar, Donald Marolf and Jose Mourao, ``Integration on the
  space of connections modulo gauge transformations'', available as
  \href{https://arxiv.org/abs/gr-qc/9403042}{\texttt{gr-qc/9403042}}.

  A. Ashtekar and R. Loll, ``New loop representations for 2+1 gravity'',
  available as
  \href{https://arxiv.org/abs/gr-qc/9405031}{\texttt{gr-qc/9405031}}.

  R. Loll, ``Independent loop invariants for 2+1 gravity'', available as
  \href{https://arxiv.org/abs/gr-qc/9408007}{\texttt{gr-qc/9408007}}.

  R. Loll, J.M. Mourão and J.N. Tavares, ``Generalized coordinates on
  the phase space of Yang--Mills theory'', available as
  \href{https://arxiv.org/abs/gr-qc/9404060}{\texttt{gr-qc/9404060}}.

  C. Di Bartolo, R. Gambini and J. Griego, ``The extended loop
  representation of quantum gravity'', available as
  \href{https://arxiv.org/abs/gr-qc/9406039}{\texttt{gr-qc/9406039}}.

  Rodolfo Gambini, Alcides Garat and Jorge Pullin, ``The constraint
  algebra of quantum gravity in the loop representation'', available as
  \href{https://arxiv.org/abs/gr-qc/9404059}{\texttt{gr-qc/9404059}}.
\end{enumerate}

\hypertarget{week39}{%
\section{September 24, 1994}\label{week39}}

I want to say a bit about Alain Connes' book, newly out in English, and
then some about Yang--Mills theory in 2 dimensions.

\begin{enumerate}
\def\labelenumi{\arabic{enumi})}
\tightlist
\item
   Alain Connes, \emph{Noncommutative Geometry}, Academic Press, Cambridge,
   Massachusetts.
\end{enumerate}

You know something is up when a prominent mathematical physicist (Daniel
Kastler) says ``Alain is great. I am just his humble prophet.'' (This
happened at a conference at Penn State I just went to.) What is
noncommutative geometry and what's so great about it?

Basically, the idea of noncommutative geometry is to generalize geometry
to ``quantum spaces''. For example, the ordinary plane has two functions
on it, the coordinate functions \(x\) and y, which commute: \(xy = yx\).
We can think of \(x\) and \(y\) as representing the position and
momentum of a classical particle. But when we consider a
quantum-mechanical particle, we must give up commutativity and instead
impose the ``canonical commutation relations'' \(xy-yx = i \hbar\),
where \(\hbar\) is Planck's constant. Now \(x\) and \(y\) are not really
functions on any space at all, but simply elements of a noncommutative
algebra. Still, we can try our best to \emph{pretend} that they are
functions on some mysterious sort of ``quantum space'' in which knowing
one coordinate of a point precisely precludes us from knowing the other
coordinate exactly, by the Heisenberg uncertainty principle.
Mathematically, noncommutative geometry consists of 1) expressing the
geometry of spaces algebraically in terms of the commutative algebra of
functions on them, and 2) then generalizing the results to classes of
noncommuative algebras.

The main trick invented by Connes was to come up with a substitute for
the ``differential forms'' on a space. Differential forms are the bread
and butter of modern geometry. If we start with a commutative algebra
\(A\) (say the algebra of smooth functions on some manifold like the
plane), we can form the algebra of differential forms over \(A\) by
introducing, for each element \(f\) in \(A\), a formal symbol \(df\),
and imposing the following rules:

\begin{itemize}
\tightlist
\item
  \(d(f+g) = df + dg\)
\item
  \(d(cf) = c df\) (\(c\) a constant)
\item
  \(d(fg) = (df)g + f dg\)
\item
  \(fdg = (dg)f\)
\item
  \(df dg = -dg df\).
\end{itemize}

More precisely, the differential forms over \(A\) are the algebra
generated by \(A\) and these differentials \(df\), modulo the above
relations. This gives a purely algebraic way of understanding what those
mysterious things like \(dx dy dz\) in integral signs are.

Now, the last two of the five rules listed above fit nicely with the
commutativity of \(A\) when it \emph{is} commutative, but they jam up the
works horribly otherwise. So: how to generalize differential forms to
the noncommutative case? There are various things one can do if \(A\) is
commutative in some generalized sense, such as ``supercommutative'' or
``braided commutative'' (which I call ``\(R\)-commutative'' in some papers
on this subject). However, if \(A\) is utterly noncommutative, it seems
that the best approach is Connes', which is first to \emph{throw out}
the last two relations, obtaining something folks call the
``differential envelope'' of \(A\) or the ``universal differential
graded algebra'' over \(A\) --- which is pleasant but quite boring by
itself --- and then to consider ``chains'' which are linear maps \(F\)
from this gadget to the complex numbers (or whatever field you're
working in) satisfying the cyclic property \[F(uv) = (-1)^{ij} F(vu)\]
where \(u\) is something that looks like \(f_0 df_1 df_2 \ldots df_i\),
and \(v\) is something like \(g_0 dg_1 dg_2 \ldots dg_j\). There are
charming things one can do with chains that wind up letting one do most
of what one could do with differential forms. More precisely, just as
differential forms allow you entry into the wonderful world of de Rham
cohomology, chains let you develop something similar called cyclic
homology (and there is a corresponding cyclic cohomology that's even
more like the de Rham theory).

Connes, being extremely inventive and ambitious, has applied
noncommutative differential geometry to many areas: index theory,
K-theory, foliations, Penrose tilings, fractals, the quantum Hall
effect, and even elementary particle physics. Perhaps the most
intriguing result is that if one develops the Yang--Mills equations using
the techniques of noncommutative geometry, but with a very simple
``commutative'' model of spacetime, namely a two-sheeted cover of
ordinary spacetime, the Higgs boson falls out rather magically on its
own. This has led Kastler and other physicists to pursue a reformulation
of the whole Standard Model in terms of noncommutative geometry, hoping
to simplify it and even make some new predictions. It is far too early
to see if this approach will get somewhere useful, but it's certainly
interesting.

I haven't read this book, just part of the French version on which it's
based (with extensive additions), but my impression is that it's quite
easy to read given the technical nature of the subject.

\begin{enumerate}
\def\labelenumi{\arabic{enumi})}
\setcounter{enumi}{1}
\tightlist
\item
   Gregory
  Moore, ``2d Yang--Mills theory and topological field theory'', available as
  \href{https://arxiv.org/abs/hep-th/9409044}{\texttt{hep-th/9409044}}.
\end{enumerate}
\noindent
This is a nice review of recent work on 2d Yang--Mills theory. While
Yang--Mills theory in 4 dimensions is the basis of our current theories
of the strong, weak, and electromagnetic forces, and mathematically
gives rise to a cornucopia of deep results about \(4\)-dimensional
topology, 2d Yang--Mills theory has traditionally been considered
``trivial'' in that one can exactly compute pretty much whatever one
wants. However, Witten, in ``On quantum gauge theories in two
dimensions'' (see \protect\hyperlink{week36}{``Week 36''}), showed that
precisely because 2d Yang--Mills theory was exactly soluble, one could
use it to study a lot of interesting mathematics problems relating to
``moduli spaces of flat connections.'' (More about those below.) And
Gross, Taylor and others have recently shown that 2d Yang--Mills theory,
at least working with gauge groups like \(\mathrm{SU}(N)\) or
\(\mathrm{SO}(N)\) and taking the ``large N limit'', could be formulated
as a string theory. So people respect 2d Yang--Mills theory more these
days; its complexities stand as a strong clue that we've just begun to
tap the depths of 4d Yang--Mills theory!

I can't help but add that Taylor and I did some work a while back in
which we formulated \(\mathrm{SU}(N)\) 2d Yang--Mills theory for
\emph{finite} N as a string theory. This was meant as evidence for my
proposal that the loop representation of quantum gravity is a kind of
string theory, a proposal described in \protect\hyperlink{week18}{``Week
18''}. For more on this sort of thing, try my paper in the book
\emph{Knots and Quantum Gravity} (see \protect\hyperlink{week23}{``Week
23''}) --- which by the way is finally out --- and also the following:

\begin{enumerate}
\def\labelenumi{\arabic{enumi})}
\setcounter{enumi}{2}
\tightlist
\item
   J. Baez and
  W. Taylor, ``Strings and two-dimensional QCD for finite \(N\)'', available as
  \href{https://arxiv.org/abs/hep-th/9401041}{\texttt{hep-th/9401041}}.
\end{enumerate}

When it comes to ``moduli spaces of flat connections'', it's hard to say
much without becoming more technical, but I certainly recommend starting
with the beautiful work of Goldman:

\begin{enumerate}
\def\labelenumi{\arabic{enumi})}
\setcounter{enumi}{3}
\item
  William Goldman, ``The symplectic nature of fundamental groups of surfaces'',  
  \emph{Adv. Math.} \textbf{54} (1984), 200--225.

   William Goldman, ``Invariant functions on Lie groups and Hamiltonian flows of surface
  group representations'', \emph{Invent. Math.}
  \textbf{83} (1986), 263--302.

  William Goldman, ``Topological components of spaces of representations'', 
  \emph{Invent. Math.} \textbf{93} (1988), 557--607.
\end{enumerate}
\noindent
The basic idea here is to take a surface \(S\) with a particular
\(G\)-bundle on it, and carefully study the space of flat connections
modulo gauge transformations, which will be a finite-dimensional
stratified space. If you fix \(G\) and \(S\), no matter what bundle you
pick, this space will appear as a subspace of a bigger space called the
moduli space of flat connections, which is the same as
\(\mathrm{Hom}(\pi_1(S),G)/\mathrm{Ad} G\). There is an open dense set
of this space, the ``top stratum'', which is a symplectic manifold.
Geometric quantization of this manifold has everything in the world to
do with Chern--Simons theory, as summarized so deftly by Atiyah:

\begin{enumerate}
\def\labelenumi{\arabic{enumi})}
\setcounter{enumi}{4}
\tightlist
\item
   Michael Atiyah, \emph{The Geometry and Physics of Knots}, Cambridge
  U. Press, Cambridge, 1990.
\end{enumerate}
\noindent
On the other hand, lately people have been using 2d Yang--Mills theory,
\(BF\) theory, and the like (see \protect\hyperlink{week36}{``Week
36''}) to get a really thorough handle on the cohomology of the moduli
space of flat connections. For a mathematical approach to this problem
that doesn't talk much about gauge theory, try:

\begin{enumerate}
\def\labelenumi{\arabic{enumi})}
\setcounter{enumi}{5}
\tightlist
\item
   Lisa C. Jeffrey, ``Group cohomology construction of the cohomology of moduli spaces of
  flat connections on 2-manifolds'', available as 
  \href{https://arxiv.org/abs/alg-geom/9404012}{\texttt{alg-geom/9404012}}.
\end{enumerate}

\hypertarget{week40}{%
\section{October 19, 1994}\label{week40}}

When I was an undergraduate I was quite interested in logic and the
foundations of mathematics --- I was always looking for the most
mind-blowing concepts I could get ahold of, and G\"odel's theorem, the
L\"owenheim--Skolem theorem, and so on were right up there with quantum
mechanics and general relativity as far as I was concerned. I did my
undergrad thesis on computability and quantum mechanics, but then I sort
of lost interest in logic and started thinking more and more about
quantum gravity. The real reason was probably that my thesis didn't turn
out as interesting as I'd hoped, but I remember feeling at the time that
logic had become less revolutionary than in it was in the early part of
the century. It seemed to me that logic had become a branch of
mathematics like any other, studying obscure properties of models of the
Zermelo--Fraenkel axioms, rather than questioning the basic presumptions
implicit in those axioms and daring to pursue new, different approaches.
I couldn't really get excited about the properties of super-huge
cardinals. Of course, I knew a bit about intuitionistic logic and
various forms of finitism, but these seemed to be the opposite of
daring; instead, they seemed to appeal mainly to grumpy people who
didn't trust abstractions and wanted to do everything as conservatively
as possible. I was pretty interested in quantum logic, too, but I tended
to think of this more as a branch of physics than `logic' proper.

Anyway, it's now quite clear to me that I just hadn't been reading the
right stuff. I think Rota has said that the really interesting work in
logic now goes under the name of `computer science', but for whatever
reason, I didn't dig into the \emph{Journal of Philosophical Logic}, other
logic journals, or proceedings of conferences on category theory,
computer science and the like and find the stuff that would have excited
me. It goes to show that one really needs to keep digging! Anyway, I
just went to a conference called the Lambda Calculus Jumelage up in
Ottawa, thanks to a kind invitation by Prakash Panangaden and Phil
Scott, who thought my ideas on category theory and physics might
interest (or at least amuse) the folks who attend this annual bash. It
became clear to me while up there that logic is alive and well!

Of course, I don't actually understand most of what these people are up
to, so take what I say with a large grain of salt. My goal here is more
to draw attention to some interesting-sounding ideas than to explain
them.

One interesting subject, which I think I'm finally beginning to get an
inkling of, is ``linear logic''. This was introduced in the following
paper (which I haven't gotten around to looking at):

\begin{enumerate}
\def\labelenumi{\arabic{enumi})}
\tightlist
\item
  Jean-Yves Girard,  ``Linear Logic'', \emph{Theoretical Computer
  Science} \textbf{50} (1987), 1--102.
\end{enumerate}

When I first heard about linear logic, it made utterly no sense. It
seemed to be a logic suitable for use in some completely different
universe than the one I inhabited! For example, there were the familiar
logical connectives ``and'' and ``or'', but they had weird alternate
versions called ``tensor'' and ``par'', the latter written with an
upside-down ampersand. There was also an alternate version of the
material implication ``\(\to\)'', and a strange operation called
``\(!\)'' (pronounced ``bang'') that somehow mediated between the
logical connectives I knew and loved and their eerie alter egos.

I understand a wee bit about these things now; one can get a certain
ways just by getting used to ``tensor'', since the rest of the weird
connectives are defined in terms of this one and the familiar ones. (I
won't worry about the ``\(!\)'' here.) One key idea, which finally
penetrated my thick skull, is that there is a good reason why ``tensor''
does not satisfy the following deduction rule so characteristic of
``and'':
\[\AxiomC{$S\vdash p$}\AxiomC{$S\vdash q$}\BinaryInfC{$S\vdash p\&q$}\DisplayProof\]
meaning: if from the set of premises \(S\) we can deduce \(p\), and
from \(S\) we can also deduce \(q\), then from \(S\) we can deduce
\(p\&q\). The point is that in linear logic one should not think of
\(S\) as a \emph{set} of premises, but rather as a \textbf{multiset},
meaning that the same premise can appear twice. The idea is that if we
use one premise in \(S\) to deduce something, we \emph{use it up}, and we can
only use it again if \(S\) has several copies of that premise in it. As
they say, linear logic is ``resource-sensitive'' (which is apparently
why computer scientists like it). So the idea is that in linear logic,
\[S\vdash p\&q\] means something like ``from the premises \(S\) one can
deduce \(p\) if one feels like it, or alternatively one can deduce \(q\)
if one feels like it, but not necessarily both at once, since there
may not be enough copies of the premises to do that.'' On the other
hand, \[S\vdash p\otimes q\] is stronger, since it means something like
``from the premises \(S\) one can deduce both \(p\) and \(q\) at once,
since there are enough copies of all the premises in \(S\) to do it.''
Thus ``\(\&\)'' satisfies the above deduction rule in linear logic just
as in classical logic, but ``tensor'' does not; instead, it satisfies
\[\AxiomC{$S\vdash p$}\AxiomC{$T\vdash q$}\BinaryInfC{$S\cup T\vdash p\otimes q$}\DisplayProof\]
where \(S \cup T\) denotes the union of the multisets \(S\) and \(T\)
(so that if both \(S\) and \(T\) have one copy of a premis,
\(S \cup T\) has two copies of it).

Well, let me leave it at that. I should add that there is a paper
available online,

\begin{enumerate}
\def\labelenumi{\arabic{enumi})}
\setcounter{enumi}{1}
\tightlist
\item
 Vaughan Pratt, ``Linear logic for generalized quantum mechanics'', 
  available at \href{http://citeseerx.ist.psu.edu/viewdoc/summary?doi=10.1.1.49.649}{http://citeseerx.ist.psu.edu/viewdoc/summary?doi=10.1.1.49.649}.
\end{enumerate}
\noindent
which relates linear logic and quantum logic, and which is part of a
body of work relating linear logic and category theory, with the key
idea being that ``linear logic is a logic of monoidal closed categories
in much the same way that intuitionistic logic is a logic of Cartesian
closed categories'' --- here I quote

\begin{enumerate}
\def\labelenumi{\arabic{enumi})}
\setcounter{enumi}{2}
\tightlist
\item
  Richard Blute, ``Hopf algebras and linear logic'', 
  \emph{Mathematical Structures in Computer Science}, \textbf{6} (1996), 189--217.
  Available as \href{https://citeseerx.ist.psu.edu/viewdoc/summary?doi=10.1.1.114.8038}{https://citeseerx.ist.psu.edu/viewdoc/} \href{https://citeseerx.ist.psu.edu/viewdoc/summary?doi=10.1.1.114.8038}{summary?doi=10.1.1.114.8038}.
\end{enumerate}
\noindent
I suppose to most people, explaining linear logic in terms of monoidal
closed categories may seem like using mud to wipe one's windshield.
However, to some of us monoidal closed categories are rather familiar
things, and in fact anyone who knows about vector spaces, linear maps,
and the vector spaces \(\mathrm{Hom}(V,W)\) and \(V\otimes W\) knows a
really good example of a monoidal closed category. Thus monoidal closed
categories can be viewed as an abstraction of linear algebra, and indeed
this is how ``linear logic'' got its name.

It seems that I should read the following papers, too, before I really
understand the connection between linear logic and category theory:

\begin{enumerate}
\def\labelenumi{\arabic{enumi})}
\setcounter{enumi}{3}
\item
  R. A. G. Seely, ``Linear logic, *-autonomous categories and cofree coalgebras'', in \emph{Categories in Computer Science and Logic}, Contemp.
  Math. \textbf{92}, American Mathematical Society, Providence, Rhode Island, 1989.
\item
  D. Yetter, ``Quantales and (noncommutative) linear logic'',
  \emph{Journal of Symbolic Logic} \textbf{55} (1990), 41--64.
\end{enumerate}
\noindent
A terse summary of linear logic in terms a categorist might like can be
found in Section 3.5 of Pratt's paper cited above. I should add that
Pratt has lots of other interesting papers available online.

\hypertarget{week41}{%
\section{October 17, 1994}\label{week41}}

In the beginning of September I went to a conference at the Center for
Gravitational Physics and Geometry at Penn State. This is the center run
by Abhay Ashtekar, and it has Jorge Pullin and Lee Smolin as faculty,
and Roger Penrose as a part-time visitor --- so it's a great place to
visit if you're interested in quantum gravity. There are a lot of good
postdocs and such there, too. I've been too busy to say much so far
about what happened at this conference, but I'd like to now.

One talk I enjoyed a lot was Steve Carlip's, on the entropy of black
holes. This has subsequently come out as a preprint, available
electronically:

\begin{enumerate}
\def\labelenumi{\arabic{enumi})}
\tightlist
\item
  Steve Carlip, ``The statistical mechanics of the (2+1)-dimensional black hole'', 
   available as
  \href{https://arxiv.org/abs/gr-qc/9409052}{\texttt{gr-qc/9409052}}.
\end{enumerate}

It's well-known by now that in certain situations it makes sense to
speak of the ``entropy'' of a black hole, but the real meaning of this
entropy is still mysterious. In particular, since the entropy of a black
hole is (often, but not always) proportional to the area of its event
horizon, it would be very satisfying if the entropy corresponded somehow
to degrees of freedom that ``lived at the event horizon''. Steve Carlip
has done a pretty credible calculation along these lines (though not
without various subtle difficulties) in the case of a black hole in
3-dimensional spacetime.

I should say a little bit about gravity in 3 dimensions and why people
are interested in it. \(3\)-dimensional gravity is drastically simpler
than 4-dimensional gravity, since in 3 dimensions the vacuum Einstein's
equations say the spacetime metric is \emph{flat}, at least if the
cosmological constant vanishes. Thus there can be no gravitational
radiation (and in quantum theory no ``gravitons''), and the metric
produced by a static point mass is not like the Schwarschild metric,
instead, on space it is just like that of a cone. Things are a bit
different if the cosmological constant is nonzero; in particular, there
are black-hole type solutions. But there is still no gravitational
radiation.

Basically, people are interested in \(3\)-dimensional quantum gravity
because it's simple enough that one can compute something and hope it
sheds some light on the \(4\)-dimensional world we live in. For some
issues this appears to be the case: primarily, conceptual issues having
to do with theories in which there is no ``background metric''.
Unfortunately, there are \emph{several different ways} to set up
\(3\)-dimensional quantum gravity, corresponding to different approaches
people have to 4-dimensional quantum gravity. For this, check out
Carlip's paper ``Six ways to quantize (2+1)-dimensional gravity,''
mentioned in \protect\hyperlink{week16}{``Week 16''}. However, I think
the ``best'' way to quantize gravity in 3 dimensions is the way
involving Chern--Simons theory, because this way is the most closely
related to Ashtekar's approach to quantizing gravity in 4 dimensions,
hence it sheds the most light on the things I'm interested in --- and I
also think it's the most beautiful. In this approach, you can compute a
lot of things, and basically what Carlip has done is to show that
associated to the event horizon there are degrees of freedom which
should give entropy proportional to its area.

I suppose I can't say how he does it much more clearly than he says it,
so I'll quote the introduction, taking the liberty of turning some of
his LaTeX into English. If you get scared by the ``Virasoro operator
\(L_0\)'' below, never fear --- in this context, it just amounts to the
angular momentum operator, which generates rotations about the origin.
So:

\begin{quote}
The basic argument is quite simple. Begin by considering general
relativity on a manifold \(M\) with boundary. We ordinarily split the
metric into true physical excitations and ``pure gauge'' degrees of
freedom that can be removed by diffeomorphisms of \(M\). But the
presence of a boundary alters the gauge invariance of general
relativity: the infinitesimal transformations {[}\ldots{]} must now be
restricted to those generated by vector fields {[}\ldots{]} with no
component normal to the boundary, that is, true diffeomorphisms that
preserve the boundary of \(M\). As a consequence, some degrees of
freedom that would naively be viewed as ``pure gauge'' become dynamical,
introducing new degrees of freedom associated with the boundary.

Now, the event horizon of a black hole is not a true boundary, although
the black hole complementarity approach of Susskind et al.~suggests that
it might be appropriately treated as such. Regardless of one's view of
that program, however, it is clear that in order to ask quantum
mechanical questions about the behavior of black holes, one must put in
``boundary conditions'' that ensure that a black hole is present. This
means requiring the existence of a hypersurface with particular metric
properties---say, those of an apparent horizon.

The simplest way to do quantum mechanics in the presence of such a
surface is to quantize fields separately on each side, imposing the
appropriate correlations as boundary conditions. In a path integral
approach, for instance, one can integrate over fields on each side,
equate the boundary values, and finally integrate over those boundary
values compatible with the existence of a black hole. But this process
again introduces boundary terms that restrict the gauge invariance of
the theory, leading once more to the appearance of new degrees of
freedom at the horizon that would otherwise be treated as unphysical.

My suggestion is that black hole entropy is determined by counting these
would-be gauge degrees of freedom. The resulting picture is similar to
Maggiore's membrane model of the black hole horizon, but with a
particular derivation and interpretation of the ``membrane'' degrees of
freedom.

The analysis of this phenomenon is fairly simple in 2+1 dimensions. It
is well known that (2+1)-dimensional gravity can be written as a
Chern--Simons theory, and it is also a standard result that a
Chern--Simons theory on a manifold with boundary induces a dynamical
Wess--Zumino--Witten (WZW) theory on the boundary. In the presence of a
cosmological constant \(\Lambda = -1/L^2\) appropriate for the
(2+1)-dimensional black hole, one obtains a slightly modified
\(\mathrm{SO}(2,1) \times \mathrm{SO}(2,1)\) WZW model, with coupling
constant \[k = \frac{L\sqrt{2}}{8G}\] This model is not completely
understood, but in the large \(k\) --- i.e., small \(\Lambda\) ---
limit, it may be approximated by a theory of six independent bosonic
oscillators. I show below that the Virasoro operator \(L_0\) for this
theory takes the form \[L_0 \sim N -\left(\frac{r}{4G}\right)^2,\] where
\(N\) is a number operator and \(r\) is the horizon radius. It is a
standard result of string theory that the number of states of such a
system behaves asymptotically as \[n(N) \sim \exp(\pi \sqrt{4} N)\] If
we demand that \(L_0\) vanish --- physically, requiring states to be
independent of the choice of origin of the angular coordinate at the
horizon --- we thus obtain \[\log n(r) \sim \frac{2 \pi r}{4G},\]
precisely the right expression for the entropy of the (2+1)-dimensional
black hole.
\end{quote}

Also, Carlo Rovelli spoke about describing the dynamics of quantum
gravity coupled to a scalar field in terms of ``spin network'' states. I
think this was based on work he did in collaboration with Lee Smolin,
and I don't think it's out yet. I'm just about to finish up a little
paper on spin network states myself, since they seem like very useful
things in quantum gravity. The simplest sort of spin network is just a
trivalent graph (i.e., 3 edges adjacent to each vertex) with edges
labelled by ``spins'' \(0,\frac{1}{2},1,\frac{3}{2},\ldots\), and satisfying the
``triangle inequality'' at each vertex:
\[j_1 + j_2 \le j_3, \quad j_2 + j_3 \le j_1, \quad j_3 + j_1 \le j_2,\]
where \(j_1\), \(j_2\), \(j_3\) are the spins labelling the edges
adjacent to the given vertex. Really, the spins should be thought of as
irreducible representations of \(\mathrm{SU}(2)\), and the triangle
inequalities are necessary for the representation \(j_3\) to appear as a
summand in the tensor product of the representations \(j_1\) and
\(j_2\). (If the last sentence was meaningless to you, reading
\protect\hyperlink{week5}{``Week 5''} will help a little, though
probably not quite enough.)

Penrose introduced spin networks as part of a purely combinatorial
approach to spacetime in the paper:

\begin{enumerate}
\def\labelenumi{\arabic{enumi})}
\setcounter{enumi}{1}
\tightlist
\item
  Roger Penrose, ``Angular momentum; an approach to combinatorial space time'', 
   in \emph{Quantum Theory and Beyond}, ed.~T. Bastin,
  Cambridge University Press, Cambridge, 1971.  Available as 
  \href{https://math.ucr.edu/home/baez/penrose/}{https://math.ucr.edu/home/baez/penrose/}.
\end{enumerate}
\noindent
It is somehow satisfying, therefore, to see that spin networks arise
naturally as a convenient description of states in the loop
representation of quantum gravity, which \emph{starts} mainly with
Einstein's equations and the principles of quantum mechanics. Certainly
there is a lot more we need to learn about them.... One place worth
reading about them is:

\begin{enumerate}
\def\labelenumi{\arabic{enumi})}
\setcounter{enumi}{2}
\tightlist
\item
 Louis Crane,  ``Conformal field theory, spin geometry, and quantum gravity'', 
  \emph{Phys. Lett.} \textbf{B259} (1991), 243--248.
\end{enumerate}
\noindent
I will be coming out with a paper on them next week if I get my act
together, and I may say a bit more about them in future
``Weeks''.

Rovelli also mentioned an interesting paper he wrote about the problem
of time in quantum gravity with the operator-algebra/noncommutative-geometry guru Alain Connes:

\begin{enumerate}
\def\labelenumi{\arabic{enumi})}
\setcounter{enumi}{3}
\tightlist
\item
  Alain Connes and Carlo Rovelli, ``Von Neumann algebra automorphisms and time-thermodynamics
   relation in general covariant quantum theories'', available as
  \href{https://arxiv.org/abs/gr-qc/9406019}{\texttt{gr-qc/9406019}}.
\end{enumerate}
\noindent
The problem of time in quantum gravity is a bit tricky to describe,
since it takes different guises in different approaches to quantum
gravity, but I have attempted to give a rough introduction to it in
\protect\hyperlink{week11}{``Week 11''} and
\protect\hyperlink{week27}{``Week 27''}. One way to get a feeling for it
is to realize that anything you are used to doing with Hamiltonians in
quantum mechanics or quantum field theory, you \emph{can't} do in quantum
gravity, at least not in any simple way, because there is no Hamiltonian
in general relativity, but only a ``Hamiltonian constraint'' --- which
in quantum gravity becomes the Wheeler--DeWitt equation 
\[H \psi = 0.\]
Now, people know there is a mystical relationship between time and
temperature that might be written \[it = \frac{1}{kT}\] where \(t\) is
time, \(T\) is temperature, and \(k\) is Boltzmann's constant. This
equation is a bit of an exaggeration! But the point is that in quantum
theory, when there is a Hamiltonian \(H\) around one evolves states
using the operator \[\exp(-itH)\] while the Gibbs state, that is, the
equilibrium state at temperature \(T\), is given by the density matrix
\[\exp(-H/kT).\] It is this fact that relates statistical mechanics and
quantum field theory so closely.

Now, in quantum gravity things aren't so simple, since there isn't a
Hamiltonian (just a Hamiltonian constraint). However, people \emph{do}
know that there are all sorts of funny relationships between statistical
mechanics and quantum gravity. For example, an accelerating observer in
Minkowski space will see the vacuum as a heat bath with temperature
proportional to her acceleration, so in curved spacetime, where there
are no truly inertial frames, there really is no well-defined notion of
a vacuum; in some vague sense, all there are is ``thermal'' states. This
fact is also somehow related to Hawking radiation, and to the notion of
black hole entropy\ldots{} but really, there is a lot that nobody
understands about all these connections!

In any event, Rovelli was prompted to use thermodynamics to
\textbf{define} time in quantum gravity as follows. Given a mixed state
with density matrix \(D\), \emph{find} some operator \(H\) such that
\(D\) is the Gibbs state \(\exp(-H/kT)\). In lots of cases this isn't
hard; it basically amounts to \[H = -kT \ln D\] Of course, \(H\) will
depend on \(T\), but this really is just saying that fixing your units
of temperature fixes your units of time!

Operator theorists have pondered this notion very carefully for a long
time and generalized it into something called the Tomita--Takesaki
theorem, which Connes and Rovelli explain. This gives a very general way
to cook up a Hamiltonian (hence a notion of time evolution) from a state
of a quantum system! For example, one can use this trick to start with a
Robertson--Walker universe full of blackbody radiation, and recover a
notion of ``time''. This is very intriguing, and it may represent some
real progress in understanding the deep relations between time,
thermodynamics, and gravity. There are, of course, lots of problems and
puzzles to deal with.

Another intriguing talk at the conference was given by Viqar Husain, on
the subject of the following paper:

\begin{enumerate}
\def\labelenumi{\arabic{enumi})}
\setcounter{enumi}{4}
\tightlist
\item
  Viqar Husain, ``The affine symmetry of self-dual gravity", available as
  \href{https://arxiv.org/abs/hep-th/9410072}{\texttt{hep-th/9410072}}.
\end{enumerate}
\noindent
Let me simply quote the abstract, since I don't feel I really understand
the essence of this business well enough to say anything useful yet:

\begin{quote}
Self-dual gravity may be reformulated as the two dimensional chiral
model with the group of area preserving diffeomorphisms as its gauge
group. Using this formulation, it is shown that self-dual gravity
contains an infinite dimensional hidden symmetry algebra, which is the
Affine (Kac--Moody) algebra associated with the Lie algebra of area
preserving diffeomorphisms. This result provides an observable algebra
and a solution generating technique for self-dual gravity.
\end{quote}

A couple more things before I wrap this up\ldots. First, in case any
mathematicians out there are wondering what this ``knots and quantum
gravity'' business is all about, here's something I wrote to review the
subject:

\begin{enumerate}
\def\labelenumi{\arabic{enumi})}
\setcounter{enumi}{5}
\tightlist
\item
  John Baez, ``Knots and quantum gravity: progress and prospects'', available as
  \href{https://arxiv.org/abs/gr-qc/9410018}{\texttt{gr-qc/9410018}}.
\end{enumerate}
\noindent
My abstract:

\begin{quote}
Recent work on the loop representation of quantum gravity has revealed
previously unsuspected connections between knot theory and quantum
gravity, or more generally, \(3\)-dimensional topology and
\(4\)-dimensional generally covariant physics. We review how some of
these relationships arise from a `ladder of field theories' including
quantum gravity and \(BF\) theory in 4 dimensions, Chern--Simons theory
in 3 dimensions, and the \(G/G\) gauged WZW model in 2 dimensions. We
also describe the relation between link (or multiloop) invariants and
generalized measures on the space of connections. In addition, we pose
some research problems and describe some new results, including a proof
(due to Sawin) that the Chern--Simons path integral is not given by a
generalized measure.
\end{quote}

Finally, let me draw people's attention to ``Matters of Gravity'', the
newsletter Jorge Pullin puts together at considerable effort, to keep
people informed about general relativity and the like, experimental and
theoretical:

\begin{enumerate}
\def\labelenumi{\arabic{enumi})}
\setcounter{enumi}{6}
\tightlist
\item
  ``Matters of Gravity'', a newsletter for the gravity community,
  Number \textbf{4}, edited by Jorge Pullin, 24 pages in Plain TeX,
  available as
  \href{https://arxiv.org/abs/gr-qc/9409004}{\texttt{gr-qc/9409004}}.
\end{enumerate}

Here's the table of contents of this issue:

\begin{itemize}
\tightlist
\item
  Editorial.
\item
  Gravity News:

  \begin{itemize}
  \tightlist
  \item
    Report on the APS topical group in gravitation, Beverly Berger.
  \end{itemize}
\item
  Research briefs:

  \begin{itemize}
  \tightlist
  \item
    Gravitational microlensing and the search for dark matter, Bohdan
    Paczynski.
  \item
    Laboratory gravity: the G mystery, Riley Newman.
  \item
    LIGO project update, Stan Whitcomb.
  \end{itemize}
\item
  Conference Reports

  \begin{itemize}
  \tightlist
  \item
    PASCOS '94, Peter Saulson.
  \item
    The Vienna Meeting, P. Aichelburg, R. Beig.
  \item
    The Pitt binary black hole grand challenge meeting, Jeff Winicour.
  \item
    International symposium on experimental gravitation at Pakistan,
    Munawar Karim.
  \item
    10th Pacific coast gravity meeting, Jim Isenberg.
  \end{itemize}
\end{itemize}

\hypertarget{week42}{%
\section{November 3, 1994}\label{week42}}

String theory means different things to different people. The original
theory of strings -- at least if I've got my history right -- was a
theory of hadrons (particles interacting via the strong force). The
strong force wasn't understood too well then, but in 1968 Veneziano
cleverly noticed when thumbing through a math book that Euler's beta
function had a lot of the properties one would expect of the formula for
how hadrons scattered (the so-called \(S\)-matrix). Later, around 1970,
Nambu and Goto noticed that this function would come out naturally if
one thought of hadrons as different vibrational modes of a relativistic
string.

This theory had problems, and eventually it was supplanted by the
current theory of the strong force, involving quarks and gluons. The
gluons are another way of talking about the strong force, which is a
gauge field. The biggest puzzle about this approach to hadrons is, ``how
come we don't see quarks?'' This is called the puzzle of confinement. In
the late 1970's, one proposed solution was that as you pulled the quark
and the antiquark in a meson apart, the strong force effectively formed
an elastic ``string'' with constant tension. This would mean that
pulling them apart took energy proportional to how far you pulled them
apart. Past a certain point, the energy would be enough to create a new
quark-antiquark pair and \emph{snap} --- the string would split into two
new strings with quark and antiquark on each end. So here the ``string''
idea is revived but as an approximation to a theory of gauge fields. One
can even try to derive approximate string equations from the equations
for the strong force: the Yang--Mills equations. In my paper on strings,
loops, knots and gauge fields (see \protect\hyperlink{week18}{``Week
18''}), I gave references to some early papers on the subject:

\begin{enumerate}
\def\labelenumi{\arabic{enumi})}
\item
  Y.\ Nambu, ``QCD and the string model'', \emph{Phys. Lett.}
  \textbf{B80} (1979), 372--376.

  A.\ Polyakov, ``Gauge fields as rings of glue'', \emph{Nucl. Phys.}
  \textbf{B164} (1979), 171--188.

  Y.\ Nambu, ``The quantum dual string wave functional in Yang--Mills theories'', 
  \emph{Phys. Lett.} \textbf{B80} (1979),
  255--258.

  F.\ Gliozzi and M.\ Virasoro, ``The interaction among dual strings as a manifestation of the gauge
  group'', \emph{Nucl. Phys.}
  \textbf{B164} (1980), 141--151.

  A.\ Jevicki, ``Loop-space representation and the large-\(N\) behavior of the
  one-plaquette Kogut-Susskind Hamiltonian'', \emph{Phys.
  Rev.} \textbf{D22} (1980), 467--471.

 Y.\ Makeenko and A.\ Migdal,  ``Quantum chromodynamics as dynamics of loops'', 
  \emph{Nucl. Phys.} \textbf{B188} (1981), 269--316.

  Y.\  Makeenko and A.\ Migdal, ``Loop dynamics: asymptotic freedom and quark confinement'', \emph{Sov. J. Nucl. Phys.} \textbf{33} (1981), 882--893.
\end{enumerate}
\noindent
These papers make very interesting reading even today. Anyone who knows
particle physics will recognize most of these names! Strings were big
back then. But then they went out of fashion, because the string models
predicted a massless spin-2 particle --- and there's no such thing in
particle physics. Later, when people were trying to cook up ``theories
of everything'' including gravity, this flaw was again seen as a plus,
since the hypothesized ``graviton'' meets that description.

The modern, more technical subject of string theory is a lot more fancy
than these early papers. In particular, the recognition that conformal
invariance was a very good thing when studying strings propagating on
fixed background metric (like that of Minkowski space) pushed string
theorists into a careful study of \(2\)-dimensional conformally invariant
quantum field theories. (Here the 2 dimensions refer to the surface the
string traces out as it moves through spacetime.) Conformal field theory
then developed a life of its own! By now it's pretty intimidating to the
outsider. Mathematicians might find the following summary handy:

\begin{enumerate}
\def\labelenumi{\arabic{enumi})}
\setcounter{enumi}{1}
\tightlist
\item
  Krzysztof Gawedzki, ``Conformal field theory'', \emph{Seminaire
  Bourbaki}, Asterisque \textbf{177--178} (1989), 95--126.
\end{enumerate}
\noindent
while physicists might try

\begin{enumerate}
\def\labelenumi{\arabic{enumi})}
\setcounter{enumi}{2}
\item
  Michio Kaku, \emph{Introduction to Superstrings},
  Springer, Berlin, 1988.

   Michio Kaku, \emph{String Fields, Conformal Fields, and Topology}, 
  New York, Springer, Berlin, 1991.
\end{enumerate}
\noindent
Kaku's books are a decent overview but rather sketchy in spots, since
they cover vast amounts of territory.

Then there is another kind of sophisticated modern string theory,
``string field theory'', which doesn't assume the strings are moving
around on a spacetime with a background geometry. This is clearly more
like what one wants to do if one is using strings to explain quantum
gravity. I don't understand this nearly as well as I'd like to, but the
guru on this subject is Barton Zwiebach, so if one was really gutsy one
would, after a suitable warmup with Kaku, plunge in and read something
like

\begin{enumerate}
\def\labelenumi{\arabic{enumi})}
\setcounter{enumi}{3}
\item
   Ashoke Sen and Barton Zwiebach,  ``Quantum background independence of closed string 
   field theory'',  available as
  \href{https://arxiv.org/abs/hep-th/9311009}{\texttt{hep-th/9311009}}.

  Ashoke Sen and Barton Zwiebach, ``Background independent algebraic structures in 
  closed string field theory'', available as
  \href{https://arxiv.org/abs/hep-th/9408053}{\texttt{hep-th/9408053}}.
\end{enumerate}
\noindent
Unfortunately I'm not quite up to it yet\ldots.

Then, in a different direction, a bunch of folks from general relativity
pursued some ideas about string and loops to the point of developing the
``loop representation of quantum gravity.'' I'm referring to

\begin{enumerate}
\def\labelenumi{\arabic{enumi})}
\setcounter{enumi}{4}
\tightlist
\item
  Carlo Rovelli and Lee Smolin,  ``Loop representation for quantum general relativity'',
  \emph{Nucl. Phys.} \textbf{B331} (1990), 80--152.
\end{enumerate}
\noindent
though it's important to credit some of the people who kept alive the
idea that one should study gauge fields as being ``loops of string'', or
more technically, ``Wilson loops'':

\begin{enumerate}
\def\labelenumi{\arabic{enumi})}
\setcounter{enumi}{5}
\tightlist
\item
  R.\ Gambini and A.\ Trias, ``Gauge dynamics in the C-representation'', \emph{Nucl. Phys.} \textbf{B278} (1986),
   436--448.
\end{enumerate}

Now what's frustrating here is that I understand the loop representation
business, but not the ``background-free closed string field theory''
business, even though they have the same historical roots and are both
trying to deal with quantum gravity (among other forces) in a way that
assumes that loops are the basic objects. Alas, the two strands speak in
different languages! Heavy-duty mathematicians like Getzler, Kapranov
and Stasheff know how to think about closed string fields in terms of
``operads'', and that stuff seems like it should be simple enough to
understand, but alas, when I read it I get snowed in detail (so far).

Let me digress to mention what an ``operad'' is. An ``operad'' is
basically a cool way to handle sets equipped with lots of \(n\)-ary
operations. These operations might be ``parametrized'' in various ways.
The operad elegantly keeps track of these parametrizations. So, for each
\(n\), an operad has a set \(X(n)\) which we think of as all the
\(n\)-ary operations. Think of something in \(X(n)\) as a black box that
has \(n\) ``input'' tubes and one ``output'' tube, or a tree-shaped
thing \[
  \begin{tikzpicture}
    \begin{knot}
      \strand[thick] (0.75,0)
        to (0.75,-2);
      \strand[thick] (0,0)
        to [out=down,in=up] (0.75,-1);
      \strand[thick] (1.5,0)
        to [out=down,in=up] (0.75,-1);
    \end{knot}
  \end{tikzpicture}
\] with \(n\) branches and one root (here \(n = 3\)). Then suppose we
have a bunch of these black boxes. Say we have something in \(X(n_1)\),
something in \(X(n_2)\), \ldots. and so on up to something in
\(X(n_k)\). Thus we've got a pile of black boxes with a total of
\(n_1 + \ldots + n_k\) input tubes and \(k\) output tubes. Now if we
\emph{also} have a guy in \(X(k)\), which has \(k\) input tubes, we can
hook up all the output tubes of all the boxes in our pile to the input
tubes of this guy, to get a monstrous machine with
\(n_1 + \ldots + n_k\) input tubes and one output. In short, there is an
operation from
\(X(n_1)\times \ldots\times X(n_k)\times X(k) \to X(n_1 + \ldots + n_k)\).
For example, if we take the tree up there, which represents something in
\(X(3)\), and another thing in \(X(3)\), we can hook up their outputs to
the inputs of something in \(X(2)\), to get something that looks like \[
  \begin{tikzpicture}[scale=0.7]
    \begin{knot}
      \strand[thick] (0.75,0)
        to (0.75,-1.5);
      \strand[thick] (0,0)
        to [out=down,in=up] (0.75,-1);
      \strand[thick] (1.5,0)
        to [out=down,in=up] (0.75,-1);
    \end{knot}
    \begin{scope}[shift={(2,0)}]
      \begin{knot}
        \strand[thick] (0.75,0)
          to (0.75,-1.5);
        \strand[thick] (0,0)
          to [out=down,in=up] (0.75,-1);
        \strand[thick] (1.5,0)
          to [out=down,in=up] (0.75,-1);
      \end{knot}
    \end{scope}
    \begin{scope}[shift={(0.75,-1.5)}]
      \begin{knot}
        \strand[thick] (0,0)
          to [out=down,in=up] (1,-1)
          to (1,-2);
        \strand[thick] (2,0)
          to [out=down,in=up] (1,-1);
      \end{knot}
    \end{scope}
  \end{tikzpicture}
\] which is in \(X(6)\). The closed string field theorists like operads
because there are lots of parametrized ways of gluing together Riemann
surfaces with punctures together. It's a handy language,
apparently\ldots{} I am a bit more familiar with operads (though not
much) in the context of homotopy theory, where they can be used to
elegantly summarize the operations one has floating around in an
infinite loop space. \emph{Very} roughly, an infinite loop space is a
space that looks like the space of loops of loops of loops of
loops\ldots{} of loops in some topological space, where you get to make
the ``dot dot dot'' part go on as long as you want! A beautifully
unpretentious and utterly readable book on these spaces, operads, and
much much more, is:

\begin{enumerate}
\def\labelenumi{\arabic{enumi})}
\setcounter{enumi}{6}
\tightlist
\item
  J. F. Adams, \emph{Infinite Loop Spaces}, Princeton U. Press,
  Princeton, 1978.
\end{enumerate}

Lest ``infinite loop spaces'' seem abstruse, I should emphasize that the
book is really a nice tour of a lot of modern homotopy theory. As he
says, ``my object has been a more elementary exposition, which I hope
may convey the basic ideas of the subject in a way as nearly
painless as I can make it. In this the Princeton audience encouraged me;
the more I found means to omit the technical details, the more they
seemed to like it.'' A lot of the general mathematical machinery he
discusses, especially in the chapter called ``Machinery'', is really too
nice to be left for only the homotopy theorists!

Anyway, once you have gotten the hang of operads you can try the work of
a reformed homotopy theorist, Jim Stasheff, on string field theory:

\begin{enumerate}
\def\labelenumi{\arabic{enumi})}
\setcounter{enumi}{7}
\tightlist
\item
  Jim Stasheff, ``Closed string field theory, strong homotopy Lie algebras and the
  operad actions of moduli spaces'', available as
  \href{https://arxiv.org/abs/hep-th/9304061}{\texttt{hep-th/9304061}}.
\end{enumerate}
\noindent
Actually Graeme Segal, another string theory guru, also used to do
homotopy theory. He's the one who's famous for:

\begin{enumerate}
\def\labelenumi{\arabic{enumi})}
\setcounter{enumi}{8}
\tightlist
\item
  Andrew Pressley and Graeme Segal, \emph{Loop groups}, Oxford
  University Press, Oxford, 1986.
\end{enumerate}
\noindent
So it's possible that these guys didn't really quit homotopy theory, but
just figured out how to get physicists interested in it. Notice all
those loops! :-)

But where was I\ldots{} romping through various approaches to string
theory, taking a detour to mention loops, but all the while sneaking up
on my goal, which is to list a few papers that lend evidence to the
thesis of my paper Strings, Loops, Knots and Gauge Fields, namely that a
profound ``string/gauge field duality'' is at work in many physical
models, and that the loop representation of quantum gravity, and string
theory, may eventually not be seen as so different after all.

Let's see what we've got here:

\begin{enumerate}
\def\labelenumi{\arabic{enumi})}
\setcounter{enumi}{9}
\tightlist
\item
  ``A reformulation of the Ponzano--Regge quantum gravity model in terms
  of surfaces'', Junichi Iwasaki, University of Pittsburgh, 11 pages in
  LaTeX format available as
  \href{https://arxiv.org/abs/gr-qc/9410010}{\texttt{gr-qc/9410010}}.
\end{enumerate}

I've discussed the Ponzano--Regge model quite a bit in
\protect\hyperlink{week16}{``Week 16''} and
\protect\hyperlink{week38}{``Week 38''}. It's an approach to quantum
gravity that is especially successful in 3 dimensions, and involves
chopping spacetime up into simplices. The exact partition function, as
they say, can be computed using this combinatorial discrete
approximation to the spacetime manifold. (In quantum field theory, when
you know enough about the partition function you can compute the
expectation values of observables to your heart's content.) Anyway, here
Iwasaki does the kind of thing I was pointing towards in my paper,
namely, to rewrite the theory, which starts out as a gauge theory, as a
theory of surfaces (``string worldsheets'') in spacetime.

Meanwhile, more work has been done on the same kind of idea for good old
quantum chromodynamics, though here there \emph{is} a background
geometry, and one approximates the spacetime manifold by a discrete
lattice not because one expects to get the \emph{exact} answers out that
way, but just because it's a decent approximation that makes things a
bit more manageable:

\begin{enumerate}
\def\labelenumi{\arabic{enumi})}
\setcounter{enumi}{10}
\item
  B. Rusakov, ``Lattice QCD as a theory of interacting surfaces'',  available as
  \href{https://arxiv.org/abs/hep-th/9410004}{\texttt{hep-th/9410004}}.

   Ivan K.  Kostov,  ``\(\mathrm{U}(N)\) gauge theory and lattice strings'', available as
  \href{https://arxiv.org/abs/hep-th/9308158}{\texttt{hep-th/9308158}}.
\end{enumerate}
\noindent
Also, if there were any gauge theory that deserved to be a string
theory, it's probably Chern--Simons theory, which has so much to do with
knots\ldots{} and indeed something like this seems to be the case,
though it's all rather subtle and mysterious so far:

\begin{enumerate}
\def\labelenumi{\arabic{enumi})}
\setcounter{enumi}{11}
\tightlist
\item
  Michael R. Douglas, ``Chern--Simons-Witten theory as a topological Fermi liquid'', 
  available as
  \href{https://arxiv.org/abs/hep-th/9403119}{\texttt{hep-th/9403119}}.
\end{enumerate}
\noindent
Frequently, when there is a whole lot of frenetic,
sophisticated-sounding activity around a certain idea, like this
relation between strings and gauge fields, there is a simple truth
yearning to be known. Sometimes it takes a while! We'll see.

\hypertarget{week43}{%
\section{November 5, 1994}\label{week43}}

It is very exciting, yet somewhat scary, as work continues on the loop
representation of quantum gravity. On the one hand, researchers are busy
making it mathematically rigorous; on the other hand, they are beginning
to understand its physical significance. The reasons for excitement are
obvious, but the scary part is that until the final touches are put on
the mathematical rigor, we don't know if the theory really exists!

Of course there is the whole separate issue of whether the theory will
find experimental confirmation. If the theory were experimentally
confirmed, questions of mathematical rigor wouldn't be quite such a big
deal. But experimental verification will probably take a long time!
Also, we don't really expect a theory of ``pure gravity'' to be
experimentally confirmed. One will need to figure out how all the other
particles and fields fit in --- except perhaps for very general,
qualitative issues. (See the paper mentioned at the very end of this
article for some of \emph{those}.) So here the suspense is of a
long-term sort. Luckily, the question of whether the theory makes
mathematical sense is already very interesting, since so many theories
of quantum gravity have already been shot down on that basis, and the
loop representation approach seems so pretty. Either it will make sense,
or we will run into some obstacle, which is bound to be enlightening.

Let me briefly review the loop representation, without too many
technical details. For more details try the original paper by Rovelli
and Smolin (see \protect\hyperlink{week42}{``Week 42''} for a
reference), the book by Ashtekar (see \protect\hyperlink{week7}{``Week
7''}), or, especially if you're a mathematician, my review article
\href{http://math.ucr.edu/home/baez/knot.ps}{``Knots and quantum
gravity: progress and prospects''}.

There are 3 basic steps in the ``canonical quantization'' of general
relativity. At each step there is a vector space of quantum states, but
only in the last do we really need a Hilbert space of states, since only
when we're done do we want to be able to compute expectation values of
observables, which takes an inner product.

In what follows I'll talk about the simplest situation, where we have
the \emph{vacuum} Einstein equations \[G = 0\] where \(G\) is the
``Einstein tensor'' cooked up from the curvature of spacetime. Say
spacetime is of the form \(\mathbb{R} \times S\), where \(\mathbb{R}\)
is the real numbers (time) and \(S\) is a \(3\)-dimensional manifold
(space). We will think of \(S\) as the ``\(t = 0\) slice'' of
\(\mathbb{R} \times S\).

\begin{enumerate}
\def\labelenumi{\Roman{enumi})}
\item
  The first stage is to get the space of ``kinematical states''. In the
  quantum mechanics of a point particle on the line, the space of
  wavefunctions is a space of functions on the real line. Similarly, in
  quantum gravity we naively expect kinematical states to be functions
  on the space of Riemannian metrics on the \(3\)-dimensional manifold
  \(S\) we're taking to be ``space''. In the loop representation one
  does something a bit more clever, but let's move on and then come back
  to that.
\item
  The second stage is getting the space of ``diffeomorphism-invariant
  states''. In fact, Einstein's equations in coordinates look like
  \[G_{\mu \nu} = 0\] where the indices \(\mu\), \(\nu\) range from 0 to
  3. It's customary to work in coordinates \(x_{\mu}\) where \(x_0\) is
  ``time'' and the other three coordinates are the ``space'' coordinates
  on \(S\). Then classically, the equations \(G_{0 \mu} = 0\) serve as
  \emph{constraints} on the initial data for Einstein's equations, while
  the remaining equations describe time evolution. I.e., only for
  certain choices of a metric and its first time derivative at \(t = 0\)
  can we get a solution of Einstein's equations. In fact, \(G_{0 \mu}\)
  can be calculated knowing only the metric and its first time
  derivative at \(t = 0\), and the equations saying they are zero are
  the constraints that this data must satisfy to get a solution of
  Einstein's equations.
\end{enumerate}

Following the usual recipes of quantum theory, we want to turn these
constraints into \emph{operators} on the kinematical Hilbert space of
stage I, and then demand that the states relevant for physics be
annihilated by these operators. The ``diffeomorphism-invariant
subspace'' is the subspace of the kinematical state space that is
annihilated by the constraints corresponding to \(G_{0i}\) where
\(i = 1, 2, 3\). Let us put off for a moment why it's called what it is!

\begin{enumerate}
\def\labelenumi{\Roman{enumi})}
\setcounter{enumi}{2}
\tightlist
\item
  The third and final stage is getting the space of ``physical states''.
  Here we look at the subspace of diffeomorphism-invariant states that
  are also annihilated by the constaint corresponding to \(G_{00}\). The
  equation saying that a diffeomorphism-invariant state is annihilated
  by this constraint is called the ``Wheeler--DeWitt equation'', and this
  is generally regarded as the fundamental equation of quantum gravity.
\end{enumerate}

Now, it should make some sense why we call the ``physical states'' what
we do. These are quantum states satisfying the quantum analogues of the
constraints that the \emph{classical} initial data must satisfy to be
initial data for a solution of Einstein's equations. But why do we
impose the constraints \(G_{\mu \nu} = 0\) in two separate stages, and
call the states in part II ``diffeomorphism-invariant states''?

This is a very important question which gives quantum gravity much of
its curious character. In classical general relativity, \(G_{0i}\) not
only gives one of Einstein's equations, namely \(G_{0i}\) = 0, it also
``generates diffeomorphisms'' of the \(3\)-dimensional manifold \(S\)
representing space. If you don't quite know what this means, let me
simply say that in classical mechanics, observables give rise to
one-parameter families of symmetries. For example, momentum gives rise
to spatial translations, while energy (aka the Hamiltonian) gives rise
to time translations. We say that the observable ``generates'' the
one-parameter family of symmetries. This is (roughly) what I mean by
saying that \(G_{0i}\) generates diffeomorphisms of \(S\). Similarly,
\(G_{00}\) generates diffeomorphisms of the spacetime
\(\mathbb{R} \times S\) corresponding to time evolution.

A similar thing happens in quantum theory. BUT: in quantum theory, if a
state is annihilated by some observable, it implies that the state is
invariant under the one-parameter family of symmetries generated by that
observable. This is not true in classical mechanics. Indeed, it's rather
odd. But what it implies is that in step II we are really restricting
ourselves to kinematical states that are invariant under diffeomorphisms
of the spatial manifold \(S\). This is why we call them
``diffeomorphism-invariant'' states. Similarly, in step III we're
further restricting ourselves to states that are invariant under time
evolution. The final ``physical states'' are, at least heuristically,
invariant under \emph{all diffeomorphisms of spacetime}. (So maybe the
physical states are the ones that really should be called
``diffeomorphism-invariant'' --- but it's too late now.) While this may
seem odd, all it really means is that in the quantum theory of gravity
--- at least when one does it this way --- the physical states describe
only those aspects of the world that are independent of any choice of
coordinate system. That has a certain charm, philosophically speaking.
It is, however, not something physicists are used to.

Now, the general scheme outlined above has been around ever since the
work of DeWitt:

\begin{enumerate}
\def\labelenumi{\arabic{enumi})}
\tightlist
\item
  Bryce S. DeWitt, ``Quantum theory of gravity, I-III'', \emph{Phys.
  Rev.} \textbf{160} (1967), 1113--1148, \textbf{162} (1967) 1195--1239,
  1239--1256.
\end{enumerate}
\noindent
However, the problem has always been making the scheme mathematically
rigorous, or else to do some kind of calculations that shed some light
on the meaning of it all! There are lots of problems. Let me not delve
into them now, but simply cut directly to the ``new variables'' idea for
handling these problems. The key idea of Ashtekar was to use as basic
variables, not the metric on \(S\) and its first time derivative, but
the ``chiral spin connection'' on \(S\) and a ``complex frame field''.
To describe these would require a digression into differential geometry
that I'm not in the mood for right now, especially since I already
explained this stuff a bit in \protect\hyperlink{week7}{``Week 7''}.
(There I call the chiral spin connection the ``right-handed''
connection.)

I do, however, want to emphasize that the new variables rely heavily
upon some of the basic group-theoretic facts about 3 and 4 dimensions.
The group of rotations in 3d space is called \(\mathrm{SO}(3)\), because
mathematically these are \(3\times3\) orthogonal matrices with
determinant \(1\). Now, a key fact in math and physics is that this
group has the group \(\mathrm{SU}(2)\) of \(2\times2\) complex unitary
matrices with determinant \(1\) as a ``double cover''. This means
roughly that there are two elements of this other group corresponding to
each element of \(\mathrm{SO}(3)\). It's this fact that allows the
existence of spin-\(\frac12\) particles!

Now, \(\mathrm{SU}(2)\) is sitting inside a bigger group,
\(\mathrm{SL}(2,\mathbb{C})\), the group of all \(2\times2\) complex
matrices with determinant \(1\), not necessarily unitary. Just as
\(\mathrm{SU}(2)\) is used to describe the symmetries of
spin-\(\frac12\) particles in space, \(\mathrm{SL}(2,\mathbb{C})\)
describes the symmetries of spin-\(\frac12\) particles in spacetime. The
reason is that \(\mathrm{SL}(2,\mathbb{C})\) is the double cover of the
group \(\mathrm{SO}(3,1)\) of Lorentz transformations.

Given a Riemannian metric on the space \(S\), there is always an
``\(\mathrm{SO}(3)\) connection'' describing how objects rotate when you
move them around a loop, due to the curvature of space. This is called
the Levi--Civita connection. With a little work we can also think of this
as an \(\mathrm{SU}(2)\) connection. However, Ashtekar works instead
with the chiral spin connection, which is an
\(\mathrm{SL}(2,\mathbb{C})\) connection cooked up from the Levi--Civita
connection and the first time derivative of the metric (which turns out
to be closely related to the ``extrinsic curvature'' of \(S\) as it sits
in the spacetime \(\mathbb{R} \times S\).)

The great advantage of Ashtekar's ``new variables'' is that the
Hamiltonian and diffeomorphism constraints are simpler in these
variables. Unfortunately, they lead to a curious new issue which at
first seemed very nasty --- the problem of ``reality conditions''. This
has a lot to do with going from \(\mathrm{SU}(2)\), which is a ``real''
group in a technical sense, to \(\mathrm{SL}(2,\mathbb{C})\), which is a
``complex'' group that's roughly twice as big. Essentially, Ashtekar's
formalism seems at first to be better suited to general relativity with
a complex-valued metric than to good old ``real'' general relativity.
For quite a while people didn't know quite what to do about this, so a
lot of work on the new variables more or less ignores this issue.
Luckily, there is now a very elegant approach to handling it, worked out
by Ashtekar and collaborators. They are coming out with a couple of
papers on this, hopefully by mid-November:

\begin{enumerate}
\def\labelenumi{\arabic{enumi})}
\setcounter{enumi}{1}
\item
  Abhay Ashtekar, Jerzy Lewandowski, Donald Marolf, Jos\'e Mour\~ao and Thomas
  Thiemann, ``Coherent state transforms for spaces of connections'',  
  \emph{Jour. Functional Analysis} \textbf{135} (1996), 519--551.
  Also available
  as \href{https://arxiv.org/abs/gr-qc/9412014}{\texttt{gr-qc/9412014}}.

  Abhay Ashtekar, Jerzy Lewandowski, Donald Marolf, Jos\'e Mour\~ao and Thomas
  Thiemann, ``Quantization of diffeomorphism invariant theories of connections with local degrees of freedom'', \emph{Jour. Math. Phys.} \textbf{36} (1995), 6456--6493.   Also available as \href{https://arxiv.org/abs/gr-qc/9504018}{\texttt{gr-qc/9504018}}.
 
\end{enumerate}
\noindent
The first paper constructs a kind of transform that takes functions on
the space of \(\mathrm{SU}(2)\) connections on \(S\) into functions on
the space of \(\mathrm{SL}(2,\mathbb{C})\) connections on \(S\). The
``kinematical states'' in Ashtekar's approach are, roughly speaking,
functions of the latter kind. (Really they are more like ``measures''.)
This is some really pretty mathematics --- it's a kind of generalization
of the Bargmann--Segal transform to the case of functions on spaces of
connections.

Physically, the transform allows us to relate Ashtekar's approach to the
traditional ``metric'' approach much more clearly, since, as I
described, \(\mathrm{SU}(2)\) connections are closely related to metrics
on \(S\). The second paper should treat the physics behind this in more
detail, and also describe a rigorous construction of ``loop states'' ---
a large class of diffeomorphism-invariant states which Rovelli and
Smolin have claimed are actually \emph{physical} states. (For more on
these, see below.) This means that to check Rovelli and Smolin's claim,
the main thing we need is a rigorous treatment of the Hamiltonian
constraint in quantum gravity.

Unfortunately, this is where it gets scary, since the Hamiltonian
constraint is a very tricky thing. For more on it, try:

\begin{enumerate}
\def\labelenumi{\arabic{enumi})}
\setcounter{enumi}{2}
\item
  M. Blencowe, ``The Hamiltonian constraint in quantum gravity'', 
  \emph{Nucl. Phys.} \textbf{B341} (1990), 213--251.

  Bernd Br\"ugmann and Jorge Pullin, ``On the constraints of quantum gravity 
  in the loop representation'', \emph{Nucl. Phys.} \textbf{B390}
  (1993), 399--438.

  Bernd Br\"ugmann, \emph{On the Constraints of Quantum General Relativity in 
  the Loop Representation},  Ph.D.~Thesis, Syracuse University, 1993.
\end{enumerate}

The loop states of Rovelli and Smolin are in one-to-one correspondence
with \emph{knots} in space, or more precisely, isotopy classes of knots.
(Roughly, two knots are isotopic if you can get one from the other by
applying a diffeomorphism of space that can be continuously deformed to
the identity.)

What is the physical meaning of these loop states? Roughly it's this.
Say you take a spin-\(\frac12\) particle and move it around in a path
that traces out a knot. When you do this using the Levi--Civita
connection, it comes back ``rotated'' by some \(\mathrm{SU}(2)\) matrix.
If you take the trace of this matrix (sum of diagonal entries) and
divide by two, you get a number between -1 and 1. This number is called
a ``Wilson loop''.

This should remind you of the Bohm--Aharonov effect where a split
electron beam takes two paths from \(A\) to \(B\). Depending on the
magnetic flux through the loop, one can have constructive or destructive
interference in the split beam experiment. Mathematically, one can
imagine moving the electron around a loop that starts at \(A\), goes to
\(B\) by one path, and goes back to \(A\) by the other path. If this
phase corresponding to going around this loop is 1 we get total constructive interference
in the split beam experiment, while if it's -1 we get total destructive
interference. So, just as the Bohm--Aharonov effect measures interference
effects due to the magnetic field, the above Wilson loop sort of
measures the interference effects due to \emph{gravity!}

In the Rovelli--Smolin loop state corresponding to a particular knot
\(K\), the expectation value of a Wilson loop around any knot \(K'\) will be 1
if \(K\) and \(K'\) are isotopic, and 0 otherwise! That's the physical meaning
of the loop states: they describe quantum states of geometry in terms of
the resulting interference effects on spin-\(\frac12\) particles.

Now, there is a more general kind of diffeomorphism-invariant state than the
loop states. These are the spin network states! Here one fancies up the
Wilson loop idea and imagines a graph embedded in space --- i.e.~a bunch
of edges and vertices --- where each edge is labelled by a spin that can
be \(0\),\(\frac{1}{2}\),\(1\),\(\frac{3}{2}\), etc. In the simplest flavor of spin
network, one only allows 3 edges to meet at each vertex, and requires
\(j_3\) to be of the form
\[j_3 = |j_1-j_2|, \; |j_1-j_2| + 1, \ldots, \; j_1+j_2-1, \; j_1+j_2.\] where
\(j_1\), \(j_2\), \(j_3\) are the spins labelling the edges adjacent to
the given vertex. For example, we can have the the three spins be
\(\frac{1}{2}\),\(3\), and \(\frac{5}{2}\), because it's possible for a spin-\(\frac12\)
particle and a spin-3 particle to interact and form a spin-\(\frac{5}{2}\)
particle. Here by ``possible'' I simply mean that it doesn't violate
conservation of angular momentum. Mathematicians would say the spins
should be thought of as irreducible representations of
\(\mathrm{SU}(2)\), and the condition above is just the condition that
the representation \(j_3\) appears as a summand in the tensor product of
the representations \(j_1\) and \(j_2\).

Just as we can compute a kind of ``Wilson loop'' number from a knot that
a spin-\(\frac12\) particle goes around, we can compute a number from a
spin network. I've thought about spin networks for quite a while, since
they are very important in topological quantum field theories. A great
introduction to how they show up in TQFTs, by the way, is:

\begin{enumerate}
\def\labelenumi{\arabic{enumi})}
\setcounter{enumi}{3}
\tightlist
\item
   Louis Crane, Louis H. Kauffman, and David N. Yetter, 
  ``State-sum invariants of manifolds, I'', 
  available as
  \href{https://arxiv.org/abs/hep-th/9409167}{\texttt{hep-th/9409167}}.
\end{enumerate}
\noindent
This explains how to cook up 3d quantum gravity (or more precisely, the
Turaev--Viro model) and a 4d TQFT field theory called the Crane--Yetter
model using spin networks.

However, Rovelli's talk on spin network states in quantum gravity (see
\protect\hyperlink{week41}{``Week 41''}), followed by some good
conversations, got me motivated to write up something on spin network
states:

\begin{enumerate}
\def\labelenumi{\arabic{enumi})}
\setcounter{enumi}{4}
\tightlist
\item
   John Baez,``Spin network states in gauge theory'', available as
  \href{https://arxiv.org/abs/gr-qc/9411007}{\texttt{gr-qc/9411007}}.
\end{enumerate}
\noindent
Basically, I show that in the loop representation of any gauge theory,
states at the kinematical level can be described by spin networks,
slightly generalized. Heck, I'll quote my abstract:

\begin{quote}
Given a real-analytic manifold \(M\), a compact connected Lie group
\(G\) and a principal \(G\)-bundle \(P \to M\), there is a canonical
`generalized measure' on the space \(A/G\) of smooth connections on
\(P\) modulo gauge transformations. This allows one to define a Hilbert
space \(L^2(A/G)\). Here we construct a set of vectors spanning
\(L^2(A/G)\). These vectors are described in terms of `spin networks':
graphs \(\varphi\) embedded in \(M\), with oriented edges labelled by
irreducible unitary representations of \(G\), and with vertices labelled
by intertwining operators from the tensor product of representations
labelling the incoming edges to the tensor product of representations
labelling the outgoing edges. We also describe an orthonormal basis of
spin networks associated to any fixed graph \(\varphi\). We conclude
with a discussion of spin networks in the loop representation of quantum
gravity, and give a category-theoretic interpretation of the spin
network states.
\end{quote}

I'm now hard at work trying to show that spin networks also give a
complete description of states at the diffeomorphism-invariant level.
Well, actually right NOW I'm goofing off by writing this darn thing, but
you know what I mean.

Rovelli and Smolin have come out with one of their papers on spin
networks and they should be coming out with another soon. These are not
about the rigorous mathematics of spin network states, but how to use
them to really understand the \emph{physics} of quantum gravity. The
first one out is:

\begin{enumerate}
\def\labelenumi{\arabic{enumi})}
\setcounter{enumi}{5}
\tightlist
\item
   Carlo Rovelli and Lee Smolin, ``Discreteness of area and volume in quantum gravity'', 
   available as
  \href{https://arxiv.org/abs/gr-qc/9411005}{\texttt{gr-qc/9411005}}.
\end{enumerate}
\noindent
This is perhaps the most careful computation so far that derives
\emph{discreteness} of geometrical quantities directly from Einstein's
equations and the principles of quantum theory! Let me quote the
abstract:

\begin{quote}
We study the operator that corresponds to the measurement of volume, in
non-perturbative quantum gravity, and we compute its spectrum. The
operator is constructed in the loop representation, via a regularization
procedure; it is finite, background independent, and
diffeomorphism-invariant, and therefore well defined on the space of
diffeomorphism invariant states (knot states). We find that the spectrum
of the volume of any physical region is discrete. A family of
eigenstates are in one to one correspondence with the spin networks,
which were introduced by Penrose in a different context. We compute the
corresponding component of the spectrum, and exhibit the eigenvalues
explicitly. The other eigenstates are related to a generalization of the
spin networks, and their eigenvalues can be computed by diagonalizing
finite dimensional matrices. Furthermore, we show that the eigenstates
of the volume diagonalize also the area operator. We argue that the
spectra of volume and area determined here can be considered as
predictions of the loop-representation formulation of quantum gravity on
the outcomes of (hypothetical) Planck-scale sensitive measurements of
the geometry of space.
\end{quote}

\hypertarget{week44}{%
\section{November 6, 1994}\label{week44}}

\begin{center}
\textbf{SPECIAL EDITION: THE END OF DONALDSON THEORY?}
\end{center}

I got some news today from Allen Knutson. Briefly, it appears that
Witten has come up with a new way of doing Donaldson theory that is far
easier than any previously known. According to Taubes, many of the main
theorems in Donaldson theory should now have proofs that are
\(1/1000\)th as long!

I suppose to find this exciting one must already have some idea of what
Donaldson theory is. Briefly, Donaldson theory is a theory born in the
1980s that revolutionized the study of smooth \(4\)-dimensional
manifolds by using an idea from physics, namely, the self-dual
Yang--Mills equations. The Yang--Mills equations describe most of the
forces we know and love (not gravity), but only in 4 dimensions can one
get solutions of them of a special form, known as self-dual solutions.
(In physics these self-dual solutions are known as instantons, and they
were used by 't Hooft to solve a problem plaguing particle physics,
called the \(\mathrm{U}(1)\) puzzle.)

Mathematically, \(4\)-dimensional manifolds are very different from
manifolds of any other dimension! For example, one can ask whether
\(\mathbb{R}^n\) admits any smooth structure other than the usual one.
(Technically, a smooth structure for a manifold is a maximal set of
coordinate charts covering the manifold which have smooth transition
functions. Loosely, it's a definition of what counts as a smooth
function.) The answer is no --- \emph{except} if \(n = 4\), where
there are uncountably many smooth structures! These ``exotic
\(\mathbb{R}^4\)'s'' were discovered in the 1980's, and their existence
was shown using the work of Donaldson using the self-dual solutions of
the Yang--Mills equation, together with work of the topologist Freedman.
More recently, a refined set of invariants of smooth 4-manifolds, the
Donaldson invariants, have been developed using closely related ideas.

Some references are:

\begin{enumerate}
\def\labelenumi{\arabic{enumi})}
\item
  Simon K. Donaldson and Peter B. Kronheimer, \emph{The Geometry of Four-Manifolds}, 
  Oxford University Press, Oxford, 1990.

  Simon K.  Donaldson, ``Polynomial invariants for smooth four-manifolds'', 
  \emph{Topology} \textbf{29} (1990), 257--315.

  Daniel S. Freed and Karen K. Uhlenbeck, \emph{Instantons and 
  Four-Manifolds}, Springer, Berlin, 1984.

  Charles Nash, \emph{Differential Topology and Quantum Field Theory}, 
  Academic Press, London, 1991.
\end{enumerate}

This is an extremely incomplete list, but it should be enough to get
started. Or, while you wait for the new, simplified treatments to come
out, you could make some microwave popcorn and watch the following
video:

\begin{enumerate}
\def\labelenumi{\arabic{enumi})}
\setcounter{enumi}{1}
\tightlist
\item
   Simon K. Donaldson, ``Geometry of four dimensional manifolds'', 
  videocassette (ca. 60 min.), color, American Mathematical
  Society, Providence, Rhode Island, 1988.
\end{enumerate}

Now, what follows is my interpretation of David Dror Ben-Zvi's comments
on a lecture by Clifford Taubes entitled ``Witten's Magical Equation'',
these comments being kindly passed on to me by Knutson. I have tried to
flesh out and make sense of what I received, and this required some
work, and I may have screwed up some things. Please take it all with a
grain of salt. I only hope it gives some of the flavor of what's going
on!

So, we start with a compact oriented 4-manifold \(X\) with \(L\) a
complex line bundle over \(X\) having first Chern class equal to
\(w_2\), the second Stiefel--Whitney class of \(TX\), modulo \(2\). If
\(X\) is spin (meaning that the \(w_2 = 0\)), take the bundle of spinors
over \(X\). Otherwise, pick a \(\mathrm{Spin}^c\) bundle and take the bundle of
complex spinors over X. Note that a \(\mathrm{Spin}^c\) structure is enough to define
complex spinors on \(X\), and it will always exist if \(w_2\) is the
\(\mod 2\) reduction of an integral characteristic class. For more on
this sort of stuff, try:

\begin{enumerate}
\def\labelenumi{\arabic{enumi})}
\setcounter{enumi}{2}
\tightlist
\item
  H. Blaine Lawson, Jr.~and Marie-Louise
  Michelson, \emph{Spin Geometry}, Princeton U. Press, Princeton, 1989.
\end{enumerate}

In either case, take our bundle of spinors, tensor it with the square
root of \(L\), and call the resulting bundle \(B\). (Perhaps someone can
explain to me why \(L\) has a square root here; it's obvious if \(X\) is
spin, but I don't understand the other case so well.) The data for our
construction are now a connection \(A\) on \(L\), and a section \(\psi\)
of the self-dual part of \(B\). (Note: I'm not sure what the ``self-dual
part of \(B\)'' is supposed to mean. I guess it is something required to
make the right-hand side of the formula below be self-dual in the
indices \(a,b\).) Consider now two equations. The first is the Dirac
equation for \(\psi\). The second is that the self-dual part \(F^+\) of
the curvature of \(A\) be given in coordinates as
\[F^+_{ab} = -\frac12 \langle\psi, e^a e^b \psi\rangle\] where the basis
\(1\)-forms \(e^a\), \(e^b\) act on \(\psi\) by Clifford multiplication.

Next form the moduli space \(M\) of solutions \((A, \psi)\) modulo the
action of the automorphisms of \(L\). The wonderful fact is that this
moduli space is always compact, and for generic metrics it's a smooth
manifold. Still more wonderfully (here I read the lines between what was
written), it is a kind of substitute for the moduli space normally used
in Donaldson theory, namely the moduli space of instantons. It is much
nicer in that it lacks the singularities characteristic of the other
space.

What this means is that everything becomes easy! Apparently Taubes,
Kronheimer, Mrowka, Fintushel, Stern and the other bigshots of Donaldson
theory are frenziedly turning out new results even as I type these
lines. On the one hand, the drastic simplifications are a bit
embarassing, since the technical complications of Donaldson theory were
the stuff of many erudite and difficult papers. On the other hand,
Donaldson invariants were always notoriously difficult to compute.
Taubes predicted that a purely combinatorial formula for them may be
around within a year. (Here it is interesting to note the work of Crane,
Frenkel, and Yetter in that direction; see
\protect\hyperlink{week2}{``Week 2''} and
\protect\hyperlink{week38}{``Week 38''}.) This is sure to lead to a
deeper understanding of \(4\)-dimensional topology, and quite possibly,
4-dimensional physics as well.

\hypertarget{week45}{%
\section{November 12, 1994}\label{week45}}

\begin{center}
\textbf{DONALDSON THEORY UPDATE}
\end{center}

In the previous edition of ``This Week's Finds'' I mentioned a burst of
recent work on Donaldson theory. I provocatively titled it ``The End of
Donaldson Theory?'', since the rumors I was hearing tended to be phrased
in such terms. But I hope I made it clear at the conclusion of the
article that this recent work should lead to a lot of \emph{new} results
in 4-dimensional topology! An example is Kronheimer and Mrowka's proof
of the Thom conjecture.

Many thanks to my network of spies for obtaining a preprint of the
following paper:

\begin{enumerate}
\def\labelenumi{\arabic{enumi})}
\tightlist
\item
  Peter B. Kronheimer and Tomasz S. Mrowka, ``The genus of 
  embedded surfaces in the projective plane''.
\end{enumerate}
\noindent
Let me simply quote the beginning of the paper:

\begin{quote}
The genus of a smooth algebraic curve of degree \(d\) in
\(\mathbb{CP}^2\) is given by the formula \(g = (d-1)(d-2)/2\). A
conjecture sometimes attributed to Thom states that the genus of the
algebraic curve is a lower bound for the genus of any smooth 2-manifold
representing the same homology class. The conjecture has previously been
proved for \(d \leqslant 4\) and for \(d = 6\), and less sharp lower
bounds for the genus are known for all degrees {[}references omitted{]}.
In this note we confirm the conjecture.

Theorem 1. Let \(S\) be an oriented 2-manifold smoothly embedded in
\(\mathbb{CP}^2\) so as to represent the same homology class as an
algebraic curve of degree \(d\). Then the genus \(g\) of \(S\) satisfies
\(g \geqslant (d-1)(d-2)/2\).

Very recently, Seiberg and Witten {[}references below{]} introduced new
invariants of 4-manifolds, closely related to Donaldson's polynomial
invariants {[}reference omitted{]}, but in many respects much simpler to
work with. The new techniques have led to more elementary proofs of many
theorems in the area. Given the monopole equation and the vanishing
theorem which holds when the scalar curvature is positive (something
which was pointed out by Witten), the rest of the argument presented
here is not hard to come by. A slightly different proof of the Theorem,
based on the same techniques, has been found by Morgan, Szabo and
Taubes.
\end{quote}

\noindent
The reference to Donaldson's polynomial invariants appears in
\protect\hyperlink{week44}{``Week 44''}. The references to the new
Seiberg--Witten invariants are:

\begin{enumerate}
\def\labelenumi{\arabic{enumi})}
\setcounter{enumi}{1}
\item
  Edward Witten,``Monopoles and four-manifolds'', available as
  \href{https://arxiv.org/abs/hep-th/9411102}{\texttt{hep-th/9411102}}.

  Edward Witten and
  Nathan Seiberg, ``Electric-magnetic duality, monopole condensation, and confinement in
  \(N=2\) supersymmetric Yang--Mills theory'', available as
  \href{https://arxiv.org/abs/hep-th/9407087}{\texttt{hep-th/9407087}}.

  Edward Witten and Nathan Seiberg, ``Monopoles, duality and chiral symmetry 
  breaking in \(N=2\) supersymmetric QCD'', 
  available as
  \href{https://arxiv.org/abs/hep-th/9408099}{\texttt{hep-th/9408099}}.
\end{enumerate}

Differential geometers attempting to read the second two papers will
find that they contain no instance of the term ``Donaldson theory'', and
they may be frustrated to find that these are very much \emph{physics}
papers. They concern the ground states of supersymmetric Yang--Mills
theory in 4 dimensions with gauge group \(\mathrm{SU}(2)\). The ``ground
states'' of a field theory are its least-energy states, which represent
candidates for the physical vacuum. In certain theories there is not a
unique ground state, but instead a ``moduli space'' of ground states.
Seiberg and Witten study these moduli spaces of ground states in both
the classical and quantum versions of \(\mathrm{SU}(2)\) supersymmetric
Yang--Mills theory in 4 dimensions. They also consider the theory coupled
to spinor fields, which they call ``quarks'', using the analogy of the
theory to quantum chromodynamics, aka ``QCD''.

I haven't had time to go through their papers, since this isn't my main
focus of interest. Perhaps the most useful thing I can do at this point
is to use Kronheimer and Mrowka's clear description of their moduli
space (which is presumably closely related to Seiberg and Witten's
moduli spaces) to simplify and fill in the holes of what I wrote in
\protect\hyperlink{week44}{``Week 44''}. I will aim my exposition to
mathematicians, but make some elementary digressions on physics to spice
things up.

We start with a compact oriented Riemannian 4-manifold \(X\), and assume
we are given a $\mathrm{Spin}^\mathbb{C}$ structure on \(X\). Recall the meaning of this.
First, the orthonormal frame bundle of \(X\) has structure group
\(\mathrm{SO}(4)\), and a spin structure would be a double cover of this
which is a principal bundle with structure group given by the double
cover of \(\mathrm{SO}(4)\), namely
\(\mathrm{SU}(2) \times \mathrm{SU}(2)\). Thus we get two principal
bundles with structure group \(\mathrm{SU}(2)\), the left-handed and
right-handed spin bundles. Using the fundamental representation of
\(\mathrm{SU}(2)\), we obtain two vector bundles called the bundles of
left-handed and right-handed spinors. This ``handedness'' or
``chirality'' phenomenon for spinors is of great importance in physics,
since neutrinos are left-handed spinors --- meaning, in down-to-earth
terms, that they always spin clockwise relative to their direction of
motion. The fact that the laws of nature lack chiral symmetry came as
quite a shock when it was first discovered, and part of Seiberg and
Witten's motivation in their second paper is to study mechanisms for
``spontaneous breaking'' of chiral symmetry. This means simply that
while the theory has chiral symmetry, its ground states need not.

A $\mathrm{Spin}^\mathbb{C}$ structure is a bit more subtle, but it allows us to define
bundles of left-handed and right-handed spinors as \(\mathrm{U}(2)\)
bundles, which Kronheimer and Mrowka denote by \(W+\) and \(W-\). The
determinant bundle \(L\) of \(W+\) is a line bundle on \(X\). The first
big ingredient of the theory is a hermitian connection \(A\) on \(L\).
In physics lingo this is the vector potential of a \(\mathrm{U}(1)\)
gauge field. This gives a Dirac operator \(D_A\) mapping sections of
\(W+\) to sections of \(W-\). The connection \(A\) has curvature \(F\),
and the self-dual part \(F^+\) of \(F\) can be identified with a section
of \(\mathfrak{sl}(W+)\). (This is just a global version of the
isomorphism between the self-dual part of \(\Lambda^2 \mathbb{C}^4\) and
\(\mathfrak{sl}(2,\mathbb{C})\).)

The second big ingredient of the theory is a section \(\psi\) of \(W+\),
i.e.~a left-handed spinor field. There is a way to pair two sections of
\(W+\) to get a section of \(\mathfrak{sl}(W+)\), which we write as
\(\sigma(.,.)\) and which is conjugate-linear in the first argument and
linear in the second. This is a global version of the similar pairing
\[\sigma(.,.)\colon \mathbb{C}^2 \times \mathbb{C}^2 \to \mathfrak{sl}(2,\mathbb{C})\]
where \(\sigma(v,w)\) given by taking the traceless part of the
\(2\times2\) matrix \(v^* \otimes w\). Here \(v^*\) is the element of
the dual of \(\mathbb{C}^2\) coming from \(v\) via the inner product on
\(\mathbb{C}^2\).

To get the magical moduli space, we consider solutions \((A,\psi)\) of
\[\begin{aligned}D_A\psi &= 0 \\ F^+ &= i\sigma(\psi,\psi).\end{aligned}\]
Here we are thinking of \(F^+\) as a section of \(\mathfrak{sl}(W+)\).
These are pretty reasonable equations for some sort of massless
left-handed spinor field coupled to a \(\mathrm{U}(1)\) gauge field. Let
\(M\) be the space of solutions modulo gauge transformations. Kronheimer
and Mrowka show the ``moduli space'' \(M\) is compact.

One can also perturb the equations above as follows. If we have any
self-dual \(2\)-form \(\delta\) on \(X\) we can consider
\[\begin{aligned}D_A\psi &= 0 \\ F^+ +i\delta &= i\sigma(\psi,\psi).\end{aligned}\]
and get a moduli space \(M(\delta)\). This will still be compact if
\(\delta\) is nice (here I gloss over issues of analysis).

Now, if \(X\) has an almost complex structure, Kronheimer and Mrowka
show that one can pick a  $\mathrm{Spin}^\mathbb{C}$  structure for \(X\) such that, for
``good'' metrics and generic small \(\delta\), \(M(\delta)\) is a
compact 0-dimensional manifold. Using this fact and some geometrical
yoga, it follows that the number \(n\) of points in \(M(\delta)\),
counted mod 2, is independent of (such) \(\delta\). (This is
essentially a glorified version of the fact that, when you look at the
multiple images of an object in a warped mirror and slowly bend the
mirror, the images generically appear or disappear in pairs.) Moreover,
if the self-dual Betti number \(b^+\) of \(X\) is \(> 1\), the space of
good metrics is path-connected, and \(n\) mod 2 is independent of the
choice of good metric. Kronheimer and Mrowka call this a ``simple
mod 2 version of the invariants of Seiberg and Witten''. It is one
ingredient of their proof of the Thom conjecture.

\hypertarget{week46}{%
\section{December 12, 1994}\label{week46}}

I will be on sabbatical during the first half of 1995. I'll be roaming
hither and thither, and also trying to get some work done on
\(n\)-categories, quantum gravity and such, so this will be the last
``This Week's Finds'' for a while. I have also taken a break from being
a co-moderator of \texttt{sci.physics.research}.

So, let me sign off with a roundup of diverse and sundry things! I'm
afraid I'll be pretty terse about describing some of them. First for
some news of general interest, then a little update on Seiberg--Witten
theory, then some neat stuff on TQFTs, \(n\)-categories, quantum gravity
and all that, and then various other goodies\ldots.

\begin{enumerate}
\def\labelenumi{\arabic{enumi})}
\item
   Gary Stix, ``The speed of write'', \emph{Scientific American},
  Dec.~1994, 106--111.

  Jacques Leslie, ``Goodbye, Gutenberg: pixilating peer review is revolutionizing 
  scholarly journals'',  \emph{Wired} 2.10,
  Oct.~1994.  Available as \href{https://www.wired.com/1994/10/ejournals/}{\texttt{https://www.wired.com/1994/10/ejournals/}}.
\end{enumerate}

Among other things, the above articles show that Paul Ginsparg is
starting to get the popular recognition he deserves for starting up
\texttt{hep-th}. In case anyone out there doesn't know yet,
\texttt{hep-th} is the ``high-energy physics --- theoretical'' preprint
archive, which revolutionized communications within this field by making
preprints easily available world-wide, thus rendering many (but not all)
aspects of traditional journals obsolete. The idea was so good it
quickly spread to other subjects. Within physics it went like this:

\begin{itemize}
\tightlist
\item
  High Energy Physics --- Theory (\texttt{hep-th}), started 8/91
\item
  High Energy Physics --- Lattice (\texttt{hep-lat}), started 2/92
\item
  High Energy Physics --- Phenomenology (\texttt{hep-ph}), started 3/92
\item
  Astrophysics (\texttt{astro-ph}), started 4/92
\item
  Condensed Matter Theory (\texttt{cond-mat}), started 4/92
\item
  General Relativity \& Quantum Cosmology (\texttt{gr-qc}), started 7/92
\item
  Nuclear Theory (\texttt{nucl-th}), started 10/92
\item
  Chemical Physics (\texttt{chem-ph}), started 3/94
\item
  High Energy Physics --- Experiment (\texttt{hep-ex}), started 4/94
\item
  Accelerator Physics (\texttt{acc-phys}), started 11/94
\item
  Nuclear Experiment (\texttt{nucl-ex}), started 11/94
\item
  Materials Theory (\texttt{mtrl-th}), started 11/94
\item
  Superconductivity (\texttt{supr-con}), started 11/94
\end{itemize}

Similar archives are sprouting up in mathematics (see below --- but also
note the existence of the American Mathematical Society preprint server,
described later in this week's finds).

There are many ways to access these preprint archives, since Ginsparg
has kept up very well with the times --- indeed, so much better than I
that I'm afraid to go into any details for fear of making a fool of
myself. The \emph{dernier cri}, I suppose, is to access the archives using the
World-Wide Web, which is conveniently done by opening the document

\begin{verbatim}
http://xxx.lanl.gov
\end{verbatim}
\noindent
If this makes no sense to you, my first and very urgent piece of advice
is to learn about the World-Wide Web (WWW), Mosaic, and the like, since
they are wonderful and very simple to use! In the meantime, however, you
can simply send mail to various addresses with subject header

\begin{verbatim}
help
\end{verbatim}

and no message body, in order to get information. Some addresses are:

\begin{itemize}
\item
  \texttt{acc-phys@xxx.lanl.gov} (accelerator physics)
\item
  \texttt{astro-ph@xxx.lanl.gov} (astrophysics)
\item
  \texttt{chem-ph@xxx.lanl.gov} (chemical physics)
\item
  \texttt{cond-mat@xxx.lanl.gov} (condensed matter)
\item
  \texttt{funct-an@xxx.lanl.gov} (functional analysis)
\item
  \texttt{gr-qc@xxx.lanl.gov} (general relativity / quantum cosmology)
\item
  \texttt{hep-lat@ftp.scri.fsu.edu} (computational and lattice physics)
\item
  \texttt{hep-ph@xxx.lanl.gov} (high energy physics phenomenological)
\item
  \texttt{hep-th@xxx.lanl.gov} (high energy physics formal)
\item
  \texttt{hep-ex@xxx.lanl.gov} (high energy physics experimental)
\item
  \texttt{nucl-th@xxx.lanl.gov} (nuclear theory)
\item
  \texttt{nucl-ex@xxx.lanl.gov} (nuclear experiment)
\item
  \texttt{mtrl-th@xxx.lanl.gov} (materials theory)
\item
  \texttt{supr-con@xxx.lanl.gov} (superconductivity)
\item
  \texttt{alg-geom@publications.math.duke.edu} (algebraic geometry)
\item
  \texttt{auto-fms@msri.org} (automorphic forms)
\item
  \texttt{cd-hg@msri.org} (complex dynamics \& hyperbolic geometry)
\item
  \texttt{dg-ga@msri.org} (differential geometry \& global analysis)
\item
  \texttt{nlin-sys@xyz.lanl.gov} (non-linear systems)
\item
  \texttt{cmp-lg@xxx.lanl.gov} (computation and language)
\item
  \texttt{e-mail@xxx.lanl.gov} (e-mail address database)
\end{itemize}

One might also want to check out the:

\begin{quote}
Directory of Electronic Journals, Newsletters, and Academic Discussion
Lists: Send e-mail to \texttt{ann@cni.org} at the Association of
Research Libraries, +1 (202) 296-2296, fax +1 (202) 872 0884.
\end{quote}

Let me say a bit about how the AMS preprint server works. Assuming you
are hip to the WWW, just go to

\begin{verbatim}
http://e-math.ams.org
\end{verbatim}
\noindent
You will then see a menu, and you can click on ``Mathematical
Preprints'', and then ``AMS Preprint Server'', where preprints are
classified by subject. Alternatively, click on ``New Items This Month
(all Subjects)''.

On a related note, you can also get some AMS stuff using telnet by doing

\begin{verbatim}
telnet e-math.ams.org
\end{verbatim}
\noindent
and using

\begin{verbatim}
e-math
\end{verbatim}
\noindent
as login and password. This doesn't seem to get you to the preprints,
though. For gopher fans,

\begin{verbatim}
gopher e-math.ams.org
\end{verbatim}
\noindent
has roughly similar effects.

\begin{enumerate}
\def\labelenumi{\arabic{enumi})}
\setcounter{enumi}{1}
\tightlist
\item
  Carlo Rovelli and Lee Smolin, ``Spin networks in quantum gravity''.
\end{enumerate}

This paper is closely related to the earlier one in which Rovelli and
Smolin argue that discreteness of area and volume arise naturally in the
loop representation of quantum gravity, and also to my own paper on spin
networks. (See \protect\hyperlink{week43}{``Week 43''} for more on
these, and a brief intro to spin networks.) Basically, while my paper
shows that spin networks give a kind of basis of states for gauge
theories with arbitrary (compact, connected) gauge group, in this paper
Rovelli and Smolin concentrate on the gauge groups
\(\mathrm{SL}(2,\mathbb{C})\) and \(\mathrm{SU}(2)\), which are relevant
to quantum gravity, and work out a lot of aspects particular to this
case, in more of a physicist's style. This makes spin networks into a
practical computational tool in quantum gravity, used to great effect in
the paper on the discreteness of area and volume.

\begin{enumerate}
\def\labelenumi{\arabic{enumi})}
\setcounter{enumi}{2}
\item
    Abhay Ashtekar, ``Recent mathematical developments in quantum general 
    relativity'', available as
  \href{https://arxiv.org/abs/gr-qc/9411055}{\texttt{gr-qc/9411055}}
  (discussed in \protect\hyperlink{week37}{``Week 37''}).

  Abhay
  Ashtekar, Jerzy Lewandowski, Donald Marolf, Jose Mourao and Thomas
  Thiemann, ``Coherent state transforms for spaces of connections'', available as
  \href{https://arxiv.org/abs/gr-qc/9412014}{\texttt{gr-qc/9412014}}
  (discussed in \protect\hyperlink{week43}{``Week 43''}).
\end{enumerate}

These are two papers on the loop representation of quantum gravity which
I talked about in earlier ``finds'', and are out now. The former is a
nice review of recent mathematically rigorous work; the latter takes a
tremendous step towards handling the infamous ``reality conditions''
problem.

\begin{enumerate}
\def\labelenumi{\arabic{enumi})}
\setcounter{enumi}{3}
\tightlist
\item
  Abhay Ashtekar and Jerzy Lewandowski, ``Differential geometry 
  on the space of connections via graphs and
  projective limits'', available as
  \href{https://arxiv.org/abs/hep-th/9412073}{\texttt{hep-th/9412073}}.
\end{enumerate}
\noindent
I've spoken quite a bit about doing rigorous functional
\emph{integration} in gauge theory using ideas from the loop
representation; this paper treats functional \emph{derivatives} and
other things that are more differential than integral in nature. This is
crucial in quantum gravity because the main remaining mystery, the
Wheeler--DeWitt equation or Hamiltonian constraint, involves a
differential operator on the space of connections. (For a wee bit more,
try \protect\hyperlink{week11}{``Week 11''} or
\protect\hyperlink{week43}{``Week 43''}, where the Hamiltonian
constraint is simply written as \(G_{00} = 0\).)

Let me quote their abstract:

\begin{quote}
In a quantum mechanical treatment of gauge theories (including general
relativity), one is led to consider a certain completion, \(A\), of the
space of gauge equivalent connections. This space serves as the quantum
configuration space, or, as the space of all Euclidean histories over
which one must integrate in the quantum theory. \(A\) is a very large
space and serves as a ``universal home'' for measures in theories in
which the Wilson loop observables are well-defined. In this paper, \(A\)
is considered as the projective limit of a projective family of compact
Hausdorff manifolds, labelled by graphs (which can be regarded as
``floating lattices'' in the physics terminology). Using this
characterization, differential geometry is developed through algebraic
methods. In particular, we are able to introduce the following notions
on \(A\): differential forms, exterior derivatives, volume forms, vector
fields and Lie brackets between them, divergence of a vector field with
respect to a volume form, Laplacians and associated heat kernels and
heat kernel measures. Thus, although \(A\) is very large, it is small
enough to be mathematically interesting and physically useful. A key
feature of this approach is that it does not require a background
metric. The geometrical framework is therefore well-suited for
diffeomorphism invariant theories such as quantum general relativity.
\end{quote}

\begin{enumerate}
\def\labelenumi{\arabic{enumi})}
\setcounter{enumi}{5}
\tightlist
\item
  A. P.  Balachandran, L. Chandar, Arshad Momen, ``Edge states in gravity and 
  black hole physics'',  available as
  \href{https://arxiv.org/abs/gr-qc/9412019}{\texttt{gr-qc/9412019}}.
\end{enumerate}

Ever since it started seeming that black holes have an entropy closely
related to (and often proportional to) the area of their event horizons,
many physicists have sought a better understanding of this entropy. In
many ways, the nicest sort of explanation would say that the entropy was
due to degrees of freedom living on the event horizon. A concrete
calculation along these lines was recently made by Steve Carlip (see
\protect\hyperlink{week41}{``Week 41''}) in the context of
2+1-dimensional gravity. The mechanism is mathematically very similar to
what happens in (a widely popular theory of) the fractional quantum Hall
effect! In both cases, Chern--Simons theory on a 3d manifold with
boundary gives rise to an interesting field theory on the boundary, or
``edge''. The above paper clarifies this, especially for those of us who
don't understand the fractional quantum Hall effect too well. Let me
quote:

\begin{quote}
Abstract: We show in the context of Einstein gravity that the removal of
a spatial region leads to the appearance of an infinite set of
observables and their associated edge states localized at its boundary.
Such a boundary occurs in certain approaches to the physics of black
holes like the one based on the membrane paradigm. The edge states can
contribute to black hole entropy in these models. A ``complementarity
principle'' is also shown to emerge whereby certain ``edge'' observables
are accessible only to certain observers. The physical significance of
edge observables and their states is discussed using their similarities
to the corresponding quantities in the quantum Hall effect. The coupling
of the edge states to the bulk gravitational field is demonstrated in
the context of (2+1) dimensional gravity.
\end{quote}

I can't resist adding that I have a personal stake in the notion that a
lot of interesting things about quantum gravity will only show up when
we consider it on manifolds with boundary, including the area-entropy
relations. The loop representation of quantum gravity has a lot to do
with knots and links, but on a manifold with boundary it has a lot to do
with \emph{tangles}, which can contain arcs that begin and end at the
boundary. I wrote a paper on this a while back:

\begin{enumerate}
\def\labelenumi{\arabic{enumi})}
\setcounter{enumi}{6}
\tightlist
\item
  John Baez, ``Quantum gravity and the algebra of tangles'', 
  \emph{Jour. Class. Quantum Grav.} \textbf{10} (1993), 673--694.
  Available as \href{https://arxiv.org/abs/hep-th/9205007}{\texttt{hep-th/9205007}}.
\end{enumerate}
\noindent
I'll be coming out with another in a while, co-authored with Javier
Muniain and Dardo Piriz. The importance of manifolds with boundary for
cutting-and-pasting constructions is also well-known in the theory of
``extended'' TQFTs (topological quantum field theories). These cutting
and pasting operations should allow one to describe extended TQFTs in n
dimensions purely algebraically using ``higher-dimensional algebra''. So
part of the plan here is to understand better the relation between
quantum gravity, TQFTs, and higher-dimensional algebra. Along these
lines, a very interesting new paper has come out:

\begin{enumerate}
\def\labelenumi{\arabic{enumi})}
\setcounter{enumi}{7}
\tightlist
\item
   Louis Crane and David Yetter,  ``On algebraic structures implicit in topological 
   quantum field theories'',  available as
  \href{https://arxiv.org/abs/hep-th/9412025}{\texttt{hep-th/9412025}}.
\end{enumerate}
\noindent
This makes more precise some of Louis Crane's ideas on
``categorification''. Nice TQFTs in 3 dimensions have a lot to do with
Hopf algebras (like quantum groups), or alternatively, their categories
of representations, which are certain braided monoidal categories. In
this paper it's shown that nice TQFTs in 4 dimensions have a lot to do
with Hopf categories, or alternatively, their categories of
representations, which are certain braided monoidal \(2\)-categories.

\begin{enumerate}
\def\labelenumi{\arabic{enumi})}
\setcounter{enumi}{8}
\tightlist
\item
  V. M. Kharlamov and V. G. Turaev, ``On the definition of 2-category of 
  2-knots'', \emph{Les Rencontres Physiciens-Math\'ematiciens de Strasbourg-RCP25}
  \textbf{45} (1993), 151-166.  Available as \href{http://www.numdam.org/item/RCP25_1993__45__151_0/}{\texttt{http://www.numdam.org/item/RCP25$\underline{\;}$1993$\underline{\;\;}$45$\underline{\;\;}$151$\underline{\;}$0/}}.
\end{enumerate}
\noindent
This preprint, which I obtained through my network of spies,
\emph{seems} to be implicitly claiming that the work of Fischer
describing 2d surfaces in \(\mathbb{R}^4\) via braided monoidal
\(2\)-categories (see \protect\hyperlink{week12}{``Week 12''}) is a bit
wrong, but they don't come out and say quite what if anything is really
wrong.

\begin{enumerate}
\def\labelenumi{\arabic{enumi})}
\setcounter{enumi}{9}
\tightlist
\item
   Greg
  Kuperberg, ``Non-involutory Hopf algebras and 3-manifold invariants'', available
  as \href{https://arxiv.org/abs/q-alg/9712047}{\texttt{q-alg/9712047}}.
\end{enumerate}
\noindent
I noted the existence of a draft of this paper, and related work, in
\protect\hyperlink{week38}{``Week 38''}. Let me quote:

\begin{quote}
Abstract: We present a definition of an invariant \(\#(M,H)\), defined
for every finite-dimensional Hopf algebra (or Hopf superalgebra or Hopf
object) \(H\) and for every closed, framed 3-manifold \(M\). When \(H\)
is a quantized universal envloping algebra, \(\#(M,H)\) is closely
related to well-known quantum link invariants such as the HOMFLY
polynomial, but it is not a topological quantum field theory.
\end{quote}

Okay, now for some miscellaneous interesting things\ldots{}

\begin{enumerate}
\def\labelenumi{\arabic{enumi})}
\setcounter{enumi}{10}
\tightlist
\item
  J. Lambek, ``If Hamilton had prevailed: quaternions in physics'', 
in \emph{Mathematical Conversations}, Springer, Berlin, pp.\ 259--274.
\end{enumerate}
\noindent
Lambek is mainly known for work in category theory, but he has a strong
side-interest in mathematical physics. This paper is, first of all, a
``nostalgic account of how certain key results in modern theoretical
physics (prior to World War II) can be expressed concisely in the
labguage of quaternions, thus suggesting how they might have been
discovered if Hamilton's views had prevailed.'' But there is a very
interesting new thing, too: a way in which the group
\(\mathrm{SU}(3) \times \mathrm{SU}(2) \times \mathrm{U}(1)\) shows up
naturally when considering Dirac's equation a la quaternions. This group
is the gauge group of the Standard Model! Lambek modestly says that
there does not appear to be any significance to this coincidence\ldots{}
but it \emph{would} be nice, wouldn't it?

\begin{enumerate}
\def\labelenumi{\arabic{enumi})}
\setcounter{enumi}{11}
\item
  Nina Byers, ``The life and times of Emmy Noether; contributions of E. 
  Noether to particle physics'', available as
  \href{https://arxiv.org/abs/hep-th/9411110}{\texttt{hep-th/9411110}}.

  Bert Schroer, ``Reminiscences about many pitfalls and some successes of QFT within
  the last three decades'', available as
  \href{https://arxiv.org/abs/hep-th/9410085}{\texttt{hep-th/9410085}}.

  Roman Jackiw, ``My encounters --- as a physicist --- with mathematics'', available as
  \href{https://arxiv.org/abs/hep-th/9410151}{\texttt{hep-th/9410151}}.
\end{enumerate}

These are some interesting historical/biographical pieces.

\begin{enumerate}
\def\labelenumi{\arabic{enumi})}
\setcounter{enumi}{12}
\tightlist
\item
   Karl Svozil, ``Speedup in quantum computation is associated with attenuation of
  processing probability'', available as
  \href{https://arxiv.org/abs/hep-th/9412046}{\texttt{hep-th/9412046}}.
\end{enumerate}
\noindent
The subject of quantum computation has become more lively recently. I
haven't had time to look at this paper, but quoting the abstract:

\begin{quote}
Quantum coherence allows the computation of an arbitrary number of
distinct computational paths in parallel. Based on quantum parallelism
it has been conjectured that exponential or even larger speedups of
computations are possible. Here it is shown that, although in principle
correct, any speedup is accompanied by an associated attenuation of
detection rates. Thus, on the average, no effective speedup is obtained
relative to classical (nondeterministic) devices.
\end{quote}

\hypertarget{week47}{%
\section{January 17, 1995}\label{week47}}

Hi. I know I'm supposed to be taking a break from ``This Week's Finds,''
but it's a habit that's hard to quit. I just want to briefly mention a
couple of new electronic venues for mathematics and physics. First,
there are a couple new preprint servers along the lines of
\texttt{hep-th}:

\begin{itemize}
\tightlist
\item
  Nuclear Experiment (\texttt{nucl-ex@xxx.lanl.gov}), started 12/94
\item
  Quantum Physics (\texttt{quant-ph@xxx.lanl.gov}), started 12/94
\end{itemize}
\noindent
To check these out, send email with subject header ``help'' and no
message body to the address listed, or use more slick modern method.

There are also a bunch of math preprint servers I have been neglecting
to mention. The newest one is devoted to ``quantum algebra'' (quantum
groups and the like) and knot theory! Here are some:

\begin{itemize}
\tightlist
\item
  Algebraic Geometry (\texttt{alg-geom@eprints.math.duke.edu}), started
  2/92
\item
  Functional Analysis (\texttt{funct-an@xxx.lanl.gov}), started 4/92
\item
  Differential Geometry and Global Analysis (\texttt{dg-ga@msri.org}),
  started 6/94
\item
  Automorphic Forms (\texttt{auto-fms@msri.org}), started 6/94
\item
  Complex Dynamics and Hyperbolic Geometry (\texttt{cd-hg@msri.org}),
  started 6/94
\item
  Quantum Algebra (Including Knot Theory)
  (\texttt{q-alg@eprints.math.duke.edu}), started 12/94
\end{itemize}

Also, there is a new electronic journal on the theory and applications
of category theory. This is a subject dear to my heart since much of the
most interesting work on topological quantum field theories uses
category theory in an interesting way, so I hope mathematicians and
mathematical physicists interested in category theory turn to this
journal:

\begin{verbatim}
               ANNOUNCEMENT AND CALL FOR PAPERS
 
       This is to announce a new refereed electronic journal
 
 
            THEORY AND APPLICATIONS OF CATEGORIES
 
 
The Editors, who are listed below, invite submission of articles 
for publication in the first volume (1995) of the journal.   
The Editorial Policy of the journal, basic information about 
subscribing and some information for authors follows. 
More details are available from the journal's WWW server at
 
URL http://www.tac.mta.ca/tac/
 
or by anonymous ftp from
 
ftp.tac.mta.ca
 
in the directory pub/tac.
 
 
EDITORIAL POLICY
 
The journal Theory and Applications of Categories will disseminate articles 
that significantly advance the study of categorical algebra or methods, or 
that make significant new contributions to mathematical science using 
categorical methods. The scope of the journal includes: all areas of 
pure category theory, including higher dimensional categories; applications 
of category theory to algebra, geometry and topology and other areas of 
mathematics; applications of category theory to computer science, physics 
and other mathematical sciences; contributions to scientific knowledge 
that make use of categorical methods.
 
Articles appearing in the journal have been carefully and critically 
refereed under the responsibility of members of the Editorial Board.
Only papers judged to be both significant and excellent are accepted for 
publication. 
 
The method of distribution of the journal is via the Internet tools 
WWW/gopher/ftp. The journal is archived electronically and in printed 
paper format. 
 
 
 
SUBSCRIPTION INFORMATION
 
Individual subscribers receive (by e-mail) abstracts of articles 
as they are published. Full text of published articles is available 
in .dvi and Postscript format. Details will be e-mailed to new 
subscribers and are available by WWW/gopher/ftp. To subscribe, 
send e-mail to tac@mta.ca including a full name and postal address. 
For institutional subscription, send enquiries to the Managing Editor, 
Robert Rosebrugh, <rrosebrugh@mta.ca>
 
 
 
 
INFORMATION FOR AUTHORS
 
The typesetting language of the journal is TeX, and LaTeX is the 
preferred flavour. TeX source of articles for publication should 
be submitted by e-mail directly to an appropriate Editor. 
Please obtain detailed information on submission format and 
style files by WWW or anonymous ftp from the sites listed above.
You may also write to tac@mta.ca to receive details by e-mail.
 
 
EDITORIAL BOARD
 
John Baez, University of California, Riverside
baez@math.ucr.edu
 
Michael Barr, McGill University
barr@triples.math.mcgill.ca
 
Lawrence Breen, Universite de Paris 13
breen@math.univ-paris.fr
 
Ronald Brown, University of North Wales 
r.brown@bangor.ac.ik
 
Jean-Luc Brylinski, Pennsylvania State University
jlb@math.psu.edu
 
Aurelio Carboni, University of Genoa
carboni@vmimat.mat.unimi.it
 
G. Max Kelly, University of Sydney
kelly_m@maths.su.oz.au
 
Anders Kock, University of Aarhus
kock@mi.aau.dk
 
F. William Lawvere, State University of New York at Buffalo
mthfwl@ubvms.cc.buffalo.edu
 
Jean-Louis Loday, Universite de Strasbourg
loday@math.u-strasbg.fr
 
Ieke Moerdijk, University of Utrecht
moerdijk@math.ruu.nl
 
Susan Niefield, Union College
niefiels@gar.union.edu
 
Robert Pare, Dalhousie University
pare@cs.dal.ca
 
Andrew Pitts, University of Cambridge
ap@cl.cam.ac.uk
 
Robert Rosebrugh, Mount Allison University
rrosebrugh@mta.ca
 
Jiri Rosicky, Masaryk University
rosicky@math.muni.cz
 
James Stasheff, University of North Carolina
jds@charlie.math.unc.edu
 
Ross Street, Macquarie University 
street@macadam.mpce.mq.edu.au
 
Walter Tholen, York University
tholen@vm1.yorku.ca
 
R. W. Thomason, Universite de Paris 7
thomason@mathp7.jussieu.fr
 
Myles Tierney, Rutgers University
tierney@math.rutgers.edu
 
Robert F. C. Walters, University of Sydney
walters_b@maths.su.oz.au
 
R. J. Wood, Dalhousie University
rjwood@cs.da.ca
\end{verbatim}

\hypertarget{week48}{%
\section{February 26, 1995}\label{week48}}

There are a few things I've bumped into that I feel I should let folks
know about, so here's a special issue from Munich, where the Weissbier
is very good. (And not at all white, but that's another subject.)

One of the most exciting aspects of mathematics over the last few years
--- in my utterly biased opinion --- has been how topological quantum
field theories have revealed the existence of deep connections between
\(3\)-dimensional topology, complex analysis, and algebra, particularly
the algebra of quantum groups.

The most interesting topological quantum field theory in this game is
Chern--Simons theory. This is a field theory in 3 dimensions, and the reason
it's called ``topological'' is that you don't need any metric or other
geometrical structure on your 3d spacetime manifold for this theory to
make sense. Thus it admits \emph{all} coordinate transformations (or
more precisely, all diffeomorphisms) as symmetries. In particular, this
means that the quantity folks like to compute whenever they see a
quantum field theory --- the partition function, which you get by doing
a path integral a la Feynman --- is an invariant of \(3\)-dimensional
manifold you happen to have taken as ``spacetime''.

Now computing path integrals is often a very dubious and tricky
business. They are integrals over the space of all possible histories of
the classical fields corresponding to your quantum field theory. This is
a big fat infinite-dimensional space, of the sort to which ordinary
integration theory doesn't really apply. If you aren't very careful,
path integrals often give infinite answers. So one very nice thing is
that, suitably interpreted, the path integrals in Chern--Simons theory
actually make rigorous sense!

A key advance here was Atiyah's axioms for topological quantum field
theories (or TQFTs). These axioms formalize exactly what properties
you'd hope path integrals would have in the case of a
diffeomorphism-invariant theory. If, by no matter what devilish tricks,
one can get a theory that satisfies these axioms, it's in many ways just
as good if one had figured out how to make sense of the path integrals
by honest labor, because all the manipulations one would normally want
to do are then permitted. The marvelous thing about Chern--Simons theory
is that one can show the TQFT axioms hold starting from some beautiful
algebraic structures called quantum groups. Corresponding to every
``semisimple Lie group'' --- examples being the groups
\(\mathrm{SU}(n)\) of unitary complex \(n \times n\) matrices with determinant
\(1\) --- there is a quantum group, which is not really a group, but a
so-called ``quasitriangular Hopf algebra.'' Amusingly these quantum
groups really depend on Planck's constant \(\hbar\), and reduce to the
ordinary groups in the ``classical limit'' \(\hbar \to 0\)!

Now, where these quantum groups come from has always been a bit of a
puzzle. They can be rigorously shown to exist, that's for sure. There
are also many good algebraic reasons why they ``should'' exist. But it
is still a bit remarkable that they have the exactly the properties
needed to get Chern--Simons theory to be a TQFT. So people have tried in
many ways to turn the tables on them, and get \emph{them} from the
\emph{path integrals}. Lots of these approaches have succeeded, but most
of them involve a few subtleties here and there, so mathematicians, who
only feel happy when everything is \emph{blindingly obvious} (to them, that
is, not you), have continued to seek the most beautiful, elegant way of
getting at them.

\begin{enumerate}
\def\labelenumi{\arabic{enumi})}
\tightlist
\item
  Daniel Freed, ``Quantum groups from path integrals'', available as
  \href{https://arxiv.org/abs/q-alg/9501025}{\texttt{q-alg/9501025}}.
\end{enumerate}

This is a nice expository treatment of the work of Free and Quinn on
topological quantum field theories, particularly Chern--Simons theory
with finite gauge group. In this case, the path integral reduces to a
finite sum and one really can get the quantum group from the path
integral very beautifully. But there are some differences in this case
of finite gauge group. For example, the resulting ``finite quantum
group'' does \emph{not} depend on \(\hbar\); it's just the ``quantum
double'' of the group algebra of the group. For more on what this has to
do with marvelous algebraic things like \(n\)-categories, the reader
should check out the paper by Freed cited in
\protect\hyperlink{week12}{``Week 12''}, which has subsequently been
published:

\begin{enumerate}
\def\labelenumi{\arabic{enumi})}
\setcounter{enumi}{1}
\tightlist
\item
  Daniel Freed, ``Higher algebraic structures and quantization'',
  \emph{Commun.\ Math.\ Phys.} \textbf{159} (1994), 343--398.
\end{enumerate}
\noindent
and also

\begin{enumerate}
\def\labelenumi{\arabic{enumi})}
\setcounter{enumi}{2}
\tightlist
\item
  Daniel Freed and
  Frank Quinn, ``Chern--Simons theory with finite gauge group'', 
  \emph{Commun.\ Math.\ Phys.} \textbf{156} (1993), 435--472.
\end{enumerate}

Now in addition to the path-integral approach to quantum field theory
there is another, the so-called Hamiltonian approach, which is very much
like the approach people usually learn in a first course on quantum
mechanics: if you know the wavefunction of your system at \(t = 0\), the
Hamiltonian tells you how it evolves in time from then on. Now this has
special subtleties in diffeomorphism-invariant theories. When there is
no unique best coordinate system, there's no unique best notion of
``time evolution''. This leads to the so-called ``problem of time'',
very important in quantum gravity, but rather easier to deal with in toy
models like Chern--Simons theory.

Now if we take a \(3\)-dimensional spacetime and look at the \(t = 0\)
slice, we will with some luck get a \(2\)-dimensional manifold, such as
a sphere, torus, or more general \(n\)-holed doughnut. This is where the
complex analysis comes in, because the complex plane is
\(2\)-dimensional, and we can cover these surfaces with coordinate
systems that look locally like the complex plane, making them into
``Riemann surfaces'' upon which we can merrily proceed with complex
analysis. In particular, the phase space of classical Chern--Simons
theory is something of which complex analysts have long been enamored,
namely the ``moduli space of flat bundles'' over our Riemann surface.
(Let me reassure physicists that this ``flat'' business is merely a
weird way of saying that in Chern--Simons theory the analog of the
magnetic field vanishes.)

Now starting from the description of the classical phase space for
Chern--Simons theory one should be able to get ahold of the quantum
theory by some ``quantization'' business just as one does in elementary
quantum mechanics, where the ``classical phase space'' is the space of
\(p\)'s and \(q\)'s, and to quantize one merely decrees that these no
longer commute: 
\[pq-qp = -i \hbar.\] 
So one should be able to get ahold
of quantum groups this way too: starting with the ``moduli space of flat
bundles'' and ``quantizing'' it. I had long wondered why nobody did this in the
way which seemed most tempting and plausible to me. To lapse into jargon
for a bit here, since the quantum group is obtained from the group by
deformation quantization, where the Poisson structure of the group
itself is described by some ``classical \(r\)-matrices'', and since
Chern--Simons theory is in some sense obtained by quantizing the moduli
space, where the Poisson structure was explicitly described by Goldman,
it seemed natural to me that somehow the classical r-matrices should be
\emph{derivable} from Goldman's formulas. But after a few wimpy tries at
figuring it out, with not much success, I gave up. Luckily it turns out
someone else succeeded nicely, as I found out in a talk by Alekseev here
at the Mathematisches Institut of the Universit\"at M\"unchen:

\begin{enumerate}
\def\labelenumi{\arabic{enumi})}
\setcounter{enumi}{3}
\tightlist
\item
   V.\ V.\ Fock and A.\ A.\ Rosly, ``Poisson structures on moduli of flat connections 
   on Riemann surfaces and $r$-matrix'', 
   \href{https://arxiv.org/abs/math/9802054}{\texttt{math/9802054}}.
\end{enumerate}
\noindent
They figured out a beautiful formula relating the classical \(r\)-matrices
and the Poisson structure on moduli space. Using this, Alekseev, Grosse,
and Schomerus were able to get at quantum groups quite directly from
deformation quantization of moduli space, though there are a few
important points left to nail down:

\begin{enumerate}
\def\labelenumi{\arabic{enumi})}
\setcounter{enumi}{4}
\tightlist
\item
  A.\ Yu. Alekseev, H.\ Grosse, and V.\ Schomerus, 
  ``Combinatorial quantization of the Hamiltonian Chern--Simons theory, I \& II''.
  Available as 
  \href{https://arxiv.org/abs/hep-th/9403066}{\texttt{hep-th/9403066}}
  and
  \href{https://arxiv.org/abs/hep-th/9408097}{\texttt{hep-th/9408097}}.
\end{enumerate}
\noindent
In fact Schomerus told me about this while I was in Cambridge Massachusetts over
Christmas, but somehow I didn't pick up on the coolest thing about it,
that it gets the quantum groups quite naturally from Chern--Simons
theory, via the relation between the Poisson brackets on moduli space
and the classical \(r\)-matrices.

Let's see. This is getting technical so I'll give in and get technical.
One other paper that I found out about while here studies rather similar
issues, but from the viewpoint of geometric quantization rather than
deformation quantization. Axelrod, Della Pietra and Witten did some very
fundamental work on geometric quantization of Chern--Simons theory:

\begin{enumerate}
\def\labelenumi{\arabic{enumi})}
\setcounter{enumi}{5}
\tightlist
\item
   Scott Axelrod,
  Steve Della Pietra and Edward Witten,  ``Geometric quantization of 
  Chern--Simons gauge theory'', \emph{Jour. Diff. Geom.} \textbf{33}
  (1991), 787--902.  Also available as \href{https://projecteuclid.org/euclid.jdg/1214446565}{\texttt{https://projecteuclid.org/euclid.jdg/1214446565}}.
\end{enumerate}
\noindent
but this only treated the case of Riemann surfaces, not the Riemann
surfaces with punctures that you need to think about when sticking
\emph{knots} in your 3-manifold. Martin Schottenloher, one of the
professors here, gave me a nice review of the work of one of his
students on the case with punctures:

\begin{enumerate}
\def\labelenumi{\arabic{enumi})}
\setcounter{enumi}{6}
\tightlist
\item
  Martin Schottenloher, ``Metaplectic quantization of the moduli space of flat and parabolic
  bundles (after Peter Scheinhost)'', in \emph{Public. I. R. M. A.
  Strasbourg}, \textbf{45} (1993), 43--70.  Also available as  \hfill \break
  \href{http://www.numdam.org/item/RCP25_1993__45__43_0/}{\texttt{http://www.numdam.org/item/RCP25$\underline{\;}$1993$\underline{\;\;}$45$\underline{\;\;}$43$\underline{\;}$0/}}.
\end{enumerate}
\noindent
This uses a lot of complex geometry that's beyond my ken, but one very
exciting remark penetrated my thick skull, namely that you really need
to take into account the so-called ``metaplectic correction'', as one
usually does in geometric quantization, and that in the case of no
punctures, this has the sole effect of accomplishing the mysterious
``level shift'' \((k \to k + N)\) that pervades \(\mathrm{SU}(N)\)
Chern--Simons theory. (Of course, I bet this must be what's going on for
other gauge groups too.) Also, when there are punctures, you apparently
really need the metaplectic correction to get the right answers from
geometric quantization; the ad hoc level shift ain't enough.

Well, I have some other things, but this issue seems more or less
devoted to Chern--Simons theory so far, and only Chern--Simons freaks are
likely to have read this far, so I'll put off the rest for later.

\hypertarget{week49}{%
\section{February 27, 1995}\label{week49}}

For the last year or so I've been really getting interested in
\(n\)-categories as a possible tool for unifying a lot of strands in
mathematics and physics. What's an \(n\)-category? Well, in a sense
that's the big question! Roughly speaking, it's a structure where there
are a bunch of ``objects'', and for any pair of objects \(x,y\) a bunch
of ``morphisms'' from \(x\) to \(y\), written \(f\colon x \to y\), and
for any pair of morphisms \(f, g\colon x \to y\) a bunch of
``2-morphisms'' from \(f\) to \(g\), written
\(F\colon f \Rightarrow g\), and for any pair of \(2\)-morphisms
\(F, G\colon f \Rightarrow g\) a bunch of ``3-morphisms'' from \(F\) to
\(G\), \ldots{} and so on, up to \(n\)-morphisms. Ordinary categories,
or \(1\)-categories, have been studied since the 1940s or so, when they
were invented by Eilenberg and Mac Lane. Anyone wanting to get going on
those could try:

\begin{enumerate}
\def\labelenumi{\arabic{enumi})}
\tightlist
\item
  Saunders Mac Lane, \emph{Categories for the Working Mathematician},
  Springer, Berlin, 1988.
\end{enumerate}

Roughly, categories give a general framework for dealing with situations
where you have ``things'' and ``ways to go between things'', like sets
and functions, vector spaces and linear maps, points in a space and
paths between points, etc. That's a pretty broad territory!
\(n\)-categories show up when you also have ``ways to go between ways'',
``ways to go between ways to go between ways'', etc. That may seem a
little weird at first. But in fact they show up in a lot of places if
you look for them. Perhaps the most obvious place is topology. If you
think of a point in a space as an object, and a path between two points
as a morphism: \[x\xrightarrow{\quad f\quad}g\] you are easily tempted
to think of a ``path of paths'' as a \(2\)-morphism. Here a ``path of
paths'' is just a continuous 1-parameter family of paths from x to y,
which you can think of as tracing out a \(2\)-dimensional surface, as
follows: \[\includegraphics[scale=0.3]{Fnatftog.pdf}\] And one
can keep on going and look at ``paths of paths of paths'', etc. In fact,
people in homotopy theory do this all the time.

There is another example, equally primordial, which is a bit more
``inbred'' in flavor. In other words, if you already know and love
\(n\)-categories, there is a wonderful example of an \((n+1)\)-category
which you should know and love too! Now this isn't so bad, actually,
because a 0-category is basically just a \emph{set}, namely the set of
objects. Since everyone knows and loves sets, everyone can start here!
Okay, there is a wonderful \(1\)-category called Set, the category of
all sets. This has sets as its objects and functions between sets as its
morphisms. So now that you know an example of a \(1\)-category, you know
and love \(1\)-categories, right? Well, it turns out there is this
wonderful thing called Cat, the \(2\)-category of all \(1\)-categories.
(Usually people restrict to ``small'' \(1\)-categories, which have a
mere \emph{set} of objects, so that set theorists don't start freaking
out at a certain point.) To understand why Cat is a \(2\)-category is a
bit of work, but as objects it has categories, as morphisms it has the
usual sort of morphisms between categories, so-called ``functors'', and
as \(2\)-morphisms it has the usual sort of things that go between
functors, so-called ``natural transformations''. These are the bread and
butter of category theory; just take my word for it if you haven't
studied them yet! Okay, so Cat is a \(2\)-category, so now you know and
love \(2\)-categories, right? (Well, I haven't even told you the
definitions, but just nod your head.) Guess what: there is this
wonderful thing called 2Cat, the \(3\)-category of all \(2\)-categories!
And so on.

So in short, the detailed theory of \(n\)-categories at each level
automatically leads one to get interested in \((n+1)\)-categories. Now
for the bad news: so far, people have only figured out the right
definition of \(n\)-category for \(n = 0, 1, 2, 3\). By the ``right''
definition I mean the ultimate, most general definition, which should be
the most useful in many ways. So far people only know about ``strict''
\(n\)-categories for all \(n\), which one can think of as a special case
of the ultimate ones; the ultimate \(1\)-categories are just categories,
the ultimate \(2\)-categories are often called bicategories (see the
reference to Benabou in \protect\hyperlink{week35}{``Week 35''}), and
ultimate \(3\)-categories are usually called tricategories (see the
reference to the paper by Gordon, Power and Street) in
\protect\hyperlink{week29}{``Week 29''}. Tricategories were just defined
last year! They have a whole lot to do with knots, Chern--Simons theory,
and other \(3\)-dimensional phenomena, as one might expect. If we could
understand the ultimate \(4\)-categories --- tetracategories? --- it
would probably help us with some of the riddles of topology and physics
in 4 dimensions. (Indeed, what little we \emph{do} understand is already
helping a bit.)

So anyway, I have been trying to learn about these things, and had the
good luck to meet James Dolan via the net, who has helped me immensely,
since he eats, lives and breathes category theory, and he is now at
Riverside hard at work figuring out the ultimate definition of
\(n\)-categories for all \(n\). (Although when he reads this, he will
not be hard at work; he will be goofing off, reading the news.)

He and I have have written one paper so far espousing our philosophy
concerning \(n\)-categories, topology, and physics:

\begin{enumerate}
\def\labelenumi{\arabic{enumi})}
\setcounter{enumi}{1}
\tightlist
\item
   John Baez and James Dolan, 
  ``Higher-dimensional algebra and topological quantum field theory'', 
  \emph{Jour. Math. Phys.} \textbf{36} (1995), 6073--6105.
  Also available as
  \href{http://arxiv.org/abs/q-alg/9503002}{\texttt{q-alg/9503002}}.
\end{enumerate}

One of the main themes of this paper is what I sometimes jokingly call
the ``periodic table''. Say you have an \((n+k)\)-category with only one
object, one morphism, one \(2\)-morphism, \ldots{} and only one
\((k-1)\)-morphism. Then all the interest lies in the \(k\)-morphisms,
the \((k+1)\)-morphisms, and so on up to the \((k+n)\)-morphisms. So
there are \(n\) interesting levels of morphism, and we can actually
think of our \((n+k)\)-category as an \(n\)-category of a special sort.
Let's call this kind a ``\(k\)-tuply monoidal \(n\)-category''. Now we
can make a chart of these:

\newpage
\begin{longtable}[]{@{}llll@{}}
\caption*{\(k\)-tuply monoidal \(n\)-categories}\tabularnewline
\toprule
\begin{minipage}[b]{0.26\columnwidth}\raggedright
\strut
\end{minipage} & \begin{minipage}[b]{0.21\columnwidth}\raggedright
\(n=0\)\strut
\end{minipage} & \begin{minipage}[b]{0.21\columnwidth}\raggedright
\(n=1\)\strut
\end{minipage} & \begin{minipage}[b]{0.21\columnwidth}\raggedright
\(n=2\)\strut
\end{minipage}\tabularnewline
\midrule
\endfirsthead
\toprule
\begin{minipage}[b]{0.26\columnwidth}\raggedright
\strut
\end{minipage} & \begin{minipage}[b]{0.21\columnwidth}\raggedright
\(n=0\)\strut
\end{minipage} & \begin{minipage}[b]{0.21\columnwidth}\raggedright
\(n=1\)\strut
\end{minipage} & \begin{minipage}[b]{0.21\columnwidth}\raggedright
\(n=2\)\strut
\end{minipage}\tabularnewline
\midrule
\endhead
\begin{minipage}[t]{0.26\columnwidth}\raggedright
\(k=0\)\strut
\end{minipage} & \begin{minipage}[t]{0.21\columnwidth}\raggedright
sets\strut
\end{minipage} & \begin{minipage}[t]{0.21\columnwidth}\raggedright
categories\strut
\end{minipage} & \begin{minipage}[t]{0.21\columnwidth}\raggedright
\(2\)-categories\strut
\end{minipage}\tabularnewline
\begin{minipage}[t]{0.26\columnwidth}\raggedright
\strut
\end{minipage} & \begin{minipage}[t]{0.21\columnwidth}\raggedright
\strut
\end{minipage} & \begin{minipage}[t]{0.21\columnwidth}\raggedright
\strut
\end{minipage} & \begin{minipage}[t]{0.21\columnwidth}\raggedright
\strut
\end{minipage}\tabularnewline
\begin{minipage}[t]{0.26\columnwidth}\raggedright
\(k=1\)\strut
\end{minipage} & \begin{minipage}[t]{0.21\columnwidth}\raggedright
monoids\strut
\end{minipage} & \begin{minipage}[t]{0.21\columnwidth}\raggedright
monoidal categories\strut
\end{minipage} & \begin{minipage}[t]{0.21\columnwidth}\raggedright
monoidal \(2\)-categories\strut
\end{minipage}\tabularnewline
\begin{minipage}[t]{0.26\columnwidth}\raggedright
\strut
\end{minipage} & \begin{minipage}[t]{0.21\columnwidth}\raggedright
\strut
\end{minipage} & \begin{minipage}[t]{0.21\columnwidth}\raggedright
\strut
\end{minipage} & \begin{minipage}[t]{0.21\columnwidth}\raggedright
\strut
\end{minipage}\tabularnewline
\begin{minipage}[t]{0.26\columnwidth}\raggedright
\(k=2\)\strut
\end{minipage} & \begin{minipage}[t]{0.21\columnwidth}\raggedright
commutative monoids\strut
\end{minipage} & \begin{minipage}[t]{0.21\columnwidth}\raggedright
braided monoidal categories\strut
\end{minipage} & \begin{minipage}[t]{0.21\columnwidth}\raggedright
braided monoidal \(2\)-categories\strut
\end{minipage}\tabularnewline
\begin{minipage}[t]{0.26\columnwidth}\raggedright
\strut
\end{minipage} & \begin{minipage}[t]{0.21\columnwidth}\raggedright
\strut
\end{minipage} & \begin{minipage}[t]{0.21\columnwidth}\raggedright
\strut
\end{minipage} & \begin{minipage}[t]{0.21\columnwidth}\raggedright
\strut
\end{minipage}\tabularnewline
\begin{minipage}[t]{0.26\columnwidth}\raggedright
\(k=3\)\strut
\end{minipage} & \begin{minipage}[t]{0.21\columnwidth}\raggedright
" "\strut
\end{minipage} & \begin{minipage}[t]{0.21\columnwidth}\raggedright
symmetric monoidal categories\strut
\end{minipage} & \begin{minipage}[t]{0.21\columnwidth}\raggedright
weakly involutory monoidal \(2\)-categories\strut
\end{minipage}\tabularnewline
\begin{minipage}[t]{0.26\columnwidth}\raggedright
\strut
\end{minipage} & \begin{minipage}[t]{0.21\columnwidth}\raggedright
\strut
\end{minipage} & \begin{minipage}[t]{0.21\columnwidth}\raggedright
\strut
\end{minipage} & \begin{minipage}[t]{0.21\columnwidth}\raggedright
\strut
\end{minipage}\tabularnewline
\begin{minipage}[t]{0.26\columnwidth}\raggedright
\(k=4\)\strut
\end{minipage} & \begin{minipage}[t]{0.21\columnwidth}\raggedright
" "\strut
\end{minipage} & \begin{minipage}[t]{0.21\columnwidth}\raggedright
" "\strut
\end{minipage} & \begin{minipage}[t]{0.21\columnwidth}\raggedright
strongly involutory monoidal \(2\)-categories\strut
\end{minipage}\tabularnewline
\begin{minipage}[t]{0.26\columnwidth}\raggedright
\strut
\end{minipage} & \begin{minipage}[t]{0.21\columnwidth}\raggedright
\strut
\end{minipage} & \begin{minipage}[t]{0.21\columnwidth}\raggedright
\strut
\end{minipage} & \begin{minipage}[t]{0.21\columnwidth}\raggedright
\strut
\end{minipage}\tabularnewline
\begin{minipage}[t]{0.26\columnwidth}\raggedright
\(k=5\)\strut
\end{minipage} & \begin{minipage}[t]{0.21\columnwidth}\raggedright
" "\strut
\end{minipage} & \begin{minipage}[t]{0.21\columnwidth}\raggedright
" "\strut
\end{minipage} & \begin{minipage}[t]{0.21\columnwidth}\raggedright
" "\strut
\end{minipage}\tabularnewline
\bottomrule
\end{longtable}

First, I should emphasize that some parts of the chart as I've drawn it
here are a bit conjectural; since we don't know what the most general
7-categories are like, for example, we don't really know for sure what
5-tuply monoidal \(2\)-categories are like. The exact status of all the
entries on the table is made more clear in the paper. For now, let me
just say, first, that these various flavors of \(n\)-categories turn out
to be of great interest in topology --- some have already been used a
lot in topological quantum field theory and knot theory, other less, so
far, but they all seem to have lot to do with generalizations of knot
theory to different dimensions. Second, it seems that the \(n\)th column
``stabilizes'' by the time you get down to the \((n+2)\)nd row. This
very interesting pattern turns out also to have a lot to do with knots
and their generalizations, and also to a subject called stable homotopy
theory.

Now it also appears that there is a nice recipe for hopping down the
columns. (Again, we only understand this \emph{perfectly} in certain
cases, but the pattern seems pretty clear.) In other words, there's a
nice recipe to get a (k+1)-tuply monoidal \(n\)-category from a
\(k\)-tuply monoidal one. It goes like this. Hang on to your seat. You
start with a \(k\)-tuply monoidal \(n\)-category \(\mathcal{C}\). It's a
special sort of \((n+k)\)-category, so its an object in
\((n+k)\mathsf{Cat}\). But \((n+k)\mathsf{Cat}\), remember, is an
\((n+k+1)\)-category. Now look at the largest sub-\((n+k+1)\)-category
of \((n+k)\mathsf{Cat}\) which has \(\mathcal{C}\) as its only object,
\(1_{\mathcal{C}}\) (the identity of \(\mathcal{C}\)) as its only
morphism, \(1_{1_{\mathcal{C}}}\) as its only \(2\)-morphism,
\(1_{1_{1_{\mathcal{C}}}}\) as its only \(3\)-morphism, and so on, up to
\(1_{1_{1_{\scriptsize\ldots}}}\) as its only \(k\)-morphism. Let's call
this \(\mathcal{C}'\). If one keeps count, this should be a
\((k+1)\)-tuply monoidal \(n\)-category. That's how it goes.

Now say we do this to an example. Say we do it to the category
\(\mathcal{C}\) of all representations of a finite group \(G\). This is
in fact a monoidal category, so the result \(\mathcal{C}'\) is a braided
monoidal category. It is, in fact, just the category of representations
of the ``quantum double'' of \(G\), which is an example of what one
might call a ``finite quantum group''. These play a big role in the
study of Chern--Simons theory with finite gauge group (see the papers by
Freed and Quinn in \protect\hyperlink{week48}{``Week 48''}). One can
also get the other quantum groups with the aid of this ``quantum
double'' trick. A good description of this case appears in:

\begin{enumerate}
\def\labelenumi{\arabic{enumi})}
\setcounter{enumi}{2}
\tightlist
\item
  Christian
  Kassel and Vladimir Turaev,  ``Double construction for monoidal categories'',
  \textsl{Acta Mathematica} \textbf{175} (1995), 1--48.  Also available as
  \href{https://projecteuclid.org/euclid.acta/1485890865}{\texttt{https://projecteuclid.org/euclid.acta/1485890865}}.
\end{enumerate}

So this is rather remarkable: starting from a finite group, and all this
\(n\)-categorical abstract nonsense, out pops precisely the raw
ingredients for a perfectly respectable \(3\)-dimensional topological
quantum field theory! Understanding \emph{why} this kind of thing works
is part of the aim of Dolan's and my paper, though there are some
important pieces of the puzzle that we don't get around to mentioning
there.

Right now I'm busily working out the details of how to get braided
monoidal \(2\)-categories from monoidal \(2\)-categories by the same
trick, with the aid of Martin Neuchl and Frank Halanke here. These
should have a lot to do with \(4\)-dimensional topological quantum field
theories (see e.g. the paper by Crane and Yetter cited in
\protect\hyperlink{week46}{``Week 46''}). And here I can't resist
mentioning a very nice paper by Neuchl and Schauenburg,

\begin{enumerate}
\def\labelenumi{\arabic{enumi})}
\setcounter{enumi}{3}
\tightlist
\item
  Martin Neuchl and Peter Schauenburg, ``Reconstruction in braided categories 
  and a notion of commutative bialgebra'',  \emph{Jour.\ Pure and Applied 
  Algebra} \textbf{124} (1998), 241--259.
\end{enumerate}

Let me conclude by describing this. I always let myself get a bit more
technical at the end of each issue, so I'll do that now. The
relationship between Hopf algebras and monoidal categories is given by
``Tannaka--Krein reconstruction theorems'', which give conditions under
which a monoidal category is equivalent to the category of
representations of a Hopf algebra, and actually constructs the Hopf
algebra for you. In physics people use related but fancier
``Doplicher--Haag--Roberts'' theorems to reconstruct the gauge group of a
quantum field theory. This paper starts with the beautiful Tannaka--Krein
theorem in

\begin{enumerate}
\def\labelenumi{\arabic{enumi})}
\setcounter{enumi}{4}
\tightlist
\item
  Peter Schauenburg, ``Tannaka duality for arbitrary Hopf algebras'', 
  \emph{Algebra-Berichte} \textbf{66} (1992).
\end{enumerate}

Leaving out a bunch of technical conditions that make the theorem
actually \emph{true}, it says roughly that when you have a braided monoidal
category \(\mathcal{B}\), a category \(\mathcal{C}\), and a functor
\(f\colon \mathcal{C} \to \mathcal{B}\), there is a coalgebra object
\(a\) in \(\mathcal{B}\), the universal one for which \(f\) factors
through the forgetful functor from \(a\mathsf{Comod}\) (the category of
\(a\)-comodule objects in \(\mathcal{B}\)) to \(\mathcal{B}\). The point
is that the ordinary Tannaka--Krein theorem is a special case of this one
where \(\mathcal{B}\) is the category of vector spaces. The point of the
new paper is as follows. Suppose \(\mathcal{C}\) is actually braided
monoidal and \(f\) preserves the braiding and monoidal structure. Then
we expect \(a\) to actually be something like commutative bialgebra
object in \(\mathcal{B}\). The paper makes this precise. There are
actually some sneaky issues involved in doing so. In particular, the
``quantum double'' trick for categories makes an appearance here. I
guess I'll leave it at that!

\hypertarget{week50}{%
\section{March 12, 1995}\label{week50}}

Your roving mathematical physics reporter is now in Milan, where (when
not busy eating various kinds of cheese) he is discussing \(BF\) theory
with his hosts, Paolo Cotta-Ramusino and Maurizio Martellini. This is
rather a long way to go to stumble on the October 1994 issue of the
\emph{Journal of Mathematical Physics}, but such is life. This issue is a
special issue on ``Topology and Physics'', a nexus dear to my heart, so
let me say a bit about some of the papers in it.

\begin{enumerate}
\def\labelenumi{\arabic{enumi})}
\tightlist
\item
  Edward Witten, ``Supersymmetric Yang--Mills theory on a four-manifold'',  
  \emph{Jour. Math. Phys.} \textbf{35} (1994), 5101--5135.
  Available as \href{https://arxiv.org/abs/hep-th/9403195}{\texttt{hep-th/9403195}}.
\end{enumerate}

This paper concerns the relation between supersymmetric Yang--Mills
theory and Donaldson theory, discovered by Seiberg and Witten, which not
too long ago hit the front page of various newspapers. (See
\protect\hyperlink{week44}{``Week 44''} and
\protect\hyperlink{week45}{``Week 45''} for my own yellow journalism on
the subject.) I don't have anything new to say about this stuff, of
which I am quite ignorant. If you are an expert on \(N = 2\)
supersymmetric Yang--Mills theory, hyper-K\"ahler manifolds, cosmic
strings and the renormalization group, the paper should be a piece of
cake. Seriously, he does seem to be making a serious effort to
communicate the ideas in simple terms to us mere mortals, so it's worth
looking at.

\begin{enumerate}
\def\labelenumi{\arabic{enumi})}
\setcounter{enumi}{1}
\tightlist
\item
  Louis Crane and Igor Frenkel, ``Four-dimensional topological quantum
  field theory, Hopf categories, and the canonical bases'', \emph{Jour.
  Math. Phys.} \textbf{35} (1994), 5136--5154. Also available as
  \href{https://arxiv.org/abs/hep-th/9405183}{\texttt{hep-th/9405183}}.
\end{enumerate}

I discussed this paper a wee bit in \protect\hyperlink{week38}{``Week
38''}. Now you can actually see the pictures. As we begin to understand
\(n\)-categories (see \protect\hyperlink{week49}{``Week 49''}) our
concept of symmetry gets deeper and deeper. This isn't surprising. When
all we knew about was 0-categories --- that is, sets! --- our concept of
symmetry revolved around the notion of a ``group''. This is a set \(G\)
where you can multiply elements in an associative way, with an identity
element 1 such that \[1g = g1 = g\] for all elements \(g\) of \(G\), and
where every element \(g\) has an inverse \(g^{-1}\) with
\[gg^{-1} = g{^-1}g = 1.\] For example, the group of rotations in
\(n\)-dimensional Euclidean space. When we started understanding
\(1\)-categories --- that is, categories! --- the real idea behind
groups and symmetry could be more clearly expressed. Sets have elements,
and they are either equal or not --- no two ways about it. Categories
have ``objects'', and even though two objects aren't equal, they can
still be ``isomorphic''. An object can be isomorphic to itself in lots
of different ways: these are its symmetries, and the symmetries form a
group. But this is really just the tip of a still mysterious iceberg.

For example, in a \(2\)-category, even though two objects aren't equal,
or even isomorphic, they can be ``equivalent'', or maybe I should say
``2-equivalent''. This is a still more general notion of ``sameness''. I
won't try to define it just now, but I'll just note that it arises from
the fact that in a \(2\)-category one can ask whether two morphisms are
isomorphic! (For people who followed \protect\hyperlink{week49}{``Week
49''} and know some category theory, let me note that the standard
notion of equivalence of categories is a good example of this
``equivalence'' business.)

As we climb up the \(n\)-categorical ladder, this keeps going. We get
ever more subtle refinements on the notion of ``sameness'', hence ever
subtler notions of symmetry. It's all rather mind-boggling at first, but
not really very hard once you get the hang of it, and since there's lots
of evidence that \(n\)-dimensional topological quantum field theories
are related to \(n\)-categories, I think these subtler notions of
symmetry are going to be quite interesting for physics.

And now, to wax technical for a bit (skip this paragraph if the last one
made you dizzy), it's beginning to seem that the symmetry groups
physicists know and love have glorious reincarnations, or avatars if one
prefers, at these higher \(n\)-categorical levels. Take your favorite
group --- mine is \(\mathrm{SU}(2)\), which describes 3d rotational
symmetry hence angular momentum in quantum mechanics. Its category of
representations isn't just any old category, it's a symmetric monoidal
category! See the chart in \protect\hyperlink{week49}{``Week 49''} if
you forget what this is. Now, there are more general things than groups
whose categories of representations are symmetric monoidal categories
--- for example, cocommutative Hopf algebras. And there are other kinds
of Hopf algebras --- ``quasitriangular'' ones, often known as ``quantum
groups'' --- whose categories of representations are \emph{braided}
monoidal categories. The cool thing is that \(\mathrm{SU}(2)\) has an
avatar called ``quantum \(\mathrm{SU}(2)\)'' which is one of these
quantum groups. Again, eyeball the chart in
\protect\hyperlink{week49}{``Week 49''}. Symmetric monoidal categories
are a special kind of \(4\)-categories, which is why they show up so
much in 4d physics, while braided monoidal categories are a special kind
of \(3\)-categories, which is why they (and quantum groups) show up in
3d physics. For example, quantum \(\mathrm{SU}(2)\) shows up in the
study of 3d quantum gravity (see \protect\hyperlink{week16}{``Week
16''}). Now the even cooler thing is that while a quaistriangular Hopf
algebra is a set with a bunch of operations, there is a souped-up
gadget, a ``quasitriangular Hopf category'', which is a \emph{category}
with an analogous bunch of operations, and these have a \emph{category}
of representations, but not just any old \(2\)-category, but in fact a
\emph{braided monoidal \(2\)-category}. If you again eyeball the chart,
you'll see this is a special kind of \(4\)-category, so it should be
related to 4d topology and --- this is the big hope --- 4d physics. Now
the \emph{really} cool thing, which is what Crane and Frenkel show here,
is that \(\mathrm{SU}(2)\) has yet another avatar which is one of these
quasitriangular Hopf categories.

Regardless of whether it has anything to do with physics, this business
about how symmetry groups have avatars living on all sorts of rungs of
the \(n\)-categorical ladder is such beautiful math that I'm sure it's
trying to tell us something. Right now I'm trying to figure out just
what.

\begin{enumerate}
\def\labelenumi{\arabic{enumi})}
\setcounter{enumi}{2}
\tightlist
\item
  Raoul Bott and Clifford Taubes, ``On the self-linking of knots'',
  \emph{Jour. Math. Phys.} \textbf{35} (1994), 5247--5287.
\end{enumerate}
\noindent
I'd need to look at this a few more times before I could say anything
intelligent about it, but it looks to be a very exciting way of
understanding what the heck is really going on as far as Vassiliev
invariants, Feynman diagrams in Chern--Simons theory, and so on are
concerned --- especially if you wanted to generalize it all to higher
dimensions.

Alas, I'm getting worn out, so let me simply \emph{list} a few more
papers, which are every bit as fun as the previous ones\ldots{} no
disrespect intended\ldots{} I just have to call it quits soon. As
always, it should be clear that what I write about and what I don't is
purely a matter of whim, caprice, chance, and my own ignorance.

\begin{enumerate}
\def\labelenumi{\arabic{enumi})}
\setcounter{enumi}{3}
\item
  Christopher King and Ambar Sengupta, ``An explicit description of the
  symplectic struture of moduli spaces of flat connections'',
  \emph{Jour. Math. Phys.} \textbf{35} (1994), 5338--5353.

  Christopher King and Ambar Sengupta, ``The semiclassical limit of the
  two-dimensional quantum Yang--Mills model'', \emph{Jour. Math. Phys.}
  \textbf{35} (1994), 5354--5363.
\item
  D. J. Thouless, ``Topological interpretations of quantum Hall
  conductance'', \emph{Jour. Math. Phys.} \textbf{35} (1994),
  5362--5372.
\item
  J. Bellisard, A. van Elst, and H. Schulz-Baldes, ``The noncommutative
  geometry of the quantum Hall effect'', \emph{Jour. Math. Phys.}
  \textbf{35} (1994), 5373--5451.
\item
  Steve Carlip and R. Cosgrove, ``Topology change in (2+1)-dimensional
  gravity'', \emph{Jour. Math. Phys.} \textbf{35} (1994), 5477--5493.
\end{enumerate}

\end{document}